\documentclass[]{article}

\usepackage{fullpage}
\usepackage{here}
\usepackage{graphicx}
\usepackage{times}
\usepackage{amscd}
\usepackage{amssymb}
\usepackage{amsmath}
\usepackage{mathrsfs}
\usepackage{verbatim}

\newcommand{\beq}{\begin{equation}}
\newcommand{\eeq}{\end{equation}}
\newcommand{\bea}{\begin{eqnarray}}
\newcommand{\eea}{\end{eqnarray}}

\newcommand{\Beq}{\begin{displaymath}}
\newcommand{\Eeq}{\end{displaymath}}
\newcommand{\Bea}{\begin{eqnarray*}}
\newcommand{\Eea}{\end{eqnarray*}}

\protect
\protect
\protect
\protect
\protect
\protect
\protect
\protect
\protect
\protect
\protect\newtheorem{*Exercise}[Exercise]{$\ast\,$Exercise}
\protect
\protect
\protect
\protect
\protect
\protect
\protect
\protect
\protect
\protect
\protect
\protect
\protect
\protect
\protect
\protect
\protect
\protect
\protect
\protect
\protect
\protect
\protect\newtheorem{corollary}{Corollary}
\protect\newtheorem{definition}{Definition}
\protect
\protect
\protect\newtheorem{*exercise}[exercise]{$\ast\,$Exercise}
\protect
\protect
\protect
\protect\newtheorem{lemma}{Lemma}
\protect
\protect
\protect
\protect
\protect
\protect
\protect
\protect\newtheorem{remark}{Remark}
\protect
\protect
\protect\newtheorem{theorem}{Theorem}
\protect

\def\bbbc{{\mathchoice {\setbox0=\hbox{$\displaystyle\rm C$}\hbox{\hbox
to0pt{\kern0.4\wd0\vrule height0.9\ht0\hss}\box0}}
{\setbox0=\hbox{$\textstyle\rm C$}\hbox{\hbox
to0pt{\kern0.4\wd0\vrule height0.9\ht0\hss}\box0}}
{\setbox0=\hbox{$\scriptstyle\rm C$}\hbox{\hbox
to0pt{\kern0.4\wd0\vrule height0.9\ht0\hss}\box0}}
{\setbox0=\hbox{$\scriptscriptstyle\rm C$}\hbox{\hbox
to0pt{\kern0.4\wd0\vrule height0.9\ht0\hss}\box0}}}}
\def\bbbq{{\mathchoice {\setbox0=\hbox{$\displaystyle\rm
Q$}\hbox{\raise
0.15\ht0\hbox to0pt{\kern0.4\wd0\vrule height0.8\ht0\hss}\box0}}
{\setbox0=\hbox{$\textstyle\rm Q$}\hbox{\raise
0.15\ht0\hbox to0pt{\kern0.4\wd0\vrule height0.8\ht0\hss}\box0}}
{\setbox0=\hbox{$\scriptstyle\rm Q$}\hbox{\raise
0.15\ht0\hbox to0pt{\kern0.4\wd0\vrule height0.7\ht0\hss}\box0}}
{\setbox0=\hbox{$\scriptscriptstyle\rm Q$}\hbox{\raise
0.15\ht0\hbox to0pt{\kern0.4\wd0\vrule height0.7\ht0\hss}\box0}}}}
\def\bbbt{{\mathchoice {\setbox0=\hbox{$\displaystyle\rm
T$}\hbox{\hbox to0pt{\kern0.3\wd0\vrule height0.9\ht0\hss}\box0}}
{\setbox0=\hbox{$\textstyle\rm T$}\hbox{\hbox
to0pt{\kern0.3\wd0\vrule height0.9\ht0\hss}\box0}}
{\setbox0=\hbox{$\scriptstyle\rm T$}\hbox{\hbox
to0pt{\kern0.3\wd0\vrule height0.9\ht0\hss}\box0}}
{\setbox0=\hbox{$\scriptscriptstyle\rm T$}\hbox{\hbox
to0pt{\kern0.3\wd0\vrule height0.9\ht0\hss}\box0}}}}
\def\bbbs{{\mathchoice
{\setbox0=\hbox{$\displaystyle     \rm S$}\hbox{\raise0.5\ht0\hbox
to0pt{\kern0.35\wd0\vrule height0.45\ht0\hss}\hbox
to0pt{\kern0.55\wd0\vrule height0.5\ht0\hss}\box0}}
{\setbox0=\hbox{$\textstyle        \rm S$}\hbox{\raise0.5\ht0\hbox
to0pt{\kern0.35\wd0\vrule height0.45\ht0\hss}\hbox
to0pt{\kern0.55\wd0\vrule height0.5\ht0\hss}\box0}}
{\setbox0=\hbox{$\scriptstyle      \rm S$}\hbox{\raise0.5\ht0\hbox
to0pt{\kern0.35\wd0\vrule height0.45\ht0\hss}\raise0.05\ht0\hbox
to0pt{\kern0.5\wd0\vrule height0.45\ht0\hss}\box0}}
{\setbox0=\hbox{$\scriptscriptstyle\rm S$}\hbox{\raise0.5\ht0\hbox
to0pt{\kern0.4\wd0\vrule height0.45\ht0\hss}\raise0.05\ht0\hbox
to0pt{\kern0.55\wd0\vrule height0.45\ht0\hss}\box0}}}}
\def\bbbz{{\mathchoice {\hbox{$\sf\textstyle Z\kern-0.4em Z$}}
{\hbox{$\sf\textstyle Z\kern-0.4em Z$}}
{\hbox{$\sf\scriptstyle Z\kern-0.3em Z$}}
{\hbox{$\sf\scriptscriptstyle Z\kern-0.2em Z$}}}}

\newcommand{\ve}{\varepsilon}

\newcommand{\ul}{\textbf}

\newcommand{\R}{\mathbb{R}}

\newcommand{\KKK}{\mbox{\boldmath$\kappa$}}
\newcommand{\TTT}{\mbox{\boldmath$\tau$}}
\newcommand{\desda}{\Leftrightarrow}

\newcommand{\ealpha}{\, \alpha}
\newcommand{\talpha}{ {\tilde{\alpha}}}
\newcommand{\ebeta}{\, \beta}
\newcommand{\tbeta}{ {\tilde{\beta}}}
\newcommand{\egamma}{\, \gamma}
\newcommand{\tgamma}{{\tilde{\gamma}}}
\newcommand{\pdx}{\partial_x}
\newcommand{\pdy}{\partial_y}
\newcommand{\pdz}{\partial_z}

\DeclareMathAlphabet\gothic{U}{euf}{m}{n}

\title{\textbf{Diffusion, Convection and Erosion on  $\R^3\rtimes S^2$ and their Application to the Enhancement of Crossing Fibers}}

\author{\mbox{Remco Duits$^{1,2}$, Eric Creusen$^{1,2}$, Arpan Ghosh$^{1,}$\footnote{The European Commission is gratefully acknowledged for supporting Arpan Ghosh under the Initial Training Network-FIRST, agreement No. PITN-GA-2009-238702.}, Tom Dela Haije$^{2}$} \\ \\
\centerline{{\small Eindhoven University of Technology, P.O. Box 513, 5600 MB, Eindhoven, The Netherlands. }} \\
\centerline{{\small Department of Mathematics and Computer Science$^{1}$, CASA Applied analysis, }} \\
\centerline{{\small Department of Biomedical Engineering$^{2}$, BMIA Biomedical image analysis.}}
\\
\centerline{{\small e-mail: R.Duits@tue.nl,
E.J.Creusen@tue.nl, A.Ghosh@tue.nl, T.C.J.Dela.Haije@student.tue.nl,}}
}
\date{\today}

\begin{document}

\maketitle
\begin{abstract}
\noindent

In this article we study both left-invariant (convection-)diffusions and left-invariant Hamilton-Jacobi equations (erosions) on the space $\R^3 \rtimes S^2$ of 3D-positions and orientations naturally embedded in the group $SE(3)$ of $3D$-rigid body movements. The general motivation for these (convection-)diffusions and erosions is to obtain crossing-preserving fiber enhancement on probability densities defined on the space of positions and orientations.

The linear left-invariant (convection-)diffusions are forward Kolmogorov equations of Brownian motions on $\R^{3}\rtimes S^2$ and can be solved by $\R^3 \rtimes S^{2}$-convolution
with the corresponding Green's functions or by a finite difference scheme. The left-invariant Hamilton-Jacobi equations are Bellman equations of cost processes on $\R^{3} \rtimes S^2$ and they are solved by a morphological
$\R^3 \rtimes S^{2}$-convolution with the corresponding Green's functions. We will reveal the remarkable analogy between these erosions/dilations and diffusions. Furthermore, we consider
pseudo-linear scale spaces on the space of positions and orientations that combines dilation and diffusion in a single evolution.

In our design and analysis for appropriate linear, non-linear, morphological and pseudo-linear scale spaces on $\R^3 \rtimes S^2$ we employ the underlying differential geometry on $SE(3)$,
where the frame of left-invariant vector fields serves as a moving frame of reference. Furthermore, we will present new and simpler finite difference schemes for our diffusions,
which are clear improvements of our previous finite difference schemes.

We apply our theory to the enhancement of fibres in
magnetic resonance imaging (MRI) techniques (HARDI and DTI) for imaging water diffusion processes in fibrous tissues such as brain white matter and muscles.
We provide experiments of our crossing-preserving (non-linear) left-invariant evolutions on neural images of a human brain containing crossing fibers.
\end{abstract}
\begin{small}
\end{small}

\newpage
\tableofcontents

\newpage

\section{Introduction}

Diffusion-Weighted Magnetic Resonance Imaging (DW-MRI) involves magnetic
resonance techniques for non-invasively measuring
local water diffusion inside tissue. Local water diffusion profiles reflect underlying
biological fiber structure of the imaged area. For instance in the brain, diffusion is less constrained
parallel to nerve fibers than
perpendicular to them. In this way, the measurement of water diffusion gives information about the fiber structures present, which
allows extraction of clinical information from these scans.

The diffusion of water molecules in tissue over time is described by a
transition density \mbox{$\ul{y} \mapsto p_{t}(Y_{t}=\ul{y} \;|\; Y_{0}=\ul{y}_{0})$}
that reflects the probability density of finding a water molecule at time $t>0$ at position
$\ul{y} \in \R^{3}$ given that it
started at $\ul{y} \in \R^{3}$ at time $t=0$. Here the family of random variables $(Y_{t})_{t \geq 0}$
(stochastic process) with joint state space $\R^3$ reflects the distribution of water molecules over time.
The function $p_{t}$ can be directly related to MRI signal attenuation of Diffusion-Weighted image sequences, so
can be estimated given enough measurements. The exact methods to do this are described by e.g. Alexander
\cite{Alexander}. Diffusion tensor imaging (DTI), introduced by Basser et al. \cite{Basser},
assumes that $p_{t}$ can be described for each
position $\ul{y} \in \R^{3}$ by an anisotropic Gaussian function. So
\[
p_{t}(Y_{t}=\ul{y}' \;|\; Y_{0}=\ul{y})= \frac{1}{ (4\pi t)^{\frac{3}{2}}
\sqrt{|\textrm{det}D(\ul{y})|}}
e^{-\frac{(\ul{y}'-\ul{y})^{T} (D(\ul{y}))^{-1} (\ul{y}'-\ul{y})}{4t}}\ ,
\]
where $D$ is a tensor field of positive definite symmetric tensors on $\R^{3}$.
In a DTI-visualization one plots the surfaces
\begin{equation} \label{DTIglyph}
\ul{y}+ \{\ul{v} \in \R^{3}\; |\; \ul{v}^{T} D^{-1}(\ul{y}) \ul{v} = \mu^2\},
\end{equation}
where $\mu>0$ is fixed and $\ul{y} \in \Omega$ with $\Omega$ some compact subset of $\R^{3}$.
From now on we refer to these ellipsoidal surfaces as DTI-glyphs.
The drawback of this anisotropic Gaussian approximation in DTI is the limited angular resolution of the
corresponding probability density $U:\R^{3} \rtimes S^{2} \to \R^{+}$ on positions and orientations
\begin{equation} \label{DTItoHARDI}
U(\ul{y},\ul{n})= \frac{3}{4\pi \, \int_{{\Omega}}\textrm{trace}\{D(\ul{y}')\}{\rm d}\ul{y}'} \; \ul{n}^{T} D(\ul{y}) \ul{n}, \ \ \ul{y}\in \R^{3}, \ul{n} \in S^{2},
\end{equation}
and thereby unprocessed DTI is not capable of representing
crossing, kissing or diverging fibers, cf. \cite{Alexander}.  

\begin{figure}[h] \vspace{-0.15cm}\mbox{}
\centerline{
\hfill
\includegraphics[width=0.7\hsize]{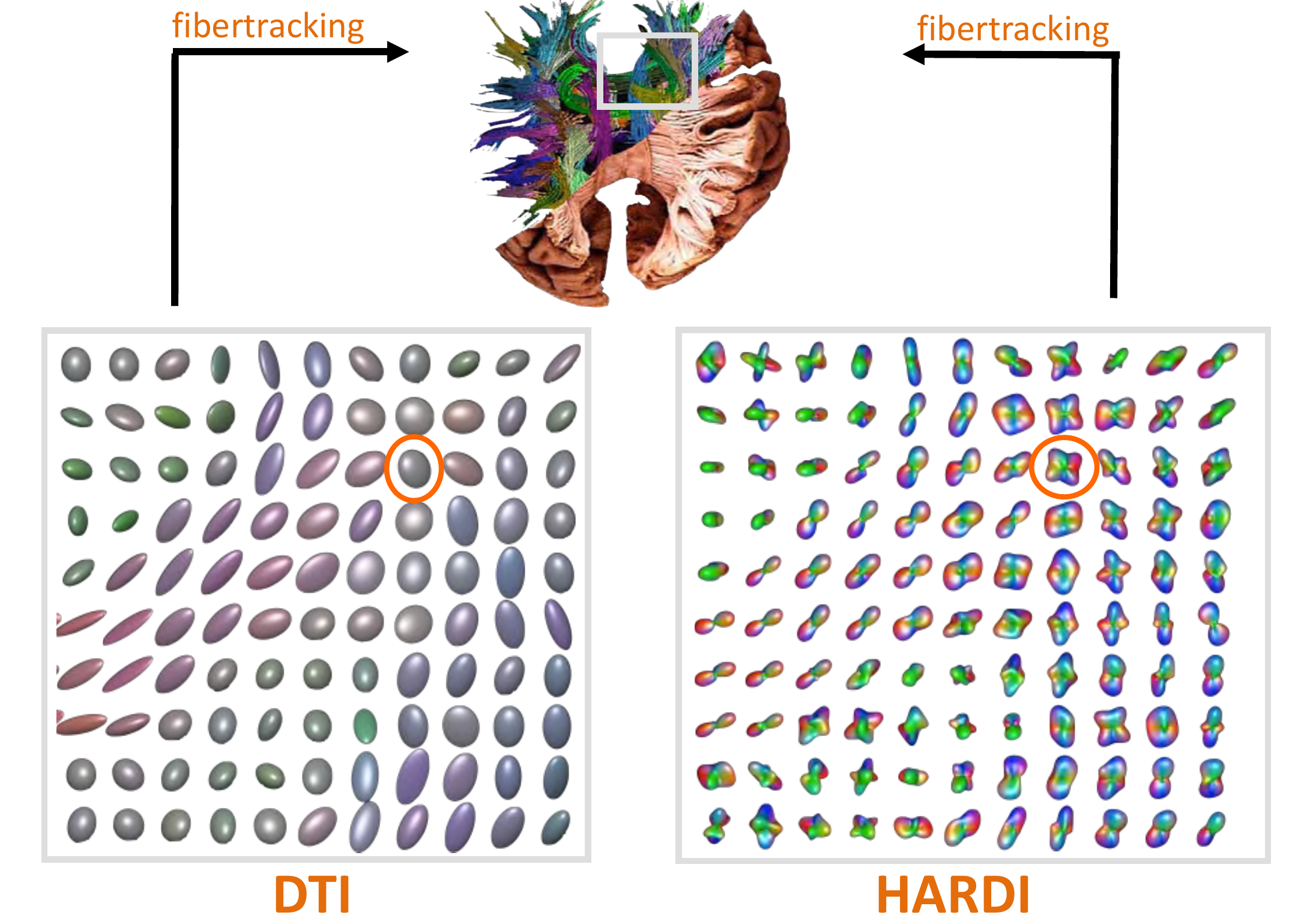}
\hfill
}
\vspace{-0.4cm}\mbox{}
\caption{This figure shows glyph visualizations of HARDI and DTI-images of a 2D-slice in the brain where neural fibers in the corona radiata cross with neural fibers in the corpus callosum.
Here DTI and HARDI are visualized differently; HARDI is visualized
according to Def.~\ref{def:viz}, whereas DTI is visualized using Eq.~(\ref{DTIglyph}).
}
\label{fig:1}
\vspace{-0.3cm}\mbox{}
\end{figure}

High Angular Resolution Diffusion Imaging (HARDI) is another recent magnetic resonance imaging technique for imaging water diffusion processes in fibrous tissues such as brain white matter and muscles.
HARDI provides for each position in $\R^3$ and for each orientation
an MRI signal attenuation profile, which can be
related to the local diffusivity of water molecules in the corresponding direction.
Typically, in HARDI modeling
the Fourier transform of the estimated transition densities is considered at a fixed characteristic radius
(generally known as the \emph{$b$-value}), \cite{Descoteaux}.
As a result, HARDI images are distributions $(\ul{y}, \ul{n}) \mapsto U(\ul{y},\ul{n})$ over positions
and orientations, which are often
normalized per position. HARDI is not restricted
to functions on the 2-sphere induced by a quadratic form and is capable of reflecting
crossing information,
see Figure~\ref{fig:1}, where we visualize HARDI
by glyph visualization as defined below.
\begin{definition}\label{def:viz}
A glyph of a distribution $U:\R^3 \times S^2 \to \R^{+}$ on positions and orientations is a surface $\mathcal{S}_{\mu}(U)(\ul{y})= \{\ul{y}+\mu\, U(\ul{y},\ul{n})\; \ul{n}\; |\; \ul{n} \in S^{2} \} \subset \R^3$ for some \, $\ul{y} \in \R^3$ and $\mu>0$.
A glyph visualization of the distribution $U:\R^3 \times S^2 \to \R^{+}$ is a visualization of a field $\ul{y} \mapsto \mathcal{S}_{\mu}(U)(\ul{y})$ of glyphs, where $\mu>0$ is a suitable
constant.
\end{definition}
For the purpose of detection and visualization of biological fibers,
DTI and HARDI data enhancement should maintain fiber junctions and crossings,
while reducing high frequency noise in the joined domain of positions
and orientations.
Promising research has been done on constructing diffusion/regularization
processes on the 2-sphere defined at each spatial locus separately
\cite{Descoteaux2007b,Florack2,Florack2008b,Hess2006a}
as an essential pre-processing step for robust fiber tracking.
In these approaches position and orientation space are decoupled,
and diffusion is only performed over the angular part, disregarding spatial context.
Consequently, these methods are inadequate for spatial denoising and enhancement, and
tend to fail precisely at the interesting locations where fibres cross or bifurcate.

In contrast to previous work on diffusion of DW-MRI
\cite{Descoteaux2007b,Florack2,Florack2008b,Hess2006a,Ozarslan},
we consider both the spatial and the orientational part to be included in the {\em domain},
so a HARDI dataset is considered as a function $U:\R^3 \times S^{2} \rightarrow \R$.
Furthermore, we explicitly employ the proper underlying group structure, that naturally arises by embedding
the coupled space of positions and orientations
\[
\R^{3}\rtimes S^{2} := SE(3)/(\{0\} \times SO(2)),
\]
as the quotient of left cosets, into the group $SE(3)=\R^{3} \rtimes SO(3)$ of 3D-rigid motions.
The relevance of group theory in DTI/HARDI (DW-MRI) imaging
has also been stressed in promising and well-founded recent works \cite{Gur,Gur2,Fletcher}.
However these works rely on bi-invariant Riemannian metrics on compact groups (such as $SO(3)$) and in our
case the group $SE(3)$ is neither compact nor does it permit a bi-invariant metric \cite{Arsigny,DuitsAMS2}.

Throughout this article we use the following identification between the DW-MRI data $(\ul{y},\ul{n}) \to U(\ul{y},\ul{n})$ and functions
$\tilde{U}:SE(3) \to \R$ given by
\begin{equation} \label{OS}
\tilde{U}(\ul{y},R)=U(\ul{y}, R\ul{e}_{z}).
\end{equation}
By definition one has $\tilde{U}(\ul{y},R R_{\ul{e}_{z},\alpha})=\tilde{U}(\ul{y},R)$
for all $\alpha \in [0,2\pi)$, where $R_{\ul{e}_{z},\alpha}$ is the
counterclockwise rotation around $\ul{e}_{z}$ by angle $\alpha$.

\begin{figure}
\centerline{\includegraphics[width=0.9\hsize]{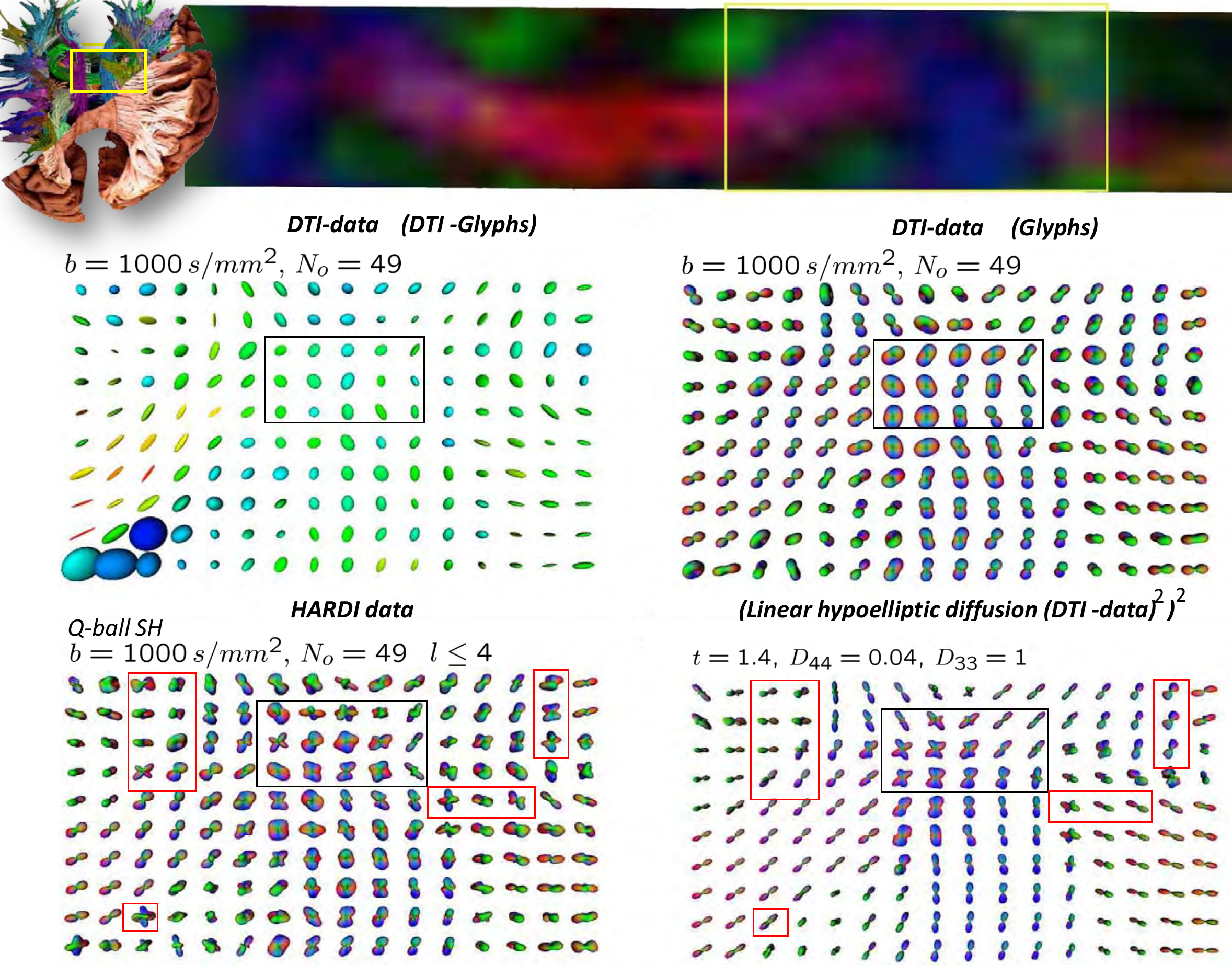}}
\caption{DTI and HARDI data containing fibers of the corpus callosum and the corona radiata in a human brain
with $b$-value $1000s/mm^2$ on voxels of $(2mm)^3$, cf.~\cite{Prckovska}. We visualize a $10 \times 16$-slice (162 samples on $S^{2}$ using icosahedron tessellations) of interest from
$104\times104\times10 \times (162\times 3)$
datasets. Top row, region of interest with fractional anisotropy
intensities with colorcoded DTI-principal directions.
Middle row, DTI data $U$
($49$ scanned orientations) vizualized according to Eq.~(\ref{DTIglyph}) respectively
Def.~\ref{def:viz}.
Bottom row: HARDI data (Q-ball with $l\leq 4$, \cite{Descoteaux}) of the same region, processed DTI data $\Phi_{t}(U)$.
For the sake of visualization we applied min-max normalization of $\ul{n} \mapsto \Phi_{t}(U)(\ul{y},\ul{n})$ for all positions $\ul{y}$.
For parameter settings of the hypo-elliptic diffusion operator $\Phi_{t}$, see Section~\ref{ch:LDSE3}, Eq.~(\ref{enhancement}).
}\label{Fig:2}
\end{figure}

\begin{figure}
\centerline{\includegraphics[width=0.9\hsize]{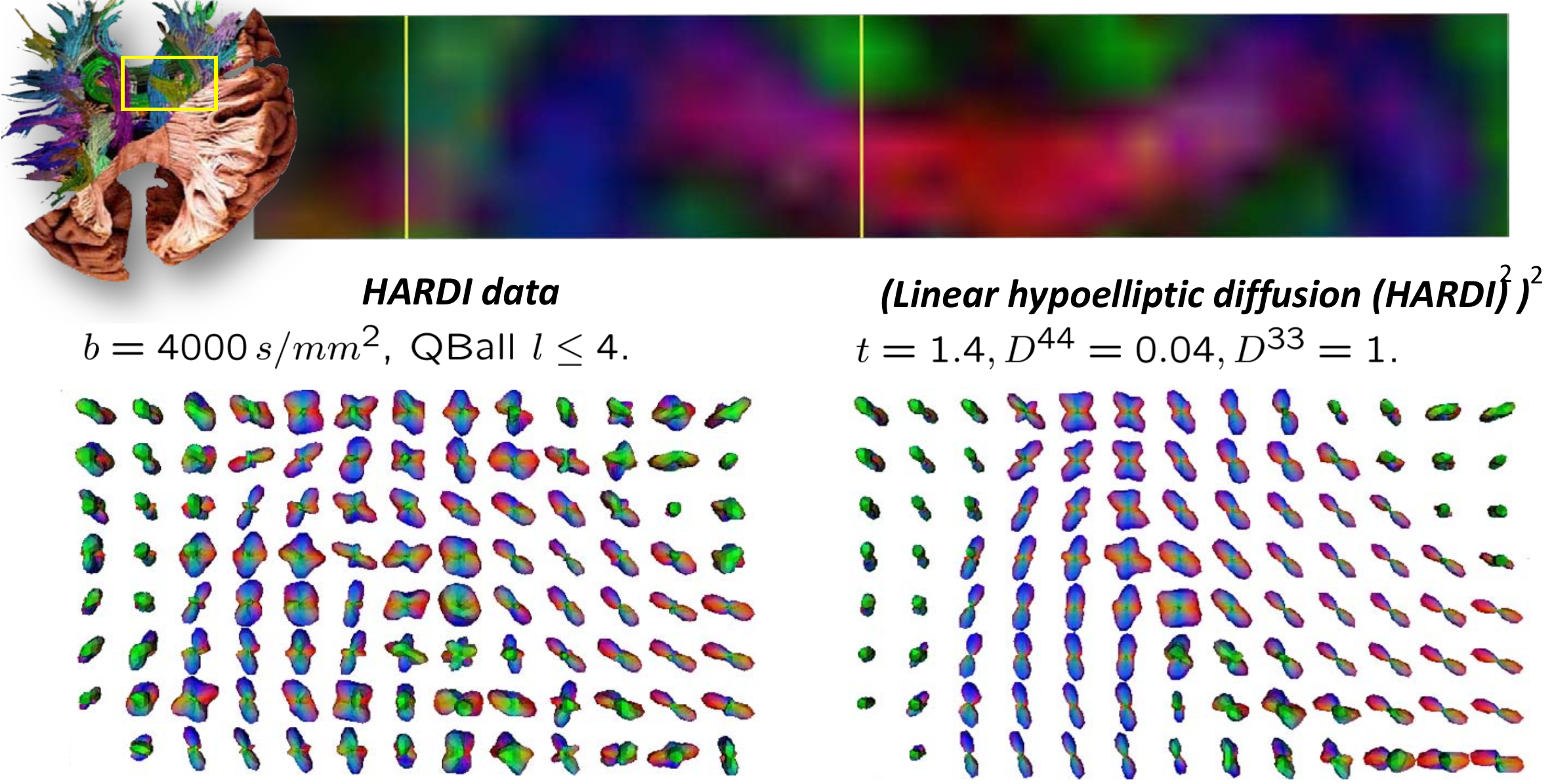}}
\caption{Same settings as Fig:2, except for different $b$-value and a different region of interest
and the hypo-elliptic diffusion, Eq.~(\ref{enhancement}), is applied to the
HARDI dataset.
}\label{Fig:3}
\end{figure}

In general the advantage of our approach on $SE(3)$ is that we can enhance the original HARDI/DTI data using
simultaneously orientational and spatial neighborhood information,
which potentially leads to improved enhancement and detection algorithms, \cite{DuitsIJCV2010,Rodrigues,Prckovska}.
See Figure~\ref{Fig:2} where fast practical implementations \cite{Rodrigues} of the theory developed in \cite[ch:8.2]{DuitsIJCV2010} have been applied.
The potential clinical impact is evident: By hypo-elliptic diffusions on $\R^{3}\rtimes S^{2}$ one
can generate distributions from DTI that are similar to HARDI-modeling as recently reported by Pr$\check{c}$kovska et al. \cite{Prckovska}. This allows a reduction of scanning
directions in areas where the random walks processes that underly hypo-elliptic diffusion \cite[ch:4.2]{DuitsIJCV2010} on $\R^{3}\rtimes S^2$ yield
reasonable fiber extrapolations. Experiments on neural DW-MRI images containing crossing fibers
of the corpus callosum and corona radiata show that extrapolation of DTI (via hypo-elliptic diffusion)
can cope with HARDI for a whole range of reasonable $b$-values, \cite{Prckovska}. See Figure \ref{Fig:2}.
However, on the locations of crossings HARDI in principle produces more detailed information
than extrapolated DTI and application of
the same hypo-elliptic diffusion on HARDI removes spurious crossings that arise in HARDI,
see Figure \ref{Fig:3} and the recent work \cite{Rodrigues}.

In this article we will build on the recent previous work \cite{Prckovska,DuitsIJCV2010,Rodrigues},
and we address the following open issues that arise immediately from these works:
\begin{itemize}
\item Can we adapt the diffusion on $\R^{3}\rtimes S^{2}$ locally to the initial HARDI/DTI image?
\item Can we apply left-invariant Hamilton-Jacobi equations (erosions) to sharpen the data without grey-scale
transformations (squaring) needed in our previous work ?
\item Can we classify the viscosity solutions of these left-invariant Hamilton-Jacobi equations?
\item Can we find analytic approximations for viscosity solutions of these left-invariant Hamilton-Jacobi equations
on $\R^{3}\rtimes S^{2}$, likewise the analytic approximations we derived for linear left-invariant
diffusions, cf. \cite[ch:6.2]{DuitsIJCV2010}?
\item Can we relate alternations of greyscale transformations and linear diffusions to alternations of
linear diffusions and erosions?
\item The resolvent Green's functions of the direction process on $\R^{3}\rtimes S^{2}$ contain singularities at the origin,
\cite[ch:~6.~1.~1,Fig.~8]{DuitsIJCV2010}.
Can we overcome this complication in our algorithms for the direction process and can we analyze iterations of
multiple time integrated direction process to control the filling of gaps?
\item Can we combine left-invariant diffusions and left-invariant dilations/erosions in a single pseudo-linear scale space
on $\R^{3}\rtimes S^{2}$, generalizing
the approach by Florack et al. \cite{Florackpseudo} for greyscale images to DW-MRI images?
\item Can we avoid spherical harmonic transforms and the sensitive regularization parameter $t_{reg}$ \cite[ch:7.1,7.2]{DuitsIJCV2010} in our finite difference schemes and obtain
both faster and simpler numerical approximations of the left-invariant vector fields ?
\end{itemize}
To address these issues, we introduce besides linear scale spaces,
morphological scale spaces and pseudo-linear scale spaces
 $(\ul{y},\ul{n},t) \mapsto W(\ul{y},\ul{n},t)$, for all $\ul{y} \in \R^{3},\ul{n} \in S^{2},t>0$,
 defined on $(\R^{3}\rtimes S^{2}) \times \R^{+}$, where we
use the
input DW-MRI image 
as initial condition $W(\ul{y},\ul{n},0)=U(\ul{y},\ul{n})$.

To get a preview of how these approaches perform on the same neural DTI dataset (different slice) considered in
\cite{Prckovska}, see Fig.~\! \ref{Fig:4},
where we used
{\small
\begin{equation}\label{Veq}
\mathcal{V}(U)(\ul{y},\ul{n})= \left(\frac{U(\ul{y},\ul{n})- U_{min}(\ul{y})}{U_{max}(\ul{y})- U_{min}(\ul{y})}\right)^2, \textrm{ with }U_{\textrm{$\underset{max}{min}$}}(\ul{y})=  \textrm{$\underset{\max}{\min}$}\{U(\ul{y},\ul{n})\;|\; \ul{n} \in S^{2}\}.
\end{equation}
}

\begin{figure}
\centerline{\includegraphics[width=0.95\hsize]{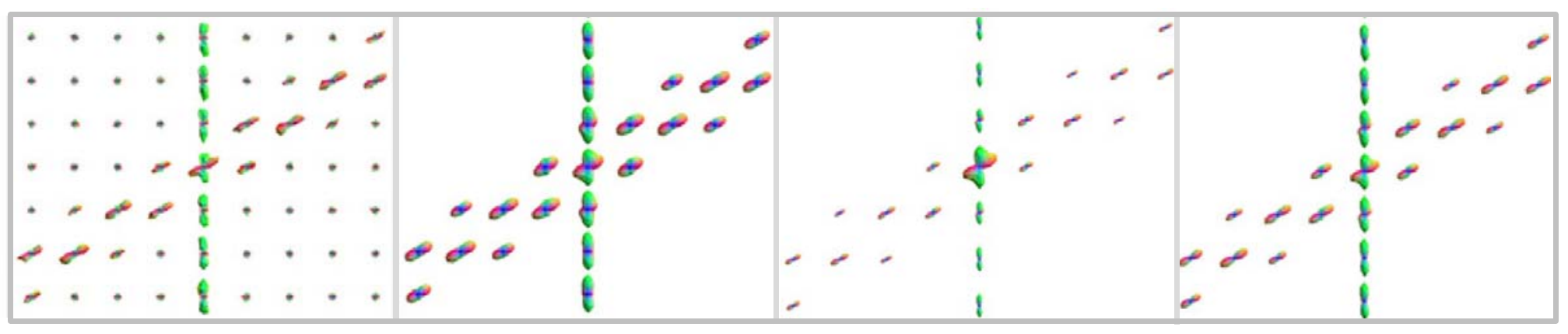}}
\caption{From left to right. Noisy artificial dataset, output diffused dataset (thresholded), squared output diffused dataset as in \cite{Prckovska,DuitsIJCV2010,Rodrigues},
$\R^{3}\rtimes S^2$-eroded output, Eq.~\!(\ref{Hamdi2}), diffused dataset, Eq.~\!(\ref{enhancement}). }\label{Fig:5}
\end{figure}

\begin{figure}
\centerline{\includegraphics[width=1.05\hsize]{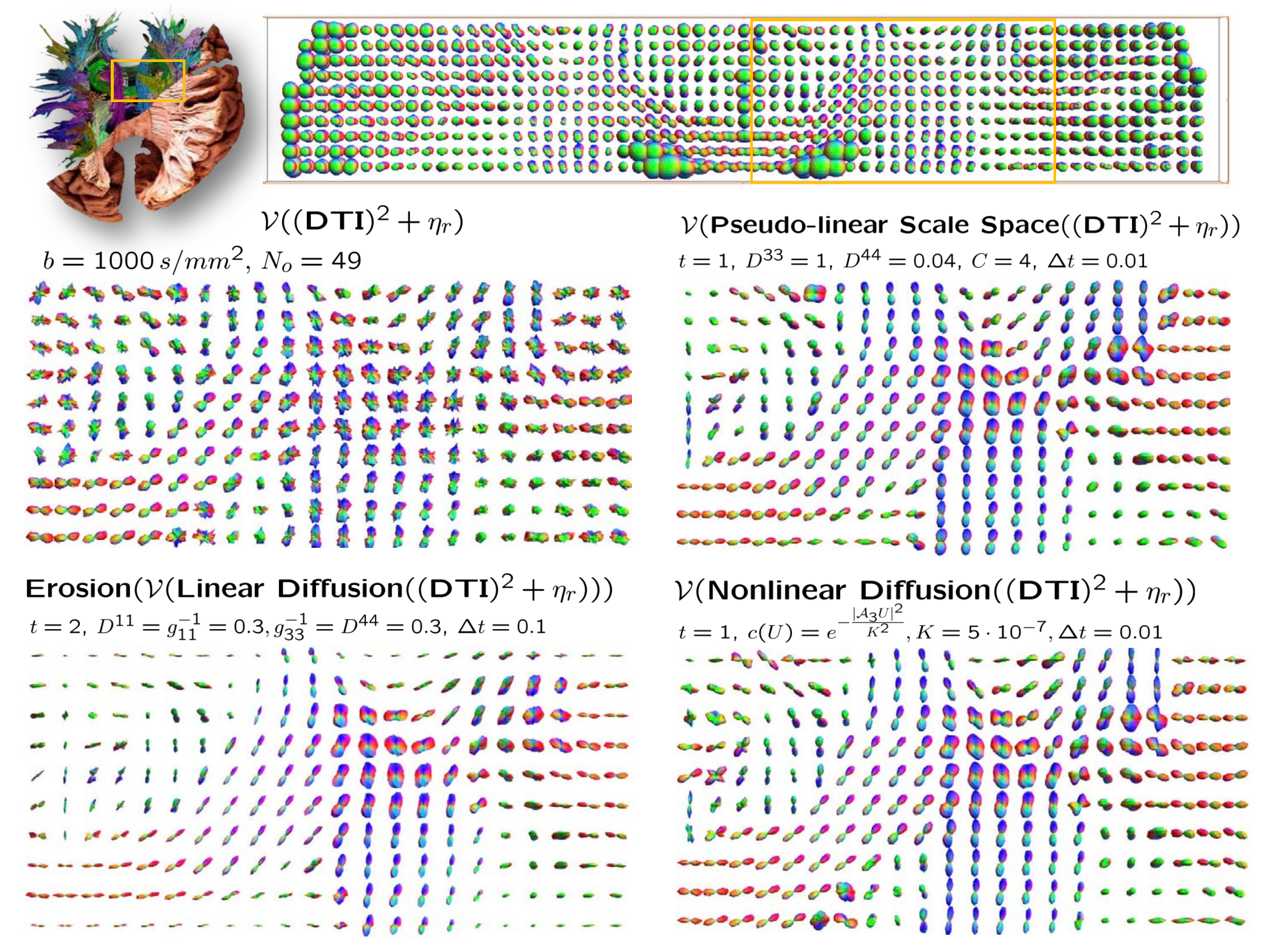}}
\caption{DTI data of corpus callosum and corona radiata fibers in a human brain
with $b$-value $1000s/mm^2$ on voxels of $(2mm)^3$. Top row: DTI-visualization according to Eq.~(\ref{DTIglyph}).
The yellow box contains $13 \times 22 \times 10$ glyphs with $162$ orientations of the input DTI-data depicted in the left image of the middle row.
This input-DTI image $U$ is visualized using Eq.~(\ref{DTItoHARDI}) and Rician noise
$\eta_{r}$ \cite[Eq.~90]{DuitsIJCV2010} with $\sigma=10^{-4}$ has been included.
Operator $\mathcal{V}$ is defined in Eq.~(\ref{Veq}). Middle row, right: output of finite difference scheme pseudo-linear scale space
(for parameters see Section \ref{ch:pseudo}).  Bottom row, left: output erosion, Eq.~(\ref{Hamdi1})  after
hypo-elliptic diffusion,
Eq.~(\ref{enhancement}) with ($D^{44}=0.04$, $D^{33}=1$, $t=1$), right: output of non-linear adaptive diffusions
discussed in \cite{Creusen} with $D^{44}=0.01$.
All evolutions commute with rotations and translations 
and are implemented by finite differences (Section \ref{ch:num}) with time $t$ and stepsize $\Delta t$.
}\label{Fig:4}
\end{figure}

\begin{figure}
\centerline{
\hfill
\includegraphics[width=0.85\hsize]{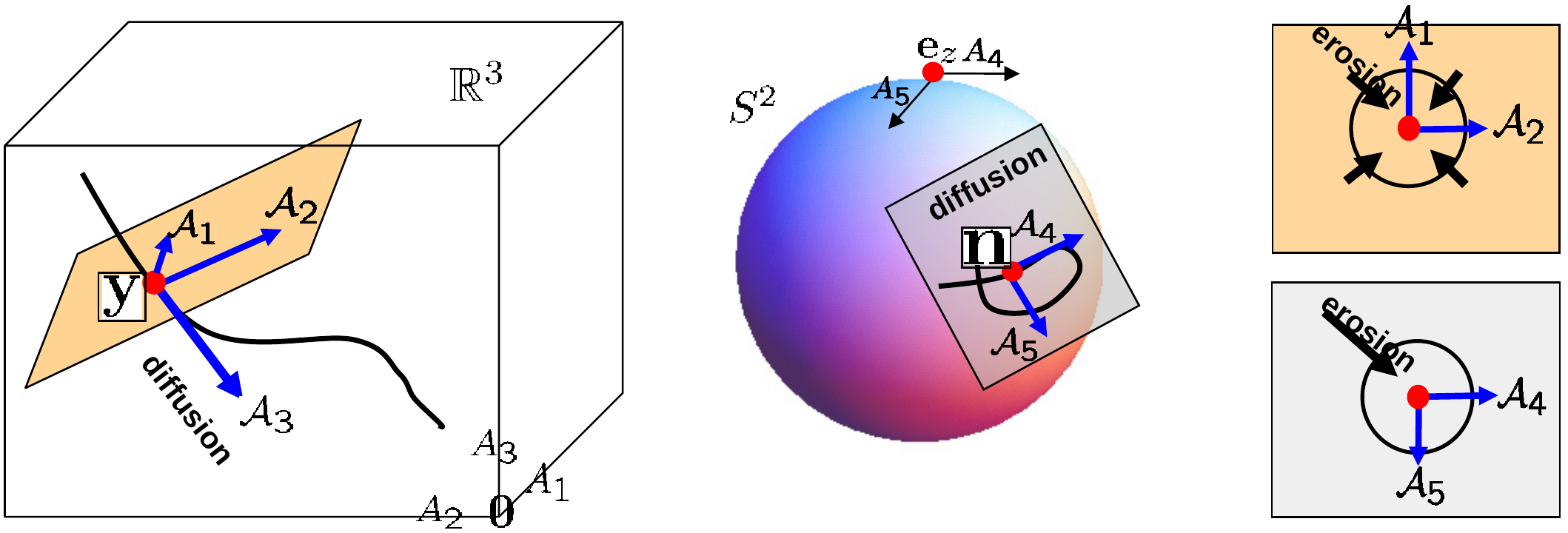}
\hfill
}
\vspace{-0.6cm}\mbox{}
\caption{{\small A curve $[0,1] \ni s \mapsto \gamma(s)=(\ul{x}(s),\ul{n}(s))\to \R^{3}\rtimes S^{2}$ consists of a spatial part $s \mapsto \ul{x}(s)$ (left) and an angular part
$s \mapsto \ul{n}(s)$ (right). Along this curve we have the moving frame of reference $\{\left.\mathcal{A}_{i}\right|_{\tilde{\gamma}(s)}\}_{i=1}^{5}$
with $\tilde{\gamma}(s)=(\ul{x}(s),R_{\ul{n}(s)})$ where $R_{\ul{n}(s)} \in SO(3)$ is \emph{any} rotation such that $R_{\ul{n}(s)}\ul{e}_{z}=\ul{n}(s) \in S^{2}$.
Here $\mathcal{A}_{i}$, with $A_{i}=\left.\mathcal{A}_{i}\right|_{(\ul{0},I)}$ denote the left-invariant vector fields in $SE(3)$.
To ensure
that the diffusions and erosions do not depend on the choice $R_{\ul{n}(s)} \in SO(3)$, Eq.~(\ref{Rn}),  these left-invariant evolution equations must be isotropic in the tangent planes
$\textrm{span}\{\mathcal{A}_{1},\mathcal{A}_{2}\}$ and $\textrm{span}\{\mathcal{A}_{4},\mathcal{A}_{5}\}$. Diffusion/convection primarily takes place along $\mathcal{A}_{3}$ in space
and (outward) in the plane $\textrm{span}\{\mathcal{A}_{4},\mathcal{A}_{5}\}$ tangent to $S^2$. Erosion takes place both inward in the tangent plane $\textrm{span}\{\mathcal{A}_{1},\mathcal{A}_{2}\}$
in space 
and inward
in the plane $\textrm{span}\{\mathcal{A}_{4},\mathcal{A}_{5}\}$.} }\label{fig:intuition}
\end{figure}

Typically, if linear diffusions are directly applied to DTI the fibers visible in DTI are propagated in too many directions.
Therefore we combined these diffusions with monotonic transformations in the codomain $\R^{+}$, such as squaring input and output cf. \cite{DuitsIJCV2010,Rodrigues,Prckovska}.
Visually, this produces anatomically plausible results, cf.~Fig.~\!\ref{Fig:2} and Fig.~\!\ref{Fig:3}, but does not allow large global
variations in the data. This is often problematic around ventricle areas in the brain,
where the glyphs are usually larger than those along the fibers, as can be seen in the top
row of Fig.~\ref{Fig:4}.
In order to achieve a better way of sharpening the data where global maxima do not dominate the sharpening of the data, cf. Fig.~\!\ref{Fig:5}, we propose morphological scale spaces
on $\R^{3}\rtimes S^{2}$ where transport takes place orthogonal to the fibers,
both spatially and spherically, see Fig.~\!\ref{fig:intuition}.
The result of such an erosion after application of a linear diffusion is depicted down left in Fig.~\!\ref{Fig:4},
where the diffusion has created crossings in the fibers and where the erosion has visually sharpened
the fibers. 
Simultaneous dilation and diffusion can be achieved in a pseudo-linear scale space that conjugates a diffusion with a specific grey-value transformation. An experiment of
applying such a left-invariant pseudo-linear scale space to DTI-data is given up-right in Figure \ref{Fig:4}. 
Regarding the numerics of the evolutions we mainly consider the
left-invariant finite difference approach in \cite[ch:7]{DuitsIJCV2010},
as an alternative to the analytic kernel implementations in \cite[ch:8.2]{DuitsIJCV2010} and
\cite{Prckovska,Rodrigues}. Here we avoid the discrete spherical Harmonic transforms \cite{DuitsIJCV2010}, but use
fast precomputed linear interpolations instead as our finite difference schemes are of first order accuracy anyway.
Regarding the Hamilton-Jacobi equations involved in morphological and pseudo-linear scale spaces we have, akin to
the linear left-invariant diffusions \cite[ch:7]{DuitsIJCV2010}, two options:
analytic morphological $\R^{3}\rtimes S^{2}$-convolutions and
finite differences. Regarding fast computation on sparse grids the second approach is preferable.
Regarding geometric analysis the first approach is preferable.

We show that our morphological $\R^{3}\rtimes S^2$-convolutions with analytical morphological
Green's functions are the unique viscosity solutions
of the corresponding Hamilton-Jacobi equations on $\R^{3}\rtimes S^{2}$.
Thereby, we generalize the results in \cite[ch:10]{Evans}, \cite{Manfredi}
(on the Heisenberg group) to Hamilton-Jacobi equations on the space
$\R^{3} \rtimes S^{2}$ of positions and orientations.

Evolutions on HARDI-DTI must commute with all rotations and translations. Therefore evolutions
on HARDI and DTI and underlying metric (tensors) are expressed in a local frame of reference
attached to fiber fragments. This frame $\{\mathcal{A}_{1},\ldots,\mathcal{A}_{6}\}$ consists
of 6 left-invariant vector fields on $SE(3)$, given by
\begin{equation}
\mathcal{A}_{i} \tilde{U} (\ul{y},R)=\lim \limits_{h \downarrow 0} \frac{U((\ul{y},R)\, e^{h A_{i}})-
U((\ul{y},R)\,  e^{-h A_{i}})}{2h}
\end{equation}
where $\{A_{1}, \ldots, A_{6}\}$ is the basis for the Lie-algebra at the unity element and
$T_{e}(SE(3)) \ni A \mapsto e^{A} \in SE(3)$ is the exponential map in $SE(3)$ and where the
group product on $SE(3)$ is given by
\begin{equation} \label{groupproduct}
(\ul{x},R)\, (\ul{x}',R')=(\ul{x}+R\ul{x}', RR'),
\end{equation}
for all positions $\ul{x},\ul{x}' \in \R^{3}$ and rotations $R,R' \in SO(3)$.
The details will follow in Section \ref{ch:VF},
see also \cite[ch:3.3,Eq.~23--25]{DuitsIJCV2010} and \cite[ch:7]{DuitsIJCV2010} for implementation.
In order to provide a relevant intuitive preview
of this moving frame of reference
we refer to Fig.~\ref{fig:intuition}.
The associated left-invariant dual frame $\{{\rm d}\mathcal{A}^{1},\ldots,{\rm d}\mathcal{A}^{6}\}$ is uniquely
determined by
\begin{equation} \label{duals}
\langle {\rm d}\mathcal{A}^{i},\mathcal{A}_{j} \rangle:={\rm d}\mathcal{A}^{i}(\mathcal{A}_{j})=\delta_{j}^{i}, \ \ i,j=1,\ldots,6,
\end{equation}
where $\delta^{i}_{j}= 1$ if $i=j$ and zero else. Then all possible left-invariant metrics are given by
\begin{equation} \label{metrictensor}
\ul{G}_{(\ul{y},R_{\ul{n}})}= \sum \limits_{i,j=1}^{6} g_{ij}
\left.{\rm d}\mathcal{A}^{i}\right|_{(\ul{y},R_{\ul{n}})} \otimes \left.{\rm d}\mathcal{A}^{j}\right|_{(\ul{y},R_{\ul{n}})}.
\end{equation}
where $\ul{y} \in \R^{3}$, $\ul{n} \in S^{2}$ and where
$R_{\ul{n}} \in SO(3)$ is any rotation that maps $\ul{e}_{z}=(0,0,1)^{T}$ onto the normal $\ul{n} \in S^2$, i.e.
\begin{equation}\label{Rn}
R_{\ul{n}}\ul{e}_{z}=\ul{n}.
\end{equation}
and where $g_{ij} \in \mathbb{C}$.
Necessary and sufficient conditions on $g_{ij}$ to induce a well-defined left-invariant metric on the quotient $\R^{3} \rtimes S^{2}= (SE(3)/(\{0\}\times SO(2)))$
can be found in Appendix \ref{app:C}. It turns out that the matrix $[g_{ij}]$ must be constant and diagonal $g_{ij}= \frac{1}{D^{ii}} \delta_{ij}$, $i,j=1\ldots,6$ with
$D^{ii} \in \R^{+} \cup \infty$, with $D^{11}=D^{22}$, $D^{44}=D^{55}$, $D^{66}=0$. Consequently, the metric is parameterized by the values $D^{11}, D^{33}, D^{44}$
and in the sequel we write
\[
\ul{G}^{{\rm D^{11}\!\!,D^{33}\!\!,D^{44}}}:= \frac{1}{D^{11}}
({\rm d}\mathcal{A}^{1} \otimes {\rm d}\mathcal{A}^{1}+
{\rm d}\mathcal{A}^{2} \otimes {\rm d}\mathcal{A}^{2}) + \frac{1}{D^{33}}
({\rm d}\mathcal{A}^{3} \otimes {\rm d}\mathcal{A}^{3})
+ \frac{1}{D^{44}} ({\rm d}\mathcal{A}^{4} \otimes {\rm d}\mathcal{A}^{4}
+{\rm d}\mathcal{A}^{5} \otimes {\rm d}\mathcal{A}^{5})
\]
The corresponding metric tensor on the quotient $\R^{3} \rtimes S^{2}= (SE(3)/(\{0\}\times SO(2)))$ is given by
\begin{equation} \label{metrictensorR3S2}
\begin{array}{ll}
\ul{G}_{(\ul{y},\ul{n})}^{{\rm D^{11}\!\!,D^{33}\!\!,D^{44}}}(\sum \limits_{i=1}^{5}
c^{i}\left. \mathcal{A}_{i}\right|_{(\ul{y},\ul{n})},
\sum \limits_{j=1}^{5}d^{j}\left. \mathcal{A}_{i}\right|_{(\ul{y},\ul{n})}) &=
\ul{G}^{{\rm D^{11}\!\!,D^{33}\!\!,D^{44}}}_{(\ul{y},R_{\ul{n}})}(\sum \limits_{i=1}^{5}c^{i}\left. \mathcal{A}_{i}\right|_{(\ul{y},\ul{n})},
\sum \limits_{j=1}^{5}d^{j}\left. \mathcal{A}_{j}\right|_{(\ul{y},\ul{n})}) \\[7pt]
 &= \frac{1}{D^{11}} (c^{1}d^{1}+c^{2}d^{2}) +
 \frac{1}{D^{33}} (c^{3}d^{3}) +
 \frac{1}{D^{44}} (c^{4}d^{4}+c^{5}d^{5}).
\end{array}
\end{equation}
It is well-defined on $\R^{3}\rtimes S^2$
as the choice of $R_{\ul{n}}$, Eq.~(\ref{Rn}), does not matter (right-multiplication
with $R_{\ul{e}_{z},\alpha}$ boils down to rotations in the isotropic planes depicted in Figure \ref{fig:intuition})
and with the differential operators on $\R^{3}\times S^2$:
\[
\begin{array}{l}
(\left.\mathcal{A}_{j}\right|_{(\ul{y},\ul{n})}U)(\ul{y},\ul{n})=
\lim \limits_{h \to 0} \frac{U(\ul{y}+h R_{\ul{n}}\ul{e}_{j},\ul{n})-U(\ul{y}-h R_{\ul{n}}\ul{e}_{j},\ul{n})}{2h}, \\
(\left.\mathcal{A}_{3+j}\right|_{(\ul{y},\ul{n})}U)(\ul{y},\ul{n})=
\lim \limits_{h \to 0} \frac{U(\ul{y}, (R_{\ul{n}}R_{\ul{e}_{j},h})\ul{e}_z)-
U(\ul{y}, (R_{\ul{n}}R_{\ul{e}_{j},-h})\ul{e}_z)}{2h}, \ \ j=1,2,3,
\end{array}
\]
where $R_{\ul{e}_{j},h}$ denotes the counter-clockwise rotation around axis $\ul{e}_{j}$ by angle $h$,
with {\small $\ul{e}_{1}=(1,0,0)$, $\ul{e}_{2}=(0,1,0)$, $\ul{e}_{3}=(0,0,1)$},
which (except for $\mathcal{A}_{3}$) do depend on the choice of $R_{\ul{n}}$ satisfying Eq.~(\ref{Rn}).

In \cite[ch:6.2]{DuitsIJCV2010} we have analytically approximated the
hypo-elliptic diffusion kernels for both the direction process and Brownian motion on the sub-Riemannian manifold
(or contact manifold \cite{Bryantbook}) $(SE(3), {\rm d}\mathcal{A}^{1}, {\rm d}\mathcal{A}^{2}, {\rm d}\mathcal{A}^{6})$
using contraction towards a nilpotent group. For the erosions we employ a similar type of technique to analytically
approximate the erosion (and dilation) kernels that describe the growth of balls in the sub-Riemannian manifold $(SE(3), {\rm d}\mathcal{A}^{3}, {\rm d}\mathcal{A}^{6})$.

A sub-Riemannian manifold is a Riemannian manifold with the extra constraint that tangent-vectors to curves in
that Riemannian manifold are not allowed to use certain subspaces of the tangent space. For example, curves in
$(SE(3), {\rm d}\mathcal{A}^{1}, {\rm d}\mathcal{A}^{2}, {\rm d}\mathcal{A}^{6})$ are curves $\tilde{\gamma}:[0,L] \to SE(3)$
with the constraint
\begin{equation} \label{hor}
\langle  \left. {\rm d}\mathcal{A}^{1} \right|_{\tilde{\gamma}(s)}, \dot{\tilde{\gamma}}(s) \rangle=
\langle  \left.{\rm d} \mathcal{A}^{2} \right|_{\tilde{\gamma}(s)}, \dot{\tilde{\gamma}}(s) \rangle=
\langle  \left. {\rm d}\mathcal{A}^{6} \right|_{\tilde{\gamma}(s)}, \dot{\tilde{\gamma}}(s) \rangle=0,
\end{equation}
for all $s \in [0,L]$, $L>0$. Note that Eq. (\ref{hor}) implies that
\[
\dot{\tilde{\gamma}}(s)=\langle  \left. {\rm d}\mathcal{A}^{3} \right|_{\tilde{\gamma}(s)}, \dot{\tilde{\gamma}}(s) \rangle
\left. \mathcal{A}_{3}\right|_{\tilde{\gamma}(s)}+ \langle  \left. {\rm d} \mathcal{A}^{4} \right|_{\tilde{\gamma}(s)}, \dot{\tilde{\gamma}}(s) \rangle
\left. \mathcal{A}_{4}\right|_{\tilde{\gamma}(s)} + \langle  \left. {\rm d}\mathcal{A}^{5} \right|_{\tilde{\gamma}(s)}, \dot{\tilde{\gamma}}(s) \rangle
\left. \mathcal{A}_{5}\right|_{\tilde{\gamma}(s)}
\]
Curves satisfying (\ref{hor}) are called \emph{horizontal curves} and we visualized a horizontal curve in Figure \ref{fig:intuition}.
For further details on differential geometry in (sub-Riemannian manifolds within) $SE(3)$ and $\R^{3}\rtimes S^2$ see Appendices \ref{app:metric}, \ref{app:viscosity} and \ref{app:meaming}, \ref{app:C} and \ref{app:E}.

\subsection{Outline of the article}

This paper is organized as follows. In section \ref{ch:group} we will explain the embedding of the coupled space of positions and orientations $\R^{3}\rtimes S^2$ into $SE(3)$.
In section \ref{ch:tools} we explain why operators on DW-MRI must be left-invariant and we consider, as an example, convolutions on $\R^{3}\rtimes S^2$. In section \ref{ch:VF} we construct the left-invariant vector fields on $SE(3)$. In general
there are two ways of considering vector fields. Either one considers them as differential operators on smooth locally defined functions, or one considers them as tangent vectors to equivalent classes of curves.
These two viewpoints are equivalent, for a formal proof see \cite[Prop. 2.4]{Aubin}. Throughout this article we will consider them as differential operators and use them as reference frame, cf. Fig \ref{fig:intuition},
in our evolution equations in later sections. Then in Section \ref{ch:morph} we consider morphological convolutions on $\R^{3}\rtimes S^2$.

In Section \ref{ch:LDSE3} we consider all possible linear left-invariant diffusions that are solved by $\R^{3}\rtimes S^2$-convolution (Section \ref{ch:tools}) with the corresponding Green's function. Subsequently,
in Section \ref{ch:HJE} we consider their morphological counter part: left-invariant Hamilton-Jacobi equations (i.e. erosions), the viscosity solutions of which are given by morphological convolution,
cf.~Section \ref{ch:morph}, with the corresponding
(morphological) Green's function. This latter result is a new fundamental mathematical result,  a detailed proof is given in Appendix \ref{app:viscosity}.

Subsequently, in Section \ref{ch:prob} and in Section \ref{ch:CA} we provide respectively the underlying probability theory and the underlying Cartan differential geometry of the evolutions considered in
Section \ref{ch:morph} and Section \ref{ch:LDSE3}.
Then in Section \ref{ch:kernels} we derive analytic approximations for the Green's functions of both the linear and the morphological evolutions on $\R^{3}\rtimes S^{2}$.
Section \ref{ch:pseudo} deals with pseudo linear scale spaces which are evolutions that combine the non-linear generator of erosions with the generator of diffusion in a single generator.

Section \ref{ch:num} deals with the numerics of the (convection)-diffusions, the Hamilton-Jacobi equations and the pseudo linear scale spaces evolutions.
Section \ref{ch:LDOS} summarizes our work on adaptive left-invariant diffusions on DW-MRI data. for more details we refer to the Master thesis by Eric Creusen \cite{Creusen}.
Finally, we provide experiments in Section
\ref{ch:experiments}.

\section{The Embedding of $\R^{3} \rtimes S^{2}$ into $SE(3)$ \label{ch:group}}

\begin{figure}
\centerline{\includegraphics[width=0.36\hsize]{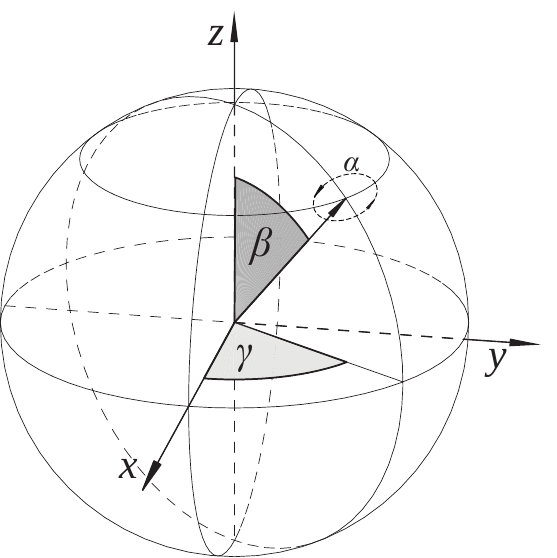}
\includegraphics[width=0.34\hsize]{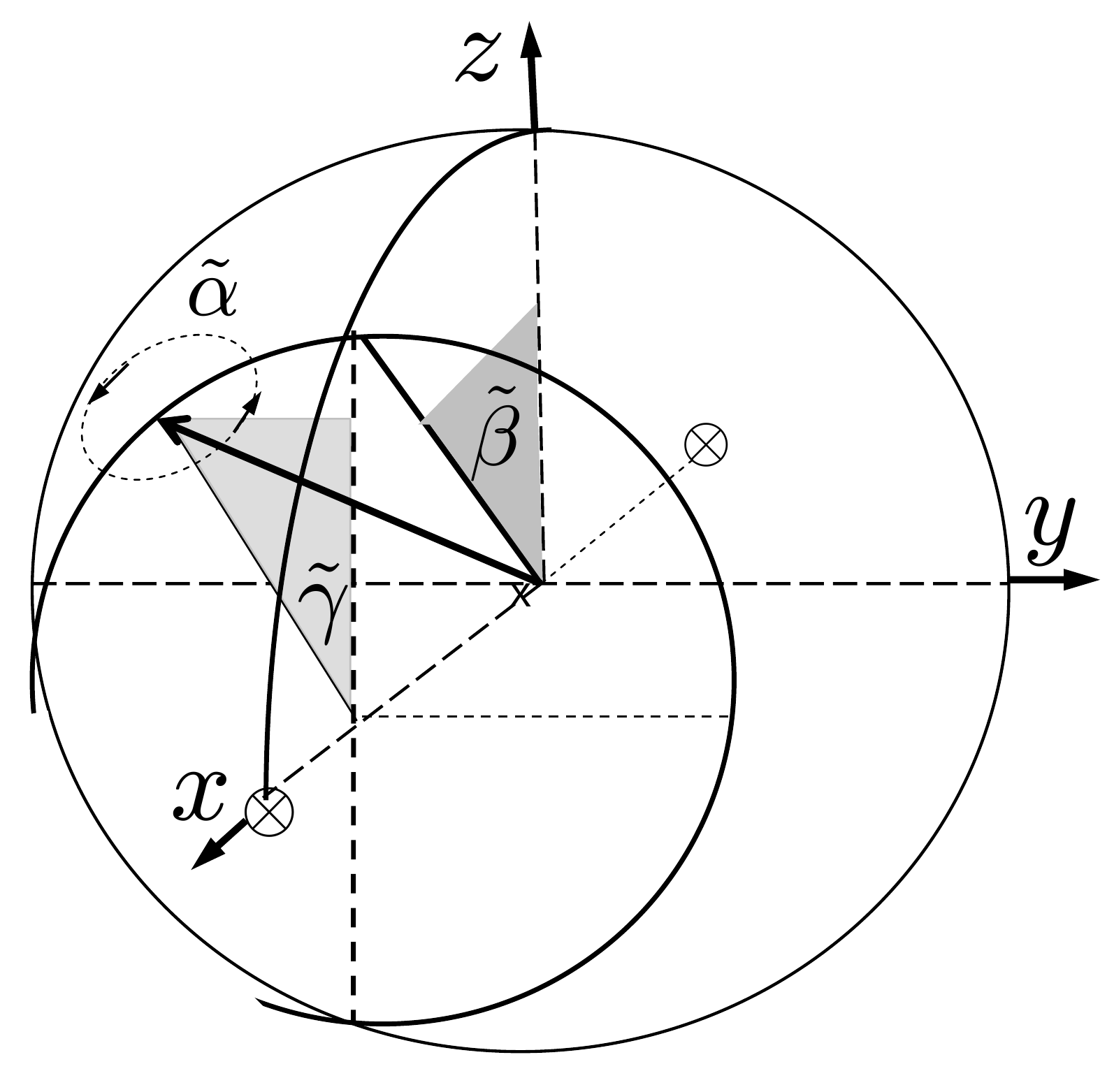}}
\caption{The two charts which together appropriately parameterize the sphere $S^{2}\equiv SO(3)/SO(2)$ where the rotation-parameters $\alpha$ and $\tilde{\alpha}$ are free. The first chart (left-image) is the common Euler-angle parametrization (\ref{firstchart}), the second chart is given by (\ref{secondchart}).
The first chart has singularities at the north and at the south-pole (inducing ill-defined parametrization of the left-invariant vector fields (\ref{eq:LeftInvVFSEthree}) at the unity element) whereas the second chart has singularities at $(\pm 1,0,0)$. }\label{fig:charts}
\end{figure}

In order to generalize our previous work on line/contour-enhancement via left-invariant diffusions on invertible orientation scores of $2D$-images we first investigate the group structure on the domain of an HARDI image. Just like orientation scores are scalar-valued functions on the
coupled space of 2D-positions and orientations, i.e. the $2D$-Euclidean motion group, HARDI images are scalar-valued functions on the coupled space of 3D-positions and orientations. This generalization involves some technicalities since the $2$-sphere
$S^{2}=\{\ul{x} \in \R^3 \; |\; \|\ul{x}\|=1\}$ is \emph{not} a Lie-group proper\footnote{If $S^{2}$ were a Lie-group then its left-invariant vector fields would be non-zero everywhere, contradicting Poincar\'{e}'s ``hairy ball theorem'' (proven by Brouwer in 1912), or more generally the Poincar\'{e}-Hopf theorem (the Euler-characteristic of an even dimensional sphere $S^{2n}$ is 2).} in contrast to the $1$-sphere $S^{1}=\{\ul{x} \in \R^2 \; |\; \|\ul{x}\|=1\}$.
To overcome this problem we will embed $\R^3 \times S^{2}$ into $SE(3)$ which is the group of $3D$-rotations and translations (i.e. the group of 3D-rigid motions).
As a concatenation of two rigid body-movements is again a rigid body movement, the product on $SE(3)$ is given by (\ref{groupproduct}).
The group $SE(3)$ is a \emph{semi-direct} product of the translation group $\R^3$ and the rotation group $SO(3)$, since it uses an isomorphism $R \mapsto (\ul{x} \mapsto R \ul{x})$ from the rotation group onto the automorphisms on $\R^3$. Therefore we write
$\R^3 \rtimes SO(3)$ rather than $\R^3 \times SO(3)$ which would yield a direct product. The groups $SE(3)$ and $SO(3)$ are not commutative. Throughout this article we will use Euler-angle parametrization for $SO(3)$,
i.e. we write every rotation as a product of a rotation around the $z$-axis, a rotation around the $y$-axis and a rotation around the $z$-axis again.
\begin{equation} \label{firstchart}
R=R_{\ul{e}_{z},\gamma}R_{\ul{e}_{y},\beta} R_{\ul{e}_{z},\alpha}\ ,
\end{equation}
where all rotations are counter-clockwise. Explicit formulas for matrices $R_{\ul{e}_{z},\gamma}, R_{\ul{e}_{y},\beta}, R_{\ul{e}_{z},\alpha}$,
are given in \cite[ch:7.3.1]{FrankenThesis}. The advantage of the Euler angle parametrization is that it directly parameterizes $SO(3)/SO(2) \equiv S^{2}$ as well. Here we recall that $SO(3)/SO(2)$ denotes the partition of all left cosets which are equivalence classes $[g]=\{h \in SO(3) \; |\; h \sim g \}=g \, SO(2)$ under the equivalence relation $g_{1}\sim g_{2} \desda g_{1}^{-1}g_{2} \in SO(2)$ where we identified $SO(2)$ with rotations around the $z$-axis and we have
\begin{equation} \label{n}
\begin{array}{l}
SO(3)/SO(2) \ni [R_{\ul{e}_{z},\gamma}R_{\ul{e}_{y},\beta}]=\{R_{\ul{e}_{z},\gamma}R_{\ul{e}_{y},\beta} R_{\ul{e}_{z},\alpha}\; |\; \alpha \in [0,2\pi)\} \leftrightarrow \\ \ul{n}(\beta,\gamma):=(\cos \gamma \sin \beta, \sin \, \gamma \sin \, \beta , \cos \, \beta)= R_{\ul{e}_{z},\gamma}R_{\ul{e}_{y},\beta} R_{\ul{e}_{z},\alpha} \ul{e}_{z} \in S^{2}.
\end{array}
\end{equation}
Like all parameterizations of $SO(3)/SO(2)$, the Euler angle parametrization suffers from the problem that there does not exist
a global diffeomorphism from a sphere to a plane. In the Euler-angle parametrization the ambiguity arises at the
north and south-poles:
{\small
\begin{equation} \label{sing1}
R_{\ul{e}_{z},\gamma}\, R_{\ul{e}_{y},\beta=0} \,R_{\ul{e}_{z},\alpha}= R_{\ul{e}_{z},\gamma-\delta}R_{\ul{e}_{y},\beta=0} R_{\ul{e}_{z},\alpha+\delta}, \textrm{ and }
R_{\ul{e}_{z},\gamma}\, R_{\ul{e}_{y},\beta=\pi}\, R_{\ul{e}_{z},\alpha}= R_{\ul{e}_{z},\gamma+\delta}\, R_{\ul{e}_{y},\beta=\pi}\, R_{\ul{e}_{z},\alpha+\delta}, \textrm{ for all }\delta \in [0,2\pi)\ .
\end{equation}
}
Consequently, we occasionally need a second chart to cover $SO(3)$;
\begin{equation}\label{secondchart}
R=R_{\ul{e}_{x},\tilde{\gamma}} R_{\ul{e}_{y},\tilde{\beta}} R_{\ul{e}_{z},\tilde{\alpha}},
\end{equation}
which again implicitly parameterizes $SO(3)/SO(2) \equiv S^{2}$ using different ball-coordinates $\tilde{\beta} \in [-\pi,\pi)$, $\tilde{\gamma} \in (-\frac{\pi}{2} ,\frac{\pi}{2})$,
\begin{equation} \label{normalOK}
\tilde{\ul{n}}(\tilde{\beta},\tilde{\gamma})=R_{\ul{e}_{x},\tilde{\gamma}} R_{\ul{e}_{y},\tilde{\beta}}\, \ul{e}_{z}=(\sin \tilde{\beta},-\cos \tilde{\beta}\, \sin \tilde{\gamma}, \cos \tilde{\beta} \, \cos \tilde{\gamma})^{T},
\end{equation}
but which has ambiguities at the intersection of the equator with the $x$-axis, \cite{DuitsIJCV2010}.
\begin{equation} \label{sing2}
R_{\ul{e}_{x},\tilde{\gamma} } R_{\ul{e}_{y},\tilde{\beta}=\pm \frac{\pi}{2}} R_{\ul{e}_{z},\tilde{\alpha}}
=R_{\ul{e}_{x},\tilde{\gamma} \mp \delta } R_{\ul{e}_{y},\tilde{\beta}=\pm \frac{\pi}{2}} R_{\ul{e}_{z},\tilde{\alpha} \pm \delta}, \textrm{ for all }\delta \in [0,2\pi)\ ,
\end{equation}
see Figure \ref{fig:charts}. Away from the intersection of the $z$ and $x$-axis with the sphere one can accomplish conversion between the two charts by solving for for either $(\tilde{\alpha},\tilde{\beta},\tilde{\gamma})$ or $(\alpha,\beta,\gamma)$ in $R_{\ul{e}_{x},\tilde{\gamma}} R_{\ul{e}_{y},\tilde{\beta}} R_{\ul{e}_{z},\tilde{\alpha}}=R_{\ul{e}_{z},\gamma}R_{\ul{e}_{y},\beta} R_{\ul{e}_{z},\alpha}$.

Now that we have explained the isomorphism $\ul{n}=R \ul{e}_{z} \in S^{2} \leftrightarrow SO(3)/SO(2) \ni [R]$ explicitly in charts, we
return to the domain of HARDI images. Considered \emph{as a set}, this domain equals the space of 3D-positions and orientations $\R^3 \times S^{2}$. However, in order to stress the fundamental embedding of the HARDI-domain in $SE(3)$ and the thereby induced (quotient) group-structure we write $\R^3 \rtimes S^{2}$, which is given by the following Lie-group quotient:
\[
\R^{3} \rtimes S^{2} := (\R^{3} \rtimes SO(3))/(\{\ul{0}\}\times SO(2)).
\]
Here the equivalence relation on the group of rigid-motions $SE(3)=\R^{3} \rtimes SO(3)$ equals
\[
(\ul{x},R) \sim (\ul{x}',R') \desda (\ul{x},R)^{-1}(\ul{x}',R') \in \{\ul{0}\} \rtimes SO(2) \desda \ul{x}=\ul{x}' \textrm{ and }R^{-1}R' \textrm{ is a rotation around z-axis }
\]
and set of equivalence classes within $SE(3)$ under this equivalence relation (i.e. left cosets) equals the space of coupled orientations and positions and is denoted by $\R^{3} \rtimes S^{2}$.

\section{Linear Convolutions on $\R^3 \rtimes S^2$ \label{ch:tools}}

In this article we will consider convection-diffusion operators on the space of HARDI images. We shall model the space of HARDI images by the space of square integrable functions on the coupled space of positions and orientations, i.e. $\mathbb{L}_{2}(\R^3 \rtimes S^{2})$.
We will first show that such operators should be left-invariant with respect to the left-action of $SE(3)$ onto the space of HARDI images. This left-action of the group $SE(3)$ onto $\R^3 \rtimes S^{2}$ is given by
\begin{equation}\label{actleft}
g \cdot (\ul{y},\ul{n})= (R \ul{y}+\ul{x}, R \ul{n}), \qquad g=(\ul{x},R) \in SE(3), \ul{x},\ul{y} \in \R^3, \ul{n} \in S^{2}, R \in SO(3)
\end{equation}
and it induces the so-called left-regular action of the same group on the space of HARDI images similar to the left-regular action on 3D-images (for example orientation-marginals of HARDI images):
\begin{definition}\label{def:action}
The left-regular actions of $SE(3)$ onto $\mathbb{L}_{2}(\R^3 \rtimes S^{2})$ respectively $\mathbb{L}_{2}(\R^3)$ are given by
\[
\begin{array}{ll}
(\gothic{L}_{g=(\ul{x}, R)}U)(\ul{y}, \ul{n})=U(g^{-1}\cdot(\ul{y},\ul{n}))=U(R^{-1}(\ul{y}-\ul{x}), R^{-1}\ul{n}), & \ul{x}, \ul{y} \in \R^3, \ul{n} \in S^{2}, U \in \mathbb{L}_{2}(\R^3 \rtimes S^{2}), \\
(\gothic{U}_{g=(\ul{x}, R)}f)(\ul{y})=f(R^{-1}(\ul{y}-\ul{x}))\ , &R \in SO(3), \ul{x},\ul{y} \in \R^3,  f \in \mathbb{L}_{2}(\R^3).
\end{array}
\]
Intuitively, $\gothic{U}_{g=(\ul{x},R)}$ represents a rigid motion operator on images, whereas $\gothic{L}_{g=(\ul{x}, R)}$ represents a rigid motion on HARDI images.
\end{definition}
Operator on HARDI-images must be left-invariant as the net operator on a HARDI-image should commute with rotations and translations. For detailed motivation see \cite{DuitsIJCV2010},
where our motivation is
similar as in our framework of invertible orientation scores \cite{MarkusThesis,Fran2009,FrankenThesis,DuitsR2006AMS,Duits2005IJCV,DuitsSS2009,DuitsRSS1,DuitsAMS2,DuitsAMS1}. 
\begin{theorem} \label{th:convolve}
Let $\mathcal{K}$ be a bounded operator from $\mathbb{L}_{2}(\R^{3}\rtimes S^{2})$ into $\mathbb{L}_{\infty}(\R^{3}\rtimes S^{2})$ then there exists an integrable kernel
$k :(\R^{3}\rtimes S^{2}) \times (\R^{3}\rtimes S^{2}) \to \mathbb{C}$ such that $\|\mathcal{K}\|^2=\sup \limits_{(\ul{y},\ul{n}) \in \R^3 \rtimes S^{2} } \int_{\R^3 \rtimes S^{2}}|k(\ul{y},\ul{n}\, ; \, \ul{y}', \ul{n}')|^{2}{\rm d}\ul{y}'{\rm d}\sigma(\ul{n}')  $ and we have
\begin{equation}\label{kernfull}
(\mathcal{K}U)(\ul{y}, \ul{n})=\int_{\R^3 \rtimes S^{2}} k(\ul{y}, \ul{n}\, ;\, \ul{y}', \ul{n}') U(\ul{y}',\ul{n}') {\rm d}\ul{y}'{\rm d}\sigma(\ul{n}')\ ,
\end{equation}
for almost every $(\ul{y},\ul{n}) \in \R^{3} \rtimes S^{2}$ and all $U \in \mathbb{L}_{2}(\R^3 \rtimes S^{2})$.
Now $\mathcal{K}_{k}:=\mathcal{K}$ is left-invariant iff $k$ is left-invariant, i.e.
\begin{equation}\label{invprop1}
\gothic{L}_{g}\circ \mathcal{K}_{k}=\mathcal{K}_{k} \circ  \gothic{L}_{g} \desda \forall_{g \in SE(3)} \forall_{\ul{y}, \ul{y}' \in \R^3} \forall_{\ul{n}, \ul{n}' \in S^2} \; :\; k(g\cdot(\ul{y}, \ul{n})\, ;\, g \cdot(\ul{y}', \ul{n}') )=k(\ul{y}, \ul{n}\, ;\, \ul{y}', \ul{n}').
\end{equation}
Then to each positive left-invariant kernel $k:\R^3 \rtimes S^{2} \times \R^3 \rtimes S^{2} \to \R^{+}$ with
\mbox{$\int_{\R^3}\int_{S^2}k(\ul{0},\ul{e}_{z}\, ;\, \ul{y}, \ul{n})  {\rm d}\sigma(\ul{n}){\rm d}\ul{y}=1$}
we can associate a unique probability density $p:\R^3 \rtimes S^{2} \to \R^{+}$ with the invariance property
\begin{equation}\label{invprop}
p(\ul{y},\ul{n})=p(R_{\ul{e}_{z},\alpha}\ul{y},R_{\ul{e}_{z},\alpha}\ul{n}), \qquad \textrm { for all } \alpha \in [0,2\pi),
\end{equation}
by means of $p(\ul{y},\ul{n})=k(\ul{y}, \ul{n} \, ; \, \ul{0},\ul{e}_{z})$. The convolution now reads
\begin{equation} \label{thetrueconvolution}
\mathcal{K}_{k}U(\ul{y},\ul{n})= (p*_{\R^3 \rtimes S^{2}}U)(\ul{y},\ul{n})= \int \limits_{\R^3} \int \limits_{S^{2}} p(R_{\ul{n}'}^{T}(\ul{y}-\ul{y}'),R_{\ul{n}'}^{T}\ul{n}) \, U(\ul{y}',\ul{n}') {\rm d}\sigma(\ul{n}') {\rm d}\ul{y}',
\end{equation}
where $\sigma$ denotes the surface measure on the sphere and where $R_{\ul{n}'}$ is any rotation such that $\ul{n}'=R_{\ul{n}'}\ul{e}_{z}$.
\end{theorem}
For details see \cite{DuitsIJCV2010}.
By the invariance property (\ref{invprop}), the convolution (\ref{thetrueconvolution}) on $\R^3 \rtimes S^{2}$ may be written as a (full) $SE(3)$-convolution.
To this end we extend our positively valued functions $U$ defined on the quotient $\R^3 \rtimes S^{2}=(\R^{3} \rtimes SO(3))/(\{0\} \times SO(2))$ to the full Euclidean motion group by means of Eq.~(\ref{OS}) which yields
\begin{equation} \label{OSEA}
\tilde{U}(\ul{x},R_{\ul{e}_{z},\gamma} R_{\ul{e}_{y},\beta} R_{\ul{e}_{z},\alpha})= U(\ul{x},\ul{n}(\beta,\gamma)).
\end{equation}
in Euler angles.
Throughout this article we will use this natural extension to the full group.
\begin{definition} \label{def:OS}
We will call $\tilde{U}:\R^3 \rtimes SO(3) \to \R$, given by Eq.~(\ref{OS}), the HARDI-orientation score corresponding to HARDI image $U: \R^{3} \rtimes S^{2} \to \R$.
\end{definition}
\begin{remark}
By the construction of a HARDI-orientation score, Eq.~(\ref{OS}), it satisfies the following invariance property \mbox{$\tilde{U}(\ul{x},R R_{\ul{e}_{z},\alpha})=U(\ul{x},R)$} for all $\ul{x} \in \R^{3}, R \in SO(3), \alpha \in [0,2\pi)$.
\end{remark}
An $SE(3)$ convolution \cite{Chirikjian} of two functions $\tilde{p}:SE(3) \to \R$, $\tilde{U}:SE(3) \to \R$ is given by:
\begin{equation} \label{SE3conv}
(\tilde{p}*_{SE(3)}\tilde{U})(g)= \int_{SE(3)} \tilde{p}(h^{-1}g) \tilde{U}(h) {\rm d}\mu_{SE(3)}(h)\ ,
\end{equation}
where Haar-measure ${\rm d}\mu_{SE(3)}(\ul{x},R)={\rm d}\ul{x} \, {\rm d}\mu_{SO(3)}(R)$ with
${\rm d}\mu_{SO(3)}(R_{\ul{e}_{z},\gamma}R_{\ul{e}_{y},\beta}R_{\ul{e}_{z},\alpha})=\sin \beta {\rm d}\alpha{\rm d}\beta {\rm d}\gamma$. It is easily verified that that the following identity holds:
\[
(\tilde{p} *_{SE(3)} \tilde{U})(\ul{x},R)\, = \, 2\pi \;(p *_{\R^3 \rtimes S^{2}}U)(\ul{x},R \ul{e}_{z})\  .
\]
Later on in this article (in Subsection \ref{ch:BM} and Subsection \ref{ch:EM}) we will relate scale spaces on HARDI data and first order Tikhonov regularization on HARDI data to Markov processes. But in order to provide a road map of how the group-convolutions will appear in the more technical remainder of this article we provide some preliminary explanations on probabilistic interpretation of $\R^3 \rtimes S^{2}$-convolutions.

\section{Left-invariant Vector Fields on $SE(3)$ and their Dual Elements \label{ch:VF}}

We will use the following basis for the tangent space $T_{e}(SE(3))$ at the unity element $e=(\ul{0},I) \in SE(3)$:
\begin{equation} \label{ebasis}
A_{1}=\partial_{x}, A_{2}=\partial_{y}, A_{3}=\partial_{z}, A_{4}=\partial_{\tilde{\gamma}}, A_{5}=\partial_{\tilde{\beta}}, A_{6}=\partial_{\tilde{\alpha}}\ ,
\end{equation}
where we stress that at the unity element $(\ul{0},R=I)$, we have $\beta=0$ and here the tangent vectors $\partial_{\beta}$ and $\partial_{\gamma}$ are not defined, which requires a description of the tangent vectors on the $SO(3)$-part by means of the second chart.

The tangent space at the unity element is a 6D Lie algebra equipped with Lie bracket
\begin{equation}\label{bracket1}
[A,B]= \lim_{t \downarrow 0}
t^{-2}\left(a(t)b(t)(a(t))^{-1}(b(t))^{-1}-e \right) ,
\end{equation}
where $t \mapsto a(t)$ resp. $t\mapsto b(t)$ are \emph{any} smooth
curves in $SE(3)$ with $a(0)=b(0)=e$ and $a'(0)=A$ and $b'(0)=B$, for explanation on the formula (\ref{bracket1}) which holds for general matrix Lie groups, see \cite[App.G]{DuitsFuehrJanssen}. The Lie-brackets of the basis given in Eq.~(\ref{ebasis}) are given by
\begin{equation} \label{Lieproduct}
[A_{i},A_{j}]=\sum_{k=1}^{6}c^{k}_{ij}A_{k}\ ,
\end{equation}
where the non-zero structure constants for all three isomorphic Lie-algebras are given by
\begin{equation}\label{structureconstants}
-c^{k}_{ji}=c^{k}_{ij}= \left\{
\begin{array}{ll}
\textrm{sgn}\, \textrm{perm}\{i-3,j-3, k-3\} &\textrm{ if }i,j,k \geq 4, i\neq j \neq k, \\
\textrm{sgn}\, \textrm{perm}\{i,j-3, k\} &\textrm{ if }i,k \leq 3, j \geq 4, i\neq j \neq k . \\
\end{array}
\right.
\end{equation}
More explicitly, we have the following table of Lie-brackets:
\begin{equation} \label{tabel}
([A_{i},A_{j}])^{i=1,\ldots 6}_{j=1,\ldots,6}=
\begin{pmatrix}
0 & 0 & 0 & 0 & A_{3} & -A_{2}  \\
0 & 0 & 0 & -A_{3} & 0 & A_{1} \\
0 & 0 & 0 & A_{2} & -A_{1} & 0 \\
0 & A_{3} & -A_{2} & 0 & A_{6} & -A_{5} \\
-A_{3} & 0 & A_{1} & -A_{6} & 0 & A_{4} \\
A_{2} & -A_{1} & 0 & A_{5} & -A_{4} & 0 \\
\end{pmatrix} ,
\end{equation}
so for example $c^{3}_{15}=1, c^{3}_{14}=c_{15}^{2}=0, c^{2}_{16}=-c^{2}_{61}=-1$.
The corresponding left-invariant vector fields $\{\mathcal{A}_{i}\}_{i=1}^{6}$ are obtained by the push-forward of the left-multiplication $L_{g}h=gh$ by
$\left. \mathcal{A}_{i} \right|_{g} \phi= (L_{g})_{*}A_{i}\phi=A_{i}(\phi \circ L_{g})$ (for all smooth $\phi:\Omega_{g} \to \R$ which are locally defined on some neighborhood $\Omega_{g}$  of $g$) and they can be obtained by the derivative of the right-regular representation:
\begin{equation} \label{dR}
\begin{array}{l}
\left. \mathcal{A}_{i}\right|_{g}\phi=({\rm d}\mathcal{R}(A_{i})\phi)(g)= \lim \limits_{t \downarrow 0} \frac{\phi(g \; e^{t A_{i}})- \phi(g)}{t}\ , \ \
\textrm{ with }\mathcal{R}_{g}\phi(h)=\phi(hg).
\end{array}
\end{equation}
Expressed in the {\em first coordinate chart}, Eq.~(\ref{firstchart}), this renders for the left-invariant derivatives at position \\ $g=(x,y,z,R_{\ul{e}_{z},\gamma} R_{\ul{e}_{y},\beta} R_{\ul{e}_{z},\alpha}) \in SE(3)$ (see also \cite[Section~~9.10]{Chirikjian})
{\small
\begin{equation}\label{eq:LeftInvVFSEthree}
\hspace{-0.2cm}\hbox{}
\begin{array}{lrcrrrr}
\mathcal{A}_{1} =& (\cos\ealpha\cos\ebeta\cos\egamma -\sin\ealpha\sin\egamma)\;\pdx & 
+& (\sin\ealpha\cos\egamma + \cos\ealpha\cos\ebeta\sin\egamma)\;\pdy & -& \cos\ealpha\sin\ebeta \;\pdz ,\\
\mathcal{A}_{2} =& (-\sin\ealpha \cos\ebeta \cos\egamma  - \cos\ealpha \sin\egamma) \;\pdx & 
+ & (\cos\ealpha \cos\egamma - \sin\ealpha \cos\ebeta \sin\egamma) \;\pdy & + & \sin\ealpha\sin\ebeta \;\pdz , \\
\mathcal{A}_{3} =& \sin\ebeta \cos\egamma \;\pdx & +& \sin\ebeta \sin\egamma \;\pdy & + & \cos\ebeta \;\pdz , \\
\mathcal{A}_{4} =& \cos\ealpha\cot\ebeta \;\partial_{\alpha}& +& \sin\ealpha \;\partial_{\beta }& - &\frac{\cos\ealpha}{\sin\ebeta} \;\partial_{\egamma}, &   \\
\mathcal{A}_{5} =& -\sin\ealpha \cot\ebeta\;\partial_{\alpha}& +& \cos\ealpha \;\partial_{\ebeta}& + & \frac{\sin\ealpha}{\sin\ebeta} \;\partial_{\gamma},&  \\
\mathcal{A}_{6} =& \partial_{\alpha} & & & & &.\\
\end{array}
\end{equation}
}
for $\beta\neq 0$ and $\beta\neq\pi$. The explicit formulae of the left-invariant vector fields in the second chart,
Eq.~(\ref{secondchart}), are :
{\small
\begin{equation}\label{eq:LeftInvVFSEthree2}
\hspace{-0.4cm}\hbox{}
\begin{array}{ll}
\begin{array}{ll}
\mathcal{A}_{1} &= \cos\talpha\cos\tbeta\;\pdx + (\cos\tgamma\sin\talpha+\cos\talpha\sin\tbeta\sin\tgamma)\;\pdy \\
&\qquad + (\sin\talpha\sin\tgamma-\cos\talpha\cos\tgamma\sin\tbeta)\;\pdz, \\
\mathcal{A}_{2} &= -\sin\talpha\cos\tbeta\;\pdx + (\cos\talpha\cos\tgamma-\sin\talpha\sin\tbeta\sin\tgamma)\;\pdy \\
&\qquad + (\sin\talpha\sin\tbeta\cos\tgamma+\cos\talpha\sin\tgamma)\;\pdz, \\
\mathcal{A}_{3} &= \sin\tbeta\;\pdx - \cos\tbeta\sin\tgamma\;\pdy + \cos\tbeta\cos\tgamma\;\pdz,
\end{array}
&
\begin{array}{ll}
\mathcal{A}_{4} &= -\cos\talpha\tan\tbeta \;\partial_{\talpha} + \sin\talpha\;\partial_{\tbeta} + \frac{\cos\talpha}{\cos\tbeta} \;\partial_{\tgamma} , \\ & \\
\mathcal{A}_{5} &= 
\sin\talpha \tan\tbeta \;\partial_{\talpha} + \cos\talpha\;\partial_{\tbeta}-\frac{\sin\talpha}{\cos\tbeta}\;\partial_{\tgamma}, \\
 & \\
\mathcal{A}_{6} &= \partial_{\talpha} ,
\end{array}
\end{array}
\end{equation}
}
for $\tbeta\neq \frac{\pi}{2}$ and $\tbeta\neq-\frac{\pi}{2}$. Note that $\rm d \mathcal{R}$ is a Lie-algebra isomorphism, i.e.
\[
[A_{i},A_{j}]=\sum \limits_{k=1}^{6}c^{k}_{ij} A_{k} \desda [{\rm d}\mathcal{R}(A_{i}),{\rm d}\mathcal{R}(A_{j})]=\sum \limits_{k=1}^{6} c^{k}_{ij} {\rm d}\mathcal{R}(A_{k}) \desda [\mathcal{A}_{i},\mathcal{A}_{j}]=\mathcal{A}_{i}\mathcal{A}_{j}-\mathcal{A}_{j}\mathcal{A}_{i}= \sum \limits_{k=1}^{6} c^{k}_{ij} \mathcal{A}_{k}\ .
\]
These vector fields form a local moving coordinate frame of reference on $SE(3)$.
The corresponding dual frame $\{{\rm d}\mathcal{A}^{1},\ldots, {\rm d}\mathcal{A}^{6}\} \in (T(SE(3)))^{*}$ is defined by duality.
A brief computation yields :
{\small
\begin{equation} \label{duals2}
\left(
\begin{array}{c}
{\rm d}\mathcal{A}^{1} \\
{\rm d}\mathcal{A}^{2} \\
{\rm d}\mathcal{A}^{3} \\
{\rm d}\mathcal{A}^{4} \\
{\rm d}\mathcal{A}^{5} \\
{\rm d}\mathcal{A}^{6} 
\end{array}
\right)
= \left(
\begin{array}{c|c}
(R_{\ul{e}_{z},\gamma}R_{\ul{e}_{y},\beta} R_{\ul{e}_{z},\alpha})^T& \ul{0}\\
\hline \\
\ul{0} & M_{\beta,\alpha}
\end{array}
\right)
\left(
\begin{array}{c}
{\rm d}x \\
{\rm d}y \\
{\rm d}z \\
{\rm d}\alpha \\
{\rm d}\beta \\
{\rm d}\gamma 
\end{array}
\right)=
\left(
\begin{array}{c|c}
(R_{\ul{e}_{x},\tilde{\gamma}}R_{\ul{e}_{y},\tilde{\beta}} R_{\ul{e}_{z},\tilde{\alpha}})^T& \ul{0}\\
\hline \\
\ul{0} & \tilde{M}_{\tilde{\beta},\tilde{\alpha}}
\end{array}
\right)
\left(
\begin{array}{c}
{\rm d}x \\
{\rm d}y \\
{\rm d}z \\
{\rm d}\tilde{\alpha} \\
{\rm d}\tilde{\beta} \\
{\rm d}\tilde{\gamma} 
\end{array}
\right)
\end{equation}
}
where the $3\times 3$-zero matrix is denoted by $\ul{0}$ and where the $3\times 3$-matrices $M_{\beta,\alpha}$, $\tilde{M}_{\tilde{\beta},\tilde{\alpha}}$ are given by
{\small
\[
\begin{array}{ll}
M_{\beta,\alpha}=
\left(
\begin{array}{ccc}
0 & \sin \alpha & -\cos \alpha \sin \beta \\
0 & \cos \alpha & \sin \alpha \sin \beta \\
1 & 0 & \cos \beta
\end{array}
\right)\ ,
&
\tilde{M}_{\tilde{\beta},\tilde{\alpha}}=
\left(
\begin{array}{ccc}
-\cos \tilde{\alpha} \tan \tilde{\beta} & \sin \tilde{\alpha} & \frac{\cos \tilde{\alpha}}{\cos \tilde{\beta}} \\
\sin \tilde{\alpha} \tan \tilde{\beta} & \cos \tilde{\alpha} & -\frac{\sin \tilde{\alpha}}{\cos \tilde{\beta}} \\
1 & 0 & 0
\end{array}
\right)^{-T}.
\end{array}
\]
}
Finally, we note that by linearity the $i$-th dual vector filters out the $i$-th component of a vector field $\sum_{j=1}^{6}v^{j}\mathcal{A}_{j}$
\[
\langle {\rm d}\mathcal{A}^{i}, \sum_{j=1}^{6} v^{j}\mathcal{A}_{j}\rangle= v^{i}\ , \qquad \textrm{ for all }i=1,\ldots,6.
\]
\begin{remark}
In our numerical schemes, we do not use the formulas (\ref{eq:LeftInvVFSEthree}) and (\ref{eq:LeftInvVFSEthree2}) for the left-invariant vector fields as we want to avoid
sampling around the inevitable singularities that arise with the the coordinate charts, given by Eq.~(\ref{firstchart}) and (\ref{secondchart}), of $S^{2}$.
Instead in our numerics we use the approach that will be described in Section \ref{ch:num}. However, we shall use formulas (\ref{eq:LeftInvVFSEthree}) and (\ref{eq:LeftInvVFSEthree2}) frequently in
our analysis and derivation of Green's functions of left-invariant diffusions and left-invariant Hamilton-Jacobi equations on $\R^{3}\rtimes S^{2}$.
These left-invariant diffusions and left-invariant erosions are similar to diffusions and erosions on $\R^{5}$, we ``only'' have to replace the fixed left-invariant vector fields
$\{\partial_{x^{1}}, \ldots, \partial_{x^{5}}\}$ by the left-invariant vector fields $\{\left.\mathcal{A}_{1}\right|_{g},\ldots,\left.\mathcal{A}_{5}\right|_{g}\}$ which serve as a moving frame of reference
(along fiber fragments) in $\R^{3}\rtimes S^{2}= SE(3)/(\{\ul{0}\}\times SO(2))$. For an a priori geometric intuition
behind our left-invariant erosions and diffusions expressed in the left-invariant vector fields see Figure \ref{fig:intuition}.
\end{remark}

\section{Morphological Convolutions on $\R^3 \rtimes  S^2$ \label{ch:morph}}

\begin{figure}
\centerline{
\hfill
\includegraphics[width=0.8\hsize]{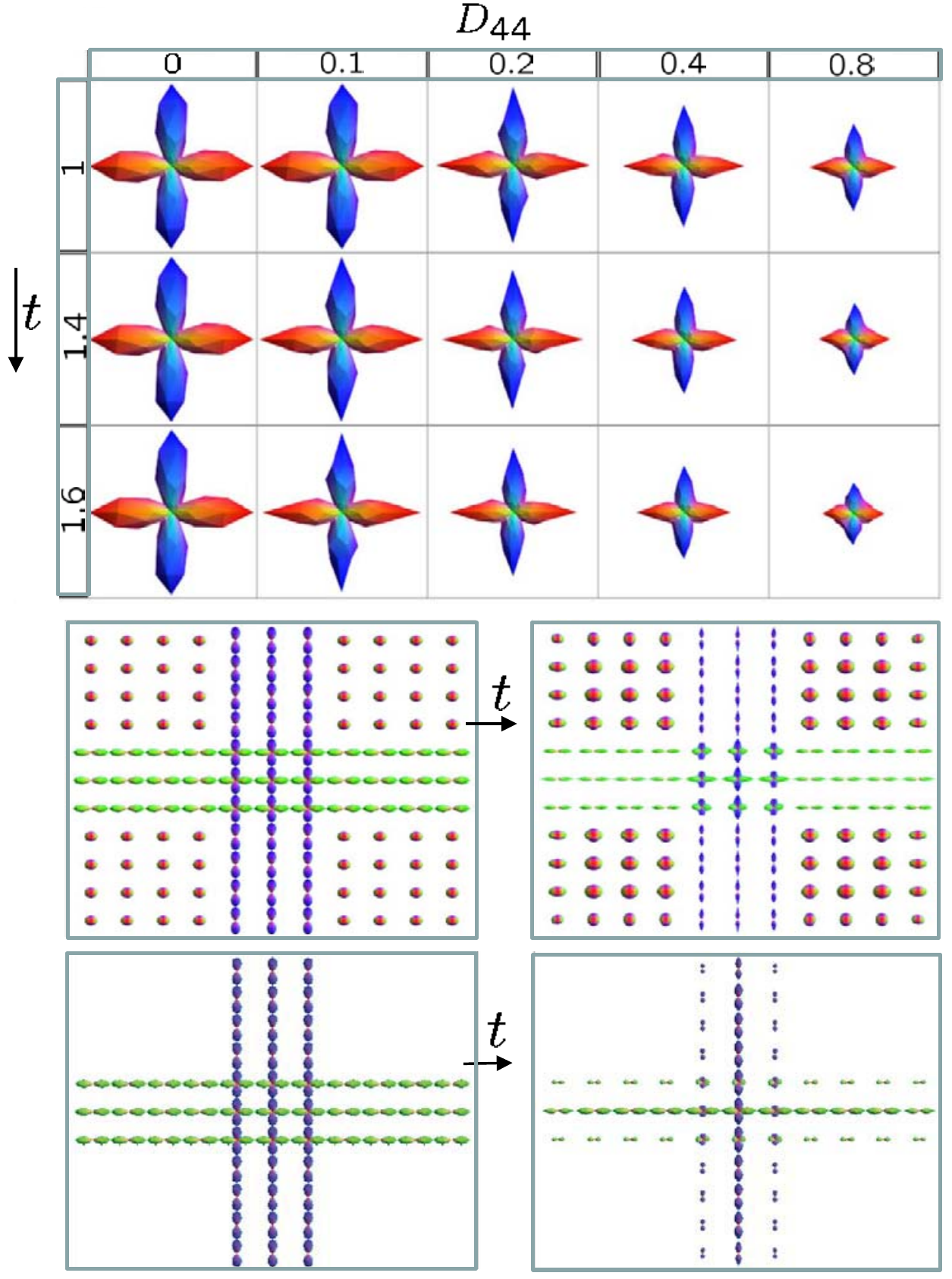}
\hfill
}
\caption{Top row: output angular erosion on a single glyph with crossings for several parameter settings.
Middle row: input and output simultaneous spatial and angular erosion (i.e. $\R^{3}\rtimes S^2$-erosion), where angular erosion dominates $t=3, D^{44}=0.6, D^{33}=0.3$.
Bottom row: input and output $\R^{3}\rtimes S^2$-erosion, where angular spatial erosion dominates $t=2, D^{44}=0.1, D^{33}=0.3$.
 }\label{fig:erode}
\end{figure}

Dilations on the joint space of positions and orientations $\R^{3}\rtimes S^2$ are obtained by replacing the $(+,\cdot)$-algebra by the $(\textrm{max},+)$-algebra in the $\R^{3}\rtimes S^{2}$-convolutions (\ref{thetrueconvolution})
\begin{equation} \label{dilation}
(k \oplus_{\R^{3} \rtimes S^2}U)(\ul{y},\ul{n})=
\sup \limits_{(\ul{y}',\ul{n}')\in \R^{3} \rtimes S^{2}} \left[k(R_{\ul{n}'}^{T}(\ul{y}-\ul{y}'),R_{\ul{n}'}^{T}\ul{n})+U(\ul{y}',\ul{n}')\right]
\end{equation}
where $k$ denotes a morphological kernel. If this morphological kernel is induced by a semigroup (or evolution)
then we write $k_{t}$ for the kernel at time $t$. Our aim is to derive suitable morphology kernels such that
\[
k_{t} \oplus_{\R^{3} \rtimes S^{2}} k_{s}= k_{s+t}\ , \textrm{ for all }s,t>0,
\]
where $t \mapsto k_{t}(\ul{y},\ul{n})$ describes the growth of balls in $\R^{3}\rtimes S^2$, i.e.
$t \mapsto k_{t}(\ul{y},\ul{n})$ is the \emph{unique viscosity solution}, see \cite{Crandall} and Appendix \ref{app:viscosity}, of
\begin{equation}\label{kernMorph}
\left\{
\begin{array}{l}
\frac{\partial W}{\partial t}(\ul{y},\ul{n},t)= \frac{1}{2}
\ul{G}^{-1}_{(\ul{y},\ul{n})}\left(\left.{\rm d}W(\cdot,\cdot,t) \right|_{\ul{y},\ul{n}},\left. {\rm d}W(\cdot,\cdot,t) \right|_{\ul{y},\ul{n}})\right) \\
W(\ul{y},\ul{n},0)= -\delta^{C}
\end{array}
\right.
\end{equation}
with the inverse $\ul{G}^{-1}$ of the metric tensor restricted to the sub-Riemannian manifold
$(SE(3),{\rm d}\mathcal{A}^{3},{\rm d}\mathcal{A}^{6})$
\begin{equation}\label{metrictensorR3S2erosion}
\ul{G} =
\frac{1}{D^{11}}(\mathcal{A}_{1} \otimes \mathcal{A}_{1}+\mathcal{A}_{2} \otimes \mathcal{A}_{2})+ \frac{1}{D^{44}}(\mathcal{A}_{4} \otimes \mathcal{A}_{4}+\mathcal{A}_{5} \otimes \mathcal{A}_{5} ),
\end{equation}
given by
\[
\ul{G}^{-1}=
D^{11}(\mathcal{A}_{1} \otimes \mathcal{A}_{1}+\mathcal{A}_{2} \otimes \mathcal{A}_{2})+ D^{44}(\mathcal{A}_{4} \otimes \mathcal{A}_{4}+\mathcal{A}_{5} \otimes \mathcal{A}_{5}).
\]
Both $\ul{G}$ and its inverse are well-defined on the cosets $\R^{3}\rtimes S^2= (\R^{3} \rtimes SO(3))/(\{\ul{0}\}\times SO(2))$, cf. Appendix \ref{app:C}.
Furthermore, in Eq.~(\ref{kernMorph}) we have used the morphological delta distribution $\delta^{C}(\ul{y},\ul{n})=+\infty $ if $(\ul{y},\ul{n}) \neq (\ul{0},\ul{e}_{z})$ and $0$ else, where we note that
\[
\lim \limits_{t \downarrow 0} k_{t} \oplus_{\R^{3} \rtimes S^{2}} U= (-\delta^{C}) \oplus_{\R^{3} \rtimes S^{2}} U = U
\]
uniformly on $\R^{3}\rtimes S^{2}$. Furthermore, we have used the
left-invariant gradient which is the co-vector field
\begin{equation} \label{leftinvgradU}
{\rm d}U(\ul{y},\ul{n})=\sum \limits_{i=1}^{6}(\mathcal{A}_{i}(U))(\ul{y},\ul{n}) \, \left. {\rm d}\mathcal{A}^{i}\right|_{(\ul{y},\ul{n})}=\sum \limits_{i=1}^{5}(\mathcal{A}_{i}(U))(\ul{y},\ul{n}) \, \left. {\rm d}\mathcal{A}^{i}\right|_{(\ul{y},\ul{n})}, \  \  (\ul{y},\ul{n}) \in \R^{3}\rtimes S^2,
\end{equation}
which we occasionally represent by a row vector given by
\begin{equation} \label{grad}
\nabla U(\ul{y},\ul{n})=(\mathcal{A}_{1}U(\ul{y},\ul{n}),\ldots,\mathcal{A}_{5}U(\ul{y},\ul{n}),0).
\end{equation}%
The dilation equation (\ref{kernMorph}) now becomes
\[
\left\{
\begin{array}{ll}
\frac{\partial W}{\partial t}(\ul{y},\ul{n},t) &= \frac{1}{2} \left( D^{11}((\mathcal{A}_{1}W(\ul{y},\ul{n},t))^2+(\mathcal{A}_{2}W(\ul{y},\ul{n},t))^2) 
+ D^{44} ((\mathcal{A}_{4}W(\ul{y},\ul{n},t))^2+(\mathcal{A}_{5}W(\ul{y},\ul{n},t))^2) \right)\\
W(\ul{y},\ul{n},0)&= -\delta^{C}(\ul{y},\ul{n})\ .
\end{array}
\right.
\]
Now we can consider either a positive definite metric (the case of dilations), or we can consider a negative definite metric (the case of erosions).
In the non-adaptive case this means; either we consider the $g_{ii}=\frac{1}{D^{ii}}>0$ or we choose them $g_{ii}=-\frac{1}{D^{ii}}<0$. Note that the erosion kernel follows from the dilation kernel by negation $k_{t} \mapsto - k_{t}$. Dilation kernels are negative and erosion kernels are positive and therefore we write
 \[
k_{t}^{+}:=-k_t \geq 0\ \  \textrm{ and }\ \  k_{t}^{-}:=k_t \leq 0\ .
\]
for the erosion kernels.
Erosions on $\R^{3}\rtimes S^{2}$ are given by:
\begin{equation}\label{Erosion}
(k_{t}^{+} \ominus_{\R^{3} \rtimes S^2}U)(\ul{y},\ul{n})=
\inf \limits_{(\ul{y}',\ul{n}')\in \R^{3} \rtimes S^{2}} \left[U(\ul{y}',\ul{n}')+ k_{t}^{+}(R_{\ul{n}'}^{T}(\ul{y}-\ul{y}'),R_{\ul{n}'}^{T}\ul{n})\right]\ .
\end{equation}
We distinguish between three types of erosions/dilations on $\R^{3}\rtimes S^2$
\begin{enumerate}
\item Angular erosion/dilation (i.e. erosion on glyphs): In case $g_{11}=g_{33}=0$ and $g_{44}<0$ the erosion kernels are given by
\[
k_{t}^{+,g_{11}=0,g_{44}}(\ul{y},\ul{n})=
\left\{
\begin{array}{ll}
\infty &\textrm{if }\ul{y} \neq \ul{0} \\
k_{t}^{+,g_{11}=0,|g_{44}|}(\ul{0},\ul{n}) & \textrm{else }
\end{array}
\right.
\]
such that the angular erosions are given by
\begin{equation}\label{angularerosion}
\begin{array}{ll}
(k_{t}^{+,g_{11}=0,|g_{44}|} \ominus_{\R^{3}\rtimes S^{2}}U)(\ul{y},\ul{n}) &= (k_{t}^{+,g_{11}=0,|g_{44}|}(\ul{0},\cdot) \ominus_{S^{2}} U(\ul{y},\cdot))(\ul{n}) \\
 &= \inf \limits_{\ul{n}' \in S^{2}} \left[k_{t}^{+,g_{11}=0,|g_{44}|}(\ul{0},R_{\ul{n}'}^{T}(\ul{n}'))+U(\ul{y}',\ul{n}')\right]
\end{array}
\end{equation}
whereas the angular dilations $g_{44}>0$ are given by
\begin{equation}\label{angulardilation}
\begin{array}{ll}
(k_{t}^{-,g_{11}=g_{33}=0,g_{44}} \oplus_{\R^{3}\rtimes S^{2}}U)(\ul{y},\ul{n}) &= (k_{t}^{-,g_{11}=g_{33}=0,g_{44}}(\ul{0},\cdot) \oplus_{S^{2}} U(\ul{y},\cdot))(\ul{n}) \\
&= \sup \limits_{\ul{n}' \in S^{2}} \left[k_{-t}^{g_{11}=g_{33}=0,g_{44}}(\ul{0},R_{\ul{n}'}^{T}\ul{n})+U(\ul{y}',\ul{n}')\right]
\end{array}
\end{equation}
\item Spatial erosion/dilation (i.e. the same spatial erosion for all orientations):  In case $g_{11}g_{33}\neq 0$, $g_{11} \leq 0$, $g_{33} \leq 0$ and $g_{44}=0$ the erosion kernels are given by
\[
k_{t}^{+,g_{11},g_{33}, g_{44}=0}(\ul{y},\ul{n})=
\left\{
\begin{array}{ll}
\infty &\textrm{if }\ul{n} \neq \ul{e}_{z} \\
k_{t}^{+,|g_{11}|=|g_{33}|,g_{44}=0}(\ul{y},\ul{e}_{z}) & \textrm{else }
\end{array}
\right.
\]
such that the spatial erosions are given by
\[
(k_{t}^{+,|g_{11}|,|g_{33}|,g_{44}=0} \ominus_{\R^{3}\rtimes S^{2}}U)(\ul{y},\ul{n})= (k_{t}^{+,|g_{11}|,|g_{33}|,g_{44}=0}(\cdot,\ul{e}_z) \ominus_{\R^3} U(\cdot,\ul{n}))(\ul{y})
\]
whereas the spatial dilations $g_{11},g_{33} \geq 0$, $g_{11}g_{33}>0$ are given by
\[
(k_{t}^{-,g_{11},g_{33},g_{44}=0} \oplus_{\R^{3}\rtimes S^{2}}U)(\ul{y},\ul{n})= (k_{t}^{-,g_{11},g_{33},g_{44}=0}(\cdot,\ul{e}_z) \oplus_{\R^3} U(\cdot,\ul{n}))(\ul{y}).
\]
\item Simultaneous spatial and angular erosion/dilations (i.e. erosions and dilations along fibers). The case $g_{44}\neq 0$ and $g_{11}\neq 0$ or $g_{33} \neq 0$.
\end{enumerate}
Similar to our previous work on $\R^{3}\rtimes S^2$-diffusion \cite{DuitsFrankenCASA} the third case is the most interesting one, simply because one would like to erode orthogonal to the fibers such that both the angular distribution and the spatial distribution of the a priori probability density
$U:\R^{3}\rtimes S^{2} \to \R^{+}$ are sharpened. See Figure \ref{fig:erode}.

\section{Left-Invariant Diffusions on $SE(3)=\R^3 \rtimes SO(3)$ and $\R^{3} \rtimes S^{2}$ \label{ch:LDSE3}}

In order to apply our general theory on diffusions on Lie groups, \cite{DuitsRSS1}, to suitable (convection-)diffusions on HARDI images,
we first extend all functions $U:\R^{3} \rtimes S^{2} \to \R^{+}$ to functions $\tilde{U}:\R^3 \rtimes SO(3) \to \R^{+}$ in the natural way, by means of (\ref{OS}).

Then we follow our general construction of scale space representations $\tilde{W}$ of functions $\tilde{U}$ (could be an image, or a score/wavelet transform of an image) defined on Lie groups, \cite{DuitsRSS1}, where we consider the special case $G=SE(3)$:
\begin{equation} \label{generaleqs}
\left\{
\begin{array}{l}
\partial_{t}\tilde{W}(g,t)=Q^{\ul{D}, \ul{a}}(\mathcal{A}_{1},\mathcal{A}_{2},\ldots,\mathcal{A}_{6}) \; \tilde{W}(g,t) \ ,\\
\lim \limits_{t\downarrow 0}\tilde{W}(g,t) =\tilde{U}(g)\ .
\end{array}
\right.
\end{equation}
which is generated by a quadratic form on the left-invariant vector fields:
\begin{equation} \label{QF}
Q^{\ul{D},\ul{a}}(\nabla):=Q^{\ul{D},\ul{a}}(\mathcal{A}_{1},\mathcal{A}_{2},\ldots,\mathcal{A}_{n})= \sum \limits_{i=1}^{6} -a_{i} \mathcal{A}_{i} +\sum \limits_{j=1}^{6}
\mathcal{A}_{i} D^{ij}\mathcal{A}_{j}
\end{equation}
Now the H\"{o}rmander requirement, \cite{Hoermander}, on the symmetric 
$\ul{D}=[D^{ij}] \in \R^{6\times 6}$, $\ul{D} \geq 0$ and $\ul{a}$, which guarantees smooth non-singular scale spaces, for $SE(3)$ tells us that $\ul{D}$ need not be strictly positive definite. The H\"{o}rmander requirement is that all included generators together with their commutators should span the full tangent space. 
To this end for diagonal $\ul{D}$ one should consider the set
\[
\mathcal{S}=\{i \in \{1,\ldots,6\} \; |\; D^{ii} \neq 0 \vee a_{i} \neq 0\}\ ,
\]
now if for example $1$ is not in here then $3$ and $5$ must be in $\mathcal{S}$, or if $4$ is not in $\mathcal{S}$ then $5$ and $6$ should be in $\mathcal{S}$.
If the H\"{o}rmander condition is satisfied the solutions of the \emph{ linear }diffusions (i.e. $D$, $\ul{a}$ are constant) are given by $SE(3)$-convolution with
a \emph{smooth probability kernel }$p_{t}^{\ul{D},\ul{a}}:SE(3) \to \R^{+}$ such that
\[
\begin{array}{l}
\tilde{W}(g,t)=(p_{t}^{D, \ul{a}}*_{SE(3)} \tilde{U})(g)= \int \limits_{SE(3)}p_{t}^{D, \ul{a}}(h^{-1}g) \tilde{U}(h){\rm d}\mu_{SE(3)}(h), \\
\lim \limits_{t \downarrow 0} \; \; p_{t}^{D, \ul{a}}*_{SE(3)} \tilde{U} = \tilde{U}\ , \qquad \textrm{ with }p_{t}^{\ul{D},\ul{a}}>0 \textrm{ and }\int_{SE(3)}p_{t}^{\ul{D},\ul{a}}(g) {\rm d}\mu_{SE(3)}(g)=1.
\end{array}
\]
where the limit is taken in $\mathbb{L}_{2}(SE(3))$-sense. On HARDI \emph{images} whose domain equals the homogeneous space $\R^3 \rtimes S^{2}$ one has the following scale space representations:
\begin{equation} \label{generalPDE}
\boxed{
\left\{
\begin{array}{l}
\partial_{t}W(\ul{y},\ul{n},t)=Q^{\ul{D}(U), \ul{a}(U)}(\mathcal{A}_{1},\mathcal{A}_{2},\ldots,\mathcal{A}_{5}) \; W(\ul{y},\ul{n},t) \ ,\\
W(\ul{y},\ul{n},0) =U(\ul{y},\ul{n})\ .
\end{array}
\right.
}
\end{equation}
with $Q^{\ul{D}(U),\ul{a}(U)}(\mathcal{A}_{1},\mathcal{A}_{2},\ldots,\mathcal{A}_{n})= \sum_{i=1}^{6} \left(a_{i} \mathcal{A}_{i} +\sum_{j=1}^{6}
\mathcal{A}_{i} D^{ij}(U)\mathcal{A}_{j}\right)$, where from now on we assume that $\ul{D}(U)$ and $\ul{a}(U)$ satisfy
\begin{equation}\label{generalrequirement}
\ul{a}(\tilde{U})(gh)= Z_{\alpha}^{T} (\ul{a})(g) \textrm{ and }
\ul{D}(\tilde{U})(gh)= Z_{\alpha} \ul{D}(\tilde{U})(g) Z_{\alpha}^{T}\ .
\end{equation}
for all $g \in SE(3)$ and all $h=(\ul{0}, R_{\ul{e}_{z},\alpha}) \in (\{\ul{0}\}\times SO(2))$ and
where
{\small
\begin{equation}\label{Zalpha}
Z_{\alpha} =
\left(\begin{array}{ccc|ccc}
\cos\alpha & -\sin\alpha  & 0 & 0 & 0 & 0 \\
\sin\alpha & \cos\alpha & 0 & 0 & 0 & 0 \\
0 & 0 & 1 & 0 & 0 & 0 \\
\hline
0 & 0 & 0 & \cos\alpha & -\sin\alpha & 0 \\
0 & 0 & 0 & \sin\alpha & \cos\alpha & 0 \\
0 & 0 & 0 & 0 & 0 & 1 \\
\end{array}\right) = R_{\ul{e}_{z},\alpha} \oplus R_{\ul{e}_{z},\alpha}, Z_{\alpha} \in SO(6), R_{\ul{e}_{z},\alpha} \in SO(3).
\end{equation}
}
Recall the grey tangent planes in Figure \ref{fig:intuition}
where we must require isotropy due to our embedding of $\R^{3}\rtimes S^{2}$ in $SE(3)$,
cf.~\cite[ch:4]{DuitsIJCV2010}. 

In the linear case the solutions of (\ref{generalPDE}) are given by the following kernel operators on $\R^3 \rtimes S^{2}$:
\begin{equation} \label{kernelsolution}
\boxed{
\begin{array}{l}
W(\ul{y},\ul{n},t)= (p_{t}^{\ul{D}, \ul{a}} *_{\R^3 \rtimes S^2} U)(\ul{y},\ul{n}) \\
= \int \limits_{0}^{\pi}
\int \limits_{0}^{2\pi} \int \limits_{\R^3} p_{t}^{\ul{D},\ul{a}}((R_{\ul{e}_{z},\gamma'}R_{\ul{e}_{y},\beta'})^{T}(\ul{y}-\ul{y}'), (R_{\ul{e}_{z},\gamma'}R_{\ul{e}_{y},\beta'})^{T}\ul{n}))\; U(\ul{y}', \ul{n}(\beta',\gamma'))\; {\rm d}\ul{y}' \,
{\rm d}\sigma(\ul{n}(\beta',\gamma')),
\end{array}
}
\end{equation}
where the surface measure on the sphere is given by ${\rm d}\sigma(\ul{n}(\beta',\gamma'))=\sin \beta' \; {\rm d}\gamma' {\rm d}\beta' \equiv {\rm d}\sigma(\tilde{\ul{n}}(\tilde{\beta},\tilde{\gamma}))=|\cos \tilde{\beta}|\; {\rm d}\tilde{\beta} {\rm d}\tilde{\gamma}$.
Now in particular in the linear case, since $(\R^3, I)$ and $(\ul{0},SO(3))$ are subgroups of $SE(3)$, we obtain the Laplace-Beltrami operators on these subgroups by means of:
 \[
 \begin{array}{l}
 \Delta_{S^{2}}=Q^{D=\textrm{diag}\{0,0,0,1,1,1\},\ul{a}=\ul{0}}=(\mathcal{A}_{4})^{2} +(\mathcal{A}_{5})^{2} +(\mathcal{A}_{6})^2= (\partial_{\beta})^2 +\cot(\beta) \partial_{\beta} + \sin^{-2}(\beta) (\partial_{\gamma})^2 \ ,\\
 \Delta_{\R^3}=Q^{D=\textrm{diag}\{1,1,1,0,0,0\},\ul{a}=\ul{0}}= (\mathcal{A}_{1})^{2} +(\mathcal{A}_{2})^{2} +(\mathcal{A}_{3})^2= (\partial_{x})^2 + (\partial_{y})^2 + (\partial_{z})^2\ .
 \end{array}
 \]
One wants to include line-models which exploit a natural coupling between position and orientation. \emph{Such a coupling is naturally included in a smooth way as long as the Hormander's condition is satisfied.}
Therefore we will consider more elaborate simple left-invariant convection, diffusions on $SE(3)$ with natural coupling between position and orientation.
To explain what we mean with natural coupling we shall need the next definitions.
\begin{definition}
A curve $\gamma:\R^{+} \to \R^3 \rtimes S^{2}$ given by $s \mapsto \gamma(s)=(\ul{y}(s),\ul{n}(s))$ is called horizontal if $\ul{n}(s)\equiv \|\dot{\ul{y}}(s)\|^{-1} \dot{\ul{y}}(s)$.
A tangent vector to a horizontal curve is called a horizontal tangent vector. A vector field $\mathcal{A}$ on $\R^3 \rtimes S^{2}$
is horizontal if for all $(\ul{y},\ul{n}) \in \R^3 \rtimes S^{2}$ the tangent vector
$\mathcal{A}_{(\ul{y},\ul{n})}$ is horizontal. The horizontal part $\mathcal{H}_{g}$ of each tangent space is the vector-subspace
of $T_{g}(SE(3))$ consisting of horizontal vector fields. Horizontal diffusion is diffusion which only takes place along
horizontal curves.
\end{definition}
It is not difficult to see that the horizontal part $\mathcal{H}_{g}$ of each tangent space $T_{g}(SE(3))$ is spanned by $\{\mathcal{A}_{3},\mathcal{A}_{4}, \mathcal{A}_{5}\}$.
So all horizontal left-invariant convection diffusions are given by Eq.~(\ref{generalPDE})
where in the linear case one must set $a_{1}=a_{2}=a_6=0$, $D_{j2}=D_{2j}=D_{1j}=D_{j1}=D_{j6}=D_{6j}=0$ for all $j=1,2,\ldots,6$. Now on a commutative group like $\R^6$ with commutative Lie-algebra $\{\partial_{x_{1}},\ldots, \partial_{x_6}\}$ omitting $3$-directions (say $\partial_{x_{1}}$, $\partial_{x_2}$, $\partial_{x_6}$) from each tangent space in the diffusion would yield no smoothing along the global $x_{1}$, $x_2$, $x_{6}$-axes. In $SE(3)$ it is different since the commutators take care of
indirect smoothing in the omitted directions $\{\mathcal{A}_{1},\mathcal{A}_{2},\mathcal{A}_{6}\}$, since \[
\textrm{span}\,\left\{\mathcal{A}_{3},\mathcal{A}_{4},\mathcal{A}_{5}, [\mathcal{A}_{3},\mathcal{A}_{5}]=\mathcal{A}_{2}, [\mathcal{A}_{4},\mathcal{A}_{5}]=\mathcal{A}_{6}, [\mathcal{A}_{5},\mathcal{A}_{3}]=\mathcal{A}_{1}\right\}= T(SE(3))
\]
Consider for example the $SE(3)$-analogues of the Forward-Kolmogorov (or Fokker-Planck) equations of the direction process for contour-completion and the stochastic process for contour enhancement which we considered in our
previous work, \cite{DuitsAMS1}, on $SE(2)$. Here we first provide the resulting PDEs and explain the underlying stochastic processes later in subsection \ref{ch:BM}. The Fokker-Planck equation for (horizontal) \emph{contour completion} on $SE(3)$ is given by
\begin{equation} \label{completion}
\left\{
\begin{array}{ll}
\partial_{t}W(\ul{y},\ul{n},t) &=(-\mathcal{A}_{3}+ D ((\mathcal{A}_{4})^{2} + (\mathcal{A}_{5})^2)) \; W(\ul{y},\ul{n},t)= (-\mathcal{A}_{3}+ D\, \Delta_{S^2}) \; W(\ul{y},\ul{n},t) \ , D=\frac{1}{2}\sigma^2 >0. \\
\lim \limits_{t\downarrow 0}W(\ul{y},\ul{n},t) &=U(\ul{y},\ul{n})\ .
\end{array}
\right.
\end{equation}
where we note that $(\mathcal{A}_{6})^2(W(\ul{y},\ul{n}(\beta,\gamma),s))=0$. This equation arises from Eq.~(\ref{generalPDE}) by setting $D^{44}=D^{55}=D$ and $a_3=1$ and all other parameters to zero. The Fokker-Planck equation for (horizontal) \emph{contour enhancement} is
\begin{equation} \label{enhancement}
\left\{
\begin{array}{ll}
\partial_{t}W(\ul{y},\ul{n},t) &=(D^{33}(\mathcal{A}_{3})^2+ D^{44} ((\mathcal{A}_{4})^{2} \! + \!(\mathcal{A}_{5})^2))\, W(\ul{y},\ul{n},t)= (D^{33}(\mathcal{A}_{3})^2+ D\, \Delta_{S^2}) \; W(\ul{y},\ul{n},t)\ , \\
\lim \limits_{t\downarrow 0}W(\ul{y},\ul{n},t) &=U(\ul{y},\ul{n})\ .
\end{array}
\right.
\end{equation}
The solutions of the left-invariant diffusions on $\R^3 \rtimes S^{2}$ given by (\ref{generalPDE}) (with in particular (\ref{completion}) and (\ref{enhancement})) are again given by convolution product (\ref{kernelsolution}) with a probability kernel $p_{t}^{\ul{D},\ul{a}}$ on $\R^3 \rtimes S^2$.
For a visualization of these probability kernels, see Figure \ref{fig:Brownianmotion}.

\section{Left-invariant Hamilton-Jacobi Equations on $\R^3 \rtimes S^2$ \label{ch:HJE}}

The unique viscosity solutions of
\begin{equation}\label{Hamdi1}
\left\{
\begin{array}{l}
\frac{\partial W}{\partial t}(\ul{y},\ul{n},t)= \frac{1}{2\eta}\ul{G}_{(\ul{y},\ul{n})}^{-1}\left(\left.{\rm d} W(\cdot,\cdot,t) \right|_{\ul{y},\ul{n}},\left. {\rm d} W(\cdot,\cdot,t) \right|_{\ul{y},\ul{n}})\right)^{\eta} \\
W(\ul{y},\ul{n},0)= U(\ul{y},\ul{n}),
\end{array}
\right.
\end{equation}
where $\eta \geq \frac{1}{2}$,
are given by dilation, Eq.~(\ref{dilation}), with the morphological Green's function $k_{t}^{D^{11}\!,D^{44}\!,\eta,-}:\R^3 \rtimes S^2 \to \R^{-}$
\begin{equation} \label{viscsol1}
W(\ul{y},\ul{n},t)=(k_{t}^{D^{11}\!,D^{44}\!,\eta,-} \oplus U)(\ul{y},\ul{n}),
\end{equation}
whereas the unique viscosity solutions of
\begin{equation}\label{Hamdi2}
\left\{
\begin{array}{l}
\frac{\partial W}{\partial t}(\ul{y},\ul{n},t)= -\frac{1}{2 \eta}\ul{G}_{(\ul{y},\ul{n})}^{-1}\left(\left. {\rm d} W(\cdot,\cdot,t) \right|_{\ul{y},\ul{n}},\left. {\rm d} W(\cdot,\cdot,t) \right|_{\ul{y},\ul{n}})\right)^{\eta} \\
W(\ul{y},\ul{n},0)= U(\ul{y},\ul{n})
\end{array}
\right.
\end{equation}
are given by erosion, Eq.~(\ref{Erosion}), with the morphological Green's function $k_{t}^{D^{11}\!,D^{44}\!,\eta,+}:\R^3 \rtimes S^2 \to \R^{+}$
\begin{equation} \label{viscsol2}
W(\ul{y},\ul{n},t)=(k_{t}^{D^{11}\!,D^{44}\!,\eta,+} \ominus U)(\ul{y},\ul{n}),
\end{equation}
This is formally shown in Appendix \ref{app:viscosity}, Theorem \ref{th:thetheorem}.
The exact morphological Green's functions are given by (where we recall (\ref{metrictensorR3S2}))
\begin{equation} \label{distR3XS2b}
-k_{t}^{D^{11}\!,D^{44}\!,\eta,-}(\ul{y},\ul{n})=
k_{t}^{D^{11}\!,D^{44}\!,\eta,+}(\ul{y},\ul{n}):=\inf \limits_{{\tiny \begin{array}{c}
\gamma=(\ul{x}(\cdot),R_{\ul{n}}(\cdot)) \in C^{\infty}((0,t), SE(3)),\\
\gamma(0)=(\ul{0},I=R_{\ul{e}_z}), \gamma(t)=(\ul{y},R_{\ul{n}}),\\
\langle \left.{d \rm}\mathcal{A}^{3}\right|_{\gamma}, \dot{\gamma}\rangle=\langle \left.{d \rm}\mathcal{A}^{6}\right|_{\gamma},\dot{\gamma}\rangle=0
\end{array}}} \int \limits_{0}^{t} \overline{\mathcal{L}}_{\eta}(\gamma(p),\dot{\gamma}(p))\, \left( \frac{dp}{ds}\right)^{\frac{1}{2\eta-1}}{\rm d}p\ ,
\end{equation}
with Lagrangian
\[
\overline{\mathcal{L}}_{\eta}(\gamma(p),\dot{\gamma}(p)):= \frac{2\eta-1}{2\eta} \left(\frac{1}{D^{11}}((\dot{\gamma}^{1}(p))^2+(\dot{\gamma}^{2}(p))^2)+
\frac{1}{D^{44}}((\dot{\gamma}^{4}(p))^2+(\dot{\gamma}^{5}(p))^2)\right)^{\frac{\eta}{2\eta-1}}=\frac{2\eta-1}{2\eta},
\]
where we applied short notation $\dot{\gamma}^{i}(p)=\langle \left.{\rm d}\mathcal{A}^{i}\right|_{\gamma(p)}, \dot{\gamma}(p) \rangle$ and with $\R^{3}\rtimes S^{2}$-``erosion arclength'' given by
\begin{equation}\label{arclengthR3XS2b}
\begin{array}{ll}
p(\tau) &= \int \limits_{0}^{\tau} \sqrt{\mathbf{G}_{\gamma(\tilde{\tau})}(\dot{\gamma}(\tilde{\tau}),\dot{\gamma}(\tilde{\tau}))}\, {\rm d}\tilde{\tau} =  \int \limits_{0}^{\tau} \sqrt{\sum \limits_{i \in \{1,2,4,5\}}\frac{1}{D^{ii}}|\langle \left. {\rm d}\mathcal{A}^{i}\right|_{\gamma(\tilde{\tau})}, \dot{\gamma}(\tilde{\tau})\rangle |^{2}  }\, {\rm d}\tilde{\tau}.
\end{array}
\end{equation}
For further explanation an details see Appendix \ref{app:viscosity}, in particular Lemma \ref{lemma:semigroup}.
As motivated in Appendix \ref{app:viscosity}, we use the following asymptotical analytical formula for the Green's function
\begin{equation} \label{assy}
k_{t}^{D^{11}\!,D^{44}\!, \eta, \pm}(\ul{y},\tilde{\ul{n}}(\tilde{\beta},\tilde{\gamma})) \equiv \pm \, \frac{2\eta-1}{2\eta} C^{\frac{2\eta}{2\eta-1}} t^{-\frac{1}{2\eta-1}} \left(\sum \limits_{i=1}^{6}
\frac{|\tilde{c}^{i}(\ul{y},\tilde{\alpha}=0,\tilde{\beta},\tilde{\gamma})|^{\frac{2}{w_{i}}}}{D^{ii}} \right)^{\frac{\eta}{2\eta-1}}
\end{equation}
for sufficiently small time $t>0$, where $C\geq 1$ is a constant that we usually set equal to $1$ and
where the constants $\tilde{c}^{i}_{\ul{y},\tilde{\alpha}=0,\tilde{\beta},\tilde{\gamma}}$, $i=1,\ldots,6$, are components
of the logarithm
\[
\sum \limits_{i=1}^{6}\tilde{c}^{i}_{\ul{y},\tilde{\alpha}=0,\tilde{\beta},\tilde{\gamma}}A_{i}=\log_{SE(3)} (\ul{y},R_{\ul{e}_{x},\tilde{\gamma}}R_{\ul{e}_{y},\tilde{\beta}})\ ,
\]
that we shall derive in Section \ref{ch:exp}, Eq.~(\ref{log2ndchart}).

Apparently, by Eq.~(\ref{distR3XS2b}) the morphological kernel describes the growth of ``erosion balls'' in $\R^{3}\rtimes S^{2}$.
In Appendix \ref{app:viscosity} we show that these erosion balls are locally equivalent to a
weighted modulus on the Lie-algebra of $SE(3)$, which explains our asymptotical formula
(\ref{assy}). This gives us a simple analytic approximation formula for balls in $\R^{3}\rtimes S^{2}$, where we do not need/use the minimizing curves (i.e. geodesics) in (\ref{distR3XS2b}).
\begin{remark}
In Appendix \ref{app:E} we provide the system of Pfaffian equations for geodesics on the sub-Riemannian manifold
$(SE(3), {\rm d}\mathcal{A}^{1}, {\rm d}\mathcal{A}^{2}, {\rm d}\mathcal{A}^{6})$ as a first step to generalize our
results on $SE(2)=\R^{2}\rtimes S^{1}$ in \cite[App. A]{DuitsAMS2},
where we generalize our results on $SE(2)=\R^{2}\rtimes S^{1}$ in \cite[App. A]{DuitsAMS2} to $(SE(3)$.
In order to compute the geodesics on the sub-Riemannian manifold
$(SE(3),{\rm d}\mathcal{A}^{3}, {\rm d}\mathcal{A}^{6})$ (used in the erosions Eq.~\ref{Hamdi2}) a similar approach
can be followed.
\end{remark}
\subsection{Data adaptive angular erosion and dilation \label{ch:adapt}}

\begin{figure}
\centerline{
\hfill
\includegraphics[width=0.8\hsize]{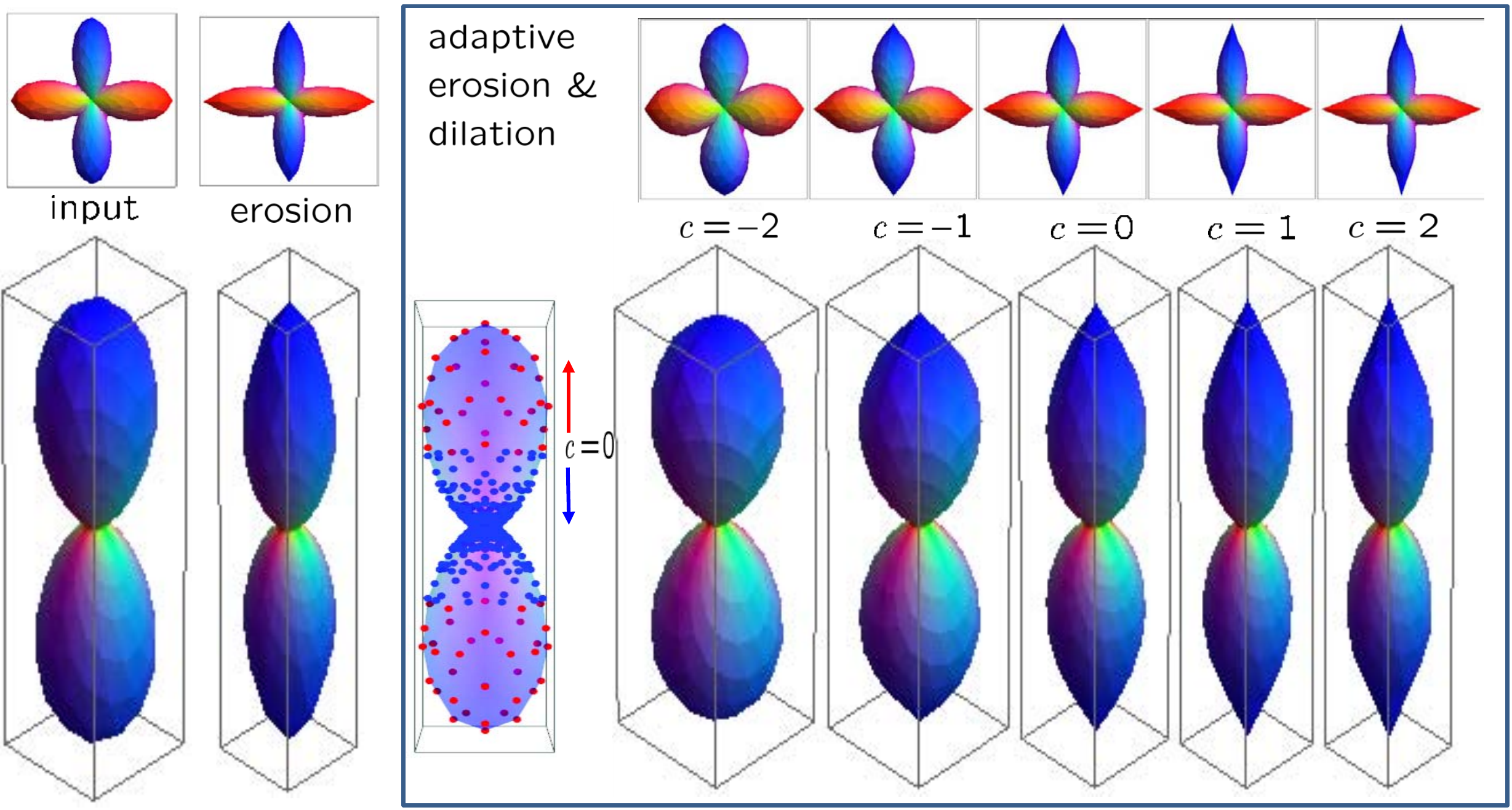}
\hfill
}\caption{From left to right, input glyph at say $\ul{y}=\ul{0}$.
Sample points where $\Delta_{LB}U(\ul{0},\ul{n})<0$ resp. $\geq 0$ are respectively indicated in red and blue.
Linearly eroded glyph ($D^{44}=0.4$, $t=0.5$, $\eta=1$, $\Delta t=0.01$), adaptive eroded glyph
($D^{44}=0.4$, $t=0.5$, $\eta=1$, $\phi(x)=D^{44} x^{\frac{1}{3}}$). }\label{Fig:adaptation}
\end{figure}

In the erosion evolution (\ref{Hamdi2}) one can include
adaptivity by making $D^{44}$ depend on the local Laplace-Beltrami-operator
\[
D^{44}(U)(\ul{y},\ul{n})= \phi( \Delta_{LB} U(\ul{y},\ul{n}) -c).
\]
with $\phi$ a non-decreasing, odd function and
$\inf \limits_{\ul{n} \in S^2, \ul{y} \in \R^{3}}\Delta_{LB} U(\ul{y},\ul{n})<c< \sup \limits_{\ul{n} \in S^2, \ul{y} \in S^{3}}\Delta_{LB} U(\ul{y},\ul{n})$.
The intuitive idea here is to dilate on points on a glyph
$\ul{x}+\mu\{U(\ul{x},\ul{x})\ul{n} \; |\; \ul{n} \in S^{2}\}$ where the Laplace-Beltrami operator is negative (usually around spherical maxima) and to erode on at locations on the glyph
where the Laplace-Beltami operator is positive. Parameter $c>0$ tunes the boundary on the glyph
where the switch between erosion and dilation takes place, See Figure \ref{Fig:adaptation} where we have set
$\phi(x)=D^{44}\, x^{\frac{1}{3}}$ while varying $c$.
%

\section{Probability Theory on $\R^3 \rtimes S^2$ \label{ch:prob}}

\subsection{Brownian Motions on $SE(3)=\R^3 \rtimes SO(3)$ and on $\R^3 \rtimes S^2$\label{ch:BM}}

Next we formulate a left-invariant discrete Brownian motion on $SE(3)$ (expressed in the moving frame of reference).
The left-invariant vector fields $\{\mathcal{A}_{1},\ldots,\mathcal{A}_{6}\}$ form a moving frame of reference to the group. Here we note that there are two ways of considering vector fields. Either one considers them as differential operators on smooth locally defined functions, or one considers them as tangent vectors to equivalent classes of curves. These two viewpoints are equivalent, for formal proof see \cite[Prop. 2.4]{Aubin}. Throughout this article we mainly use the first way of considering vector fields, but in this section we prefer to use the second way. We will write $\{\ul{e}_{1}(g), \ldots, \ul{e}_{6}(g)\}$ for the left-invariant vector fields (as tangent vectors to equivalence classes of curves) rather than the differential operators $\{\left.\mathcal{A}_{1}\right|_{g},\ldots, \left.\mathcal{A}_{6}\right|_{g}\}$. We obtain the tangent vector $\ul{e}_{i}$
from $\mathcal{A}_{i}$ by replacing
\begin{equation} \label{replace}
\begin{array}{ll}
\begin{array}{l}
\partial_{x} \leftrightarrow (1,0,0,0,0,0), \\
\partial_{y} \leftrightarrow
(0,1,0,0,0,0), \\
\partial_{z} \leftrightarrow (0,0,1,0,0,0), \\
\end{array}
&
\begin{array}{l}
\partial_{\beta} \leftrightarrow (0,0,0,\alpha \cos \, \beta \; \cos \, \gamma,\alpha \, \cos\, \beta \; \sin \, \gamma, -\alpha \, \sin \, \beta), \\
\partial_{\gamma} \leftrightarrow (0,0,0,\alpha \cos \, \gamma, \alpha \sin \gamma,0), \\
\partial_{\alpha} \leftrightarrow (0,0,0,\cos \gamma \sin \beta, \sin \, \gamma \sin \, \beta , \cos \, \beta),
\end{array}
\end{array}
\end{equation}
where we identified $SO(3)$ with a ball with radius $2\pi$ whose outer-sphere is identified with
the origin, using Euler angles $R_{\ul{e}_{z},\gamma} R_{\ul{e}_{y},\beta} R_{\ul{e}_{z},\alpha} \leftrightarrow \alpha \ul{n}(\beta,\gamma) \in B_{\ul{0},2\pi}$.
Next we formulate left-invariant discrete random walks on $SE(3)$ expressed in the moving frame of reference $\{\ul{e}_{i}\}_{i=1}^{6}$ given by (\ref{eq:LeftInvVFSEthree}) and (\ref{replace}):
%
\[
\begin{array}{l}
(\ul{Y}_{n+1}, \ul{N}_{n+1})= (\ul{Y}_{n}, \ul{N}_{n})+ \Delta s \sum \limits_{i=1}^{5} a_{i} \left.\ul{e}_{i}\right|_{(\ul{Y}_{n}, \ul{N}_{n})} + \sqrt{\Delta s} \sum \limits_{i=1}^{5} \ve_{i,n+1} \sum \limits_{j=1}^{5}
\sigma_{ji} \left.\ul{e}_{j}\right|_{(\ul{Y}_{n}, \ul{n}_{n})} \textrm{ for all }n=0,\ldots, N-1, \\
(\ul{Y}_{0}, \ul{n}_{0}) \sim U^{D},
\end{array}
\]
with random variable $(\ul{Y}_{0}, \ul{n}_{0})$ is distributed by $U^{D}$, where $U^{D}$ are the discretely sampled HARDI data (equidistant sampling in position and second order tessalation of the sphere) and where the random variables $(\ul{Y}_{n}, \ul{N}_{n})$ are recursively determined using the independently normally distributed random variables $\{\ve_{i,n+1}\}_{i=1,\ldots,5}^{n=1,\ldots, N-1}$, $\ve_{i,n+1} \sim \mathcal{N}(0,1)$ and where the stepsize equals
$\Delta s=\frac{s}{N}$ and where $\ul{a}:=\sum_{i=1}^{5} a_{i}\ul{e}_{i}$ denotes an apriori spatial velocity vector having \emph{constant} coefficients $a_{i}$ with respect to the moving frame of reference $\{\ul{e}_{i}\}_{i=1}^{5}$ (just like in (\ref{QF})). Now if we apply recursion and let $N\to \infty$ we get the following continuous Brownian motion processes on $SE(3)$:
\begin{equation} \label{stochproc}
\begin{array}{l}
Y(t)= Y(0) + \int \limits_{0}^{t} \left(\sum \limits_{i=1}^{3} a_{i} \left.\ul{e}_{i}\right|_{(Y(\tau),N(\tau))} + \frac{1}{2}\tau^{-\frac{1}{2}} \ve_{i} \sum \limits_{j=1}^{3} \sigma_{ji} \left.\ul{e}_{j}\right|_{(Y(\tau),N(\tau))} \right) {\rm d}\tau\ , \\
N(t)= N(0) + \int \limits_{0}^{t} \left(\sum \limits_{i=4}^{5} a_{i} \left.\ul{e}_{i}\right|_{(Y(\tau),N(\tau))} + \frac{1}{2}\tau^{-\frac{1}{2}} \ve_{i} \sum \limits_{j=4}^{5} \sigma_{ji} \left.\ul{e}_{j}\right|_{(Y(\tau),N(\tau))} \right) {\rm d}\tau\ ,
\end{array}
\end{equation}
with $\ve_{i} \sim \mathcal{N}(0,1)$ and $(X(0),N(0)) \sim U$ and where $\sigma= \sqrt{2D} \in \R^{6\times 6}$, $\sigma>0$. Note that ${\rm d}\sqrt{\tau}= \frac{1}{2}\tau^{-\frac{1}{2}} {\rm d}\tau$.

Now if we set $U=\delta_{\ul{0},\ul{e}_{z}}$ (i.e. at time zero ) then suitable averaging of infinitely many random walks of this process yields the transition probability $(\ul{y},\ul{n}) \mapsto p_{t}^{\ul{D},\ul{a}}(\ul{y},\ul{n})$ which is the Green's function
of the left-invariant evolution equations (\ref{generalPDE}) on $\R^{3}\rtimes S^{2}$. In general the PDE's (\ref{generalPDE}) are the Forward Kolmogorov equation of the general stochastic process (\ref{stochproc}). This follows by Ito-calculus and in particular Ito's formula for formulas on a stochastic process, for details see \cite[app.A]{MarkusThesis} where one should consistently replace the left-invariant vector fields of $\R^n$ by the left-invariant vector fields on $\R^{3}\rtimes S^{2}$.

In particular we have now formulated the \emph{direction process for contour completion} in $\R^3 \rtimes S^{2}$  (i.e. non-zero parameters in (\ref{stochproc}) are $D^{44}=D^{55}>0, a_{3}>0$ with Fokker-Planck equation (\ref{completion})) and
the \emph{(horizontal) Brownian motion for contour-enhancement} in $\R^3 \rtimes S^{2}$ (i.e. non-zero parameters in (\ref{stochproc}) are $D^{33}>0$, $D^{44}=D^{55}>0$ with Fokker-Planck equation (\ref{enhancement})).

\begin{figure}
\centerline{
\hfill
\includegraphics[width=0.99\hsize]{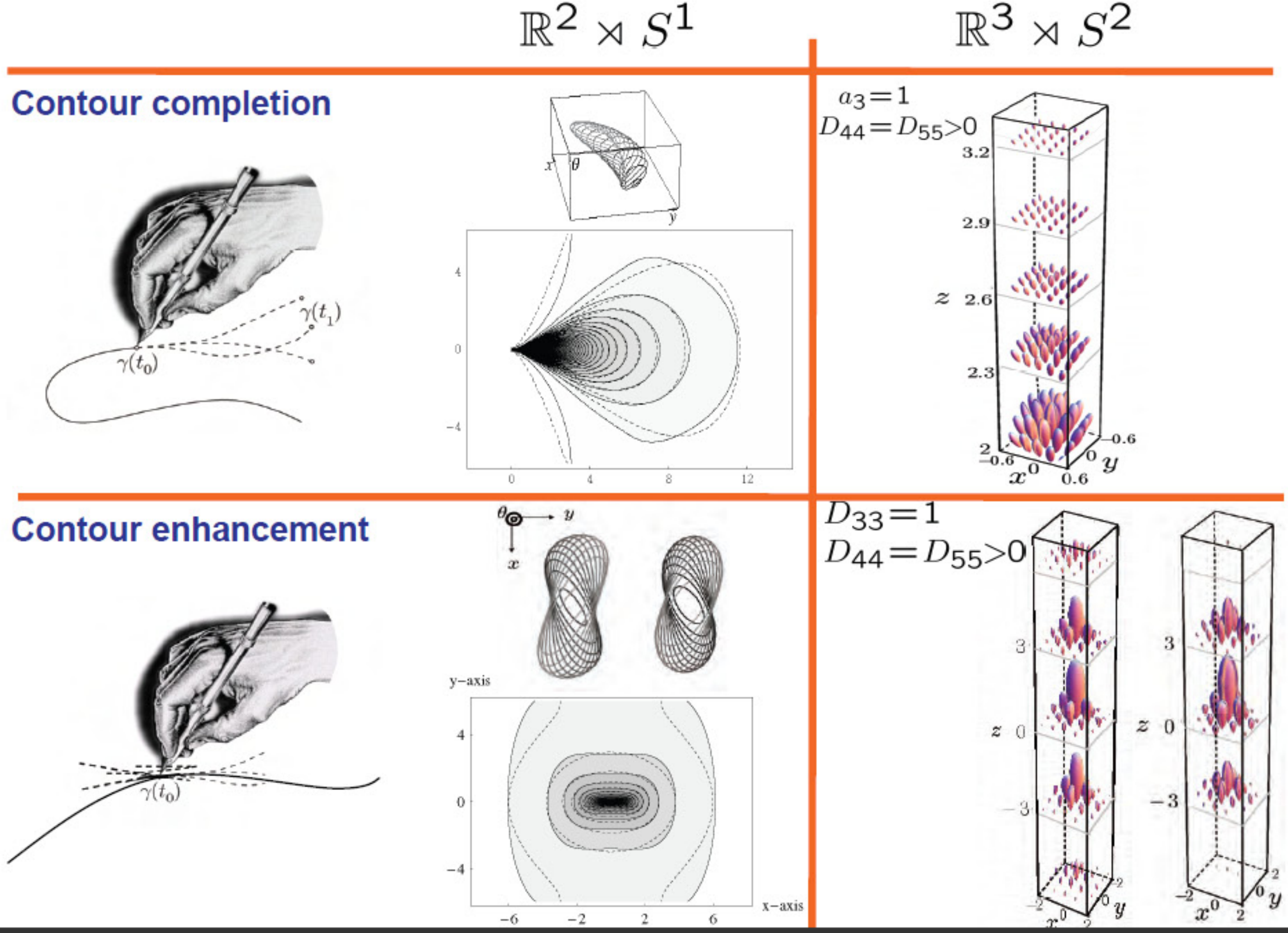}
\hfill
}
\caption{The Green's function of Forward Kolmogorov equation of (the horizontal) direction process and the Green's function of the Forward Kolmogorov equation of
(horizontal) Brownian motion in $\R^{d}\rtimes S^{d-1}$, for $d=3$ as considered
in this article and for $d=2$ as considered in our previous work \cite{DuitsR2006AMS} (completion) and \cite{DuitsAMS1,DuitsAMS2} (enhancement).}\label{fig:Brownianmotion}
\end{figure}

See Figure \ref{fig:Brownianmotion} for a visualization of typical Green's functions of contour completion and contour enhancement in $\R^{d}\rtimes S^{d-1}$, $d=2,3$.

\subsection{Time Integrated Brownian Motions \label{ch:EM}}

In the previous subsection we have formulated the Brownian-motions (\ref{stochproc}) underlying all linear left-invariant convection-diffusion equations on HARDI data, with in particular the direction process for contour completion and (horizontal) Brownian motion for contour-enhancement. However, we only considered the time dependent stochastic processes and as mentioned before in Markov-processes traveling time is memoryless and thereby negatively exponentially distributed $T \sim NE(\lambda)$, i.e. $P(T=t)=\lambda e^{-\lambda t}$ with expectation $E(T)=\lambda^{-1}$, for some $\lambda>0$. By means of Laplace-transform with respect to time we relate the time-dependent Fokker-Planck equations to their resolvent equations, as at least formally we have
\[
W(\ul{y},\ul{n},t)= (e^{t(Q^{\ul{D},\ul{a}}(\nabla))} U)(\ul{y},\ul{n}) \textrm{ and } P_{\gamma}(\ul{y},\ul{n},t)=\lambda \int_{0}^{\infty}e^{-t \lambda}(e^{t(Q^{\ul{D},\ul{a}}(\nabla))} U)(\ul{y},\ul{n})= \lambda (\lambda I - Q^{\ul{D},\ul{a}}(\nabla))^{-1}U (\ul{y},\ul{n}),
\]
for $t,\lambda>0$ and all $\ul{y} \in \R^3$, $\ul{n} \in S^{2}$, where the negative definite generator $Q^{\ul{D},\ul{a}}$ is given by (\ref{QF}) and again with $\nabla U=(\mathcal{A}_{1}U, \ldots, \mathcal{A}_{6}U)^{T}$. The resolvent operator \mbox{$\lambda (\lambda I - Q^{D=\textrm{diag}(D^{ii}),\ul{a}=0}( \nabla \,))^{-1}$} occurs in a first order Tikhonov regularization as we show in the next theorem.
\begin{theorem} \label{th:varprob}
Let $U \in \mathbb{L}_{2}(\R^3 \rtimes S^2)$ and $\lambda$, $D^{33} >0$, $D^{44}=D^{55}>0$. Then the unique solution of the variational problem
{\small
\begin{equation}\label{varprob}
\arg \mbox{}\hspace{-0.6cm}\mbox{} \min \limits_{P \in \mathbb{H}^{1}(\R^3 \rtimes S^{2})} \int \limits_{\R^3 \rtimes S^{2})}
\frac{\lambda}{2}(P(\ul{y},\ul{n})-U(\ul{y},\ul{n}))^{2} + \sum \limits_{k=3}^{5} D^{kk}|\mathcal{A}_{k} P(\ul{y},\ul{n})|^{2} {\rm d}\ul{y} {\rm d}\sigma(\ul{n})
\end{equation}
}
is given by
$
P_{U}^{\lambda}(\ul{y},\ul{n})=(R_{\lambda}^{D}*_{\R^3 \rtimes S^{2}}U)(\ul{y},\ul{n})$, where the Green's function $R_{\lambda}^{D}:\R^3 \rtimes S^{2} \to \R^{+}$ is the Laplace-transform of the heat-kernel with respect to time: $R_{\lambda}^{D}(\ul{y},\ul{n})= \lambda \, \int \limits_{0}^{\infty}p_{t}^{D, \ul{a}=\ul{0}}(\ul{y},n) e^{-t \lambda}\, {\rm d}t$ with $D=\textrm{diag}\{D^{11}, \ldots, D^{55}\}$. $P_{U}^{\lambda}(\ul{y},\ul{n})$ equals the probability of finding a random walker in $\R^3 \rtimes S^2$ regardless its traveling time at position $\ul{y} \in \R^3$ with orientation $\ul{n} \in S^{2}$ starting from initial distribution $U$ at time $t=0$.
\end{theorem}
For a proof see \cite{DuitsFrankenCASA}.
Basically, $P_{U}^{\lambda}(\ul{y},\ul{n})=(R_{\lambda}^{D}*_{\R^3 \rtimes S^{2}}U)(\ul{y},\ul{n})$ represents the probability density of finding
a random walker at position $\ul{y}$ with orientation $\ul{n}$ given that it started from the initial distribution $U$ regardless its traveling time, under the assumption that traveling time is memoryless and thereby negatively exponentially distributed $T \sim NE(\lambda)$.
There is however, a practical drawback due to the latter assumption: Both the time-integrated resolvent kernel of the direction process and the time-integrated resolvent kernel of the enhancement process suffer from a serious singularity at the unity element. In fact by some asymptotics one has
\[
\begin{array}{l}
R_{\lambda}^{\ul{a}=\ul{0}\, , \,  D=\textrm{diag}\{0,0,D^{33},D^{44},D^{55}\}}(\ul{y},\ul{n}) \sim \frac{1}{|g|^{6}} \textrm{ with }|g|=|(\ul{y},R_{\ul{n}})| \textrm{ for }|g|_{D^{33},D^{44},D^{55}} <<1 , \\
\textrm{ and }
R_{\lambda}^{\ul{a}=(0,0,1)\, ,\, D=\textrm{diag}\{0,0,0,D^{44},D^{55}\}}(0,0,z,\ul{n}) \sim \frac{1}{z^4}
\textrm{ for }0 < z << 1,
\end{array}
\]
where $|g|_{D^{33},D^{44},D^{55}}$ is the weighted modulus on $SE(3)=\R^3 \rtimes SO(3)$. For details see Appendix \ref{ch:asymptotics}.
These kernels can not be sampled using an ordinary mid-point rule. But even if the kernels are analytically integrated spatially and then numerically differentiated
the kernels are too much concentrated around the singularity for visually appealing results.

\subsubsection{A $k$-step Approach: Temporal Gamma Distributions and the Iteration of Resolvents}

The sum $T$ of $k$ independent negatively exponentially distributed random variables $T_{i} \sim NE(\lambda)$ (all with expectation $E(T_{i})=\lambda^{-1}$) is Gamma distributed:
\[
T=\sum \limits_{i=1}^{k} T_{i} \textrm{ with pdf :  } p(T=t)= p(T_{1}=t)*^{k-1}_{\R^{+}}p(T_{k}=t)=\Gamma(t\, ;\, k,\lambda):= \frac{\lambda^{k} t^{k-1}}{(k-1)!} e^{-\lambda t}\ , k \in \mathbb{N},
\]
where we recall that temporal convolutions are given by
$(f*_{\R^{+}}g)(t)=\int_{0}^{t}f(t-\tau) g(\tau)\; {\rm d}\tau$ and note that application of the laplace transform $\mathcal{L}$, given by $\mathcal{L}f(\lambda)=\int_{0}^{\infty}f(t)e^{-t \lambda} {\rm d}t$ yields $\mathcal{L}(f*_{\R^{+}}g)=\mathcal{L}(f) \cdot \mathcal{L}(g)$ .
Now for the sake of illustration we set $k=2$ for the moment and we compute the probability density of finding a random walker with traveling time $T=T_{1}+T_{2}$ at position $\ul{y}$ with orientation $\ul{n}$ given that it at $t=0$ started at $(\ul{0})$ with orientation $\ul{e}_{z}$. Basicly this means that the path of the random walker has two stages, first one with time $T_{1} \sim NE(\lambda)$ and subsequently one
with traveling time $T_{2}\sim NE(\lambda)$ and straightforward computations yield
\begin{equation}
\begin{array}{ll}
R_{\lambda,k=2}^{\ul{D},\ul{a}}(\ul{y},\ul{n}) &= \int \limits_{0}^{\infty} p(G_{T}=(\ul{y},\ul{n}) | T=t \textrm{ and }G_{0}=e)\; p(T=t)\, {\rm d}t \\
 &= \int \limits_{0}^{\infty} p(G_{T}=(\ul{y},\ul{n}) \; |\; T=T_{1}+T_{2}=t \textrm{ and }G_{0}=e)\; p(T_{1}+T_{2}=t) \, {\rm d}t \\
 &=\int \limits_{0}^{\infty} \int \limits_{0}^{t} p(G_{T_{1}+T_2}=(\ul{y},\ul{n}) \; |\; T_{1}=t-s \textrm{ and } T_{2}=s \textrm{ and }G_{0}=e)\; p(T_{1}=t-s)\; p(T_{2}=s) \, {\rm d}s {\rm d}t \\
 &=\lambda^2 \, \mathcal{L}\left(t \mapsto \int \limits_{0}^{t} (K_{t-s}*_{\R^{3}\rtimes S^{2}}K_{s} *_{\R^{3}\rtimes S^{2}} \delta_{e})(\ul{y},\ul{n}) {\rm d}s\right)(\lambda) \\
 &= \lambda^2 \, \mathcal{L}\left(t \mapsto \int \limits_{0}^{t} (K_{t-s}*_{\R^{3}\rtimes S^{2}}K_{s} )(\ul{y},\ul{n}) {\rm d}s\right)(\lambda) \\
 &= \lambda^2 \, \mathcal{L}\left(t \mapsto K_{t}(\cdot)\right)(\lambda) *_{\R^{3}\rtimes S^{2}}\mathcal{L}\left(t \mapsto K_{t}(\cdot)\right)(\lambda)(\ul{y},\ul{n}) = (R_{\lambda,k=1}^{\ul{D},\ul{a}}*_{\R^{3}\rtimes S^{2}}R_{\lambda,k=1}^{\ul{D},\ul{a}})(\ul{y},\ul{n})\ .
\end{array}
\end{equation}
By induction this can easily be generalized to the general case where we have
\[
\begin{array}{l}
R_{\lambda,k=2}^{\ul{D},\ul{a}}=R_{\lambda}^{\ul{D},\ul{a}} *^{k-1}_{\R^{3}\rtimes S^{2}}R_{\lambda}^{\ul{D},\ul{a}} \textrm{ with }R_{\lambda}^{\ul{D},\ul{a}}=R_{\lambda, k=1}^{\ul{D},\ul{a}} \textrm{ and }
p(T=t)= (p(T_{1}=\cdot) *^{k-1}_{\R^{+}} p(T_{k}=\cdot))(t)\ .
\end{array}
\]
As an alternative to our probabilistic derivation one has the following derivation (which holds in distributional sense):
\[
\begin{array}{ll}
R_{\lambda,k}^{\ul{D},\ul{a}} &= \int \limits_{0}^{\infty} (e^{t Q^{\ul{D},\ul{a}}( \nabla \,)} \delta_{e}) \; \Gamma(t \, ; k,\lambda)\, {\rm d}t \\
 &=\lambda^{k} (-Q^{\ul{D},\ul{a}}( \nabla \,)+\lambda I)^{-k} \delta_{e} = (\lambda \, (-Q^{\ul{D},\ul{a}}(( \nabla \,)+\lambda I)^{-1})^{k} \delta_{e} \\
 & = R_{\lambda}^{\ul{D},\ul{a}} *^{k-1}_{\R^{3}\rtimes S^{2}}R_{\lambda}^{\ul{D},\ul{a}}.
\end{array}
\]
where we note that the Laplace transformation of a Gamma distribution equals
\[
\mathcal{L} \Gamma(\cdot, k,\lambda)(s)=(1+\lambda^{-1}s)^{-k}\ .
\]
\begin{figure}[H] \vspace{-0.15cm}\mbox{}
\centerline{
\hfill
\includegraphics[width=0.51\hsize]{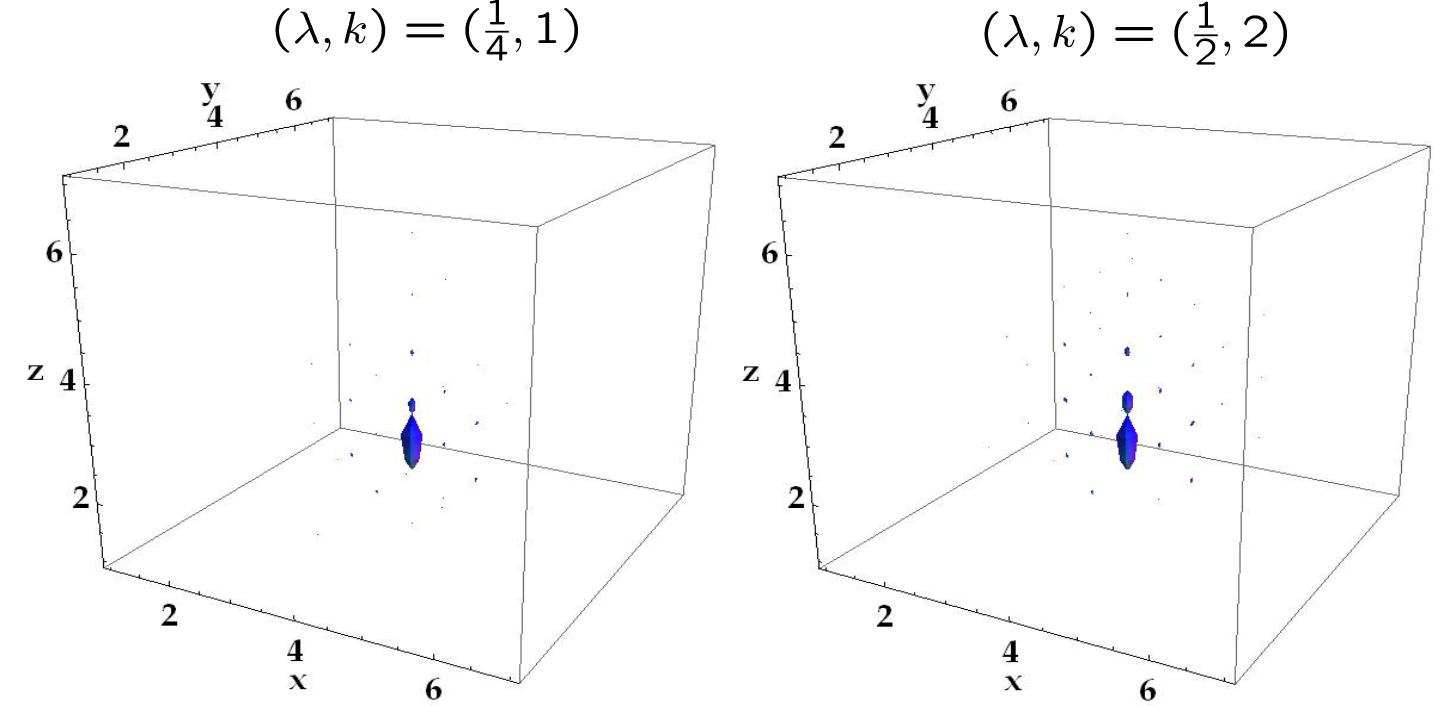}
\includegraphics[width=0.51\hsize]{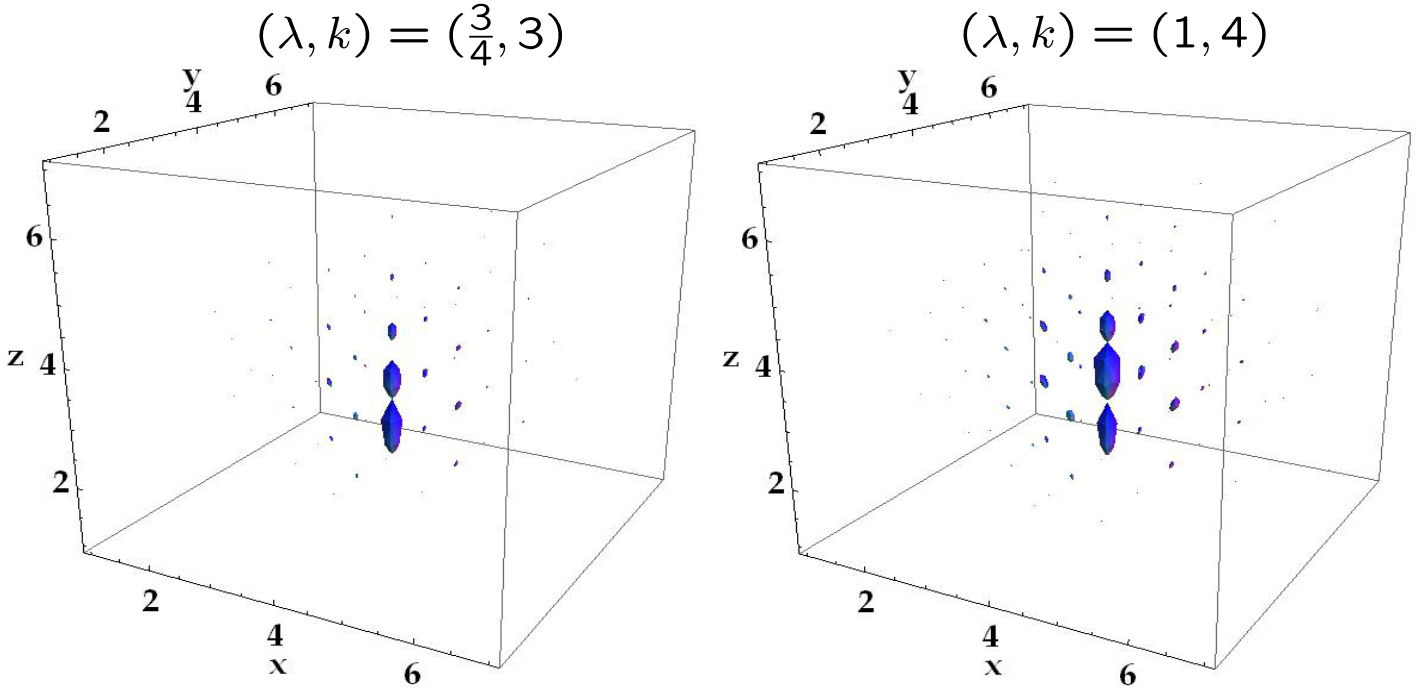}
\hfill
}
\vspace{-0.3cm}\mbox{}
\caption{Glyph visualization, recall Def. \ref{def:viz}, of the kernels $R_{\lambda,k=2}^{\ul{D},\ul{a}}:\R^{3}\rtimes S^{2} \to
\R^{+}$, with diffusion matrix $\ul{D}=\textrm{diag}\{0,0,0,D^{44},D^{44},0\}$ and convection vector $\ul{a}=(0,0,1,0,0,0)$, for several parameter-settings $(\lambda,k)$ for $D^{44}=0.005$. Kernels are sampled and computed on a
spatial $8\times 8 \times 8$-grid and on an a 162-point tessellation of the icosahedron using a left-invariant
finite difference scheme, cf.~\cite{Creusen}. For the sake of comparison we
fixed the expected value of the Gamma-distribution to $\frac{k}{\lambda} = 4$. The glyph-visualization parameter
$\mu$ (which determines the global scaling of the glyphs) has been set manually in all cases.
\label{fig:HARDI-DTI}}\label{fig:completionk}
\vspace{-0.3cm}\mbox{}
\end{figure}
\begin{figure}[H] \vspace{-0.15cm}\mbox{}
\centerline{
\hfill
\includegraphics[width=0.85\hsize]{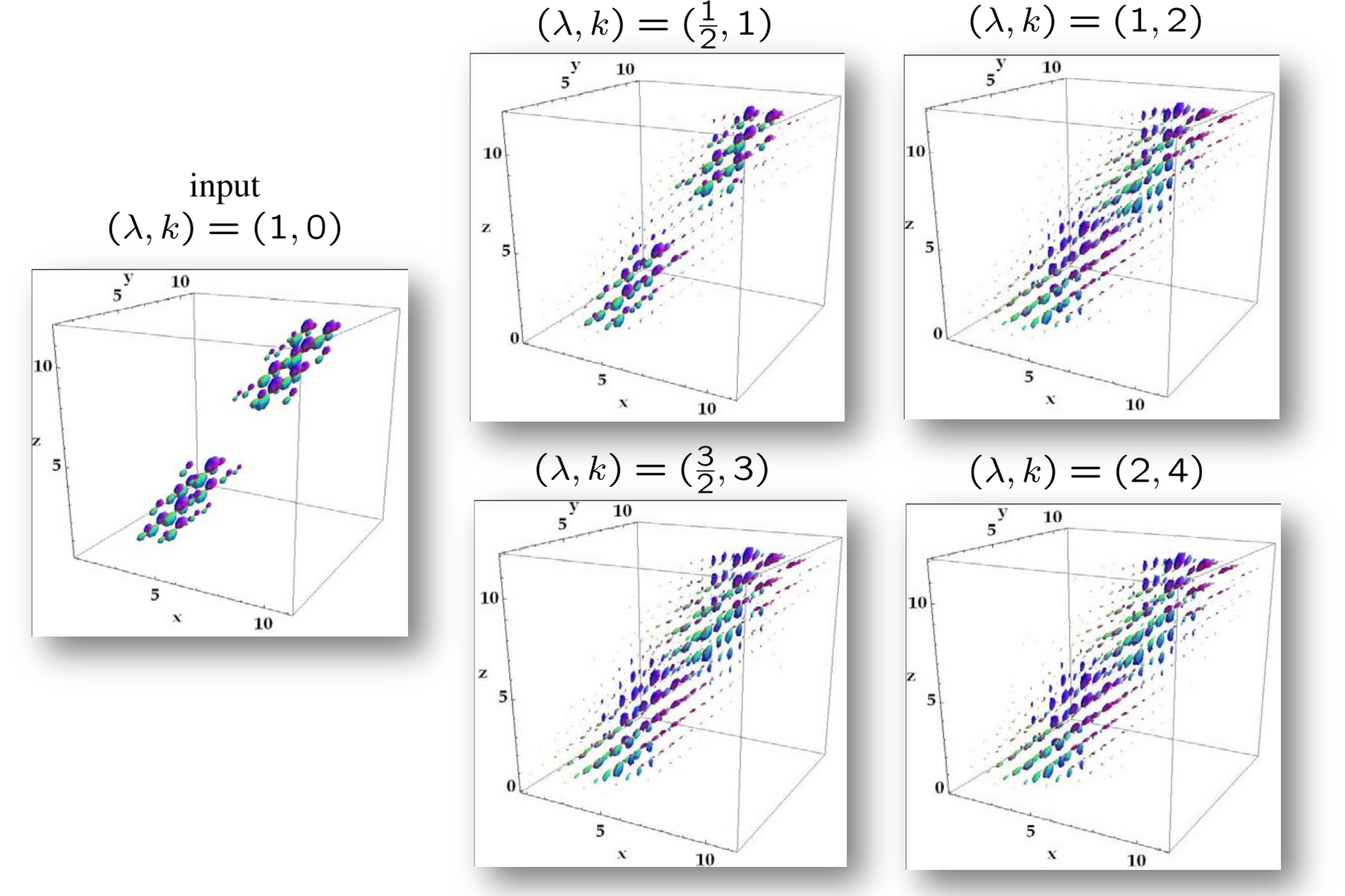}
\hfill
}
\vspace{-0.3cm}\mbox{}
\caption{Glyph visualization, recall Def. \ref{def:viz}, of left-invariant finite difference scheme, on a $12\times 12 \times 12$-grid and 162-point tessellation of the icosahedron, applied to the input data set, according to our time integrated
contour-completion process, ref. Eq.~(\ref{stochproc}), with diffusion matrix $\ul{D}=\textrm{diag}\{0,0,0,D^{44},D^{44}\}$, $D^{44}=0.005$ and convection vector $\ul{a}=(0,0,1,0,0)$, iterated $k$-times with
$E(T_{i})=\lambda^{-1}$, $i=1,\ldots,k$.
\label{fig:HARDI-DTI}}\label{fig:completion}
\vspace{-0.3cm}\mbox{}
\end{figure}
Figure \ref{fig:completion} shows some experiments of contour completion on an artificial data set containing fibers with a gap we would like to complete, for various parameter settings of
$(\lambda,k)$, where $T_{i} \sim \textrm{NE}(\lambda)$ and $k$ denotes the iteration index of the time integrated
contour-completion process. In principle this scheme boils down to $\R^{3}\rtimes S^{2}$-convolution with the Green's functions depicted in Fig. \ref{fig:erosion} (this is only approximately the case due to discretization).
For a fair comparison we kept the expected value constant, that is
$E(T)=E(T_{1}+\ldots T_{k})=\frac{k}{\lambda}=4$ in Figure \ref{fig:completionk} and
$E(T)=E(T_{1}+\ldots T_{k})=\frac{k}{\lambda}=2$ in Figure \ref{fig:completion},
which roughly coincides with half the size of the gap. Note that the results hardly change after two iterations, as
the graphs of the scaled Gamma distributions $\Gamma(\cdot\, ;\, \frac{k}{2},k)/\Gamma((2(k-1)/k) ;\, \frac{k}{2},k)$ are similar for $k=2,3,4,5$.

\subsection{Cost Processes on $SE(3)$ \label{ch:cost}}

In this subsection we present a short overview of cost processes on $SE(3)$. 
The mapping
$I \mapsto \lim_{\varepsilon \to 0} \log_{\epsilon}I$ defines a morphism of the $(+,\cdot)$-algebra to the $(\min,+)$-algebra on $\R^{+}$ (the co-domain of our HARDI orientation scores, recall Def. \ref{def:OS}) and it is indeed
readily verified that $\lim_{\varepsilon \to 0} \log_{\epsilon} \epsilon^{a}\cdot \epsilon^{b} =a+b$ and $\lim_{\varepsilon \to 0} \log_{\epsilon} (\epsilon^{a}+\epsilon^{b}) = \min(a,b)$.
Using this map, various notions of probability calculus can be mapped to their counterparts in optimization problems.
Next we mention the following definitions as given by Akian, Quadrat and Viot in \cite{Akian} adapted to our case.
\begin{definition}
The \textit{decision space} is the triplet $(SE(3),\mathcal{O}(SE(3)),C)$ with $\mathcal{O}(SE(3))$
denoting the set of all open subsets of the topological space $SE(3)$. The function
$C:\mathcal{O}(SE(3))\rightarrow \overline{\mathbb{R}}^{+}:= \R^{+} \cup \{\infty\}$ is such that
\begin{enumerate}
\item $C(SE(3))=0$
\item $C(\phi)=+\infty$
\item $C(\bigcup_{n}{A_{n}})=\inf_{n}{C(A_{n})}$ for any $A_{n}\in\mathcal{O}(SE(3))$
\end{enumerate}
$C$ is called a \textit{cost measure} on $SE(3)$.
The map $c:SE(3)\rightarrow\overline{\mathbb{R}}^{+}$ given by $g\mapsto c(g)$ such that $C(A)=\inf_{g\in A}c(g)$ for all $A\subseteq SE(3)$, is called the \textit{cost density} of the cost measure $C$.
\end{definition}
\begin{definition}
Analogous to the random variables of probability theory, a \textit{decision variable} $G$
(on $(SE(3),\mathcal{O}(SE(3)),C)$) is a mapping from $SE(3)$ to $\mathbb{R}$.
It induces a cost measure $C_G$ on $\mathbb{R}$ given by $C_G(I)=C(G^{-1}(I))$ for all
$I\in\mathcal{O}(\mathbb{R})$. The associated cost density is denoted by $c_G$.
\end{definition}
One can formulate related concepts such as independent decision variables, conditional cost, mean of a decision variable, characteristic function of a decision variable etc.
in the same way as in probability theory keeping in mind the morphism of the $(+,\cdot)$-algebra to the $(\min,+)$-algebra.
The Laplace or Fourier transform in $(+,\cdot)$-algebra corresponds to the Frenchel transform in the $(\min,+)$-algebra.
Now we present the decision counterpart of the Markov processes, namely Bellman processes.
\begin{definition}
A continuous time Bellman process $G_t$ on $(SE(3),\mathcal{O}(SE(3)),C)$ is a function from $SE(3)$ to $\mathcal{C}(\mathbb{R}^+)$ (set of continuous functions on non negative reals) with cost density
\begin{equation}
c_G(\tilde{\gamma}(\cdot))=c_0(\tilde{\gamma}(0))+\int^{\infty}_0{c(t,\tilde{\gamma}(t),\tilde{\gamma}'(t))}\, {\rm d}t
\end{equation}
where $c$ is called the transition cost which is a map from $\mathbb{R}\times SE(3)\times\mathcal{T}(SE(3))$ to $\overline{\mathbb{R}}^{+}$ such that
\[
\inf_{\mathcal{A}\in\mathcal{L}(SE(3))}c(t,g,\left.\mathcal{A}\right|_{g})\equiv 0 \textrm{ for all }(t,g) \in \R^{+}\times SE(3),
\]
and where $c_0$ is some cost density on $SE(3)$.
\end{definition}
We set
\[
c(t,\tilde{\gamma}(t),\tilde{\gamma}'(t))= \sum \limits_{i \in \{1,2,4,5\}}
|\langle \left.{\rm d}\mathcal{A}^{i}\right|_{\tilde{\gamma}(t)}, \tilde{\gamma}'(t) \rangle|^2.
\]
Then the marginal cost for a Bellman process $G_t$ on
$(SE(3),\mathcal{O}(SE(3)),C)$ to be in a state $g$ at time $t$ given an initial cost $c_0$,
$\tilde{m}(g,t):=C(G_t=g)$ satisfies the following relation known as the Bellman equation
\begin{equation}\label{Bellman}
\begin{array}{l}
\partial_{t}\tilde{m}+(\mathcal{F}_{\mathcal{L}(SE(3))}(L))({\rm d}\tilde{m})=0, \\
\tilde{m}(g,0)=c_0(g)
\end{array}
\end{equation}
where $\mathcal{F}_{\mathcal{L}(SE(3))}$ denotes the Fenchel transform, see Appendix \ref{app:viscosity}, Definition \ref{def:LT},
on the Lie-Algebra $\mathcal{L}(SE(3))$ of
left-invariant vector fields and where
\[
\begin{array}{l}
L(\mathcal{A})=\!\!\sum \limits_{i\in \{1,2,4,5\}}\!\! \frac{1}{D^{ii}}
\left|\langle {\rm d}\mathcal{A}^{i}, \mathcal{A}\rangle \right|^2,   \\
{\rm d}\tilde{m}= \sum \limits_{i\in \{1,2,4,5\}} \mathcal{A}_{i}\tilde{m} \; {\rm d}\mathcal{A}^{i}
\end{array}
\]
for all $\mathcal{A}$ left-invariant vector fields within contact-manifold $(SE(3), {\rm d}\mathcal{A}^{3},{\rm d}\mathcal{A}^{6})$
and with left-invariant gradient of $\tilde{m}$ within $(SE(3), {\rm d}\mathcal{A}^{3},{\rm}\mathcal{A}^{6})$. This Bellman equation for the
cost process coincides with the Hamilton-Jacobi equation on $SE(3)$ (\ref{Hamdi1})
whose viscosity solution is given by morphological convolution with the corresponding morphological Green's function as proven in Appendix \ref{app:viscosity}.

\section{Differential Geometry: The underlying Cartan-Connection on $SE(3)$ and the Auto-Parallels in $SE(3)$\label{ch:CA}}

Now that we have constructed all left-invariant scale space representations on HARDI images, generated by means of a quadratic form (\ref{QF}) on the
left-invariant vector fields on $SE(3)$. The question rises what is the underlying differential geometry for these evolutions ?

For example, as the left-invariant vector fields clearly vary per position in the group yielding a moving frame of reference attached to luminosity particles (random walkers in $\R^3 \rtimes S^{2}$ embedded in $SE(3)$) with both a position and an orientation, the question rises along which trajectories in \mbox{$\R^3 \rtimes S^2$} do these particles move ?
Furthermore, as the left-invariant vector fields are obtained by the push-forward of the left-multiplication on the group,
\[
\mathcal{A}_{g}=(L_{g})_{*} \mathcal{A}_{e}, \textrm{    i.e. } \mathcal{A}_{g}\tilde{\phi}= \mathcal{A}_{e}(\tilde{\phi} \circ L_{g}), \textrm{ where }L_gh=gh, \ \ g,h \in SE(3), \tilde{\phi}:SE(3) \to \R \textrm{ smooth },
\]
the question rises whether this defines a connection between all tangent spaces, such that these trajectories are auto-parallel with respect to this connection ? Finally, we need a connection to rigid body mechanics described in a moving frame of reference, to get some physical intuition in the choice of the fundamental constants\footnote{Or later in Subsection \ref{ch:LDOS} to get some intuition in the choice of \emph{functions} $\{a_{i}\}_{i=1}^{6}$ and $\{D^{ij}\}_{i,j=1}^{6}$.} $\{a_{i}\}_{i=1}^{6}$ and $\{D^{ij}\}_{i,j=1}^{6}$ within our generators (\ref{QF}).

In order to get some first physical intuition on analysis and differential geometry along the moving frame $\{\mathcal{A}_{1},\ldots, \mathcal{A}_{6}\}$ and its dual frame $\{{\rm d}\mathcal{A}^{1},\ldots, {\rm d}\mathcal{A}^{6}\}$, we will make some preliminary remarks on the well-known theory of rigid body movements described in moving coordinate systems.
Imagine a curve in $\R^3$ described in the moving frame of reference (embedded in the spatial part of the group $SE(3)$), describing a rigid body movement \emph{with constant spatial velocity $\hat{\ul{c}}^{(1)}$ and constant angular velocity $\hat{\ul{c}}^{(2)}$} and parameterized by arc-length $s>0$. Suppose the curve is given by
\[
\ul{y}(s)=\sum \limits_{i=1}^{3}\alpha^{i}(s) \left. \mathcal{A}_{i} \right|_{\ul{y}(s)}\  \textrm{ where }\alpha^{i} \in C^{2}([0,L],\R),
\]
such that $\hat{\ul{c}}^{(1)}= \sum_{i=1}^{3}\dot{\alpha}^{i}(s) \left. \mathcal{A}_{i} \right|_{\ul{y}(s)}$ for all $s>0$.
Now if we differentiate twice with respect to the arc-length parameter and keep in mind that $\frac{d}{ds} \left. \mathcal{A}_{i}\right|_{\ul{y}(s)}= \hat{\ul{c}}^{(2)} \times \left. \mathcal{A}_{i} \right|_{\ul{y}(s)}$, we get
\[
\ddot{\ul{y}}(s)=  0+2 \hat{\ul{c}}^{(2)} \times \hat{\ul{c}}^{(1)} + \hat{\ul{c}}^{(2)} \times (\hat{\ul{c}}^{(2)} \times \ul{y}(s))\ .
\]
In words: The absolute acceleration equals the relative acceleration (which is zero, since $\hat{\ul{c}}^{(1)}$ is constant) plus the Coriolis acceleration $2\hat{\ul{c}}^{(2)} \times \hat{\ul{c}}^{(1)}$ and the centrifugal acceleration $\hat{\ul{c}}^{(2)} \times (\hat{\ul{c}}^{(2)}\times \ul{y}(s))$. Now in case of uniform circular motion the speed is constant but the velocity is always tangent to the orbit of acceleration and the acceleration has constant magnitude and always points to the center of rotation. In this case, the total sum of Coriolis acceleration and \emph{centrifugal acceleration} add up to the well-known centripetal acceleration,
\[
\ddot{\ul{y}}(s)= 2 \hat{\ul{c}}^{(2)} \times (-\hat{\ul{c}}^{(2)} \times R\ul{r}(s))+ \hat{\ul{c}}^{2} \times (\hat{\ul{c}}^{(2)}\times R \ul{r}(s))= - \|\hat{\ul{c}}^{(2)}\|^{2} R \ul{r}(s)= - \frac{\|\hat{\ul{c}}^{1}\|^{2}}{R} \ul{r}(s),
\]
where $R$ is the radius of the circular orbit $\ul{y}(s)= \ul{m}+ R \, \ul{r}(s), \quad \|\ul{r}(s)\|=1$).
The centripetal acceleration equals half the Coriolis acceleration, i.e. $\ddot{\ul{y}}(s)=\hat{\ul{c}}^{(2)} \times \hat{\ul{c}}^{(1)}$.

In our previous work \cite[part II]{DuitsAMS1} on contour-enhancement and completion via left-invariant diffusions on invertible orientation scores (complex-valued functions on $SE(2)$) we put a lot of emphasis on the underlying differential geometry in $SE(2)$. All results straightforwardly generalize to the case of HARDI images, which can be considered as functions on $\R^3 \rtimes S^{2}$ embedded in $SE(3)$. These rather technical results are summarized in Theorem \ref{th:CartanSE2}, which answers all questions raised in the beginning of this section. Unfortunately, this theorem requires general differential geometrical concepts such as principal fiber bundles, associated vector bundles, tangent bundles, frame-bundles and the Cartan-Ehresmann connection defined on them. These concepts are explained in full detail in \cite{Spiv75b} (with a very nice overview on p.386 ).

The reader who is not familiar with these technicalities from differential geometry can skip the first part of the theorem while accepting the formula of the covariant derivatives given in Eq.~(\ref{codDV}), where the \emph{ anti-symmetric } Christoffel symbols are equal to minus the structure constants $c^{k}_{ij}=-c^{k}_{ji}$ (recall Eq.~(\ref{structureconstants})) of the Lie-algebra. Here we stress that we follow the Cartan viewpoint on differential geometry, where connections are expressed in moving coordinate frames (we use the frame of left-invariant vector fields $\{\mathcal{A}_{1}, \ldots, \mathcal{A}_{6}\}$ derived in Subsection \ref{ch:VF} for this purpose) and thereby we have non-vanishing torsion.\footnote{The torsion tensor $T_{\nabla}$ of a connection $\nabla$ is given by $T_{\nabla}[X,Y]=\nabla_{X}Y-\nabla_{Y}X-[X,Y]$. The torsion-tensor $T_{\overline{\nabla}}$ of a Levi-Civita connection vanishes, whereas the torsion-tensor of our Cartan connection $\nabla$ on $SE(3)$ is given by
$T_{\nabla}=3 \sum_{i,j,k=1}^{6}c_{ij}^{k} {\rm d}\mathcal{A}^{i} \otimes {\rm d}\mathcal{A}^{j} \otimes \mathcal{A}_{k}$.} This is different from the Levi-Civita connection for differential geometry on Riemannian manifolds, which is much more common in image analysis. The Levi-Civita connection is the unique torsion free metric compatible connection on a Riemannian manifold and because of this vanishing torsion of the Levi-Civita connection $\overline{\nabla}$ there is a 1-to-1 relation\footnote{In a Levi-Civita connection one has $\Gamma^{i}_{kl}=\Gamma^{i}_{lk}=\frac{1}{2}\sum_{m}g^{im}(g_{mk,l}+g_{ml,k}-g_{kl,m})$ with respect to a holonomic basis.} to the Christoffel symbols (required for covariant derivatives $\overline{\nabla}_{i}v^{j}=\partial_{i}v^{j} + \Gamma^{k}_{ij} \partial_{k}v^{j}$
) and the derivatives of the metric tensor. In the more general Cartan connection outlined below, however, one can have non-vanishing torsion and the Christoffels are not necessarily related to a metric tensor, nor need they be symmetric.
\begin{theorem} \label{th:CartanSE2}
The Maurer-Cartan form $\omega$ on $SE(3)$ is given by
\begin{equation} \label{MCC}
\omega_{g}(X_{g})= \sum \limits_{i=1}^{6} \langle \left.{\rm d}\mathcal{A}^{i}\right|_{g}, X_{g} \rangle A_{i}, \qquad X_{g}\in T_{g}(SE(3)),
 \end{equation}
where the dual vectors $\{{\rm d}\mathcal{A}^{i}\}_{i=1}^{6}$ are given by (\ref{duals}) and $A_{i}=\left.\mathcal{A}_{i}\right|_{e}$. 
It is a Cartan Ehresmann connection form on the principal fiber bundle
$P=(SE(3), \pi: SE(3) \to e \equiv SE(3)/SE(3), SE(3))$, where $\pi(g)=e$, $R_{g}u=ug$, $u,g \in SE(3)$. Let $\textrm{Ad}$ denote the adjoint action of $SE(3)$ on its own Lie-algebra $T_{e}(SE(3))$, i.e.
$\textrm{Ad}(g)=(R_{g^{-1}}L_{g})_{*}$, i.e. the push-forward of conjugation. Then the adjoint representation of $SE(3)$ on the vector space $\mathcal{L}(SE(3))$ of left-invariant vector fields is given by
\begin{equation} \label{Ads}
\widetilde{\textrm{Ad}}(g)={\rm d}\mathcal{R} \circ \textrm{Ad}(g) \circ \omega.
\end{equation}
This adjoint representation gives rise to the \emph{associated vector bundle } $SE(3) \times_{\widetilde{\textrm{Ad}}} \mathcal{L}(SE(3))$.
The corresponding connection form
on this vector bundle
is given by
\begin{equation} \label{tildeomega}
\tilde{\omega} = \sum \limits_{j=1}^{6} \widetilde{\textrm{ad}}(\mathcal{A}_{j}) \otimes {\rm d}\mathcal{A}^{j}= \sum \limits_{i,j,k=1}^{6} c^{k}_{ij} \; \mathcal{A}_{k} \otimes {\rm d}\mathcal{A}^{i} \otimes {\rm d}\mathcal{A}^{j},
\end{equation}
\mbox{with $\widetilde{\textrm{ad}}=(\widetilde{\textrm{Ad}})_{*}$, i.e. $\widetilde{\textrm{ad}}(\mathcal{A}_{j})=\sum \limits_{i=1}^{6}\![\mathcal{A}_{i},\mathcal{A}_{j}] \otimes {\rm d}\mathcal{A}^{i}\!$,\mbox{\cite[p.265]{Jost}}.}
Then $\tilde{\omega}$ yields the following $6\times 6$-matrix valued 1-form
\begin{equation} \label{matrix1form}
\tilde{\omega}^{k}_{j}(\cdot) := -\tilde{\omega}({\rm d}\mathcal{A}^{k}, \cdot, \mathcal{A}_{j}) \qquad k,j=1,2,\ldots,6.
\end{equation}
on the frame bundle, \cite[p.353,p.359]{Spiv75b}, where the sections are moving frames \cite[p.354]{Spiv75b}.
Let $\{\mu_{k}\}_{k=1}^{6}$ denote the sections in the tangent bundle $E:=(SE(3),T(SE(3)))$ which coincide with the left-invariant vector fields $\{\mathcal{A}_{k}\}_{k=1}^{6}$.
Then the matrix-valued 1-form given by Eq.~(\ref{matrix1form}) yields the Cartan connection
given by the covariant derivatives
\begin{equation}\label{codDV}
\begin{array}{ll}
D_{\left. X \right|_{\gamma(t)}}(\mu(\gamma(t)))&:= D \mu(\gamma(t))(\left. X \right|_{\gamma(t)}) \\
 &= \sum \limits_{k=1}^{6} \dot{a}^{k}(t) \mu_{k}(\gamma(t))+ \sum \limits_{k=1}^{6} a^{k}(\gamma(t)) \sum \limits_{j=1}^{6} \tilde{\omega}^{j}_{k}(\left. X \right|_{\gamma(t)}) \; \mu_{j}(\gamma(t)) \\
 &= \sum \limits_{k=1}^{6} \dot{a}^{k}(t) \mu_{k}(\gamma(t)) + \sum \limits_{i,j,k=1}^{6} \dot{\gamma}^{i}(t) \, a^{k}(\gamma(t))\; \Gamma^{j}_{ik} \; \mu_{j}(\gamma(t))
 \end{array}
\end{equation}
with $\dot{a}^{k}(t)= \sum \limits_{i=1}^{6}\dot{\gamma}^{i}(t) \, (\left.\mathcal{A}_{i}\right|_{\gamma(t)} a^{k})$,
for all tangent vectors
$ \left. X \right|_{\gamma(t)}= \sum \limits_{i=1}^{6}\dot{\gamma}^{i}(t)\, \left. \mathcal{A}_{i}\right|_{\gamma(t)}$
along a curve $t \mapsto \gamma(t)\in SE(2)$ and all sections
$\mu(\gamma(t))=\sum \limits_{k=1}^{6}a^{k}(\gamma(t))\, \mu_{k}(\gamma(t))$. The Christoffel symbols in (\ref{codDV}) are constant $\Gamma^{j}_{ik}=-c^{j}_{ik}$, with $c^{j}_{ik}$ the structure constants of Lie-algebra
$T_{e}(SE(3))$.
Consequently, the connection $D$ has constant curvature and constant torsion and
the left-invariant evolution equations given in Eq.~(\ref{generaleqs}) can be rewritten in covariant derivatives
(using short notation
$\nabla_{j}:=D_{\mathcal{A}_{j}}$):
{\small
\begin{equation} \label{diffcov}
\hspace{-0.068cm}\mbox{}
\!\!\left\{ \! \!
\begin{array}{l}
\partial_{t} W(g,t)=\sum \limits_{i=1}^{6}-a^{i}(W)\mathcal{A}_{i}W(g,t)+ \!\sum \limits_{i,j=1}^{6}\! \mathcal{A}_{i}\, (\, (D^{ij}(W))(g,t) \mathcal{A}_{j} W)(g,t)\\
=\sum \limits_{i=1}^{6}-a^{i}(W)\nabla_{i}W(g,t)+\! \sum \limits_{i,j=1}^{6}\! \nabla_{i}\, ((D^{ij}(W))(g,t) \nabla_{j} W)(g,t)  \\
W(g,0) =\tilde{U}(g)\ , \qquad \textrm{ for all }g \in SE(3),t>0.
\end{array}
\right.
\end{equation}
}
Both convection and diffusion in the left-invariant evolution equations (\ref{generaleqs}) take place along the exponential curves $\gamma_{\ul{c},g}(t)=g\cdot e^{t \sum \limits_{i=1}^{6} c^{i}A_{i}}$ in $SE(3)$ which are the covariantly constant curves (i.e. auto-parallels) with respect to the Cartan connection. In particular, if $a^{i}(W)=c^{i}$ constant and if $D^{ij}(W)=0$ (convection case) then the solutions are
\begin{equation} \label{LIconv}
W(g,t)= \tilde{U}(g \cdot e^{-t \sum \limits_{i=1}^{6} c^{i}A_{i}})\ .
\end{equation}
The spatial projections $\mathbb{P}_{\R^3}\gamma$ of these of the auto-parallel/exponential curves $\gamma$ are circular spirals with constant curvature and constant torsion.
The curvature magnitude equals $\|\hat{\ul{c}}^{(1)}\|^{-1} \|\hat{\ul{c}}^{(2)} \times \hat{\ul{c}}^{(1)}\|$ and the curvature vector equals
\begin{equation} \label{curvtorsion}
\KKK(t)=\frac{1}{\|\hat{\ul{c}}^{(1)}\|}\left(\cos ( t\; \|\hat{\ul{c}}^{(2)}\|)\; \hat{\ul{c}}^{(2)}\times \hat{\ul{c}}^{(1)}+  \frac{\sin(t \; \|\hat{\ul{c}}^{(2)}\|)}{ \|\hat{\ul{c}}^{(2)}\|}  \hat{\ul{c}}^{(2)} \times \hat{\ul{c}}^{(2)} \times \hat{\ul{c}}^{(1)}\right) \ ,
\end{equation}
where $\ul{c}=(c^{1},c^{2}, c^{3}\, ;\, c^{4}, c^{5}, c^{6})=(\hat{\ul{c}}^{(1)}\, ; \, \hat{\ul{c}}^{(2)})$.
The torsion vector equals $\TTT(t)=|\hat{\ul{c}}_{1}\cdot \hat{\ul{c}}_{2}|\;  \KKK(t)$.
\end{theorem}
\textbf{Proof } The proof is a straightforward generalization from our previous results \cite[Part II, Thm 3.8 and Thm 3.9]{DuitsAMS1} on the $SE(2)$-case to the case $SE(3)$. The formulas of the constant torsion and curvature of the spatial part of the auto-parallel curves (which are the exponential curves) follow by the formula (\ref{exp}) for (the spatial part $\ul{x}(s)$ of) the exponential curves, which we will derive in Section \ref{ch:exp}. Here we stress that $s(t)=t\, \sqrt{(c^{1})^2+(c^{2})^2 + (c^{3})^2}$ is the arc-length of the spatial part of the exponential curve and where we recall that $\KKK(s)=\ddot{\ul{x}}(s)$ and $\TTT(s)=\frac{d}{ds}\;(\dot{\ul{x}}(s)\times \ddot{\ul{x}}(s)  )$. Note that both the formula (\ref{exp}) for the exponential curves in the next section and the formulas for torsion and curvature are simplifications of our earlier formulas \cite[p.175-177]{FrankenThesis}. In the special case of only convection the solution (\ref{LIconv}) follows by
$e^{t {\rm d}\mathcal{R}(A)}\tilde{U}(g)=
\mathcal{R}_{e^{t A}}\tilde{U}(g)$, with $A=-\sum_{i=1}^{6}c^{i}A_{i}$ and ${\rm d}\mathcal{R}(A)=-\sum_{i=1}^{6}c^{i}\mathcal{A}_{i}$ with $\mathcal{A}_{i}={\rm d}\mathcal{R}(A_{i})$. 

\subsection{The Exponential Curves and the Logarithmic Map explicitly in Euler Angles \label{ch:exp}}

The surjective exponent mapping is given by
\begin{equation}\label{exp}
\gamma_{\ul{c}}(t)=e^{t \sum \limits_{i=1}^{6} c^{i} A_{i}}=
\left\{
\begin{array}{l}
(c_{1}t, c_{2}t, c_{3}t \; , \;  I)  \textrm{ if }\hat{\ul{c}}^{(2)}=\ul{0}\ , \\
 (t \hat{\ul{c}}^{(1)} + \frac{1-\cos(\tilde{q}t)}{\tilde{q}^2}\Omega \hat{\ul{c}}^{(1)} +( t \tilde{q}^{-2}- \frac{\sin (\tilde{q}t)}{\tilde{q}^3} )\; \Omega^{2}\hat{\ul{c}}^{(1)}\; , \;I + \frac{\sin(\tilde{q}t)}{\tilde{q}}\Omega + \frac{(1-\cos (\tilde{q}t))}{\tilde{q}^2}\Omega^2 ) \textrm{ else.} \\
\end{array}
\right.
\end{equation}
where $\tilde{q}=\|\hat{\ul{c}}^{(2)}\|=\sqrt{(c^{4})^{2}+(c^{5})^{2}+(c^{6})^{2}}$ and $\Omega \ul{x}=\hat{\ul{c}}^{(2)}\times \ul{x}$ for all $\ul{x} \in \R^{3}$.

The logarithmic mapping on $SE(3)$ is given by:
\begin{equation}\label{logSE3}
\begin{array}{ll}
\log_{SE(3)}(\ul{x}, R_{\gamma, \beta,\alpha}) &= \sum \limits_{i=1}^{3} c^{i}_{\ul{x},\gamma, \beta,\alpha} A_{i}
\; + \; \sum \limits_{i=4}^{6} c^{i}_{\gamma, \beta,\alpha} A_{i} \\
 &= \sum \limits_{i=1}^{3} c^{i}_{\ul{x},\tilde{\gamma}, \tilde{\beta},\tilde{\alpha}} A_{i}
\; + \; \sum \limits_{i=4}^{6} c^{i}_{\tilde{\gamma}, \tilde{\beta},\tilde{\alpha}} A_{i}
\end{array}
\end{equation}
Expressed in the first chart (using short notation $c^{i}:=c^{i}_{\ul{x},\gamma, \beta,\alpha}$) we have
\begin{equation}\label{one}
\begin{array}{l}
\tilde{q}=\arcsin \sqrt{\cos^{2}\left( \frac{\alpha+\gamma}{2}\right) \sin^{2}\beta + \cos^{4}\left(\frac{\beta}{2}\right)\sin^{2}(\alpha+\gamma)} \\
\ul{c}^{(2)}:=(c^{4},c^{5},c^{6})^T  =\frac{\tilde{q}}{2 \sin \tilde{q}}\; ( \sin \beta (\sin \alpha - \sin \gamma)\ , \
 \sin \beta (\cos \alpha + \cos \gamma)\ , \ 2 \cos^{2}\left(\frac{\beta}{2}\right)\sin (\alpha +\gamma))
\end{array}
\end{equation} and
\begin{equation}\label{two}
\begin{array}{l}
\ul{c}^{(1)}=(c^{1},c^{2},c^{3})^{T}= \ul{x} -\frac{1}{2}\, \ul{c}^{(2)} \times \ul{x} + \tilde{q}^{-2}(1-(\frac{\tilde{q}}{2})\, \cot(\frac{\tilde{q}}{2}))\;
\ul{c}^{(2)} \times (\ul{c}^{(2)} \times \ul{x}).
\end{array}
\end{equation}
with $\ul{c}^{(1)}=(c^{1},c^{2},c^{3})$, $\ul{c}^{(2)}=(c^{4},c^{5},c^{6})$.

Throughout this article we will take the section $\alpha=0$ (this is just a choice, we could have taken another section)
in the partition $\R^{3}\rtimes S^{2}/(\{\ul{0}\}\times SO(2))$, which means that we will
only consider the case
\[
\left.( \ul{c}^{(1)}\; ;\; \ul{c}^{(2)} )\right|_{\tilde{\alpha}=0}=
(c^{1}_{\ul{x},\gamma, \beta,\alpha=0},
c^{2}_{\ul{x},\gamma, \beta,0},c^{3}_{\ul{x},\gamma, \beta,0},c^{4}_{\ul{x},\gamma, \beta,0},c^{5}_{\ul{x},\gamma, \beta,0},c^{6}_{\ul{x},\gamma, \beta,0}).
\]
Expressed in the second chart the section $\tilde{\alpha}=0$ coincides with the section $\alpha=0$
and along this section we again have (\ref{two}) but now with
\begin{equation} \label{log2ndchart}
\begin{array}{ll}
\tilde{q}&=\arcsin \sqrt{\cos^{4}(\tilde{\gamma}/2)\, \sin^{2}(\tilde{\beta})+\cos^{2}(\tilde{\beta}/2)\, \sin^{2}(\tilde{\gamma}) }\ , \\
\ul{c}^{(2)}&=(\tilde{c}^{4},\tilde{c}^{5},\tilde{c}^{6})^{T}=
\frac{\tilde{q}}{\sin(\tilde{q})}\; (\, \sin \tilde{\gamma} \, \cos^{2}(\frac{\tilde{\beta}}{2})\, ,\, \sin \tilde{\beta}\, \cos^{2}(\frac{\tilde{\gamma}}{2})\, , \, \frac{1}{2}\sin{\tilde{\gamma}}\, \sin{\tilde{\beta}}\, )^{T}\ ,
\end{array}
\end{equation}
where we again used short notation $\tilde{c}^{i}:=\tilde{c}^{i}_{\ul{x},\tilde{\gamma}, \tilde{\beta},\tilde{\alpha}=0}$. Roughly speaking, $\ul{c}^{(1)}$ is the spatial velocity of the exponential curve (fiber) and $\ul{c}^{(2)}$ is the angular velocity of the exponential curve (fiber).

\section{Analysis of the Convolution Kernels of Scale Spaces on HARDI images \label{ch:kernels}}

It is notorious problem to find explicit formulas for the exact Green's functions $p_{t}^{\ul{D},\ul{a}}:\R^3 \rtimes S^{2}$ of the left-invariant diffusions (\ref{generalPDE}) on $\R^{3}\rtimes S^{2}$. Explicit, tangible and exact formulas for heat-kernels on $SE(3)$ do not seem to exist in literature. Nevertheless, there does exist a nice general theory overlapping the fields of functional analysis and group theory, see for example \cite{TerElst3,Nagel}, which at least provides Gaussian estimates for Green's functions of left-invariant diffusions on Lie groups, generated by subcoercive operators. In the remainder of this section we will employ this general theory to our special case where $\R^3 \rtimes S^{2}$ is embedded into $SE(3)$ and we will derive new explicit and useful approximation formulas for these Green's functions. Within this section
we 
use the second coordinate chart (\ref{secondchart}), as it is highly preferable over the more common Euler angle parametrization (\ref{firstchart}) because of the much more suitable singularity locations on the sphere.


We shall first carry out the method of contraction. This method typically relates the group of positions and rotations to a (nilpotent) group positions and velocities and serves as an essential pre-requisite for our Gaussian estimates and approximation kernels later on.
The reader who is not so much interested in the detailed analysis can skip this section and continue with the numerics explained in Chapter \ref{ch:num}.

\subsection{Local Approximation of $SE(3)$ by a Nilpotent Group via Contraction \label{ch:contraction}}

The group $SE(3)$ is not nilpotent. This makes it hard to get tangible explicit formulae for the heat-kernels.
Therefore we shall generalize our Heisenberg approximations of the Green's functions on $SE(2)$, \cite{DuitsR2006AMS}, \cite{Thornber2},\cite{MarkusThesis}, to the case $SE(3)$. Again we will follow the general work by ter Elst and Robinson \cite{TerElst3} on semigroups on Lie groups generated by weighted subcoercive operators. 
In their general work we consider a particular case by setting the Hilbert space $\mathbb{L}_{2}(SE(3))$, the group $SE(3)$ and the right-regular representation $\mathcal{R}$. Furthermore we consider the algebraic basis
\mbox{$\{\mathcal{A}_{3},\mathcal{A}_{4},\mathcal{A}_{5}\}$} leading to the following filtration of the Lie algebra
\begin{equation} \label{filtration}
\gothic{g}_{1}=\textrm{span}\{\mathcal{A}_{3},\mathcal{A}_{4},\mathcal{A}_{5}\} \subset \gothic{g}_{2}=\textrm{span}\{\mathcal{A}_{1},\mathcal{A}_{2},\mathcal{A}_{3},\mathcal{A}_{4},\mathcal{A}_{5},\mathcal{A}_{6}\}=\mathcal{L}(SE(3))\ .
\end{equation}
Now that we have this filtration we have to assign weights to the generators
\begin{equation} \label{weights}
w_{3}=w_{4}=w_{5}=1 \textrm{ and }w_{1}=w_{2}=w_{6}=2.
\end{equation}
For example $w_{3}=1$ since $\mathcal{A}_{3}$ already occurs in $\gothic{g}_{1}$, $w_{6}=2$ since $\mathcal{A}_{6}$ is within in $\gothic{g}_{2}$ and not in $\gothic{g}_{1}$.

Now that we have these weights we define the following dilations on the Lie-algebra $T_{e}(SE(3))$ (recall $A_{i}=\left.\mathcal{A}_{i}\right|_{e}$):
\[
\begin{array}{l}
\gamma_{q}(\sum \limits_{i=1}^{6} c^{i}\, A_{i})=\sum \limits_{i=1}^{6} q^{w_{i}}\, c^{i}\, A_{i}, \ \textrm{ for all } c^{i} \in \R , \\
\tilde{\gamma}_{q}(x,y,z,R_{\tilde{\gamma},\tilde{\beta},\tilde{\alpha}})=\left(\frac{x}{q^{w_1}}, \frac{y}{q^{w_2}},\frac{z}{q^{w_3}}, R_{\frac{\tilde{\gamma}}{q^{w_4}}, \frac{\tilde{\beta}}{q^{w_5}}, \frac{\tilde{\alpha}}{q^{w_6}}} \right), q>0,
\end{array}
\]
and for $0<q\leq 1$ we define the Lie product $[A,B]_{q}=\gamma^{-1}_{q}[\gamma_{q}(A),\gamma_{q}(B)]$. Now let $(SE(3))_{q}$ be the simply connected Lie group generated by the Lie algebra $(T_{e}(SE(3)),[\cdot,\cdot]_{q})$. This Lie group is isomorphic to the matrix group
with group product:
\begin{equation} \label{matrixgroup}
(\ul{x}, R_{\tilde{\gamma}\tilde{\beta}\tilde{\alpha}}) \cdot_{q} (\ul{x}', R_{\tilde{\gamma}'\tilde{\beta}'\tilde{\alpha}'} )=(\, \ul{x}+S_q \cdot R_{\tilde{\gamma}q, \tilde{\beta}q, \tilde{\alpha}q^2} \cdot S_{q^{-1}}\; \ul{x}' \; , \; R_{\tilde{\gamma}\tilde{\beta}\tilde{\alpha}} \cdot R_{\tilde{\gamma}'\tilde{\beta}'\tilde{\alpha}'}\,)
\end{equation}
where the diagonal $3\times 3$-matrix is defined by $S_{q}:=\textrm{diag}\{1,1,q\}$ and we used short-notation $R_{\tilde{\gamma}\tilde{\beta}\tilde{\alpha}}=R_{\ul{e}_{x},\tilde{\gamma}} R_{\ul{e}_{y}, \tilde{\beta}} R_{\ul{e}_{z}, \tilde{\alpha}}$, i.e. our elements of $SO(3)$ are expressed in the second coordinate chart (\ref{secondchart}).
Now the left-invariant vector fields on the group $(SE(3))_{q}$ are given by
\[
\left.\mathcal{A}_{i}^{q} \right|_{g}= (\tilde{\gamma}_{q}^{-1} \circ L_{g} \circ \tilde{\gamma}_{q})_{*}A_{i}, \qquad i=1,\ldots,6.
\]
Straightforward (but intense) calculations yield (for each $g=(\ul{x}, R_{\tilde{\gamma}\tilde{\beta}\tilde{\alpha}}) \in (SE(3))_{q}$ ):
{\small
\[
\begin{array}{ll}
\left.\mathcal{A}_{1}^{q} \right|_{g}&= \cos (q^2 \tilde{\alpha}) \cos (q \tilde{\beta})\, \partial_{x} + (\cos(\tilde{\gamma}q)\sin(\tilde{\alpha}q^2)+ \cos(\tilde{\alpha}q^2)\sin(\tilde{\beta}q) \sin(\tilde{\gamma} q))\, \partial_{y} + \\ & \qquad + q(\sin(\tilde{\alpha}q^2)\sin(\tilde{\gamma}q)-\cos(\tilde{\alpha}q^2) \cos(\tilde{\gamma}q)\sin(\tilde{\beta}q) )\, \partial_{z} \\
\left.\mathcal{A}_{2}^{q} \right|_{g} &= -\sin(\tilde{\alpha}q^2)\cos(\tilde{\beta}q)\, \partial_{x} +(\cos(q^2 \tilde{\alpha})\cos(\tilde{\gamma}q)-\sin(\tilde{\alpha}q^2)\sin(\tilde{\beta}q)\sin(\tilde{\gamma}q))\, \partial_{y} + \\ &\qquad + q(\sin(\tilde{\alpha}q^2)\sin(\tilde{\beta}q) \cos(\tilde{\beta}q)+ \cos(\tilde{\alpha}q^2)\sin(\tilde{\gamma}q))\, \partial_{z} \\
\left.\mathcal{A}_{3}^{q} \right|_{g}&=q^{-1} \sin(\tilde{\beta}q) \, \partial_{x} -q^{-1}\cos(\tilde{\beta}q)\sin(\tilde{\gamma}q)\, \partial_{y} + \cos(\tilde{\beta}q)\cos(\tilde{\gamma}q)\, \partial_{z}\, \\
\left.\mathcal{A}_{4}^{q} \right|_{g}&= -q^{-1} \cos(\tilde{\alpha}q^{2} ) \tan(\tilde{\beta}q)\, \partial_{\tilde{\alpha}} + \sin(\tilde{\alpha}q)\, \partial_{\tilde{\beta}}+ \frac{\cos(\tilde{\alpha}q^2)}{\cos(\tilde{\beta}q)}\, \partial_{\tilde{\gamma}} \\
\left.\mathcal{A}_{5}^{q} \right|_{g}&= q^{-1} \sin(\tilde{\alpha}q^2)\tan(\tilde{\beta}q)\, \partial_{\tilde{\alpha}}+ \cos(\tilde{\alpha}q^2) \, \partial_{\tilde{\beta}}- \frac{\sin(q^2 \tilde{\alpha})}{\cos(q \tilde{\beta})}\, \partial_{\tilde{\gamma}} \\
\left.\mathcal{A}_{6}^{q} \right|_{g}&=\partial_{\tilde{\alpha}}.
\end{array}
\]
}
Now note that
$
[A_{i}, A_{j}]_{q}= \gamma_{q}^{-1}[\gamma_{q}(A_i), \gamma_{q}(A_{j})]= \gamma_{q}^{-1} q^{w_{i}+w_{j}} [A_{i}, A_{j}]= \sum \limits_{k=1}^{6} q^{w_{i}+w_{j}-w_{k}} c^{k}_{ij} A_{k}$\ and thereby
we have
{\small
\begin{equation} \label{nzcomm}
\begin{array}{llll}
[A_{4},A_{5}]_{q}= A_{6},  & [A_{4},A_{6}]_{q}= -q^2 A_{5}, & \,
[A_{5},A_{6}]_{q}=q^2 A_{4}, & [A_{4}, A_{3}]_{q}=-A_{2},  \, \\
\, [A_{4},A_{2}]_{q}= q^2 A_{3}, & [A_{5}, A_{1}]_{q}=-q^2 A_{3},  &
[A_{5}, A_{3}]_{q}= A_{1}, & [A_{6}, A_{1}]_{q}= q^2 A_{2} \textrm{ and } [A_{6}, A_{2}]_{q}= - q^2 A_{1}  .
\end{array}
\end{equation}
}
Analogously to the case $q=1$, $(SE(3))_{q=1}=SE(3)$ we have an isomorphism of the common Lie-algebra at the unity element $T_{e}(SE(3))=T_{e}((SE(3))_{q})$ and left-invariant vector fields on the group $(SE(3))_{q}$ :
\[
(A_{i} \leftrightarrow \mathcal{A}_{i}^{q} \textrm{ and }A_{j} \leftrightarrow \mathcal{A}_{j}^{q}) \Rightarrow [A_{i},A_{j}]_{q} \leftrightarrow [\mathcal{A}_{i}^{q},\mathcal{A}_{j}^{q}]\ .
\]
It can be verified that the left-invariant vector fields $\mathcal{A}_{i}^{q}$ satisfy the same commutation relations (\ref{nzcomm}).

Now let us consider the case $q \downarrow 0$, then we get a nilpotent-group $(SE(3))_{0}$ with left-invariant vector fields
\begin{equation} \label{genheis}
\begin{array}{l}
\mathcal{A}_{1}^{0}=\partial_{x} ,\
\mathcal{A}_{2}^{0}= \partial_{y} , \
\mathcal{A}_{3}^{0}= \tilde{\beta}\partial_{x}-\tilde{\gamma}\partial_{y} +\partial_{z}, \
\mathcal{A}_{4}^{0}=-\tilde{\beta} \partial_{\tilde{\alpha}} +\partial_{\tilde{\gamma}} , \
\mathcal{A}_{5}^{0}=\partial_{\tilde{\beta}} , \
\mathcal{A}_{6}^{0}=\partial_{\tilde{\alpha}}\ .
\end{array}
\end{equation}

\subsubsection{The Heisenberg-approximation of the Time-integrated $k$-step Contour Completion Kernel}

Recall that the generator of contour completion diffusion equals $\mathcal{A}_{3} + D^{44}((\mathcal{A}_{4})^{2}+(\mathcal{A}_{5})^2)$. So let us replace the true left-invariant vector fields $\{\mathcal{A}_{i}\}_{i=3}^{5}$ on $SE(3)=(SE(3))_{q=1}$ by their Heisenberg-approximations $\{\mathcal{A}_{i}^{0}\}_{i=3}^{5}$ that are given by (\ref{genheis}) and compute the Green's function
$\overline{p}_{t}^{a_{3}=1,D^{44}=D^{55}}$ on $(SE(3))_{0}$ (i.e. the convolution kernel which yields the solutions of contour completion on $(SE(3))_{0}$ by group convolution on $(SE(3))_{0}$). For $0<D^{44}<<1$ this kernel is a local approximation of the true contour completion kernel\footnote{The superscript for the kernel is actually $p_{t}^{D, \ul{a}}$ so in the superscript-labels, for the sake of simplicity, we only mention the \emph{non-zero} coefficients $D^{ij}$, $a_{i}$ of (\ref{generalPDE}).} $p_{t}^{a_{3}=1,D^{44}=D^{55}}$, on $\R^3 \rtimes S^{2}$:
\begin{equation} \label{formule}
\begin{array}{ll}
\overline{p}_{t}^{a_{3}=1,D^{44}=D^{55}} &:=e^{t(\mathcal{A}_{3}^{0}+D^{44}((\mathcal{A}_{4}^{0})^2+(\mathcal{A}_{5}^{0})^2))} \delta_{0}^{x} \otimes \delta_{0}^{y} \otimes \delta_{0}^{z}
 \otimes \delta_{0}^{\tilde{\gamma}} \otimes \delta_{0}^{\tilde{\beta}} \Rightarrow \\
\overline{p}_{t}^{a_{3}=1,D^{44}=D^{55}}(x,y,z,\tilde{\ul{n}}(\tilde{\beta},\tilde{\gamma})) &=\delta(t-z) \; (e^{t(\tilde{\beta}\partial_{x}+D^{44}(\partial_{\tilde{\beta}})^2)} \delta_{0}^{x}\otimes \delta_{0}^{\tilde{\beta}})(x, \tilde{\beta}) \; (e^{t(-\tilde{\gamma}\partial_{y}+D^{44}(\partial_{\tilde{\gamma}})^2)} \delta_{0}^{y}\otimes \delta_{0}^{\tilde{\gamma}})(y, \tilde{\gamma}) \\
 &=\delta(t-z) \frac{3}{4 (D^{44} \pi z^2)^2} \, e^{-\frac{12(x- (1/2)z\tilde{\beta})^2 + z^2 \tilde{\beta}^{2}}{4 z^3 D^{44}}}\, e^{-\frac{12(y+ (1/2)z\tilde{\gamma})^2 + z^2 \tilde{\gamma}^{2}}{4 z^3 D^{44}}}\ ,
 \end{array}
 \end{equation}
where $\tilde{\ul{n}}(\tilde{\beta},\tilde{\gamma})=R_{\ul{e}_{x}, \tilde{\gamma}} R_{\ul{e}_{y}, \tilde{\beta}} \ul{e}_{z}=(\sin \tilde{\beta},-\sin  \tilde{\gamma} \; \cos  \tilde{\beta}, \cos \tilde{\gamma} \cos \tilde{\beta})^{T}$. The corresponding $k$-step resolvent kernel on the group $(SE(3))_{0}$ is now directly obtained by conditional integration over time\footnote{Note that the delta distribution $\delta(s-z)$ allowed us to replace all $s$ by $z$ in the remaining factor in (\ref{formule}) which makes it easy to apply the integration $R_{\gamma,k}^{a_{3}=1,D^{33},D^{44}}= \int_{\R^{+}} p_{t}^{a_{3}=1,D^{33},D^{44}} \Gamma(t\, ;\, k,\lambda) {\rm d t}$.}
{\small
\begin{equation} \label{kresolvent}
\overline{R}_{\lambda,k}^{a_{3}=1,D^{44}=D^{55}}(x,y,z,\tilde{\ul{n}}(\tilde{\beta},\tilde{\gamma}))=
\left\{
\begin{array}{l}
\frac{3}{4 (D^{44} \pi)^2}\, \frac{(\lambda)^k z^{k-5}}{(k-1)!} \, e^{-\lambda z} \, e^{-\frac{12(x- (1/2)z\tilde{\beta})^2 + z^2 \tilde{\beta}^{2}}{4 z^3 D^{44}}}\, e^{-\frac{12(y+ (1/2)z\tilde{\gamma})^2 + z^2 \tilde{\gamma}^{2}}{4 z^3 D^{44}}}\ \textrm{if } z>0 \\
0 \qquad \textrm{ if }z\leq 0 \textrm{ and }(x,y) \neq (0,0).
\end{array}
\right.
\end{equation}
}

\subsubsection{Approximations of the Contour Enhancement Kernel}

Recall that the generator of contour completion diffusion equals $D^{33}(\mathcal{A}_{3})^2 + D^{44}((\mathcal{A}_{4})^{2}+(\mathcal{A}_{5})^2)$. So let us replace the true left-invariant vector fields $\{\mathcal{A}_{i}\}_{i=3}^{5}$ on $SE(3)=(SE(3))_{q=1}$ by their Heisenberg-approximations $\{\mathcal{A}_{i}^{0}\}_{i=3}^{5}$ given by (\ref{genheis}) and consider the Green's function $\overline{p}_{t}^{D^{33},D^{44}=D^{55}}$ on $(SE(3))_{0}$:

Now since the Heisenberg approximation kernel $\overline{p}_{t}^{D^{33},D^{44}\; ;\; (SE(2))_{0}}$ is for
reasonable parameter settings (that is {\small $0<\frac{D^{44}}{D^{33}}<<1$}) close to the
exact kernel $\overline{p}_{t}^{D^{33},D^{44}\; ;\; (SE(2))}$ we \emph{heuristically propose} for
these reasonable parameter settings the same direct-product approximation for the exact contour-enhancement kernels on $\R^{3} \rtimes S^{2}$:
{\small
\begin{equation} \label{heur2}
p_{t}^{D^{33},D^{44}=D^{55}\; ;\; \R^{3} \rtimes S^{2}}(x,y,z,\tilde{\ul{n}}(\tilde{\beta},\tilde{\gamma})) \approx
N(D^{33},D^{44},t)\, p_{t}^{D^{33},D^{44}\; ;\; (SE(2))}(z/2,x,\tilde{\beta})\cdot
p_{t}^{D^{33},D^{44}\; ;\; (SE(2))}(z/2,-y,\tilde{\gamma})\ ,
\end{equation}
}
where
\begin{equation} \label{l1norm}
N(D^{33},D^{44},t)\approx \frac{1}{8 \sqrt{2}} \sqrt{\pi} t \sqrt{t D^{33}} \sqrt{D^{33}D^{44}}
\end{equation}
takes care of $\mathbb{L}_{1}(\R^3 \rtimes S^2)$-normalization and
with
\begin{equation}\label{estimates}
\begin{array}{l}
\overline{p}_{t}^{D^{33},D^{44}}(x,y,\theta) \equiv \frac{1}{4\pi t^{2} D^{44}D^{33}} e^{- \frac{1}{4t\, c^2} \sqrt{ \left(\frac{x^2}{D^{33}}+\frac{\theta^2}{D^{44}}\right)^{2} + \frac{|y-\frac{x\theta}{2}|^{2}}{D^{44}D^{33}}} }, \\
p_{t}^{D^{33},D^{44}}(x,y,\theta) \equiv
\left\{
\begin{array}{l}
\frac{1}{4\pi t^{2} D^{44}D^{33}} e^{-\frac{1}{4t\, c^2} \sqrt{\left( \frac{\theta^2}{D^{44}} + \frac{\theta^{2}(y-(- x \, \sin \theta + y \, \cos \theta))^2}{4(1-\cos(\theta))^2 D^{33}}\right)^2 + \frac{1}{D^{44}D^{33}}\left|\frac{\theta((x\, \cos \theta  +y\, \sin \theta) -x)}{2(1-\cos \theta)}\right|^2}}, \textrm{if }\theta \neq 0, \\
\frac{1}{4\pi t^{2} D^{44}D^{33}} e^{-\frac{1}{4t\, c^2} \sqrt{\left( \frac{x^2}{D^{33}}\right)^2+\frac{|y|^{2}}{D^{44}D^{33}}}} , \textrm{ if } \theta=0
 \end{array}
 \right.
\end{array}
\end{equation}
which are reasonably sharp estimates of hypoelliptic diffusion on $SE(2)$, with $\frac{1}{2} \leq c \leq \sqrt[4]{2}$,
for details see \cite[ch 5.4]{DuitsAMS1}. For the purpose of numerical computation, we simplify $p_{t}^{D^{33},D^{44}}(x,y,\theta)$ in (\ref{estimates}) to
{\small
\[
p_{t}^{D^{33},D^{44}}(x,y,\theta)=
\frac{1}{4\pi t^{2} D^{44}D^{33}} e^{-\frac{1}{4t\, c^2} \sqrt{\left( \frac{\theta^2}{D^{44}} +
\frac{\left(\frac{\theta y}{2}+ \frac{\theta/2}{\tan (\theta/2)}\, x \right)^2}{D^{33}} \right)^2
+ \frac{1}{D^{44}D^{33}}\left(\frac{-x \theta}{2}+ \frac{\theta/2}{\tan (\theta/2)}\, y \right)^2 }}
\]
}
where one can use the estimate
$\frac{\theta/2}{\tan (\theta/2)} \approx \frac{\cos (\theta/2)}{1-(\theta^2/24)}$ for $|\theta|< \frac{\pi}{10}$ to avoid numerical errors.


\subsection{Gaussian Estimates for the Heat-kernels on $SE(3)$ \label{ch:GE}}

In \cite[ch:6.2]{DuitsFrankenCASA} it is shown that the constants $\tilde{C}_{3},\tilde{C}_{4}$ are very close and that a reasonably sharp approximation and upperbound of the horizontal diffusion kernel
on $\R^{3} \rtimes S^{2}$ is given by
\begin{equation} \label{estfin}
\begin{array}{l}
p_{t}^{D=\textrm{diag}\{0,0,D^{33},D^{44},D^{55},0\}}(\ul{y},\tilde{n}(\tilde{\beta},\tilde{\gamma}))
 \approx
\frac{1}{(4\pi t^2 D^{33} D^{44})^2}\,  e^{-\frac{|(\ul{y},R_{\ul{n}})|^{2}_{D^{33},D^{44}}}{4t}},
\end{array}
\end{equation}
with weighted modulus
\begin{equation} \label{diffusionmodulus}
|(\ul{y},R_{\ul{n}})|_{D^{33},D^{44}}^{2}:=\, \sqrt{\frac{|c^{1}|^2 + |c^{2}|^2}{D^{33}D^{44}}+ \frac{|c^{6}|^2}{D^{44}} +  \left( \frac{(c^{3})^2}{D^{33}}+
\frac{|c^{4}|^2 + |c^{5}|^2}{D^{44}} \right)^2}
\end{equation}
where $\ul{n}=\tilde{\ul{n}}(\tilde{\beta},\tilde{\gamma}) \in S^2$ and where we again use \emph{short notation} $c^{k}:=c^{k}_{q=1}(\ul{x}, R_{\tilde{\gamma},\tilde{\beta},0})$, $k=1,\ldots,6$. Recall from Section \ref{ch:exp} that these constants are computed by the logarithm (\ref{logSE3}) on $SE(3)$ or more explicitly by (\ref{log2ndchart}).
\subsection{Analytic estimates for the Green's functions of the Cost Processes on $\R^3 \rtimes S^2$}

In Appendix \ref{app:viscosity} we have derived the following analytic approximation for the Green's function
$k_{t}^{D^{11},D^{44}, \eta, \pm}$
of the Hamilton-Jacobi equation (\ref{Hamdi1}) and corresponding cost-process explained in Section \ref{ch:cost},
Eq.~(\ref{Bellman}) on $\R^{3}\rtimes S^2$:
\begin{equation}
k_{t}^{D^{11},D^{44}, \eta, \pm}(\ul{y},\tilde{\ul{n}}(\tilde{\beta},\tilde{\gamma})) \equiv \pm \,
\frac{2\eta-1}{2\eta} C^{\frac{2\eta}{2\eta-1}} t^{-\frac{1}{2\eta-1}}
\left(\sum \limits_{i=1}^{6}
\frac{|\tilde{c}^{i}(\ul{y},\tilde{\alpha}=0,\tilde{\beta},\tilde{\gamma})|^{\frac{2}{w_{i}}}}{D^{ii}} \right)^{\frac{\eta}{2\eta-1}}
\end{equation}
for sufficiently small time $t>0$, with weights $w_{1}=w_{2}=w_{4}=w_{5}=1$, $w_{3}=w_{6}=2$ and $2>C>0$ and
where the $(c^{1},\ldots,c^{6})$ are given by (\ref{log2ndchart}) and where we have set $\ul{n}=\tilde{\ul{n}}(\tilde{\beta},\tilde{\gamma})$ and $c^{k}:=c^{k}_{q=1}(\ul{y}, R_{\tilde{\gamma},\tilde{\beta},0})$.
Here we recall that $D^{11}$ tunes the spatial erosion orthogonal to fibers,
$D^{44}$ tunes angular erosion and $\eta \in (\frac{1}{2},1]$ relates homogeneous erosion $\eta=\frac{1}{2}$ to
standard quadratic erosion $\eta=1$.

\begin{figure}[h] \vspace{-0.15cm}\mbox{}
\centerline{
\hfill
\includegraphics[width=0.6\hsize]{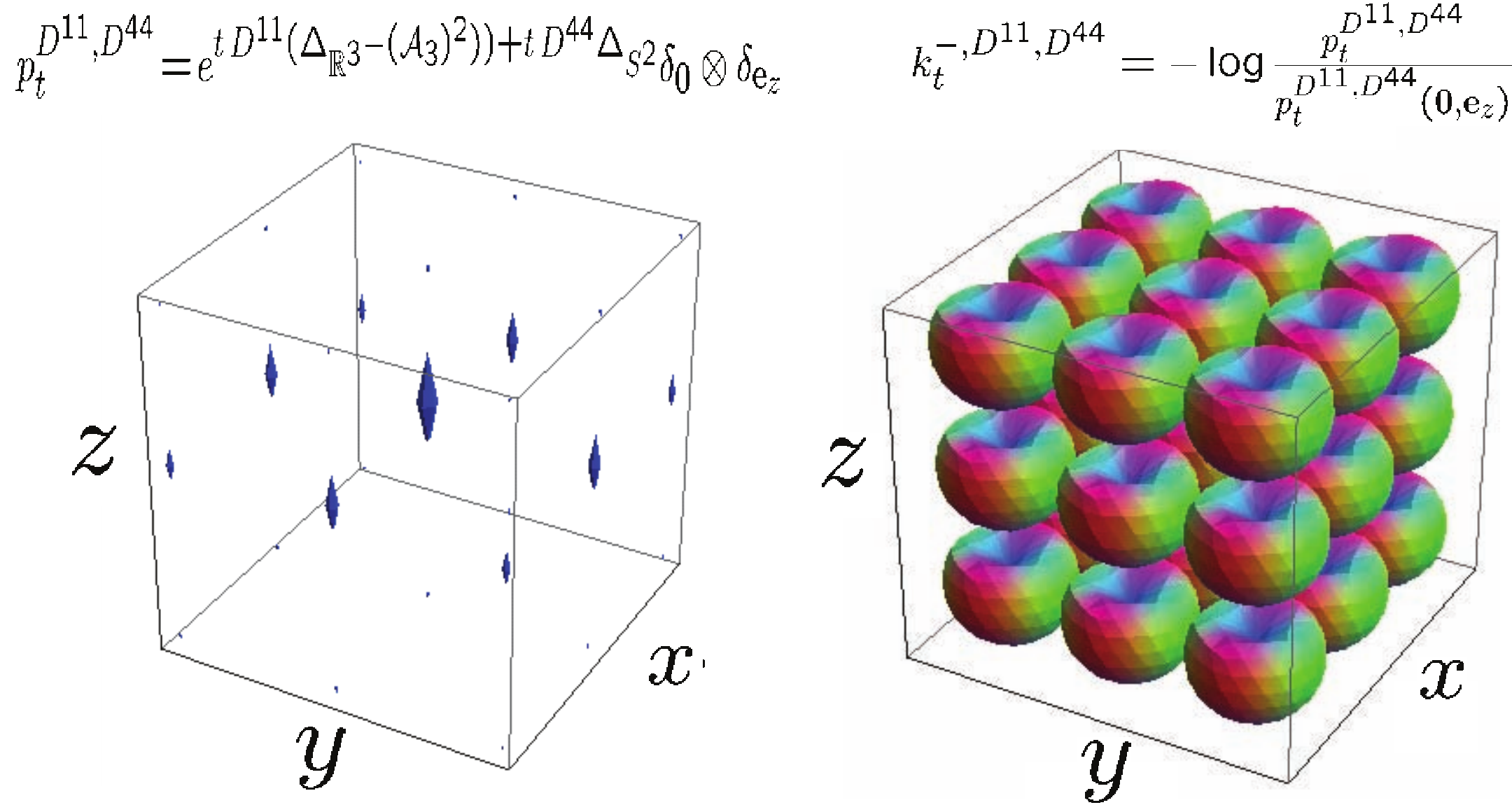}
\hfill
}
\vspace{-0.3cm}\mbox{}
\caption{Left: Diffusion kernel in Eq.~(\ref{corrDiff}), where angular diffusion takes place over the full $S^{2}$ sphere and spatial diffusion simultaneously takes place in the plane
$\ul{G}^{{\rm D^{11}\!,D^{33}\!, D^{44}}}$-orthogonal to $\mathcal{A}_{3}$. Right: Erosion kernel in (\ref{approx}).
\label{fig:HARDI-DTI}}\label{fig:erosion}
\vspace{-0.3cm}\mbox{}
\end{figure}

\begin{figure}[h] \vspace{-0.15cm}\mbox{}
\centerline{
\hfill
\includegraphics[width=0.5\hsize]{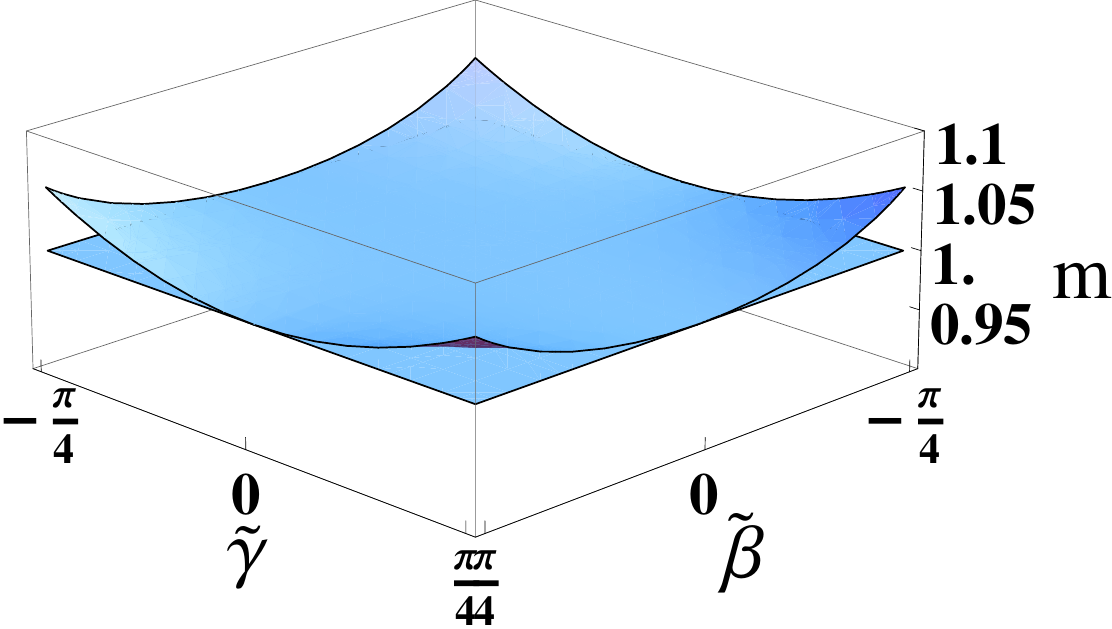}
\hfill
}
\vspace{-0.3cm}\mbox{}
\caption{This figure shows that our approximation (\ref{approx}) is rather accurate for relevant parameter settings
$\tilde{\beta},\tilde{\gamma} \in [-\frac{\pi}{4},\frac{\pi}{4}]$. We have computed the surface $m=m(\tilde{\beta},\tilde{\gamma})$
where $m(\tilde{\beta},\tilde{\gamma})$ is given by (\ref{mformule}). The approximation is exact iff $m=1$, which is the case along the lines $\tilde{\beta}=0$ and $\tilde{\gamma}=0$. Relative errors are smaller than 0.05 percent.
\label{fig:goodapprox}}
\vspace{-0.3cm}\mbox{}
\end{figure}

Now analogously to our approach on $SE(2)$ in \cite[Ch:5.4]{DuitsAMS1} we use the estimate
\[
|a|+|b| \geq \sqrt{|a|^2+|b|^2} \geq \frac{1}{\sqrt{2}} (|a|+|b|)
\]
$a,b \in \R$ to obtain differentiable analytic local approximations:
\begin{equation}
\label{approx}
\begin{array}{l}
k_{t}^{D^{11},D^{44}, \eta,+}(\ul{y},\ul{n}) \approx
\frac{2\eta-1}{2\eta} C^{\frac{2\eta}{2\eta-1}} t^{-\frac{1}{2\eta-1}} \left(\left(\frac{|c^1|^2 +|c^2|^2}{D^{11}} +\frac{|c^4|^2+|c^{5}|^2}{D^{44}}\right)^2 + \frac{|c^3|^2}{D^{11}D^{44}}\right)^{\frac{\eta}{2\eta-1}} , \eta \in (\frac{1}{2},1]
\end{array}
\end{equation}
and for $\eta=\frac{1}{2}$ we obtain the flat analytic local approximation (that arises by taking the limit $\eta \downarrow \frac{1}{2}$):
\begin{equation}
\label{approxetahalf}
\begin{array}{l}
k_{t}^{D^{11},D^{44}, \eta, +}(\ul{y},\ul{n}) \approx
\{
\begin{array}{l}
\infty \textrm{ if } \sqrt{\left(\frac{|c^1|^2 +|c^2|^2}{D^{11}} +\frac{|c^4|^2+|c^{5}|^2}{D^{44}}\right)^2 + \frac{|c^3|^2}{D^{11}D^{44}}} \geq t^2 \\
0      \textrm{ else }
\end{array}
\end{array}
\end{equation}
See Figure \ref{fig:erosion} for glyph visualizations (recall Definition \ref{def:viz}) of the erosion and corresponding diffusion kernel
\begin{equation} \label{corrDiff}
p_{t}^{D=\textrm{diag}\{D^{11},D^{11},0,D^{44},D^{55},0\}}=e^{t(D^{11}(\Delta_{\R^3}-(\mathcal{A}_{3})^2)) + D^{44} \Delta_{LB}} \delta_{0} \otimes \delta_{\ul{e}_{z}}
\end{equation}
on the contact manifold $(\R^{3}\rtimes S^{2}, {\rm d}\mathcal{A}^{3})$.
In order to study the accuracy of this approximation formula for the case of angular erosion only (i.e. $D^{11}=0$)
we have analytically computed
\begin{equation}\label{mformule}
\begin{array}{ll}
m(\tilde{\beta},\tilde{\gamma}) &= \frac{|\mathcal{A}_{4}k_{t}^{D^{11}=0,D^{44}=1}|^2+|\mathcal{A}_{5}k_{t}^{D^{11}=0,D^{44}=1}|^2}{\partial_{t}k_{t}^{D^{11}=0,D^{44}=1}}=
\frac{|\partial_{\tilde{\beta}}k_{t}^{D^{11}=0,D^{44}=1}|^2 +(\cos \tilde{\beta})^{-2} |\partial_{\tilde{\gamma}}k_{t}^{D^{11}=0,D^{44}=1}|^2}{\partial_{t}k_{t}^{D^{11}=0,D^{44}=1}} \\
 &=\frac{1}{4} \frac{(\partial_{\tilde{\beta}}((c^4)^2+(c^{5})^2))^2+ \cos^{-2} \tilde{\beta} (\partial_{\tilde{\gamma}}((c^4)^2+(c^{5})^2))^2}{(c^4)^2+(c^{5})^2},
\end{array}
\end{equation}
where we used the following identities
\[
\begin{array}{l}
|\mathcal{A}_{4}U|^2+|\mathcal{A}_{5}U|^2= |\partial_{\tilde{\beta}}^{2}U|^{2}+(\cos \tilde{\beta})^{-2} |\partial_{\tilde{\gamma}}^{2}U|^{2} \\
|\mathcal{A}_{1}U|^2+|\mathcal{A}_{2}U|^2=|\partial_{x}U|^2+|\partial_{y}U|^2+|\partial{z}U|^2- |\mathcal{A}_{3}U|^2 \textrm{ for all }U \in C^{1}(\R^{3}\rtimes S^{2}).
\end{array}
\]
Ideally $m=1$ since then the approximation is exact. For relevant parameter settings we indeed have $m \approx 1$ up to $5$-percent $\ell_{\infty}$-errors as can be seen in Figure \ref{fig:goodapprox}.

\section{Pseudo Linear Scale Spaces on $\R^{3}\rtimes S^2$ \label{ch:pseudo}}

So far we have considered anisotropic diffusions aligned with fibers and erosions orthogonal to the fibers. As these two types of left-invariant evolutions are supposed to be alternated,
\[
\Phi_{t/2}^{\ul{G}=-\textrm{diag}\{g^{11},g^{11},0,g^{44},g^{44},0\}} \circ e^{s Q^{\ul{D}=\textrm{diag}\{D^{11},D^{11},D^{33},D^{44},D^{44},0\},\ul{a}=\ul{0}}(\nabla)} \circ \Phi_{t/2}^{\ul{G}=-\textrm{diag}\{g^{11},g^{11},0,g^{44},g^{44},0\}}
\]
where $g^{ij}$ are the components of the inverse metric, which in case of a diagonal metric tensor simply reads $g^{ii}=g_{ii}^{-1}=D^{ii}$, like in (\ref{Hamdi2})
where $D^{33}\gg D^{11}$ and where
$\Phi_{t}$ denotes the viscosity solution operator $U \mapsto W(\cdot,\cdot,t)$ for the erosion Hamilton-Jacobi equation (\ref{Hamdi2}).
the natural question arises is there a single evolution process that combines erosion/dilation and diffusion.
Moreover, a from a practical point of view quite satisfactory alternative to visually sharpen distributions on positions and orientations is to apply monotonic greyvalue transformations (instead of an erosion) such as for example the power operator $\tilde{\chi}_{p}$ 
\[
\begin{array}{l}
(\tilde{\chi}_{p}(U))(\ul{y},\ul{n})=(U(\ul{y},\ul{n}))^{p}\ , 
\end{array}
\]
where $p\geq 1$, where we also recall the drawback illustrated in Figure \ref{Fig:5}
Conjugation of the diffusion operator with a monotonically increasing grey-value transformation $\chi:\R^{+} \to \R^{+}$
\begin{equation} \label{conjGT}
\chi^{-1} \circ e^{t\, Q^{\ul{D},\ul{a}}(\nabla)} \circ \chi
\end{equation}
is related to simultaneous erosion and diffusion. For a specific choice of grey-value transformations this is indeed the case we will show next, where we extend the theory for pseudo linear scale space representations
of greyscale images \cite{Florackpseudo} to DW-MRI (HARDI and DTI).

\begin{figure}
\centerline{
\includegraphics[width= 0.45 \hsize]{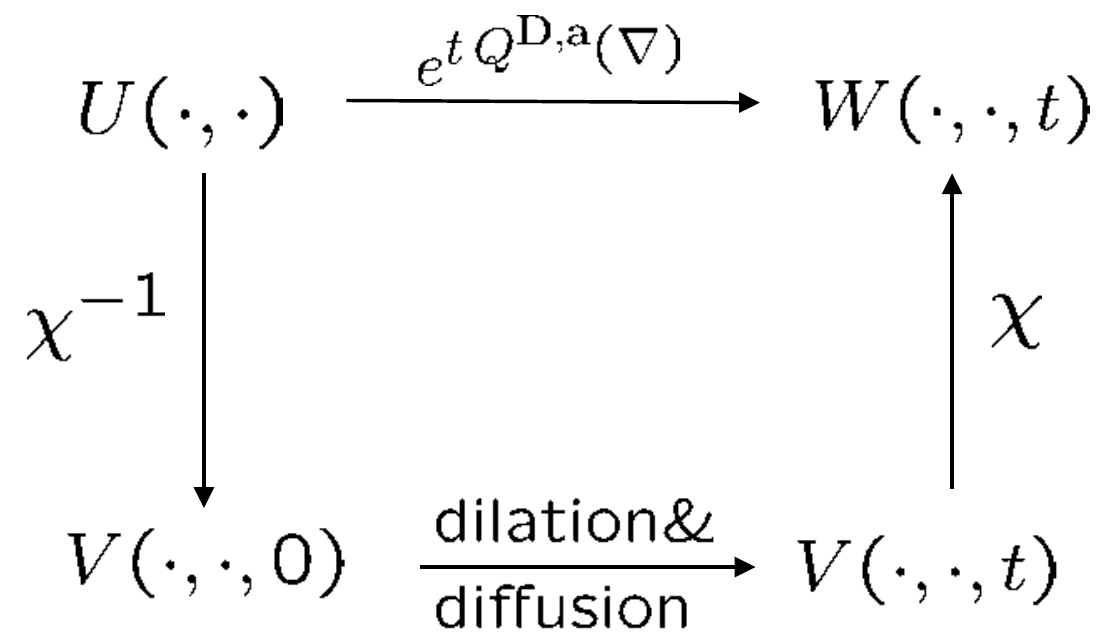}
}
\caption{The commutative diagram of simultaneous left-invariant dilation and diffusion along the fibers ($D^{33}\gg D^{11}$)
given by Eq.~(\ref{generalPDE2}) with $\tilde{\mu}(V(\ul{y},\ul{n},t))=C$ illustrates that conjugation of
specific grey-value transformations (\ref{GT}) with left-invariant diffusion is equivalent.\label{fig:Comm}}
\end{figure}

Next we derive the operator $\chi^{-1} \circ e^{t \, Q^{\ul{D},\ul{a}}(\nabla)} \circ \chi$ more explicitly using the chain-law for differentiation:
\[
\begin{array}{rl}
\mathcal{A}_{j}(\chi(V(\cdot,\cdot,t)))&=\chi'(V(\cdot,\cdot,t))\mathcal{A}_{j}(V(\cdot,\cdot,t))\ ,\\
\partial_{t}\chi(V(\cdot,\cdot,t))&= \chi'(V(\cdot,\cdot,t))\partial_{t}(V(\cdot,\cdot,t))\ , \\
\mathcal{A}_{i}\mathcal{A}_{j}(\chi(V(\cdot,\cdot,t)))&=\chi''(V(\cdot,\cdot,t))\mathcal{A}_{i}(V(\cdot,\cdot,t)) \, \mathcal{A}_{j}(V(\cdot,\cdot,t)) + \chi'(V(\cdot,\cdot,t)) \mathcal{A}_{i}\mathcal{A}_{j}(V(\cdot,\cdot,t)).
\end{array}
\]
Consequently, if we set $W(\cdot,\cdot,t)=\chi(V(\cdot,\cdot,t)) \desda V(\cdot,\cdot,t)=\chi^{-1}(W(\cdot,\cdot,t))$ and
\[
W(\cdot,\cdot,t)=e^{t Q^{\ul{D}=\textrm{diag}\{D^{11},D^{11},D^{33},D^{44},D^{44},0\}},\ul{a}=\ul{0}}W(\cdot,\cdot,0),
\]
see Figure \ref{fig:Comm}, then $V$ satisfies
\begin{equation} \label{generalPDE2}
\boxed{
\left\{
\begin{array}{l}
\partial_{t}V(\ul{y},\ul{n},t)= \sum \limits_{i=1}^{5}D^{ii}\left((\left.\mathcal{A}^{i}\right|_{(\ul{y},\ul{n})})^{2}V\right)(\ul{y},\ul{n},t)+
\tilde{\mu}(V(\ul{y},\ul{n},t))\sum \limits_{i=1}^{5}D^{ii}\left(\left.\mathcal{A}_{i}\right|_{(\ul{y},\ul{n})}V(\ul{y},\ul{n},t)\right)^2  \, \ ,\\
V(\ul{y},\ul{n},0) =\chi^{-1}(U(\ul{y},\ul{n}))\ .
\end{array}
\right.
}
\end{equation}
where $\tilde{\mu}(V(\ul{y},\ul{n},t))=\frac{\chi''(\chi(V(\ul{y},\ul{n},t)))}{\chi'(\chi(V(\ul{y},\ul{n},t)))}$.
So if we set $\tilde{\mu}(V(\ul{y},\ul{n},t))=C$ constant we achieved that (\ref{conjGT}) coincides with a
simultaneous erosion/dilation and diffusion with
\[
\ul{G}^{-1}=\sum \limits_{i,j=1}^{5} g^{ij}\mathcal{A}_{i}\otimes \mathcal{A}_{j}= \sum \limits_{i,j=1}^{5} D^{ij} \mathcal{A}_{i}\otimes \mathcal{A}_{j}\ .
\]
This means we have to solve the following ODE-system
\begin{equation} \label{ODE}
\begin{array}{l}
\chi''(I)-C\;\chi'(I)=0 \qquad I \in [0,1]\ , \\
\chi(0)=0 \textrm{ and }\chi(1)=1\ ,
\end{array}
\end{equation}
where $I \in [0,1]$ stands for intensity where in particular we set
$I:=
\frac{V(\ul{y},\ul{n},0)-\min_{(\ul{y},\ul{n})}V(\ul{y},\ul{n},0) }{
\max_{\ul{y},\ul{n}} V(\ul{y},\ul{n},0)- \min_{(\ul{y},\ul{n})}V(\ul{y},\ul{n},0)}$
 $(\ul{y},\ul{n}) \in \R^{3}\rtimes S^2$.
The unique solutions of (\ref{ODE}) are given by
\begin{equation}\label{GT}
\chi_{C}(I)=
\left\{
\begin{array}{ll}
\frac{e^{C \, I}-1}{e^{C}-1} &\textrm{ if }C \neq 0\ , \\
I &\textrm{ if }C=0\ ,
\end{array}
\right.
\end{equation}
so that $V(\ul{y},\ul{n},t)$ is the solution of an evolution (\ref{generalPDE}) where the generator is a weighted sum of a diffusion and erosion/dilation operator. The inverse of $\chi$ is given by
\begin{equation}\label{GTinv}
\chi_{C}^{-1}(I)=
\left\{
\begin{array}{ll}
\frac{1}{C} \ln (1+ (e^{C}-1)I) &\textrm{ if }C \neq 0\ , \\
I &\textrm{ if }C=0\ .
\end{array}
\right.
\end{equation}
\begin{figure}
\centerline{
\includegraphics[width= 0.241 \hsize]{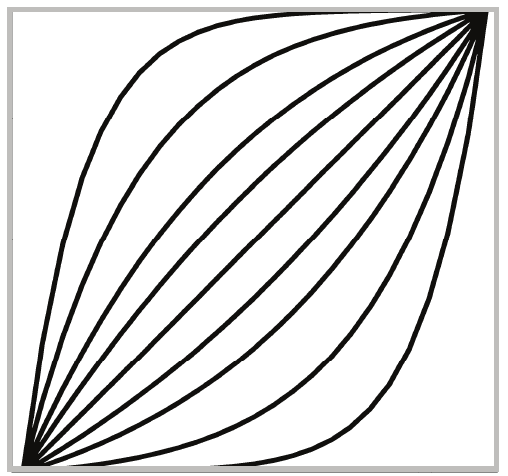}
\includegraphics[width= 0.24\hsize]{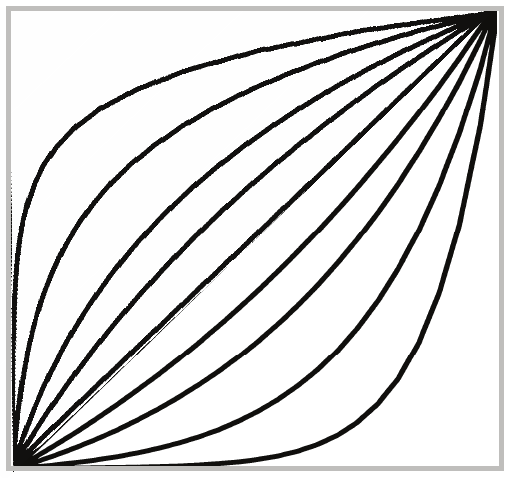}
}
\vspace{-0.5cm}\mbox{}
\caption{Left: Graphs of the grey-value transformations, Eq.~(\ref{GT}),
that are solutions $\chi_{C}:[0,1] \to [0,1]$ of Eq.~(\ref{ODE}), for $C=-8,-4,-2,-1,0,1,2,4,8$, \cite{Florackpseudo} (depicted from left to right). The concave solutions correspond to diffusion and erosion whereas the convex solutions yield convection and diffusion.
Right: Graphs of the convex $\chi_{C}$ of Eq.~(\ref{ODE})
depicted together with their inverse $\chi_{C}^{-1}$, Eq (\ref{GTinv})  , for $C=-8,-4,-2,-1$
\protect\label{fig:GT}}
\vspace{-0.3cm} \mbox{}
\end{figure}
The drawback of this intriguing correspondence, is that our diffusions primarily take place along the fibers, whereas our erosions take place orthogonal to the fibers. Therefore, at this point the correspondence
between pseudolinear scale spaces and hypo-elliptic diffusion conjugated with $\xi$ is primarily useful for \emph{simultaneous dilation and diffusion along the fibers}\footnote{In general one does not want to
erode and diffuse in the same direction.}, i.e. the case
\[
D=\textrm{diag}\{0,0,D^{11},D^{22}, D^{33},0\}
\textrm{ and }C>0.
\]
\begin{remark}
Finally, for $C\approx 2$ the operator in Eq.~(\ref{conjGT}) is close to the operators applied in Figure~\ref{fig:1} and Figure~\ref{Fig:2} on DTI and HARDI data of the brain. 
The difference though is that (\ref{conjGT}) conjugates a diffusion operator with greyvalue transformations, whereas in Figure~\ref{fig:1}~and~\ref{Fig:2} we applied
\mbox{$\chi_{p} \circ e^{t\, Q^{\ul{D},\ul{a}}(\nabla)} \circ \chi_p$} with $p=2$.
\end{remark}

\section{Implementation of the Left-Invariant Derivatives and $\R^{3} \rtimes S^{2}$-Evolutions \label{ch:num}}

In our implementations we do \emph{not} use the two charts (among which the Euler-angles
parametrization) of $S^{2}$ because this would involve cumbersome and expensive bookkeeping of mapping the coordinates from one chart to the other (which becomes necessary each time the singularities (\ref{sing1}) and (\ref{sing2}) are reached).
Instead we recall that the left-invariant vector fields on HARDI-orientation scores $\tilde{U}:SE(3) \to \R$, which by
definition (recall Definition \ref{def:OS}) automatically satisfy
\begin{equation} \label{alphainv}
\tilde{U}(R R_{\ul{e}_{z,\alpha}})=\tilde{U}(R),
\end{equation}
are constructed by the derivative of the right-regular
representation
\[
\mathcal{A}_{i}\tilde{U}(g)=({\rm d}\mathcal{R}(A_{i})\tilde{U})(g)= \lim \limits_{h \downarrow 0}\frac{\tilde{U}(g\, e^{ h\, A_{i}})-\tilde{U}(g)}{h}=
\lim \limits_{h \downarrow 0}\frac{\tilde{U}(g\, e^{h \, A_{i}})-\tilde{U}(g\, e^{-h \, A_{i}})}{2h}\ ,
\]
where in the numerics we can take finite step-sizes in the righthand side. Now in order to avoid a redundant computation we can also avoid taking the de-tour via HARDI-orientation scores and actually work with the left-invariant vector fields on the HARDI data itself. To this end we need the
consistent right-action $\; \gothic{R} \; $ of $SE(3)$ acting on the space of HARDI images $\mathbb{L}_{2}(\R^3 \rtimes S^{2})$.
To construct this consistent right-action we formally define $S:\mathbb{L}_{2}(\R^{3} \rtimes S^{2}) \to H$, where $H$ denotes the space of HARDI-orientation scores,
that equals the space of quadratic integrable functions on the group $SE(3)$ which are $\alpha$ right-\emph{invariant}, i.e. satisfying (\ref{alphainv}) by
\[
(\mathcal{S} U)(\ul{x},R)=\tilde{U}(\ul{x},R)=U(\ul{x}, R \ul{e}_{z}).
\]
This mapping is injective and its left-inverse is given by $(S^{-1}\tilde{U})(\ul{x},\ul{n})=\tilde{U}(\ul{x},R_{\ul{n}})$, where again
$R_{\ul{n}} \in SO(3)$ is some rotation such that $R_{\ul{n}}\ul{e}_{z}=\ul{n}$. Now the consistent right-action $\gothic{R}: SE(3) \to B(\mathbb{L}_{2}(\R^3 \rtimes S^{2}))$, where $B(\mathbb{L}_{2}(\R^3 \rtimes S^{2}))$ stands for all bounded linear operators on the space of HARDI images, is (almost everywhere) given by
\[
(\gothic{R}_{(\ul{x},R)}U)(\ul{y},\ul{n})= (S^{-1} \circ \mathcal{R}_{(\ul{x},R)}\circ S \, U)(\ul{y},\ul{n})= U(R_{\ul{n}}\ul{x}+\ul{y}, R_{\ul{n}}R \ul{e}_{z}).
\]
This yields the left-invariant vector fields (directly) on sufficiently smooth HARDI images:
\[
\mathcal{A}_{i}U(\ul{y},\ul{n})=({\rm d}\mathcal{R}(A_{i})U)(\ul{y},\ul{n})=
\lim \limits_{h \downarrow 0}\frac{(\gothic{R}_{e^{ h\, A_{i}}}U)(\ul{y},\ul{n})-U(\ul{y},\ul{n})}{h}=
\lim \limits_{h \downarrow 0}\frac{(\gothic{R}_{e^{ h\, A_{i}}}U)(\ul{y},\ul{n})-(\gothic{R}_{e^{ -h\, A_{i}}}U)(\ul{y},\ul{n})}{2h}.
\]
Now in our algorithms we take finite step-sizes and elementary computations (using the exponent (\ref{exp})) yield the following simple expressions for the discrete left-invariant vector fields, for respectively
central,
{\small
\begin{equation} \label{LI}
\begin{array}{l}
\mathcal{A}_{1}U(\ul{y},\ul{n}) \approx \mathcal{A}_{1}^{c}U(\ul{y},\ul{n}):=\frac{U(\ul{y}+ h \, R_{\ul{n}} \ul{e}_{x} \, , \, \ul{n})-U(\ul{y}- h \, R_{\ul{n}} \ul{e}_{x} \, ,\, \ul{n})}{2h}\ , \\
\mathcal{A}_{2}U(\ul{y},\ul{n}) \approx \mathcal{A}_{2}^{c}U(\ul{y},\ul{n}):=\frac{U(\ul{y}+ h \, R_{\ul{n}} \ul{e}_{y} \, , \, \ul{n})-U(\ul{y}- h \, R_{\ul{n}} \ul{e}_{y} \, ,\, \ul{n})}{2h}\ , \\
\mathcal{A}_{3}U(\ul{y},\ul{n}) \approx \mathcal{A}_{3}^{c}U(\ul{y},\ul{n}):=\frac{U(\ul{y}+ h \, R_{\ul{n}} \ul{e}_{z} \, , \, \ul{n})-U(\ul{y}- h \, R_{\ul{n}} \ul{e}_{z} \, ,\, \ul{n})}{2h} \ ,
\end{array}
\begin{array}{l}
\mathcal{A}_{4}U(\ul{y},\ul{n}) \approx \mathcal{A}_{4}^{c}U(\ul{y},\ul{n}):=\frac{U(\ul{y}\, , \, R_{\ul{n}} \; R_{\ul{e}_{x},h} \; \ul{e}_{z})-U(\ul{y}\, , \, R_{\ul{n}} \; R_{\ul{e}_{x},-h}\; \ul{e}_{z})}{2h}\ , \\
\mathcal{A}_{5}U(\ul{y},\ul{n}) \approx \mathcal{A}_{5}^{c}U(\ul{y},\ul{n}):=\frac{U(\ul{y}\, , \, R_{\ul{n}} \; R_{\ul{e}_{y},h} \; \ul{e}_{z})-U(\ul{y}\, ,\,  R_{\ul{n}} \; R_{\ul{e}_{y},-h}\;  \ul{e}_{z}) }{2h} \ .\\
\end{array}
\end{equation}
}
forward,
{\small
\begin{equation} \label{LIforward}
\begin{array}{l}
\mathcal{A}_{1}U(\ul{y},\ul{n}) \approx \mathcal{A}_{1}^{f}U(\ul{y},\ul{n}):=\frac{U(\ul{y}+ h \, R_{\ul{n}} \ul{e}_{x} \, , \, \ul{n})-U(\ul{y}  \, ,\, \ul{n})}{h}\ , \\
\mathcal{A}_{2}U(\ul{y},\ul{n}) \approx \mathcal{A}_{2}^{f}U(\ul{y},\ul{n}):=\frac{U(\ul{y}+ h \, R_{\ul{n}} \ul{e}_{y} \, , \, \ul{n})-U(\ul{y}  \, ,\, \ul{n})}{h}\ , \\
\mathcal{A}_{3}U(\ul{y},\ul{n}) \approx \mathcal{A}_{3}^{f}U(\ul{y},\ul{n}):=\frac{U(\ul{y}+ h \, R_{\ul{n}} \ul{e}_{z} \, , \, \ul{n})-U(\ul{y}  \, ,\, \ul{n})}{h} \ ,
\end{array}
\begin{array}{l}
\mathcal{A}_{4}U(\ul{y},\ul{n}) \approx \mathcal{A}_{4}^{f}U(\ul{y},\ul{n}):=\frac{U(\ul{y}\, , \, R_{\ul{n}} \; R_{\ul{e}_{x},h} \; \ul{e}_{z})-U(\ul{y}\, , \,\ul{n})}{h}\ , \\
\mathcal{A}_{5}U(\ul{y},\ul{n}) \approx \mathcal{A}_{5}^{f}U(\ul{y},\ul{n}):=\frac{U(\ul{y}\, , \, R_{\ul{n}} \; R_{\ul{e}_{y},h} \; \ul{e}_{z})-U(\ul{y}\, ,\,  \ul{n})}{h} \ .\\
\end{array}
\end{equation}
}
and backward,
{\small
\begin{equation} \label{LIbackward}
\begin{array}{l}
\mathcal{A}_{1}U(\ul{y},\ul{n}) \approx \mathcal{A}_{1}^{b}U(\ul{y},\ul{n}):=\frac{U(\ul{y} \, , \, \ul{n})-U(\ul{y}- h \, R_{\ul{n}} \ul{e}_{x} \, ,\, \ul{n})}{h}\ , \\
\mathcal{A}_{2}U(\ul{y},\ul{n}) \approx \mathcal{A}_{2}^{b}U(\ul{y},\ul{n}):=\frac{U(\ul{y} \, , \, \ul{n})-U(\ul{y}- h \, R_{\ul{n}} \ul{e}_{y} \, ,\, \ul{n})}{h}\ , \\
\mathcal{A}_{3}U(\ul{y},\ul{n}) \approx \mathcal{A}_{3}^{b}U(\ul{y},\ul{n}):=\frac{U(\ul{y} \, , \, \ul{n})-U(\ul{y}- h \, R_{\ul{n}} \ul{e}_{z} \, ,\, \ul{n})}{h} \ ,
\end{array}
\begin{array}{l}
\mathcal{A}_{4}U(\ul{y},\ul{n}) \approx \mathcal{A}_{4}^{b}U(\ul{y},\ul{n}):=\frac{U(\ul{y}\, , \, \ul{n})-U(\ul{y}\, , \, R_{\ul{n}} \; R_{\ul{e}_{x},-h}\; \ul{e}_{z})}{h}\ , \\
\mathcal{A}_{5}U(\ul{y},\ul{n}) \approx \mathcal{A}_{5}^{b}U(\ul{y},\ul{n}):=\frac{U(\ul{y}\, , \, \ul{n})-U(\ul{y}\, ,\,  R_{\ul{n}} \; R_{\ul{e}_{y},-h}\;  \ul{e}_{z}) }{h} \ .\\
\end{array}
\end{equation}
}
left-invariant finite differences. The left-invariant vector fields $\{\mathcal{A}_{1}, \mathcal{A}_{2},\mathcal{A}_{4},\mathcal{A}_{5}\}$ clearly depend on the choice of
$R_{\ul{n}} \in SO(3)$ which maps $R_{\ul{n}}\ul{e}_{z}=\ul{n}$. Now functions in the space $H$ are $\alpha$-right \emph{invariant}, so thereby we may assume that $R$ can be written as $R=R_{\ul{e}_{x}, \gamma} R_{\ul{e}_{y}, \beta}$, now if we choose
$R_{\ul{n}}$ again such that $R_{\ul{n}(\beta,\gamma)}=R_{\ul{e}_{x}, \gamma} R_{\ul{e}_{y}, \beta} R_{\ul{e}_{y}, \alpha=\alpha_{0}=0}$ then we take consistent sections in $SO(3)/SO(2)$ and we get full invertibility $S^{-1} \circ S =S \circ S^{-1}=\mathcal{I}$. 

Our evolution schemes, however, the choice of representant $R_{\ul{n}}$ is irrelevant, because they are well-defined
on the quotient $\R^{3}\rtimes S^{2}=SE(3)/(\{0\}\times SO(2))$.

In the computation of (\ref{LI}) one would have liked to work with discrete subgroups of $SO(3)$ acting on $S^{2}$ in order to avoid interpolations, but unfortunately the platonic solid with the
largest amount of vertices (only $20$) is the dodecahedron and the platonic solid with the largest amount of faces (again only $20$)
is the icosahedron. Nevertheless, we would like to sample the
$2$-sphere such that the distance between sampling points should be as equal as possible and simultaneously the area around each sample point should be as equal as possible. Therefore we follow the common approach by regular triangulations (i.e.
each triangle is regularly divided into $(o+1)^2$ triangles) of the
icosahedron, followed by a projection on the sphere. This leads to
$N_{o}= 2+ 10(o+1)^2$ vertices. We typically considered $o=1,2,3$, for further motivation regarding uniform spherical sampling,
see \cite[ch.7.8.1]{FrankenThesis}.

For the required interpolations to compute (\ref{LI}) within our spherical sampling there are two simple options. Either one uses a triangular interpolation of using the three closest sampling points, or one uses a discrete spherical harmonic interpolation. The disadvantage of the first
and simplest approach is that it introduces additional blurring, whereas the second approach can lead to
overshoots and undershoots. In both cases it is for computational efficiency
advisable to pre-compute the interpolation matrix, cf.~\cite[ch:2.1]{Creusen}, \cite[p.193]{FrankenPhDThesis}. See Appendix \ref{app:G}.

\subsection{Finite difference schemes for diffusion and pseudo-linear scale spaces on $\R^{3}\rtimes S^2$ \label{ch:FDDiff}}

The linear diffusion system on $\R^3 \rtimes S^{2}$ can be rewritten as
\begin{equation} \label{goal}
\left\{
\begin{array}{l}
\partial_{t} W(\ul{y},\ul{n},t) = \left(D^{11}((\mathcal{A}_{1})^{2} + (\mathcal{A}_{2})^2) + D^{33}(\mathcal{A}_{3})^2 + D^{44} \Delta_{S^{2}}\right) W(\ul{y},\ul{n},t) \\
W(\ul{y},\ul{n},0)=U(\ul{y},\ul{n})\ ,
\end{array}
\right.
\end{equation}
This system is the Fokker-Planck equation of horizontal Brownian motion on $\R^3 \rtimes S^{2}$ if $D^{11}=0$.
Spatially, we take second order centered finite differences for $(\mathcal{A}_{1})^2$, $(\mathcal{A}_{2})^{2}$
and $(\mathcal{A}_{3})^2$, i.e. we applied the discrete operators in the righthand side of (\ref{LI}) twice (where we replaced $2h \mapsto h$ to ensure direct-neighbors interaction), so we have
\begin{equation} \label{A32}
((\mathcal{A}_{3})^2W)(\ul{y},\ul{n},t) \approx \frac{W(\ul{y}+ h R_{\ul{n}} \ul{e}_{z},\ul{n},t)-2\, W(\ul{y},\ul{n},t)+ W(\ul{y}- h R_{\ul{n}} \ul{e}_{z},\ul{n},t)}{h^2}\ ,
\end{equation}
For each $q \in \mathbb{N} \cup \{0\}$, we define the vector $\ul{w}^{q}=(w^{q}_{\ul{y},k})_{k=1 \ldots N_o, \ul{y} \in \mathcal{I}} \in \R^{N^3 N_{o}}$, where $k$ enumerates the number of samples on the sphere, where $\ul{y}$ enumerates the samples
on the discrete spatial grid $\mathcal{I}=\{1,\ldots,N\}\times \{1,\ldots,N\} \times \{1,\ldots,N\}$ and where $q$ enumerates the discrete time frames.
Rewrite (\ref{A32}) and (\ref{goal}) in vector form
using Euler-forward first order approximation in time:
 \[
\ul{w}^{q+1}= (I -\Delta t (J_{S_2}+J_{\R^3}) ) \ul{w}^{q}\ ,
\]
where $J_{S^{2}} \in \R^{N^{3}N_{o}\times N^{3}N_{o}}$ denotes the angular increments block-matrix and where
$J_{\R^{3}} \in \in \R^{N^{3}N_{o}\times N^{3}N_{o}}$ denotes the spatial increments block-matrix.

\subsubsection{Angular increments block-matrix}

The angular increments block-matrix equals in the basic (more practical) approach, cf.\cite{Creusen}, of (tri-)linear spherical interpolation
\begin{equation} \label{simple}
J_{S_{2}}= h_{a}^{-2}D^{44}(B_{4}+B_{5}),
\end{equation}
with angular stepsize $h_{a}>0$ and with $B_{4}$ and $B_{5}$ are matrices with $2$ on the diagonal
and where for each column the of off-diagonal elements is also equal to 2. See Appendix \ref{app:G}, Eq.~(\ref{Js}), for details.

In the more complicated discrete spherical harmonic transform interpolation approach, cf. \cite[Ch:7]{DuitsIJCV2010}
the angular increments block-matrix is given by
\begin{equation} \label{difficult}
J_{S_2}=  I \otimes D^{44} \Lambda^{-1} M^{T} Q \overline{M} \Lambda \in \R^{(N^{3} N_o) \times (N^{3} N_{o})},
\textrm{ with }Q=\textrm{diag}_{j=1,\ldots,N_{SH}}\{l(j)(l(j)+1)e^{-t_{reg}l(j)(l(j)+1)}\},
\end{equation}
with regularization parameter $t_{reg}>0$ and
where $N_{SH}$ represents the number of spherical harmonics and with
$n_{SH}\times N_{o}$-matrix
\[
M=[M^{j}_{k}]=[\frac{1}{\sqrt{C}} Y^{l(j)}_{m(j)}(\ul{n}_{k})], \textrm{ with }l(j)=\lfloor \sqrt{j-1} \rfloor \textrm{ and }m(j)=j-(l(j))^2-l(j) -1
\]
with $C=\sum_{j=1}^{n_{SH}}|Y_{m(j)}^{l(j)}(0,0)|^2$ and
where the diagonal matrix $\mbox{\boldmath$\Lambda$}=\textrm{diag}\{\delta_{S^2}(\ul{n}_{1}), \ldots, \delta_{S^2}(\ul{n}_{N_{o}}) \}$ \
contains discrete surface measures $\delta_{S^2}(\ul{n}_{k})$ (for spherical sampling by means of higher order tessellations of the icosahedron)
given by
\begin{equation}\label{sm}
\delta_{S^2}(\ul{n}_{k})= \frac{1}{6} \sum \limits_{i\neq k,j \neq k, i\neq j, i\sim j\sim k} A(\ul{n}_{i},\ul{n}_{j},\ul{n}_{k})\ ,
\end{equation}
where $i \sim j$ means that $\ul{n}_{i}$ and $\ul{n}_{j}$ are part of a locally smallest triangle in the tessellation and where the surface measure of the spherical projection of such a triangle is given by
\[
\begin{array}{l}
A(\ul{n}_i,\ul{n}_{j},\ul{n}_{k})= 4 \arctan(\sqrt{\tan(s_{ijk}/2)\tan((s_{ijk}-s_{ij})/2)
\tan((s_{ijk}-s_{ik})/2)\tan((s_{ijk}-s_{jk})/2)})\ , \\
\textrm{with } s_{ijk}= \frac{1}{2}(s_{ij}+s_{ik}+s_{jk}) \textrm{ and } s_{ij}=\arccos(\ul{n}_{i}\cdot \ul{n}_{j}). 
\end{array}
\]

\subsubsection{Spatial increments matrix}

In both approaches the spatial increments block matrix equals $J_{\R^{3}}$ given by  Eq.~(\ref{Js}) in Appendix \ref{app:G}.


\subsubsection{Stability bounds on the step-size}

We have guaranteed stability iff
\[
\|I-\Delta t (J_{S^2}+J_{\R^3})\|= \textrm{sup}_{\|\ul{w}\|=1} \|(I-\Delta t (J_{S^2}+J_{\R^3}))\ul{w}\| < 1,
\]
which is by $I-\Delta t (J_{S^2}+J_{\R^3})=(\nu I-\Delta t J_{\R^3})+((1-\nu)I-\Delta t J_{S^2})$, with $\nu \in (0,1)$, the case if
\begin{equation} \label{suff}
\|I-\frac{\Delta t}{\nu} J_{\R^3}\|<1 \textrm{ and }\|I-\frac{\Delta t}{1-\nu}\; J_{S^2}\|<1\ \Rightarrow \|(\nu I-\Delta t J_{\R^3})+((1-\nu)I-\Delta t J_{S^2})\|< \nu+(1-\nu)=1.
\end{equation}
Sufficient (and sharp) conditions for the first inequality are obtained by means of the Gerschgorin Theorem \cite{Gershgorin}\; :
\begin{equation} \label{spatest}
\Delta t \leq \frac{\nu h^2}{4D^{11}+2D^{33}}.
\end{equation}

Application of the Gerschgorin Theorem \cite{Gershgorin} on the second equality gives after optimization over $\nu$, the following global guaranteed-stability bound
\begin{equation} \label{final1}
\Delta t \leq \frac{1}{\frac{4D^{11}+2 D^{33}}{h^2}+ \frac{4D^{44}}{h_{a}^2}}
\end{equation}
in the linear interpolation case and
\begin{equation} \label{final2}
\begin{array}{ll}
\Delta t \leq \frac{h^{2}}{4D^{11}+2D^{33}+D^{44} h^2 \frac{L(L+1)}{2\, e^{t_{reg}L(L+1)}}} &\textrm{ if }t_{reg} \cdot L(L+1) \leq 1 \ ,\\
\Delta t \leq \frac{h^{2}}{4D^{11}+2D^{33}+D^{44} h^2 \frac{1}{2 e\, t_{reg}}} &\textrm{ if }t_{reg} \cdot L(L+1) >1
\end{array}
\end{equation}
in the spherical harmonic interpolation case. For details, see respectively \cite{Creusen} and \cite{DuitsIJCV2010}.

The pseudo-linear scale spaces are implemented using the above together with Eq.~(\ref{conjGT}), and Eq.~(\ref{GT}), Eq.~(\ref{GTinv}).

\subsection{Finite difference schemes for Hamilton-Jacobi equations on $\R^{3}\rtimes S^2$ \label{ch:FDerode}}

Similar to the previous section we use a Euler-forward first order time integration scheme,
but now we use a so-called upwind-scheme, where the sign of the central difference determines whether a forward
or central difference is taken.
{\small
\[
\begin{array}{l}
W(\ul{y},\ul{n},t+\Delta t)= W(\ul{y},\ul{n})+ \frac{\Delta t }{2\eta}
\left(\; \; D^{11}\left(\left(a^{-,1}(\ul{y},\ul{n})\frac{U(\ul{y} \, , \, \ul{n})-U(\ul{y}- h \, R_{\ul{n}} \ul{e}_{x} \, ,\, \ul{n})}{h}+ a^{+,1}(\ul{y},\ul{n})\frac{U(\ul{y}+ h \, R_{\ul{n}} \ul{e}_{x} \, , \, \ul{n})-U(\ul{y}  \, ,\, \ul{n})}{h}\right)^{2}
\right. \right.
\\
\left. +
\left(a^{-,2}(\ul{y},\ul{n})\frac{U(\ul{y} \, , \, \ul{n})-U(\ul{y}- h \, R_{\ul{n}} \ul{e}_{y} \, ,\, \ul{n})}{h}+
a^{+,2}(\ul{y},\ul{n})\frac{U(\ul{y}+ h \, R_{\ul{n}} \ul{e}_{y} \, , \, \ul{n})-U(\ul{y}\, ,\, \ul{n})}{h}\right)^{2}\right)+
\\
+D^{44} \left( \left( a^{-,4}(\ul{y},\ul{n})
\frac{U(\ul{y}\, , \, \ul{n})-U(\ul{y}\, , \, R_{\ul{n}} \; R_{\ul{e}_{x},-h}\; \ul{e}_{z})}{h}+a^{+,4}(\ul{y},\ul{n})
\frac{U(\ul{y}\, , \, R_{\ul{n}} \; R_{\ul{e}_{x},h})-U(\ul{y}\, , \, \ul{n})}{h}\right)^2\ \right.  \\
+\left.\left.\left( a^{-,5}(\ul{y},\ul{n})
\frac{U(\ul{y}\, , \, \ul{n})-U(\ul{y}\, ,\,  R_{\ul{n}} \; R_{\ul{e}_{y},-h}\;  \ul{e}_{z}) }{h} +
a^{+,5}(\ul{y},\ul{n})
\frac{U(\ul{y}\, , \, R_{\ul{n}} \; R_{\ul{e}_{x},h})-U(\ul{y}\, , \, \ul{n})}{h}
\right)^2 \right)
\; \;\right)^{\eta},
\end{array}
\]
}
where the functions $a^{-,k}$ are given by
{\small
\[
\begin{array}{ll}
a^{-,1}(\ul{y},\ul{n})= \max\left\{0, \frac{U(\ul{y}+ h \, R_{\ul{n}} \ul{e}_{x} \, , \, \ul{n})-U(\ul{y}- h \, R_{\ul{n}} \ul{e}_{x} \, ,\, \ul{n})}{2h}\right\} ,&
a^{+,1}(\ul{y},\ul{n})= \min\left\{0, \frac{U(\ul{y}+ h \, R_{\ul{n}} \ul{e}_{x} \, , \, \ul{n})-U(\ul{y}- h \, R_{\ul{n}} \ul{e}_{x} \, ,\, \ul{n})}{2h}\right\}, \\
a^{-,2}(\ul{y},\ul{n})= \max\left\{0, \frac{U(\ul{y}+ h \, R_{\ul{n}} \ul{e}_{y} \, , \, \ul{n})-U(\ul{y}- h \, R_{\ul{n}} \ul{e}_{y} \, ,\, \ul{n})}{2h}\right\} ,&
a^{+,2}(\ul{y},\ul{n})= \min\left\{0, \frac{U(\ul{y}+ h \, R_{\ul{n}} \ul{e}_{y} \, , \, \ul{n})-U(\ul{y}- h \, R_{\ul{n}} \ul{e}_{y} \, ,\, \ul{n})}{2h}\right\},  \\
a^{+,4}(\ul{y},\ul{n})= \min\left\{0,\frac{U(\ul{y}\, , \, R_{\ul{n}} \; R_{\ul{e}_{x},h} \; \ul{e}_{z})-U(\ul{y}\, ,\,  R_{\ul{n}} \; R_{\ul{e}_{x},-h}\;  \ul{e}_{z}) }{2h}\right\},&
a^{-,4}(\ul{y},\ul{n})= \max\left\{0,\frac{U(\ul{y}\, , \, R_{\ul{n}} \; R_{\ul{e}_{x},h} \; \ul{e}_{z})-U(\ul{y}\, ,\,  R_{\ul{n}} \; R_{\ul{e}_{x},-h}\;  \ul{e}_{z})}{2h}\right\}. \\
a^{+,5}(\ul{y},\ul{n})= \min\left\{0,\frac{U(\ul{y}\, , \, R_{\ul{n}} \; R_{\ul{e}_{y},h} \; \ul{e}_{z})-U(\ul{y}\, ,\,  R_{\ul{n}} \; R_{\ul{e}_{y},-h}\;  \ul{e}_{z}) }{2h}\right\}. &
a^{-,5}(\ul{y},\ul{n})= \max\left\{0,\frac{U(\ul{y}\, , \, R_{\ul{n}} \; R_{\ul{e}_{y},h} \; \ul{e}_{z})-U(\ul{y}\, ,\,  R_{\ul{n}} \; R_{\ul{e}_{y},-h}\;  \ul{e}_{z})}{2h}\right\}.
\end{array}
\]
}

\subsection{Convolution implementations}

The algorithm for solving Hamilton-Jacobi Equations Eq.~\ref{Hamdi1}, (\ref{Hamdi2})
by $\R^{3}\rtimes S^{2}$ dilation/erosion with the corresponding analytic Green's function (\ref{approx})
boils down to the same algorithm as $\R^{3}\rtimes S^2$-convolutions with analytic Green's functions for diffusion.
Basically, the algebra $(+,\cdot)$ is to be replaced by respectively the $(\max,+)$ and $(\min,+)$-algebra.
So the results on fast efficient computation, using lookup-tables and/or parallelization cf.~\cite{Rodrigues},
apply also to the morphological convolutions. The difference though is that hypo-elliptic diffusion kernels
are mainly supported near the $z$-axis
(thereby within the $\R^{3}\rtimes S^2$ convolution after translation and rotation near the $\ul{A}_{3}$-axis), whereas
the erosion kernels are mainly supported near the $xy$-plane, but the efficiency
principles are easily carried over to the morphological convolutions.

Next we address two minor issues that arise when implementing the convolution algorithms \cite[ch:8.2]{DuitsIJCV2010} and \cite{Rodrigues}:
\begin{enumerate}
\item In the implementation with analytical kernels expressed in the second chart, such as (\ref{heur2}) one needs to extract the second chart Euler angles $(\tilde{\beta},\tilde{\gamma})$ from each normal $\ul{n}$
in the spherical sampling, i.e. one must solve $\ul{n}(\tilde{\beta},\tilde{\gamma})=\ul{n}$. The solutions are
\begin{equation} \label{anglesfromn}
(\tilde{\beta},\tilde{\gamma})=(\textrm{sign}(n_{1})\arccos (\widetilde{\textrm{sign}}(n_{2} \sqrt{n_{2}^2+n_{3}^2})), -\widetilde{\textrm{sign}}(n_{3}) \arcsin \left( \frac{n_{2}}{\sqrt{n_{2}^2 +n_{3}^{2}}}\right) )
\end{equation}
with $\textrm{sign}(x)= 1$ if $x >0$ and zero else, and with $\widetilde{\textrm{sign}}(x)= 1$ if $x \geq 0$ and zero else.
\item $\R^{3}\rtimes S^{2}$-convolution requires computation of $R_{\ul{n}}^{T}\ul{v}$ , with $\ul{v}=\ul{y}-\ul{y}'$, Theorem \ref{th:convolve}.
For the sake of computation speed this can be done without goniometric formulas:
\[
R_{\ul{n}}^{T}(\ul{v})=\frac{1}{(n^{1})^{2}+(n^{2})^{2}}
\left(
\begin{array}{c}
(n^{2})^2 v^{1} + n^{1}n^{2}(n^{3}-1)v^{2} + n^{1}(n^{1}n^{3}v^{1}- \sqrt{(n^{1})^2+(n^{2})^2} \sqrt{1-(n^{3})^2} v^{3}) \\
(n^{1})^2 v^{2} + n^{1}n^{2}(n^{3}-1)v^{1} + n^{2}(n^{2}n^{3}v^{2}- \sqrt{(n^{1})^2+(n^{2})^2} \sqrt{1-(n^{3})^2} v^{3}) \\
\sqrt{1-(n^{3})^2} \sqrt{(n^{1})^2+(n^{2})^2} (n^{1}v^{1}+n^{2}v^{2})+ n^{3}v^{3}((n^{1})^2+(n^{2})^2)
\end{array}
\right)
\]
$\textrm{if } (n^{1},n^{2}) \neq (0,0)$,
where $R_{\ul{n}}$ is indeed \emph{a} rotation given by
\[
R_{\ul{n}}^{T}(\ul{v})= \frac{1}{(n^{1})^{2}+(n^{2})^{2}}
\left(
\begin{array}{c}
(n^{2})^2 v^{1} + n^{1}n^{2}(n^{3}-1)v^{2} + n^{1}(n^{1}n^{3}v^{1}+ \sqrt{(n^{1})^2+(n^{2})^2} \sqrt{1-(n^{3})^2} v^{3}) \\
(n^{1})^2 v^{2} + n^{1}n^{2}(n^{3}-1)v^{1} + n^{2}(n^{2}n^{3}v^{2}+ \sqrt{(n^{1})^2+(n^{2})^2} \sqrt{1-(n^{3})^2} v^{3}) \\
-\sqrt{1-(n^{3})^2} \sqrt{(n^{1})^2+(n^{2})^2} (n^{1}v^{1}+n^{2}v^{2})+ n^{3}v^{3}((n^{1})^2+(n^{2})^2)
\end{array}
\right)
\]
$\textrm{if } (n^{1},n^{2}) \neq (0,0)$
that maps $(0,0,1)^{T}$ onto $(n^1, n^{2}, n^{3}) \in S^{2}$. If $(n^{1},n^{2})=(0,0)$ then we set $R_{\ul{n}=(0,0,n^{3})}= \textrm{sign}(n^{3})\; I$.
\end{enumerate}

\section{Adaptive, Left-Invariant Diffusions on HARDI images \label{ch:LDOS}}

\subsection{Scalar Valued Adaptive Conductivity \label{ch:PM}}

In order to avoid mutual influence of anisotropic regions (areas with fibers) and isotropic regions (ventricles) one can replace the constant diffusivity/conductivity $D^{33}$ by
\begin{equation} \label{replace}
D^{33}\mathcal{A}_{3}\mathcal{A}_{3} \mapsto  D^{33}\mathcal{A}_{3} \circ
e^{-\frac{|\mathcal{A}_{3}W(\cdot,t)|^2}{K^2}} \circ \mathcal{A}_{3}\ ,
\end{equation}
in the generator of the left-invariant diffusion (\ref{generalPDE}) (where $\ul{a}=\ul{0}$ and $D=\textrm{diag}\{0,0,D^{33},0,0\}$).
This now yields the following hypoelliptic diffusion system
\begin{equation} \label{adaptive}
\begin{array}{ll}
\frac{\partial W}{\partial t}(\ul{y},\ul{n},t) &
=D^{33}\left.\mathcal{A}_{3}\right|_{\ul{y},\ul{n}}(e^{-\frac{|\left.\mathcal{A}_{3}\right|_{\ul{y},\ul{n}}W(\ul{y},\ul{n},t)|^2}{K^2}} \mathcal{A}_{3}W)(\ul{y},\ul{n},t), \\
W(\ul{y},\ul{n},0) &=U(\ul{y},\ul{n}),
\end{array}
\end{equation}
Here could also choose to adapt the diffusivity by the original DW-MRI data $U:\R^{3}\rtimes S^{2} \to \R^{+}$ at time $0$, so that the diffusion equation itself is linear, whereas the mapping
$U \mapsto \Phi_{t}(U):=W(\cdot,t)$ in Eq.~(\ref{adaptive}) included is well-posed and
nonlinear. In the experiments, however we extended the standard approach by Perona\and Malik \cite{Perona} to $\R^{3}\rtimes S^2$ and we used (\ref{adaptive}) where both the diffusion equation and
 the  mapping $U \mapsto W(\cdot,t)$ are nonlinear.
The idea is simple: the replacement sets a soft-threshold on the diffusion in $\mathcal{A}_{3}$-direction, at fiber locations one expects $|\mathcal{A}_{3}U(\ul{y},\ul{n})|$ to be small,
whereas in the transition areas between ventricles and white matter, where one needs to block the fiber propagation by hypo-elliptic diffusion, one expects a large $|\mathcal{A}_{3}U(\ul{y},\ul{n})|$.
For further details see the second author's master thesis \cite{Creusen}. Regarding discretization of (\ref{adaptive}) in the finite difference schemes of Subsection \ref{ch:FDDiff} we propose
{\small
\[
(\mathcal{A}_{3}(D^{33}\,
e^{-\frac{|\mathcal{A}_{3}W(\cdot,t)|^2}{K^2}} \,\mathcal{A}_{3}W(\cdot,t)))(\ul{y},\ul{n})
\approx
\frac{\tilde{D}^{33}(\ul{y}\!+\!\frac{1}{2}\ul{h},\ul{n},t) \cdot
\mathcal{A}_{3}^{c}W(\ul{y}\!+\!\frac{1}{2}\ul{h},\ul{n},t) - \tilde{D}^{33}(\ul{y}\!-\!\frac{1}{2}\ul{h},\ul{n},t) \cdot
\mathcal{A}_{3}^{c}W(\ul{y}\!-\!\frac{1}{2}\ul{h},\ul{n},t)}{h}
\]
}
where $h$ is the spatial stepsize and where $\ul{h}=h R_{\ul{n}}\ul{e}_{z}$
for notational convenience and where
\[
\tilde{D}_{33}(\ul{y},\ul{n},t)= D^{33} \, e^{-\left( \frac{\max\{|\mathcal{A}_{3}^{f}W(\ul{y},\ul{n},t)|,|\mathcal{A}_{3}^{b}W(\ul{y},\ul{n},t)|\}}{K} \right)^{2}}.
\]
and where we recall Eq.~(\ref{LI}), Eq.~(\ref{LIforward}) and Eq.~(\ref{LIbackward}) for central,
forward and backward finite differences.

\subsection{Tensor-valued adaptive conductivity}

For details see \cite{DuitsIJCV2010} and \cite[Ch:9]{DuitsFrankenCASA}, where we generalized the results in \cite{DuitsAMS2,Fran2009}.
Current implementations do not produce the expected results. Here the following problems arise. Firstly, the logarithm and exponential curves in Section \ref{ch:exp} are not well-defined on the quotient $\R^{3}\rtimes S^2$
and we need to take an appropriate section through the partition of cosets in $\R^{3}\rtimes S^2$.  Secondly, the. $6\times6$ Hessian in \cite[Eq.103,Ch:9]{DuitsFrankenCASA} has a nilspace in $\mathcal{A}_{6}$-direction.

\section{Experiments and Evaluation \label{ch:experiments}}

First we provide a chronological evaluation of the experiments depicted in various Figures so far.

In Figure~\ref{Fig:2} and
Figure~\ref{Fig:3}
one can find experiments of $\R^{3}\rtimes S^{2}$-convolution implementations of
hypo-elliptic diffusion (``contour enhancement'',cf.~Figure \ref{fig:Brownianmotion}),
Eq.~(\ref{enhancement}), using the analytic approximation of the Green's function
Eq.~(\ref{heur2}). Typically, these $\R^{3}\rtimes S^2$-convolution kernel operators
propagate fibers slightly better at crossings, as they do not suffer from numerical blur
(spherical interpolation) in finite difference
schemes explained in Section \ref{ch:FDDiff}. Compare Figure~\ref{Fig:2} to Figure~\ref{Fig:4}, where
the same type of fiber crossings of the corona radiata and corpus callosum occur. Here we note that all methods
in Figure~\ref{Fig:4} have been evaluated by finite difference schemes, based on linear spherical interpolation,
with sufficiently small time steps
$\Delta t$, recall Eq.~(\ref{final1}) using $162$ spherical samples
from a higher order tessellation of the icosahedron. In Figure~\ref{Fig:4} the visually
most appealing results are obtained by a single concatenation of diffusion and erosion. Then in Figure~\ref{Fig:5} we
have experimentally illustrated the short-coming of sharpening the data by means of squaring in
comparison to the erosions, Eq.~(\ref{Hamdi2}). Experiments of the latter evolution, implemented
by the upwind finite difference method described in Section \ref{ch:FDerode},
has been illustrated for various parameter settings in Figure~\ref{fig:erode} (also in Figure~\ref{Fig:4}).
The possibilities of a further modification, where erosion is locally adapted to the data,
can be found in Figure~\ref{Fig:adaptation}. Then in Figure~\ref{fig:Brownianmotion}
and Figure~\ref{fig:completionk} we have depicted the Green's functions of the (iterated) direction process
(hypo-elliptic convection and diffusion) that we applied in a fiber completion experiment in
Figure~\ref{fig:completion}.

\begin{figure}
\centerline{
\includegraphics[width=0.4 \hsize]{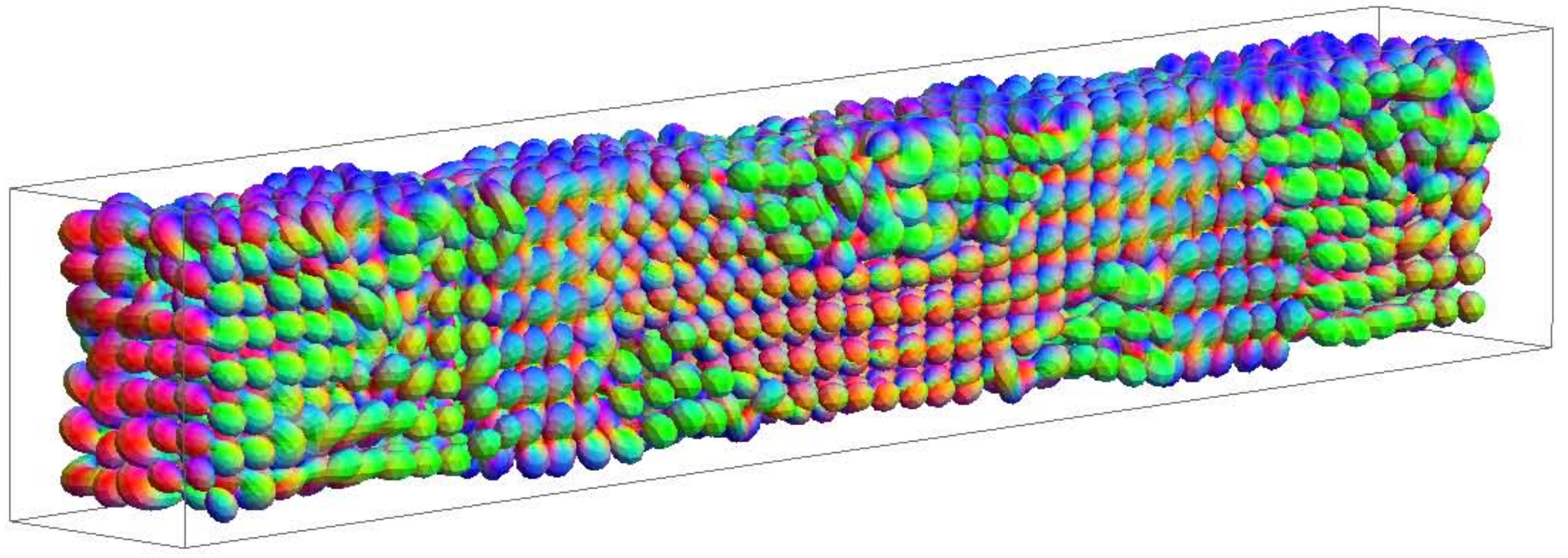}
\includegraphics[width=0.6 \hsize]{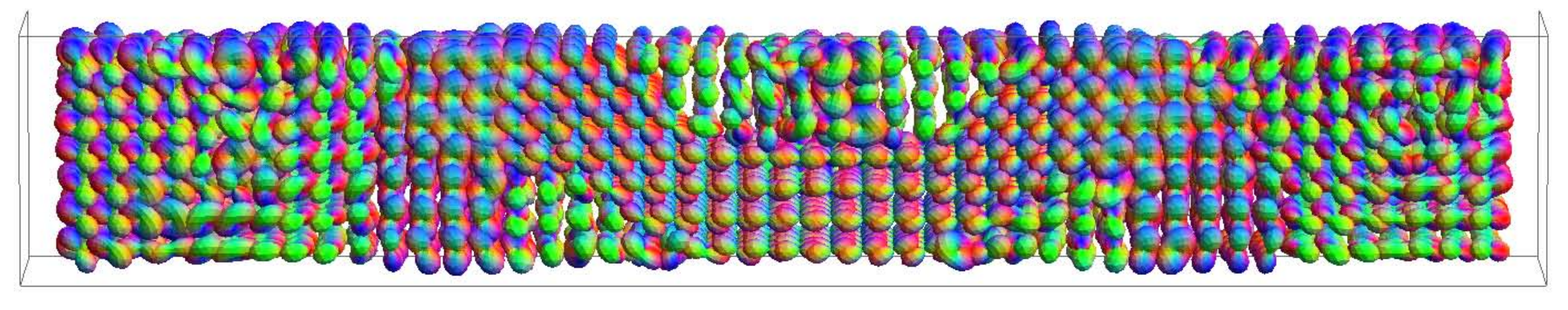}
}
\centerline{
\includegraphics[width=0.4 \hsize]{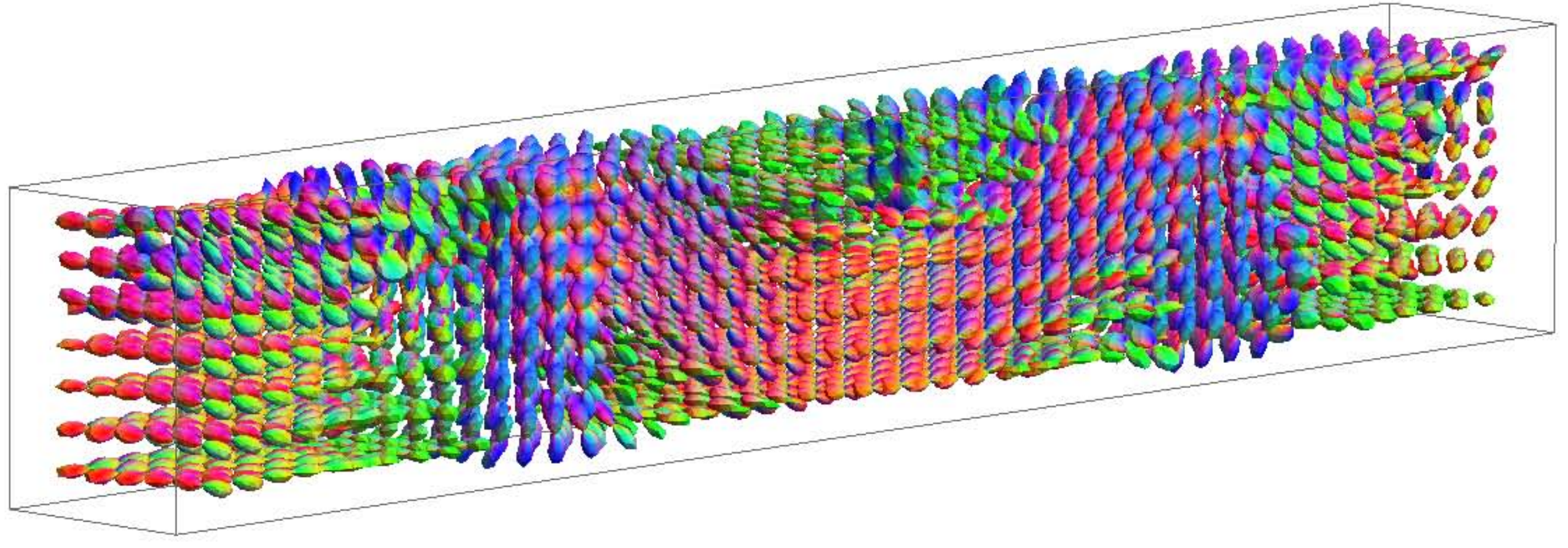}
\includegraphics[width=0.6 \hsize]{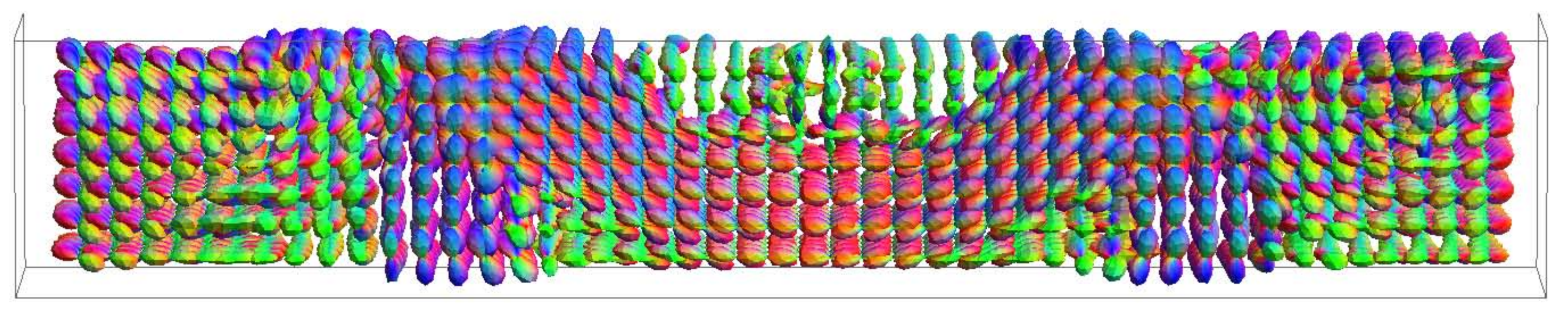}
}
\centerline{
\includegraphics[width=0.4 \hsize]{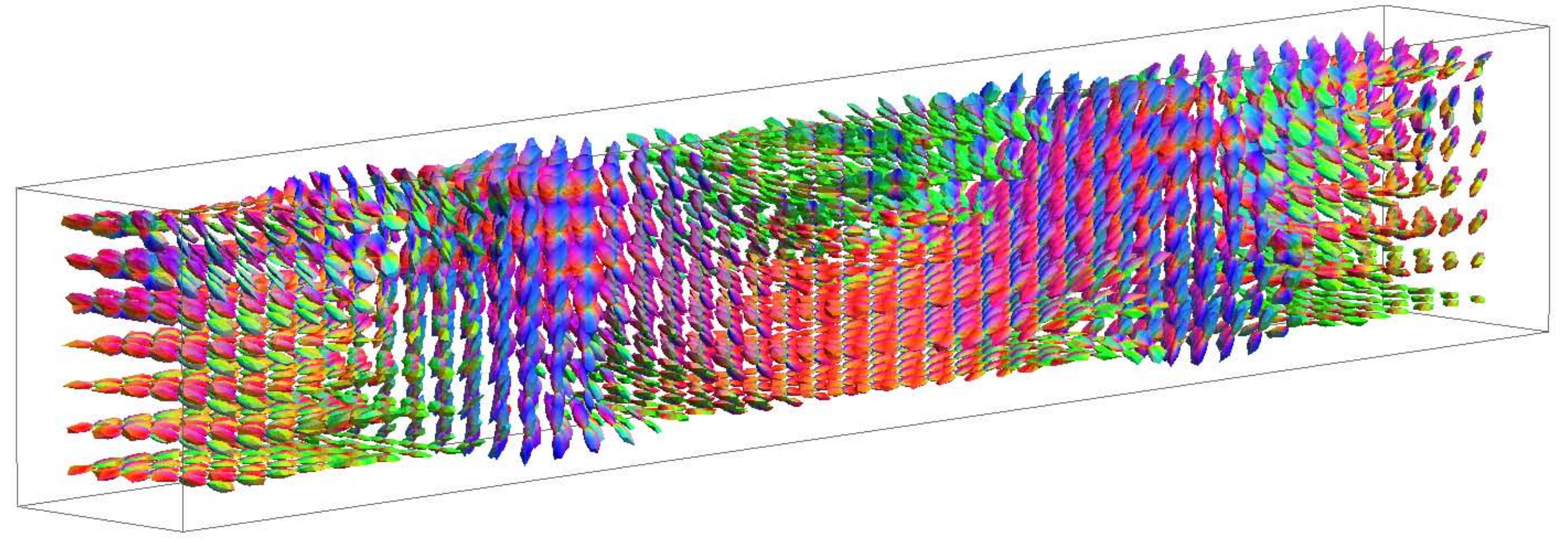}
\includegraphics[width=0.6 \hsize]{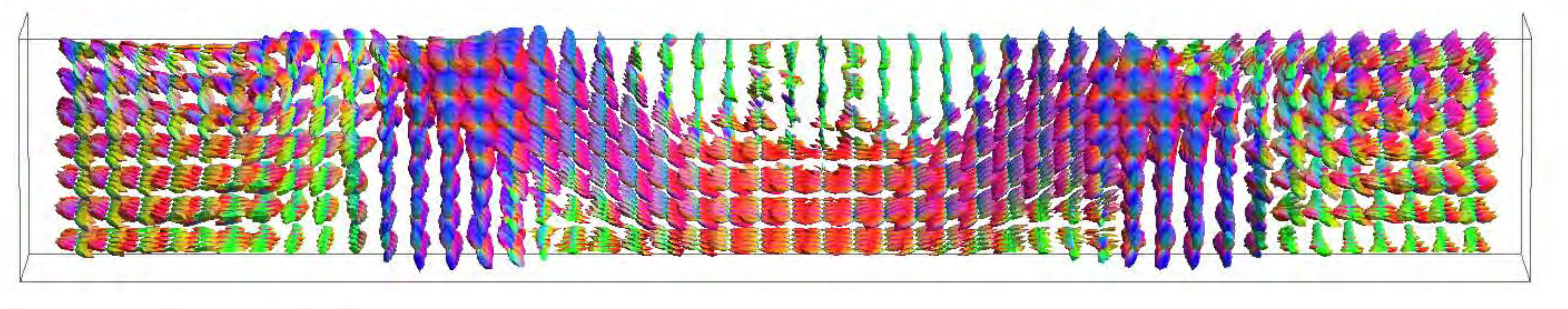}
}
\caption{Same experiments as in Figure~\ref{Fig:4}, but now visualized with several slices in $3D$.
First row: Input DTI-data set. Second row: Output squared linear diffusion on squared data-set. Third row: Output erosion applied to the diffused dataset depicted in the second row. For parameter settings, see Figure~\ref{Fig:4}.}
\label{Fig:10}
\end{figure}

The experiments in Figure~\ref{Fig:4} (implemented in \emph{Mathematica 7}) take place on volumetric DTI-data-sets $104 \times 104 \times 10$ where we applied (\ref{DTItoHARDI}) using
162 samples on the sphere (tessellation icosahedron). In Figure \ref{Fig:10} we illustrate the enhancement and eroding effect in $3D$, of the same experiments as in Figure~\ref{Fig:4}
where we took two different viewpoints on a larger 3D-part of the total data-set.

In the remaining subsections we will consider some further experiments on respectively diffusion, erosion and pseudo-linear scale spaces.

\subsection{Further experiments diffusions}

\begin{figure}
\centerline{
\includegraphics[width=0.9 \hsize]{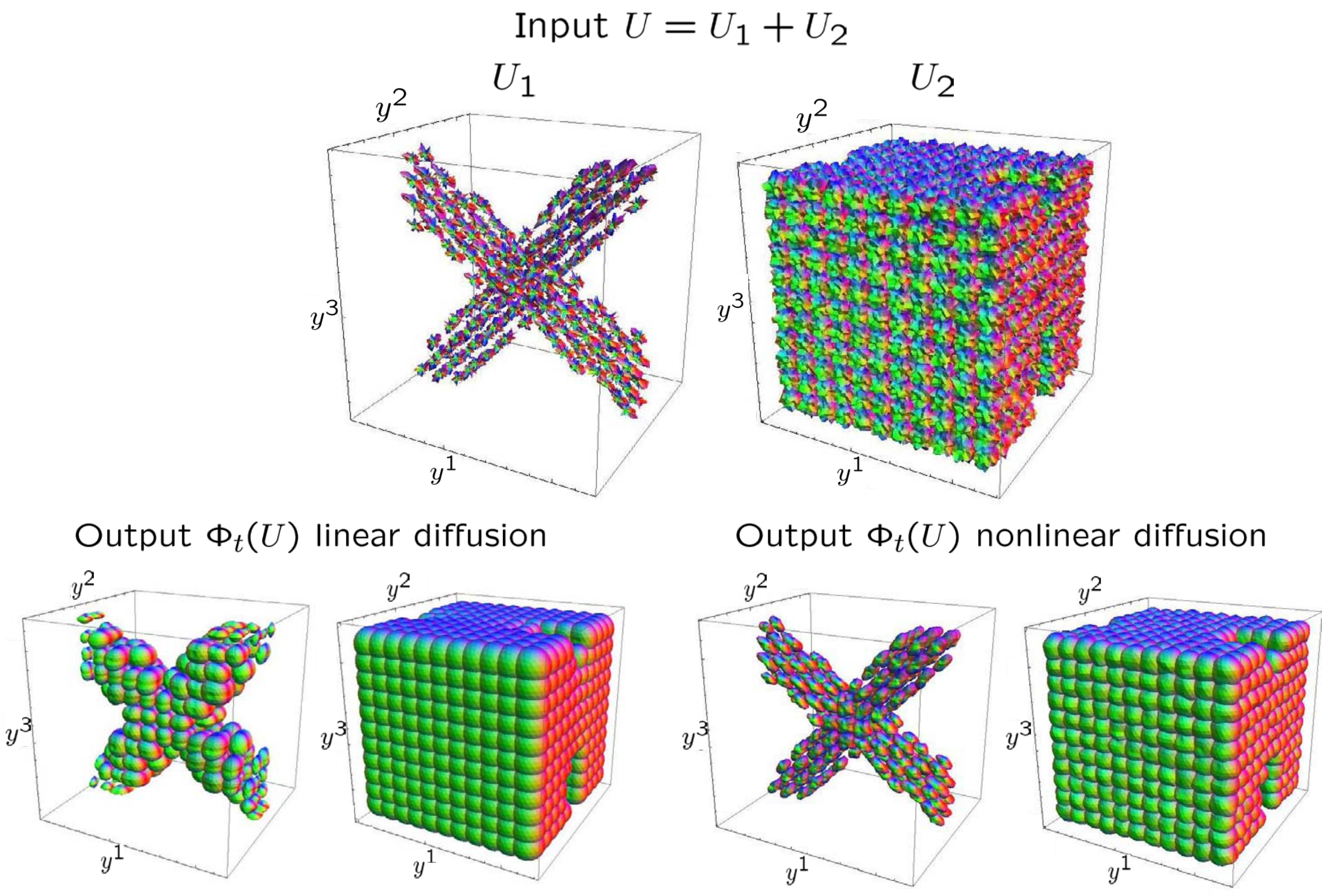}}
\caption{Adaptive diffusivity/conductivity based on the data, Section \ref{ch:PM}. Top row: Artificial
$15\times 15 \times 15 \times 162$-input data
that is a sum of a noisy fiber part and a noisy isotropic part. For the sake of visualization we depicted these parts
separately. Bottom row: Output of just linear diffusion (without greyvalue transformations, nor Eq.~(\ref{Veq})), $t=1$, $D^{33}=1$,
$D^{44}=0.04$.
Output nonlinear diffusion with $t=1$, $D^{33}=1$, $D^{44}=0.015$, $K=0.05$, $\Delta t=0.01$.}\label{Fig:adaptdiff}
\end{figure}

In Figure \ref{Fig:adaptdiff} we have applied both linear and nonlinear diffusion on an noisy artificial data set
which both contains an anisotropic part (modeling neural fibers) and containing a large isotropic part
(modeling ventricles, compare to Fig.~\ref{Fig:4}). When applying linear diffusion on the data-set the isotropic
part propagates and destroys the fiber structure in the anisotropic part, when applying nonlinear diffusion with
the same parameter settings, but now with $D^{33}= e^{-\frac{|\mathcal{A}_{3}U|^{2}}{K^{2}}}$ and suitable choice of
$K>0$ one turns of the propagation (diffusion-flow) through the boundary between isotropic and non-isotropic areas
where $(\ul{y},\ul{n} \mapsto |\mathcal{A}_{3}(\ul{y},\ul{n})|^{2}$ is relatively high. Consequently, both the fiber-part and the ventricle part are both
appropriately diffused without too much influence on each other, yielding both smooth/aligned fibers and
smooth isotropic glyphs.

\subsection{Further experiments erosions}

\begin{figure}
\centerline{
\includegraphics[width= 0.99 \hsize]{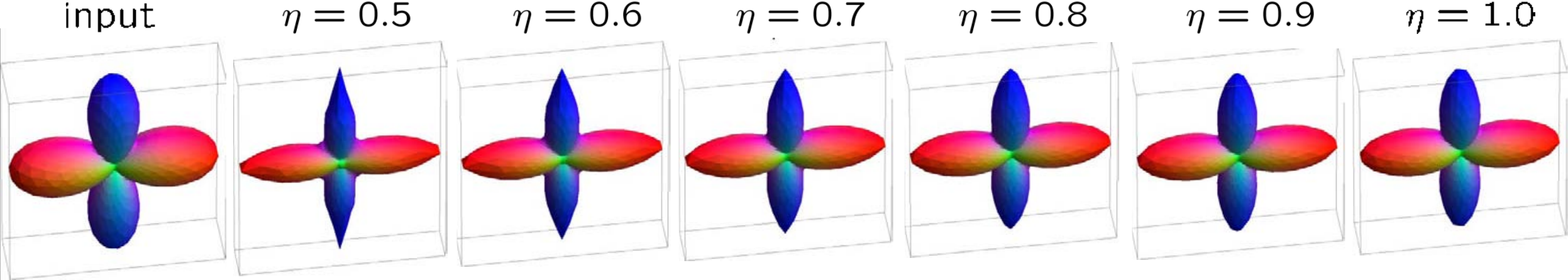}
}
\caption{The effect on $\eta \in [\frac{1}{2},1]$ on angular erosion Eq.~(\ref{Hamdi2}), $D^{44}\!=0.4$, $D^{11}\!=0$ and $t=0.4$, $\Delta t=0.02$. Left: original glyph, right eroded glyphs for $\eta=0.5, \ldots, 1.0$.}\label{Fig:eta}
\end{figure}

Erosions yield a geometrical alternative to more ad-hoc data-sharpening by squaring, see Fig:\ref{Fig:5} for a
basic illustration. We compared angular erosions for several values of $\eta \in [0,1]$, in Figure \ref{Fig:eta}, where we fixed $D^{44}=0.4$, $D^{11}=0$ and $t=0.4$ in Eq.~(\ref{Hamdi2}).
Typically, for $\eta$ close to a half the structure elements are more flat, Eq.~(\ref{approx}) and  Eq.~(\ref{approxetahalf}) leading
to a more radical erosion effect if time is fixed. Note that this difference is also due to the fact that
$D^{11}$ in Eq.~(\ref{Hamdi2}) has physical dimension $\frac{[\textrm{Intensity}][\textrm{Length}]^{2\eta}}{[\textrm{Time}]}$ and
$D^{44}$ has physical dimension $\frac{[\textrm{Intensity}]}{[\textrm{Time}]}$.

\begin{figure}
\centerline{
\hfill
\includegraphics[width= 0.25 \hsize]{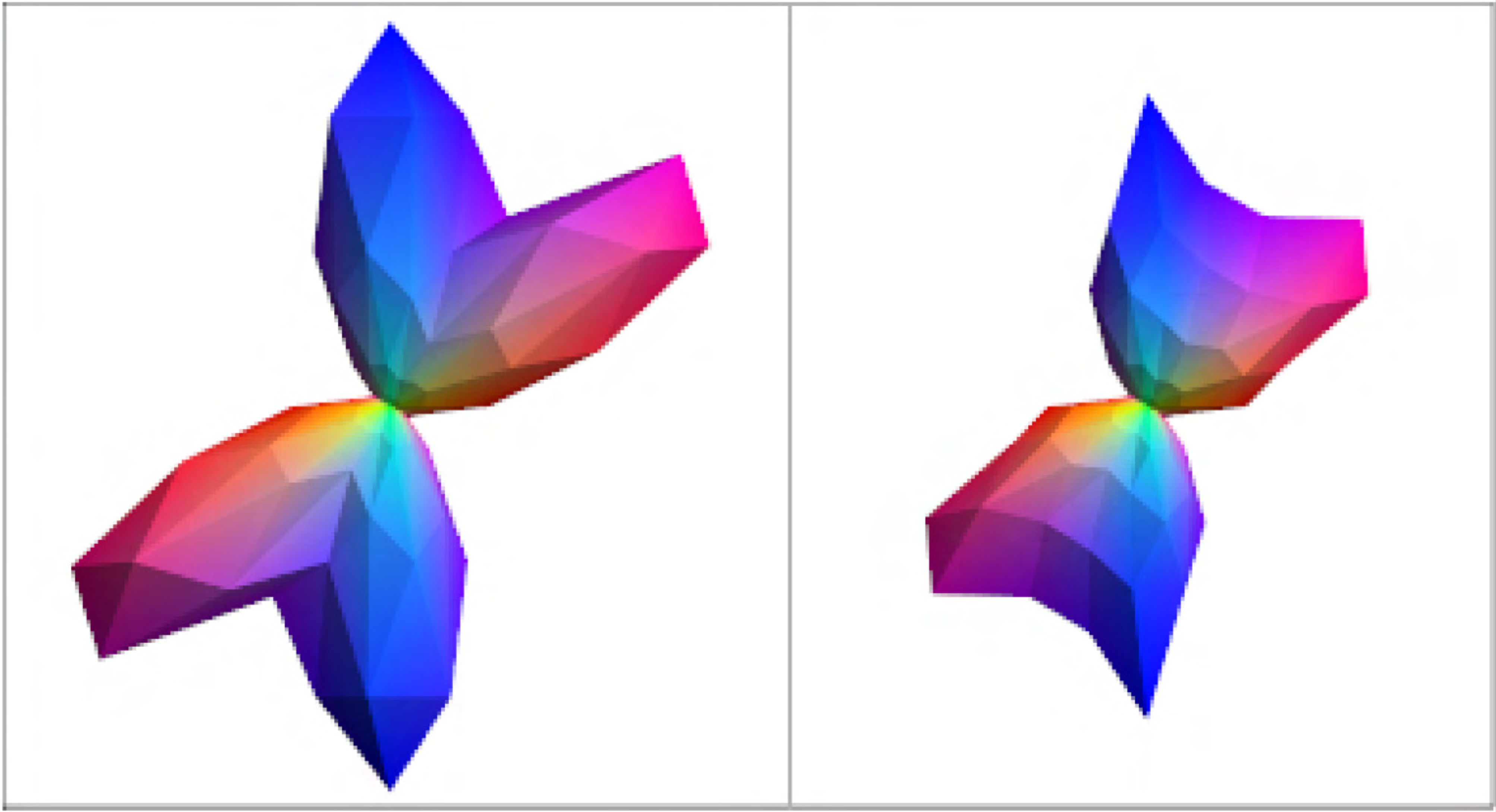} \hfill
\includegraphics[width= 0.25 \hsize]{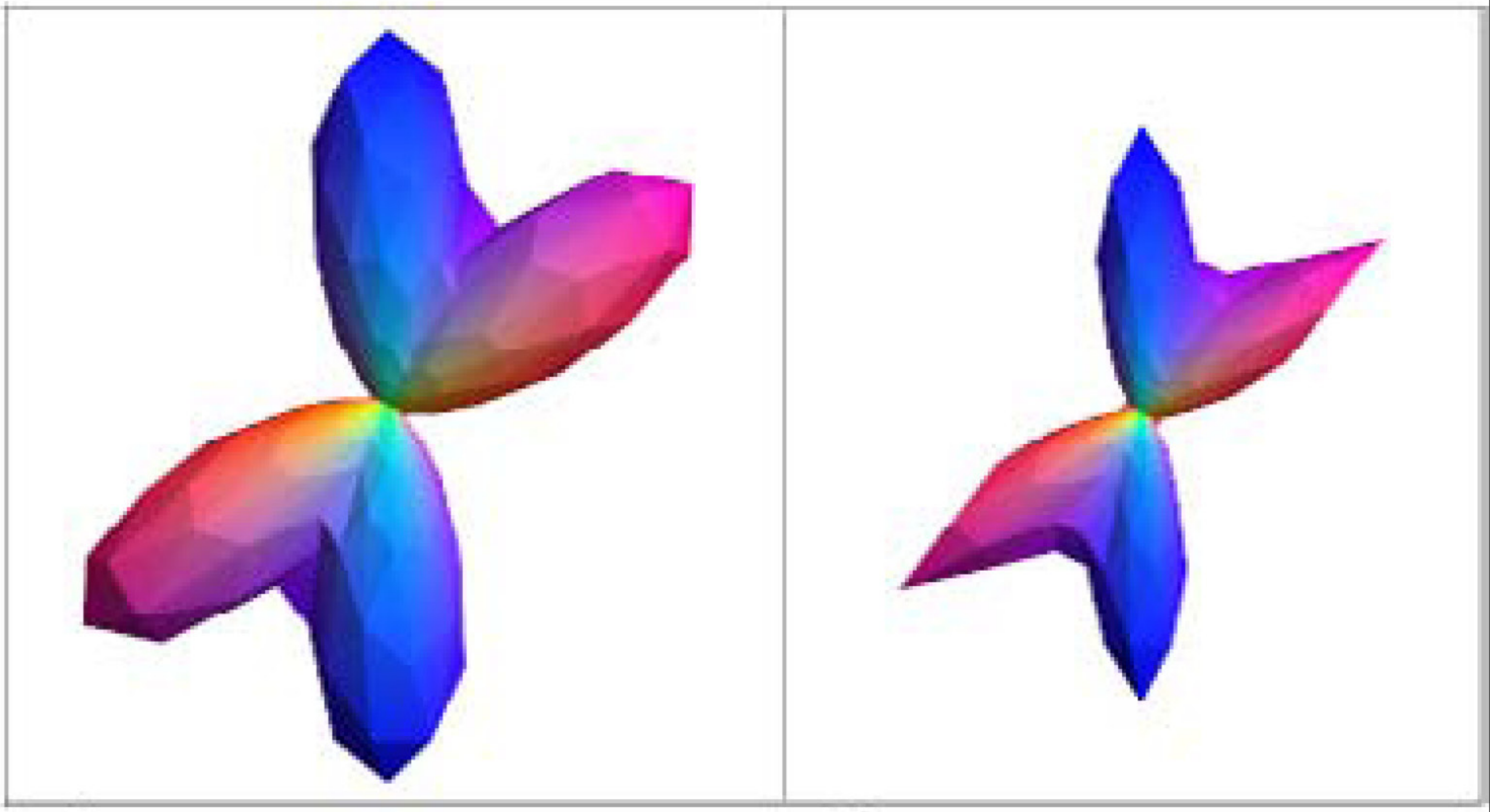}
\hfill
}
\caption{Morphological kernel implementation of the Hamilton-Jacobi equations
 Eq.~(\ref{Hamdi2}) by $\R^{3}\rtimes S^2$-erosion, Eq.~(\ref{Erosion}), $\eta=1$, $D^{11}\!=0$, $D^{44}\!=0.4$, $t=1$,
 Left: input and output for $N_{o}=162$, Right: same for $N_{o}=362$. }\label{Fig:No}
\end{figure}

We have implemented two algorithms for erosion, by morphological convolution in $\R^{3}\rtimes S^2$ erosions Eq.~(\ref{Erosion}) with analytical kernel approximations Eq.~(\ref{approx})
and finite difference upwind schemes with stepsize $t>0$.
A drawback of the morphological convolutions is that it typically requires high sampling on the sphere,
as can be seen in Figure~\ref{Fig:No},
whereas the finite difference schemes also work out well for low sampling rates
(in the experiments on real medical data sets in \ref{Fig:4} we used
a 162-point tessellation of the icosahedron). Nevertheless, for high sampling on the sphere the analytical kernel implementations
were close to the numerical finite difference schemes and apparently our analytical approximation Eq.~(\ref{approxetahalf})
approximates the propagation of balls in $\R^{3}\rtimes S^2$ close enough for reasonable parameter-settings in
practice. Upwind finite difference schemes are better suited for low sampling rates on the
sphere, but for such sampling rates they do involve inevitable numerical blur due to spherical grid interpolations.
Typically, erosions should be applied to strongly smoothed data as it enhances local extrema.
In Figure \ref{Fig:4}, lower left corner one can see the result on a DTI-data set containing fibers of the corpus
callosum and corona radiata, where the effect of both spatial erosion
(glyphs in the middle of a ``fiber bundle'' are larger than the glyphs at the boundary) and
angular erosion glyphs are much sharper so that it is easier to visually track the smoothly varying fibers in the output dataset.
Typically, as can be seen in Fig.~\ref{Fig:4}, applying first diffusion and then erosion is stable with respect to the parameter settings and
produces visually more appealing sharper results then simultaneous diffusion and dilation
in pseudo-linear scale spaces. 

\subsection{Further experiments pseudo-linear scale spaces}

Pseudo-linear scale spaces, Eq.\ref{generalPDE2} combine diffusion and dilation along fibers in a single evolution. The
drawback is a relatively sensitive parameter $C>0$
that balances infinitesimally between dilation and diffusion.
\begin{figure}
\centerline{
\includegraphics[width= 0.7 \hsize]{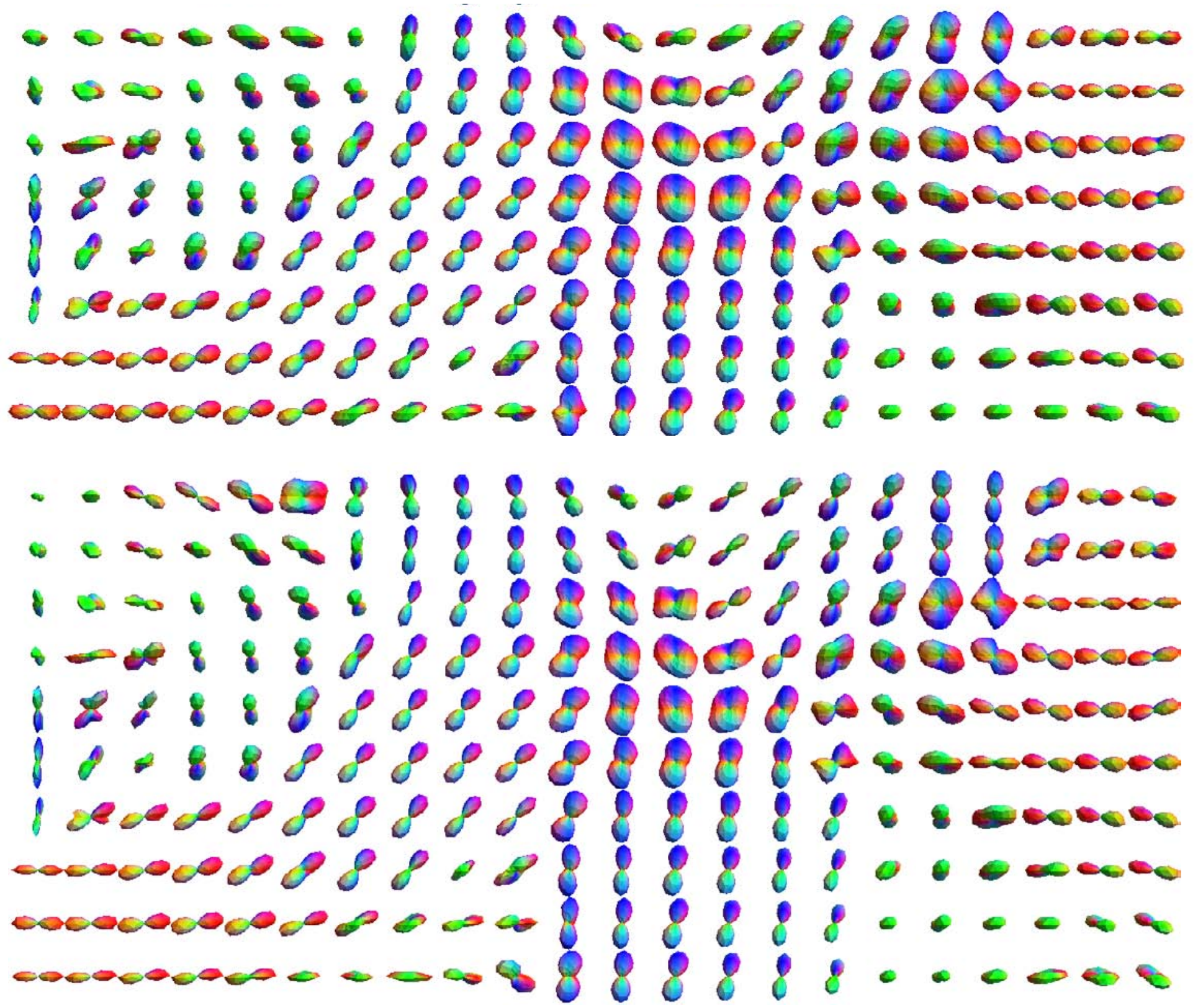} 
}
\caption{Pseudo-linear scale space on $\R^{3}\rtimes S^{2}$ applied on
the initial DTI-dataset in Figure \ref{Fig:4}. $t=D^{33}=1$, $D^{44}=0.04$,
$\Delta t=0.01$, top: $C=2$, bottom: $C=4$.}\label{Fig:pseudolinear}
\end{figure}
In Figure \ref{Fig:pseudolinear} we applied some experiments where we kept parameters fixed except
for $C>0$. The generator of a pseudo-linear scale space consists of a diffusion generator plus $C$ times
an erosion generator and indeed for $C=4$ we see more elongated glyphs and more apparent crossings than for the
case $C=2$. 

\section{Conclusion}

For the purpose of tractography detection and visualization of biological fibers,
Diffusion-Weighted MRI-images, such as DTI and HARDI, should be enhanced such that fiber junctions are maintained,
while reducing high frequency noise in the joint domain $\R^3 \rtimes S^2$ of positions
and orientations.
Therefore we have considered diffusions and Hamilton-Jacobi equations on DW-MRI HARDI and DTI induced by fundamental stochastic processes on $\R^3 \rtimes S^2$ embedded in the group manifold $SE(3)$ of $3D$-rigid body motions.
In order to achieve rotation and translation covariance,
processing must be left-invariant. We have classified all
left-invariant linear and morphological scale spaces on DW-MRI.
This classifications yield novel approaches, that in contrast to the existing methods,
employ fiber extensions simultaneously over positions and orientations,
while preserving crossings and/or bifurcations.
Application on DTI creates fiber-crossings by extrapolation,
comparable to HARDI, allowing to reduce the number of scanning directions,
whereas applying the same method to HARDI removes
spurious non-aligned crossings.
In order to sharpen the DW-MRI images, we apply
morphological scale spaces (Hamilton-Jacobi equations) related to probabilistic
cost-processes on $\R^{3}\rtimes S^2$.
Moreover, we provide pseudo-linear scale spaces combining morphological
scale spaces and linear scale spaces in a single evolution.
All evolutions can be implemented by stable,
left-invariant finite difference methods, or by convolution with
analytic kernels.
For example, we have shown by extending existing results on $\R^{n}$ (and
recently on the Heisenberg group to $SE(3)$, cf.\cite{Manfredi}),
that morphological $\R^{3}\rtimes S^2$-convolutions are the unique
viscosity solutions of Hamilton-Jacobi equations on $\R^{3}\rtimes S^2$ and moreover 
we have provided analytical approximations of the
morphological kernels (Green's functions). Finally, we have extended anisotropic linear diffusions systems on $\R^{3}\rtimes S^2$ to data-adaptive diffusion systems
of Perona \and Malik type that are useful
to avoid undesired interaction between more isotropic regions (at the neural fibers)
and less isotropic regions (at the ventricles).

In future work we aim to extend our differential geometrical approach on $\R^{3}\rtimes S^2$
to the actual tracking of fibers.
\section*{Acknowledgements.}
The authors wish to thank Vesna Pr\v{c}kovska and Paulo Rodrigues, Biomedical image analysis group \emph{TU/e, Eindhoven University of Technology}, who
are gratefully acknowledged for both providing the medical data sets and for their previous efforts, cf.\cite{Rodrigues,Prckovska}, on fast kernel implementations of linear hypoelliptic diffusion on HARDI and DTI
(Figure \ref{Fig:2} and Figure \ref{Fig:3}).
Furthermore, the authors wish to thank Mark Bruurmijn, master student Biomedical Engineering at TU/e for his DW-MRI visualization programme
\emph{ListHardiPlot},
which enabled us to do all experiments in \emph{Mathematica}.
Furthermore, the authors wish to thank Ronald Rook from the Department of Mathematics and Computer Science (CASA) TU/e for his suggestions on the finite difference schemes for adaptive diffusions.

\appendix

\section{Metric and Lagrangians
on $(\R^{3}\rtimes S^2, {\rm d}\mathcal{A}^{1}, {\rm d}\mathcal{A}^{2})$ \label{app:metric}}
\vspace{-0.6 cm} \mbox{}
Recall from Figure \ref{fig:intuition} that our left-invariant diffusions on $\R^{3}\rtimes S^2$ take place on the contact/sub-Riemannian manifold
$(SE(3), {\rm d}\mathcal{A}^{1}, {\rm d}\mathcal{A}^{2}, {\rm d}\mathcal{A}^{6})$.
On this contact manifold $(SE(3), {\rm d}\mathcal{A}^{1}, {\rm d}\mathcal{A}^{2}, {\rm d}\mathcal{A}^{6})$ we set the
Lagrangian
\begin{equation}
\gothic{L}_{\eta}(\tilde{\gamma}(s),\dot{\tilde{\gamma}}(s))= \frac{2\eta-1}{2\eta} \left(\frac{|\dot{\tilde{\gamma}}^{3}(s)|^2}{D^{33}} + \frac{|\dot{\tilde{\gamma}}^{4}(s)|^2 + |\dot{\tilde{\gamma}}^{5}(s)|^2}{D^{44}}\right)^{\frac{\eta}{2\eta-1}},
\end{equation}
where $\eta > \frac{1}{2}$ and where $D^{33}, D^{44}>0$ and where
$\dot{\tilde{\gamma}}^{i}(s)= \langle \left.{\rm d} \mathcal{A}^{i}\right|_{\tilde{\gamma}(s)}, \dot{\tilde{\gamma}}(s) \rangle$ are the components
w.r.t. the moving frame of reference attached to the curve, cf. Figure \ref{fig:intuition}. For $\eta \to \infty$
we have a homogeneous Lagrangian and a left-invariant
semi-metric on $(SE(3), {\rm d}\mathcal{A}^{1}, {\rm d}\mathcal{A}^{2}, {\rm d}\mathcal{A}^{6})$ is set by
\begin{equation} \label{semimetric}
\begin{array}{l}
\gothic{d}(g_{1},g_{2}):= \gothic{d}(g_{2}^{-1}g_{1},e=(\ul{0},I))= \\
\inf \limits_{h_{1}, h_2 \in \{0\} \times SO(2) \subset SE(3)}\! \! \! \! \! \!\inf \limits_{{\tiny \begin{array}{c}
\gamma=(\ul{x}(\cdot),R_{\ul{n}}(\cdot)) \in C^{\infty}((0,L), SE(3)),\\
\gamma(0)=e h_{2}, \gamma(1)=g_{2}^{-1}g_{1}h_{1},\\
\langle \left.{d \rm}\mathcal{A}^{1}\right|_{\gamma}, \dot{\gamma}\rangle =
\langle \left.{d \rm}\mathcal{A}^{2}\right|_{\gamma}, \dot{\gamma}\rangle =
\langle \left.{d \rm}\mathcal{A}^{6}\right|_{\gamma}, \dot{\gamma}\rangle =0
\end{array}
}}\! \! \! \! \! \! \!
\int_{0}^{1} \sqrt{\sum \limits_{i \in \{3,4,5\}}\frac{1}{D^{ii}}|{\rm d}\mathcal{A}^{i}(\dot{\gamma}(s))|^{2}}\, {\rm d}s
\end{array}
\end{equation}
By left-invariance of the Lagrangian and the fact that
$e h_{2}=h_{2} e$ we may as well set $\gamma(0)=e$ and
$\gamma(L)=h_{2}^{-1}g_{2}^{-1}g_{1}h_{1}= (g_{2}h_{2})^{-1}g_{1}h_{1}$ in Eq.~(\ref{semimetric})
so that $\gothic{d}(g_{1},g_{2})=\gothic{d}(g_{1}h_{1},g_{2}h_{2})$ for all $h_{1},h_{2} \in \{\ul{0}\}\times SO(2)$
and thereby inducing a well-defined metric on $(\R^{3}\rtimes S^2, {\rm d}\mathcal{A}^{1}, {\rm d}\mathcal{A}^{2})$:
\[
d^{\textrm{hor}}((\ul{y},\ul{n}), (\ul{y}', \ul{n}'))=\gothic{d}((\ul{y},R_{\ul{n}}), (\ul{y}', R_{\ul{n}'}))\]
on $\R^{3}\rtimes S^2$. Along horizontal curves $s \mapsto \tilde{\gamma}(s)=(\ul{x}(s),R(s))$ in $SE(3)$ we have
$\ul{n}(s)=R(s)\ul{e}_{z}=\dot{\ul{x}}(s)$ and
\[
\begin{array}{l}
\|\KKK(s)\|^2= |\langle \left.{\rm d}\mathcal{A}^{4}\right|_{\gamma(s)}, \dot{\gamma}(s)\rangle|^2+
|\langle\left.{\rm d}\mathcal{A}^{5}\right|_{\gamma(s)}, \dot{\gamma}(s)\rangle|^2 , \\
\langle \left.{\rm d}\mathcal{A}^{3}\right|_{\gamma(s)}, \dot{\gamma}(s)\rangle=\|\dot{x}(s)\|=1 ,
\end{array}
\]
and thereby the metric reduces to
\begin{equation}\label{metricred}
d^{\textrm{hor}}((\ul{y},\ul{n}), (\ul{y}', \ul{n}'))=
\inf \limits_{{\scriptsize \begin{array}{c}
\ul{x}(\cdot) \in C^{\infty}((0,L), \R^3),\\
\ul{x}(0)=\ul{0}, \; \dot{\ul{x}}(0)=\ul{e}_{z}, \\
\ul{x}(L)=R_{\ul{n}}^{T}(\ul{y}'-\ul{y}), \\ \dot{\ul{x}}(L)=R_{\ul{n}}^{T} R_{\ul{n}'}\ul{e}_{z}=R_{\ul{n}}^{T}\ul{n}', \\
\end{array}
}}\!
\int \limits_{0}^{L} \sqrt{\frac{1}{D^{44}}\|\KKK(s)\|^2 + \frac{1}{D^{33}}}\; {\rm d}s
\end{equation}
where $s$, $L>0$, and $\KKK(s)$ are respectively
spatial arclength, total length, and curvature of the
spatial part of the curve. Note that the distance (\ref{metricred}) is indeed invariant under
the choice of $R_{\ul{n}}$ such that Eq.~(\ref{Rn}) holds, due to left-invariance of (\ref{semimetric})
and rotation covariance of $(\ref{metricred})$ and
$(R_{\ul{n}}R_{\ul{e}_{z},\alpha})^T=R_{\ul{e}_{z},\alpha}^{T}R_{\ul{n}}$. Furthermore, we note that torsion is
implicitly penalized by simultaneous penalization of length and curvature.

\section{Viscosity Solutions of Hamilton-Jacobi Equations on $(\R^{3}\rtimes S^{2}, {\rm d}\mathcal{A}^{3})$ \label{app:viscosity}}

After introducing a positive semi-definite\footnote{It need not be strictly positive definite, but the H\"{o}rmander condition \cite{Hoermander} should not be violated.} left-invariant metric tensor
$\ul{G}: \R^{3}\rtimes S^{2} \times T( \R^{3}\rtimes S^{2}) \times T( \R^{3}\rtimes S^{2}) \to \R$ given by Eq.~(\ref{metrictensorR3S2erosion}),
we considered the solutions (\ref{viscsol1}) and (\ref{viscsol2}) of the Hamilton-Jacobi equations on $\R^{3}\rtimes S^{2}$:
\begin{equation}\label{kernMorph2}
\left\{
\begin{array}{l}
\frac{\partial W}{\partial t}(\ul{y},\ul{n},t)= \pm (H( {\rm d}W(\cdot,\cdot,t)))(\ul{y},\ul{n}) :=\frac{1}{2}(\ul{G}_{(\ul{y},\ul{n})})^{-1}\left(\left. {\rm d}W(\cdot,\cdot,t) \right|_{\ul{y},\ul{n}},\left. {\rm d}W(\cdot,\cdot,t) \right|_{\ul{y},\ul{n}})\right) \\
W(\ul{y},\ul{n},0)= U(\ul{y},\ul{n})
\end{array}
\right.
\end{equation}
where $H$ stands for the Hamiltonian. It is well-known \cite{Crandall} that Hamilton-Jacobi equations in general do not have a unique solution, unless extra requirements such as the practically reasonable viscosity requirement are added.\footnote{Consider for example the initial value problem $\frac{\partial U}{\partial t}(x,t)=(\frac{\partial U}{\partial x}(x,t))^2$ and $U(x,0)=0$, $x \in \R, t \geq 0$, with continuous (non-sensible) solution $U(x,t)=|x|-t$ if $t>|x|$ and zero else. This ``nonsense'' solution is not a viscosity solution, namely consider $\phi(x,t)=e^{-t}e^{-x^2}$ and $k=-1$ in the Definition of viscosity solutions according to Crandall \cite{Crandall}.}

Our claim was that the solutions (\ref{viscsol1}) and (\ref{viscsol2}) are the unique \cite{Crandall} viscosity solutions of the Hamilton-Jacobi equations (\ref{Hamdi1}) and (\ref{Hamdi2}).
\begin{definition}
Suppose that the Hamiltonian is convex and $H(p) \to \infty$ as $p\to \infty$. Then a viscosity solution is a bounded and continuous (not necessarily differentiable) weak solution $W: (\R^{3} \rtimes S^2)\times \R^{+} \to \R$ of (\ref{kernMorph2}) such that
\begin{enumerate}
\item for any smooth function $V:(\R^{3} \rtimes S^2)\times \R^{+} \to \R$ such that $W-V$ attains a local maximum at $(\ul{y}_{0},\ul{n}_{0},t_{0})$ one has
\begin{equation} \label{visc1}
\frac{\partial V}{\partial t}(\ul{y}_{0},\ul{n}_{0},t_{0})\mp (H( {\rm d}V(\cdot,\cdot,t)))(\ul{y}_{0},\ul{n}_{0}) \leq 0  \ .
\end{equation}
\item for any smooth function $V:(\R^{3} \rtimes S^2)\times \R^{+} \to \R$ such that $W-V$ attains a local minimum at $(\ul{y}_{0},\ul{n}_{0},t_{0})$ one has
\begin{equation} \label{visc2}
\frac{\partial V}{\partial t}(\ul{y}_{0},\ul{n}_{0},t_{0})\mp (H( {\rm d}V(\cdot,\cdot,t)))(\ul{y}_{0},\ul{n}_{0}) \geq 0  \ .
\end{equation}
\end{enumerate}
\end{definition}
We will first provide a quick review of Hamilton-Jacobi theory, the Hopf-Lax formula and erosions on respectively $\R^{4}$ and on the Heisenberg group $(SE(3))_{0}$ that we obtain by contraction from $SE(3)$ and which serves as a local approximation of $SE(3)$, as explained in Section \ref{ch:contraction}.
\begin{definition} \label{def:LT}
Let $X$ be a normed space, then the Legendre-Fenchel transform
$L \mapsto \gothic{F}L$ on $X$ is given by
\[
(\gothic{F}_{X}L)(x)= \sup \limits_{y \in X}
\left\{ \langle x, y\rangle - L(y)\right\}\
\]
where $L:X \to \R \cup {\infty}$ and $\langle x, y\rangle=x(y)$ and $x \in X^{*}$.
In particular if $X=\R^{n}$ one gets 
\[
(\gothic{F}_{\R^{n}}L)(\ul{x})= \sup \limits_{\ul{y} \in \R^{n}} \left\{\ul{x} \cdot \ul{y} - L(\ul{y})\right\}\ \]
and in case of the Lie-algebra $\mathcal{L}(SE(3))$ of left-invariant vector fields on $SE(3)$
we get
\[
(\gothic{F}_{\mathcal{L}(SE(3))}L)( \sum \limits_{i=1}^{6} p_{i}{\rm d}\mathcal{A}^{i})=
\sup \limits_{\sum \limits_{i=1}^{6} c^{i}\mathcal{A}_{i}}\{
\langle \sum \limits_{i=1}^{6} p_{i}{\rm d}\mathcal{A}^{i},
\sum \limits_{j=1}^{6} c^{j} \mathcal{A}_{j}  \rangle -L(\sum \limits_{i=1}^{6} c^{i}\mathcal{A}_{i})\}=
((\gothic{F}_{\R^{6}}l)(p_{1},\ldots,p_{6})\ ,
\]
with $l(c^{1},\ldots,c^{6})=L(\sum \limits_{i=1}^{6} c^{i}\mathcal{A}_{i})$.
\end{definition}

\begin{lemma}\label{lemma:LF}
Let $\eta>\frac{1}{2}$. The Legendre-Fenchel transform of the function
\[
\R^{4} \ni (c^{1},c^{2},c^{4},c^{5}) \mapsto \frac{2\eta-1}{2\eta} \left(\frac{(c^{1})^2+ (c^{2})^2 }{D^{11}}+\frac{(c^{4})^2+ (c^{5})^2 }{D^{44}} \right)^{\frac{\eta}{2\eta-1}} \in \R
\]
is given by
\[
\R^{4} \ni (p_{1},p_{2},p_{4},p_{5}) \mapsto \frac{1}{2\eta} \left(D^{11}((p_{1})^2+ (p^{2})^2) +D^{44}((p^{4})^2+ (p^{5})^2  \right)^{\eta} \in \R
\]
\end{lemma}
\textbf{Proof }It is well-known that the Fenchel transform of $\ul{c} \mapsto \frac{1}{a}\|\ul{c}\|^a$ equals
$\ul{p} \mapsto \frac{1}{b}\|\ul{p}\|^{b}$ with $\frac{1}{a}+\frac{1}{b}=1$, furthermore $\gothic{F}\circ \mathcal{D}_{\lambda}= \mathcal{D}_{\lambda^{-1}} \circ \gothic{F}$ where $\mathcal{D}_{\lambda}$ denotes the
scaling operator given by $\mathcal{D}_{\lambda}f(\ul{x})=f(\lambda^{-1} \ul{x})$ from which the result follows by setting
$b=2\eta$, $a=\frac{2\eta}{2\eta-1}$ and respectively $\lambda=\sqrt{D^{11}}, \sqrt{D^{44}}$. $\hfill \Box$
\begin{corollary} Let $\eta >\frac{1}{2}$, $t>0$, then
on $\R^{4}$ the Hopf-Lax formula is given by
\begin{equation}\label{Hopf}
w(\ul{x},t):= \inf \limits_{\ul{y} \in \R^{4}} \inf \limits_{{\tiny \begin{array}{l}
\gamma\in C^{1}(0,t) \\
\gamma(t)=\ul{x} \\
\gamma(0)=\ul{y}
\end{array}}}\{ \int \limits_{0}^{t} \mathcal{L}_{\eta}(\dot{\gamma}(s))\, {\rm d}s + u(\ul{y} )\}= \min \limits_{\ul{y} \in \R^{4}} \{ t \mathcal{L}_{\eta}\left(\frac{1}{t}(\ul{x}-\ul{y})\right)+ u(\ul{y})\},
\end{equation}
with $w(\cdot,0)=u$ a given Lipschitz continuous function on $\R^{4}$, with
\begin{equation}\label{Leta}
\begin{array}{l}
\mathcal{L}_{\eta}(\dot{\gamma}(s))=\mathcal{L}_{\eta}(\dot{\gamma}^{1}(s),\dot{\gamma}^{2}(s),\dot{\gamma}^{4}(s),\dot{\gamma}^{5}(s)) \\ := \frac{2\eta-1}{2\eta} \left(\frac{1}{D^{11}}((\dot{\gamma}^{1}(s))^2+(\dot{\gamma}^{2}(s))^2)+
\frac{1}{D^{44}}((\dot{\gamma}^{4}(s))^2+(\dot{\gamma}^{5}(s))^2)\right)^{\frac{\eta}{2\eta-1}}
\end{array}
\end{equation}
and $w(\ul{x},t)$ given by (\ref{Hopf}) is the unique viscosity solution of the Hamilton-Jacobi-Bellman
system
\[
\left\{
\begin{array}{l}
\frac{\partial w}{\partial t}(\ul{x},t)= -\frac{1}{2\eta} \left(D^{11}\left(
\left(
\frac{\partial w}{\partial x^{1}}(\ul{x},t)\right)^2+\left(
\frac{\partial w}{\partial x^{2}}(\ul{x},t)\right)^2\right)+ D^{44}\left(
\left(
\frac{\partial w}{\partial x^{4}}(\ul{x},t)\right)^2+\left(
\frac{\partial w}{\partial x^{5}}(\ul{x},t)\right)^2\right)\right)^{\eta}, \\
w(\ul{x},0)=u(\ul{x})
\end{array}
\right.
\]
with $\ul{x}=(x^{1},x^{2},x^{4},x^{5}) \in \R^{4}, t \geq 0$. The morphological Green's function of this PDE is given by
\[
k_{t}^{D^{11},D^{44},\eta,+}(\ul{x})= t \mathcal{L}_{\eta}(t^{-1}\ul{x})= \frac{2\eta-1}{2\eta}\, t^{-\frac{1}{2\eta-1}} \, \left(\frac{(x^{1})^2+ (x^{2})^2 }{D^{11}}+\frac{(x^{4})^2+ (x^{5})^2 }{D^{44}} \right)^{\frac{\eta}{2\eta-1}}.
\]
In the limiting case $\eta \to \infty$ the
Lagrangian is homogeneous\footnote{
The corresponding Hamiltonian
$\mathcal{H}(\tilde{\ul{p}})=\mathcal{L}(\dot{\gamma})$ is homogeneous as well, but now with respect to
$\tilde{p}_{i}:=\frac{1}{2}\frac{\partial^{2}(\mathcal{L}_{\infty}(\dot{\gamma}))}{\partial \dot{\gamma}^{i}\dot{\gamma}^{j}}$, [ch:3]\cite{Rund}.}
and
\[
w(\ul{x},t)=\min \limits_{\ul{y}}\{u(\ul{y})+\left(\frac{1}{D^{11}}\left((x^{1}-y^{1})^2+ (x^{2}-y^{2})^2\right)
+ \frac{1}{D^{44}}\left( (x^{4}-y^{4})^2+ (x^{5}-y^{5})^2\right)\right)^{\frac{1}{2}}\}
\] and we arrive at the time-independent Hamilton-Jacobi equation
\[
\begin{array}{l}
1= D^{11}\left(
\left(
\frac{\partial w}{\partial x^{1}}(\ul{x},t)\right)^2+\left(
\frac{\partial w}{\partial x^{2}}(\ul{x},t)\right)^2\right)+ D^{44}\left(
\left(
\frac{\partial w}{\partial x^{4}}(\ul{x},t)\right)^2+\left(
\frac{\partial w}{\partial x^{5}}(\ul{x},t)\right)^2\right), \\
\end{array}
\]
In the limiting case $\eta \downarrow \frac{1}{2}$ the Hamiltonian is homogeneous
and the Hopf-Lax formula reads
\[
w(x,t)= \min_{\ul{y}\in \R^{4}}\{k_{t}^{D^{11},D^{44},1/2,+}(\ul{x}-\ul{y})+u(\ul{y})\}
\]
where the flat morphological erosion kernel is given by
\[
k_{t}^{D^{11},D^{44},1/2,+}(\ul{y})= \left\{
\begin{array}{ll}
0 &\textrm{ if }\frac{(y^{1})^2+ (y^{2})^2 }{D^{11}}+\frac{(y^{4})^2+ (y^{5})^2 }{D^{44}}\leq t^2,  \\
\infty &\textrm{ else }.
\end{array}
\right.
\]
\end{corollary}
\textbf{Proof } Directly follows by Lemma \ref{lemma:LF} and \cite[ch:3(Thm 4), ch:10 (Thm 1, Thm 3)]{Evans}. For the special case of a homogeneous Lagrangian
($\eta \to \infty$), see \cite[Ch:3]{Rund}.
\begin{remark}
We parameterized elements in $\R^{4}$ by means of $\ul{y}=(y^{1},y^{2},y^{4},y^{5})$ since we want to stress the analogy to Hamilton-Jacobi equations on $\R^{3}\rtimes S^2$.
\end{remark}

\begin{lemma}\label{lemma:semigroup} (semigroup property of erosion operator on $\R^{3}\rtimes S^{2}$)\\
Let $\eta>\frac{1}{2}$ and let $D^{11}>0,D^{44}>0,t>0$. Let the morphological Green's function be defined as\footnote{We assume that $\gamma(t)$ is chosen such that a minimizer exists, which is the case if $\gamma(t)$ is close enough to the unity $(\ul{0},I)$. However, one can show that for $\eta \to \infty$ and $\gamma(t)$ far away from the unity a minimizer may not exist, similar to the 2D-case \cite{Boscain}.}
\begin{equation} \label{distR3XS2}
\begin{array}{ll}
k_{t}^{D^{11},D^{44},\eta}(\ul{y},\ul{n})&:=\inf \limits_{{\tiny \begin{array}{c}
\gamma=(\ul{x}(\cdot),R(\cdot)) \in C^{\infty}((0,t), SE(3)),\\
\gamma(0)=(\ul{0},I=R_{\ul{e}_z}), \gamma(t)=(\ul{y},R_{\ul{n}}),\\
\langle \left.{d \rm}\mathcal{A}^{3}\right|_{\gamma}, \dot{\gamma}\rangle=\langle \left.{d \rm}\mathcal{A}^{6}\right|_{\gamma},\dot{\gamma}\rangle=0
\end{array}}} \int \limits_{0}^{t} \overline{\mathcal{L}}_{\eta}(\gamma(p),\dot{\gamma}(p))\, \left(\frac{dp}{ds}\right)^{\frac{1}{2\eta-1}}\,{\rm d}p\ , \\
 & =\inf \limits_{{\tiny \begin{array}{c}
\gamma=(\ul{x}(\cdot),R(\cdot)) \in C^{\infty}((0,\psi^{-1}(t)), SE(3)),\\
\gamma(0)=(\ul{0},I=R_{\ul{e}_z}), \gamma(\psi^{-1}(t))=(\ul{y},R_{\ul{n}}),\\
\langle \left.{d \rm}\mathcal{A}^{3}\right|_{\gamma}, \dot{\gamma}\rangle=\langle \left.{d \rm}\mathcal{A}^{6}\right|_{\gamma},\dot{\gamma}\rangle=0
\end{array}}} \int \limits_{0}^{\psi^{-1}(t)} \tilde{\mathcal{L}}_{\eta}(\gamma(s),\dot{\gamma}(s))\,{\rm d}s\
\end{array}
\end{equation}
with
\begin{equation} \label{2eta}
\begin{array}{l}
\overline{\mathcal{L}}_{\eta}(\gamma(p),\dot{\gamma}(p)):=\mathcal{L}_{\eta}( \langle
\left.{\rm d}\mathcal{A}^{1}\right|_{\gamma(p)}, \dot{\gamma}(p)\rangle, \ldots, \langle
\left.{\rm d}\mathcal{A}^{5}\right|_{\gamma(p)}, \dot{\gamma}(p)\rangle)=\frac{2\eta-1}{2\eta}, \\
\tilde{\mathcal{L}}_{\eta}(\gamma(s),\dot{\gamma}(s))= \mathcal{L}_{\eta}( \langle
\left.{\rm d}\mathcal{A}^{1}\right|_{\gamma(s)}, \dot{\gamma}(s)\rangle, \ldots, \langle
\left.{\rm d}\mathcal{A}^{5}\right|_{\gamma(s)}, \dot{\gamma}(s)\rangle),
\end{array}
\end{equation}
where we recall Eq.~(\ref{Leta}) and with $\R^{3}\rtimes S^{2}$-``erosion arclength'' given by
\begin{equation}\label{arclengthR3XS2}
\begin{array}{ll}
p(\tau) &= \int \limits_{0}^{\tau} \sqrt{\mathbf{G}_{\gamma(\tilde{\tau})}(\dot{\gamma}(\tilde{\tau}),\dot{\gamma}(\tilde{\tau}))}\, {\rm d}\tilde{\tau} =  \int \limits_{0}^{\tau} \sqrt{\sum \limits_{i \in \{1,2,4,5\}}\frac{1}{D^{ii}}|\langle \left. {\rm d}\mathcal{A}^{i}\right|_{\gamma(\tilde{\tau})}, \dot{\gamma}(\tilde{\tau})\rangle |^{2}  }\, {\rm d}\tilde{\tau}
\end{array}
\end{equation}
with $D^{11}=D^{22}$ and $D^{44}=D^{55}$ and with left-invariant metric tensor $\mathbf{G}(\cdot,\cdot)=\lim \limits_{\eta \to \infty}(\overline{\mathcal{L}}_{\eta}(\cdot,\cdot))^2$ and where $\psi(s)=p$ given by
\begin{equation} \label{psis}
\psi(s)=\int_{0}^{s} \sqrt{\kappa^{2}(\tilde{s})+\beta^2}\, {\rm d}\tilde{s}.
\end{equation}
%
Then we have the following identity
\begin{equation} \label{semigroup}
k_{t-\tau}^{D^{11},D^{44},\eta,+} \ominus (k_{\tau}^{D^{11},D^{44},\eta,+} \ominus U)= k_{t}^{D^{11},D^{44},\eta,+} \ominus U\ ,
\end{equation}
for all $\tau \in [0,t]$ and all Lipschitz continuous functions $U:\R^{3}\rtimes S^{2} \to \R$,
where we recall that our erosion operator $\ominus$ is defined in Eq.~(\ref{Erosion}).
\end{lemma}
\begin{remark}
We apply the following convention regarding curves within $SE(3)$:
For concise notation we shall write $\gamma(p)$ instead of $\gamma(s(p))$. When writing $\dot{\gamma}(p)$ we mean
$\frac{d}{dp}\gamma(\psi^{-1}(p))$, when writing $\dot{\gamma}(s)$ we mean $\frac{d}{ds}\gamma(s)$ !
\end{remark}
\textbf{Proof }
First of all we note that $D^{11}=D^{22}$ and $D^{44}=D^{55}$ are necessary requirements for a
well-defined morphological kernel $k_{t}^{D^{11},D^{44},\eta,+}:\R^{3}\rtimes S^{2} \to \R^{+}$ and in particular (as $\eta \to \infty$) a well-defined metric
on the group quotient \mbox{$\R^{3}\rtimes S^{2}=(\R^{3}\rtimes SO(3))/(\{0\}\times SO(2))$}.

Secondly, we will explain and derive the identities Eq.~(\ref{distR3XS2}) and Eq.~\!(\ref{2eta}).
The limit $\lim \limits_{\eta \to \infty}\overline{\mathcal{L}}_{\eta}(\cdot,\cdot)$ is homogeneous in the second entry
so that
\begin{equation} \label{p}
\overline{\mathcal{L}}_{\infty}(\gamma(p), \frac{d\gamma}{dp})= \overline{\mathcal{L}}_{\infty}\left(\gamma(p(\tau)),(\overline{\mathcal{L}}_{\infty}(\gamma(p(\tau)), \frac{d\gamma}{d\tau}))^{-1} \frac{d\gamma}{d\tau}\right)=1 \textrm{ and }\lim \limits_{\eta \to \infty}k_{t}^{D^{11},D^{44},\eta,+}(\ul{y},\ul{n})=t,
\end{equation}
so $t>0$ is the total length of the curve in the sub-Riemannian manifold $(\R^{3}\rtimes S^{2}, {\rm d}\mathcal{A}^{3})$ and where $p$ is the
arc-length parameter in  $(\R^{3}\rtimes S^{2}, {\rm d}\mathcal{A}^{3})$
since connecting curves over which is optimized are not allowed to use the $\mathcal{A}_{3}$-direction in the tangent bundle.
Furthermore, by Eq.~\!(\ref{p}) we have that
$\frac{2\eta-1}{2\eta}\left(\overline{\mathcal{L}}_{\infty}(\gamma(p), \frac{d\gamma}{dp})\right)^{\frac{\eta}{\eta-\frac{1}{2}}}=\frac{2\eta-1}{2\eta}$, which explains Eq.~\!(\ref{2eta}).
Regarding Eq.~\!(\ref{distR3XS2}) it is relevant to keep $t>0$ fixed during the curve optimization in (\ref{distR3XS2}) for $\frac{1}{2}<\eta<\infty$, otherwise the minimizing curve may not exist, and direct computation yields
\[
\begin{array}{l}
\int \limits_{0}^{t} \overline{\mathcal{L}}_{\eta}(\gamma(p),\dot{\gamma}(p))\, \left(\frac{dp}{ds}\right)^{\frac{1}{2\eta-1}}\,{\rm d}p =
 \int \limits_{0}^{t} \overline{\mathcal{L}}_{\eta}(\gamma(p),\frac{d\gamma}{dp} \frac{dp}{ds})\,
\frac{ds}{dp}\,{\rm d}p
= \int \limits_{0}^{\psi^{-1}(t)} \tilde{\mathcal{L}}_{\eta}(\gamma(s),\frac{d\gamma}{ds})\, {\rm d}s,
\end{array}
\]
where we used the homogeneity of $\overline{\mathcal{L}}_{\eta}$:
\[
\left(\frac{dp}{ds} \right)^{\frac{1}{2\eta-1}} \overline{\mathcal{L}}_{\eta}(\gamma, \frac{d\gamma}{dp})= \left(\frac{dp}{ds}\right)^{-1}\left(\frac{dp}{ds} \right)^{\frac{2\eta}{2\eta-1}} \overline{\mathcal{L}}_{\eta}(\gamma, \frac{d\gamma}{dp})=
\frac{ds}{dp} \overline{\mathcal{L}}_{\eta}(\gamma, \frac{d\gamma}{dp} \frac{dp}{ds}).
\]
Thirdly, we recall that functions $U:\R^{3}\rtimes S^{2} \to \R$ are related to functions $\tilde{U}:SE(3) \to \R$
with invariance property ($\forall_{g \in SE(3)} \forall_{h=(\ul{0},R_{\ul{e}_{z},\alpha}) \in SE(3)}
\tilde{U}(g h)=\tilde{U}(g)$)  by means of
\[
\tilde{U}(g)=\tilde{U}(\ul{x},R)=U(\ul{x},R\ul{e}_{z})\ ,
\]
with $g=(\ul{x},R) \in SE(3), \ul{y} \in \R^{3}, \ul{e}_{z}=(0,0,1)^{T}$. We apply this identification by setting
\[
\tilde{W}(g,t):= \inf \limits_{q \in SE(3)} \{\tilde{k}_{t}(q^{-1}g) + \tilde{U}(q) \}
\]
for all $g \in SE(3)$ and $t \geq 0$. Consequently, we may rewrite our Eq.~(\ref{semigroup}) as
\begin{equation} \label{aap}
\tilde{W}(g,t)= \inf \limits_{v \in SE(3)}\{\tilde{k}_{t-\tau}(v^{-1}g) +\tilde{W}(v,\tau)\}\ ,
\end{equation}
where we use short notation $\tilde{k}_{t}(\ul{y},R):=k_{t}^{D^{11},D^{44},\eta,+}(\ul{y},R\ul{e}_{z})$.

Fourthly, we apply basically the same approach as in \cite[Thm.1]{Manfredi} where we replace the Heisenberg group by the 3D-Euclidean motion group $SE(3)$ where in this lemma we carefully omit scaling properties of the Lagrangian
that do hold for the Heisenberg group, but not on the $SE(3)$ (even though one has the contraction $\lim_{q \downarrow 0}(SE(3))_{q}=(SE(3))_{0}$ which serves as a local nilpotent \emph{approximation}, \cite{TerElst3}). We shall return to these local approximations later when we derive approximations of $k_{t}^{D^{11},D^{44},\eta,+}(\ul{y},\ul{n})$. For, now let us proceed with proving Eq.~(\ref{aap}).

In general one has
\begin{equation} \label{estkern}
\tilde{k}_{t}(v^{-1}g) \leq \tilde{k}_{\tau}(v^{-1}q) + \tilde{k}_{t-\tau}(q^{-1}g)
\end{equation}
for all $g,h,v \in SE(3)$ and all $\tau \in [0,t]$ as the concatenation of two optimal curves yields an admissible curve over which is optimized in the left hand side. Here equality is obtained if $v^{-1}q$ is on the minimizing curve between $e$ and $v^{-1}g$, i.e. by left-invariance $g$ is on the minimizing curve between $e$ and $v^{-1}g$.  Now due to continuity and convexity of $(c^{1},c^{2},c^{4},c^{5}) \mapsto \mathcal{L}_{\eta}(c^{1},c^{2},c^{4},c^{5})$ the infimum in (\ref{aap}) and (\ref{distR3XS2}) is actually a minimum and therefore we can choose $v \in SE(3)$ such that
\[
\tilde{W}(q,\tau)= \tilde{k}_{\tau}(v^{-1}q) + \tilde{U}(v).
\]
Then by Eq.~(\ref{estkern}) one has
\[
\begin{array}{ll}
 \tilde{W}(g,t) &\leq \tilde{k}_{t}(v^{-1}g)+\tilde{U}(v) \\
  & \leq \tilde{k}_{\tau}(v^{-1}q) + \tilde{k}_{t-\tau}(q^{-1}g) +\tilde{U}(v) \\
  & = \tilde{W}(q,\tau)+ \tilde{k}_{t-\tau}(q^{-1}g)
\end{array}
\]
for all $q \in SE(3)$ and thereby (by taking the infimum over all $q$) we obtain
\begin{equation}\label{onesideest}
k_{t-\tau}^{D^{11},D^{44},\eta,+} \ominus (k_{\tau}^{D^{11},D^{44},\eta,+} \ominus U) \geq  k_{t}^{D^{11},D^{44},\eta,+} \ominus U\ .
\end{equation}
So in order to prove (\ref{semigroup}) it remains to be shown that $k_{t-\tau}^{D^{11},D^{44},\eta,+} \ominus (k_{\tau}^{D^{11},D^{44},\eta,+} \ominus U) \leq  k_{t}^{D^{11},D^{44},\eta,+} \ominus U$. Let $w \in SE(3)$ such that
\[
\tilde{W}(g,t)=\min \limits_{v \in SE(3)} \tilde{k}_{t}(v^{-1}g)+\tilde{U}(v) = \tilde{k}_{t}(w^{-1}g)+\tilde{U}(w)\ ,
\]
now take $q$ as a point along the minimizing curve between $g$ and $w$ where the vertical $\R^{3}\rtimes S^{2}$ arc-length equals $p=\tau$.
Then we have
\[
\begin{array}{ll}
\tilde{W}(g,t)= \tilde{k}_{\tau}(w^{-1}q)+\tilde{k}_{t-\tau}(q^{-1}g)+\tilde{U}(w) &
\geq \tilde{k}_{t-\tau}(q^{-1}g) +\tilde{W}(q,\tau)
 \\&\geq
\min \limits_{r \in SE(3)}\{\tilde{k}_{t-\tau}(r^{-1}g) +\tilde{W}(r,\tau)\}
\end{array}
\]
from which the result follows. $\hfill \Box$

\begin{lemma}\label{lemmagen}
Let $\{A_{i}\}_{i=1}^{6}$ be a basis of the Lie-algebra of tangent vectors at the unity element and let $\{\mathcal{A}_{i}\}$ be the corresponding left-invariant vector fields and let $\{{\rm d}\mathcal{A}^{i}\}$ be the corresponding dual vector fields, e.g. Eq.~(\ref{duals2}). Then the exponential curves in $SE(3)$
\begin{equation} \label{gh}
\gamma(h)=(\ul{x}(h),R(h))=\textrm{exp}(h \sum \limits_{i=1}^{6} c^{i}A_{i})
\end{equation}
are auto-parallel w.r.t. Cartan connection. Their tangent vectors have constant components
w.r.t. the attached left-invariant moving frame of reference, i.e.
\begin{equation} \label{tuss}
\langle \left.{\rm d}\mathcal{A}^{i}\right|_{\gamma(h)}, \frac{d\gamma}{dh}(h)\rangle =c^{i}
\end{equation}
for $i=1,\ldots,6$ and all $h \in \R$. The spatial arclength of the spatial part $\ul{x}$ of the curve is given by
\[
s(h)=h\, \sqrt{\sum \limits_{i=1}^{3}|c^{i}|^2}.
\]
Let $c^{3}=c^{6}=0$, then the sub-Riemannian $(SE(3),{\rm d}\mathcal{A}^{3},{\rm d}\mathcal{A}^{6})$-arclength (
if $c^{1}=c^{2}$, and $c^{4}=c^{5}$ this corresponds to
$\R^{3}\rtimes S^{2}$-``erosion arclength'', Eq.\!~(\ref{arclengthR3XS2b}))
is given by
\[
p(h)= h \, \sqrt{\sum \limits_{i \in \{1,2,4,5\}}\frac{|c^i|^2}{D^{ii}}}
\]
Let $c^{1}=c^{2}=c^{6}=0$, then the sub-Riemannian $(SE(3),{\rm d}\mathcal{A}^{1}, {\rm d}\mathcal{A}^{2},{\rm d}\mathcal{A}^{6})$-arclength is given by
$
q(h)= h \, \sqrt{\sum \limits_{i \in \{3,4,5\}}\frac{|c^i|^2}{D^{ii}}}$.
\end{lemma}
\textbf{Proof }
By definition, $\gamma$ is given by Eq.~(\ref{gh}) the unique 1-parameter subgroup of $SE(3)$ whose tangent vector at the unity element
$e$ equals $\sum \limits_{i=1}^{6} c^{i}A_{i}$. It is thereby the integral curve of the corresponding left-invariant vector field given by $(L_{g})_{*}\sum \limits_{i=1}^{6} c^{i}A_{i}=\sum \limits_{i=1}^{6} c^{i}(L_{g})_{*}A_{i}= \sum \limits_{i=1}^{6} c^{i} \left. \mathcal{A}_{i}\right|_{g}$, obtained by push-forward of the left-multiplication, yielding
$
\dot{\gamma}(h)= \sum \limits_{i=1}^{6} c^{i} \left.\mathcal{A}_{i}\right|_{\gamma(h)}$,
from which Eq.~\!(\ref{tuss}) follows. In the Cartan connection on the tangent bundle $T(SE(3))$, the antisymmetric structure constants serve as Christoffel symbols, cf.~\!\cite{DuitsAMS2} and we have
\[
\nabla_{\dot{\gamma}} \dot{\gamma}=0 \desda \forall_{i=1,\ldots,6}\, :\,\ddot{\gamma}^{i} + \sum_{j,k} c_{jk}^{i}  \dot{\gamma}^{j} \dot{\gamma}^{k}= \ddot{\gamma}^{i}=0 \desda \forall_{i=1,\ldots,6}\, :\,\dot{\gamma}^{i}:=\langle\left.{\rm d}\mathcal{A}^{i}\right|_{\gamma}, \dot{\gamma} \rangle=c^{i}.
\]
For spatial arclength we have {\small $\|\dot{\ul{x}}(s)\|=1 \desda \|\ul{x}'(h)\|^2=\sum \limits_{i \in \{1,2,3\}}(c^{1})^2$}.
Regarding arclength in the sub-Riemannian manifold $(SE(3),{\rm d}\mathcal{A}^{3},{\rm d}\mathcal{A}^{6})$
we have by Eq.~\!(\ref{tuss}) that along a horizontal exponential curve \\
{\small $
p'(h)= \sqrt{\sum \limits_{i \in \{1,2,4,5\}}\frac{|\langle \left.{\rm d}\mathcal{A}^{i}\right|_{\gamma(h)},\dot{\gamma}(h) \rangle|^2}{D^{ii}}}=
 \sqrt{\sum \limits_{i \in \{1,2,4,5\}}\frac{|c^i|^2}{D^{ii}}}
$}
and the other sub-Riemannian case can be derived analogously. $\hfill \Box$ \\
\\

Let $k_{h}^{D^{11},D^{33},\eta, \pm}$ respectively denote the viscosity solution of
Eq.~(\ref{Hamdi1}), Eq.~(\ref{Hamdi2}), with initial condition $\pm\delta_{C}$, $\eta>\frac{1}{2}, D^{11},D^{33}>0$.
This kernel describes the propagation of ``balls'' in $\R^{3}\rtimes S^{2}=SE(3)/(\{0\}\times SO(2))$ centered around unity element $e=(\ul{0},\ul{e}_{z})\equiv \{(\ul{0},R_{\ul{e}_{z},\alpha})\;|\; \alpha \in [0,2\pi)\} \in \R^{3}\rtimes S^{2}$. Next we derive a few estimates for this morphological kernel in order to
\begin{enumerate}
\item motivate our analytic approximations in our erosion algorithms described in Section \ref{ch:HJE},
\item serve as ingredients in Theorem \ref{th:thetheorem}, which is the main result of this section.
\end{enumerate}
\begin{lemma} \label{lemma:1}
Let $\eta>\frac{1}{2}$, $D^{11},D^{33}>0$, $h>0$.
The morphological kernel $k_{h}^{D^{11},D^{33},\eta, \pm}$ satisfies the following estimate
\begin{equation}\label{estimatelemma1}
\begin{array}{ll}
\frac{2\eta-1}{2\eta} C^{- \frac{2\eta}{2\eta-1}}|\sum \limits_{i \in \{1,2,4,5\}} c^{i}A_{i}|^{\frac{2\eta}{2\eta-1}} &\leq \frac{k_{h}^{D^{11},D^{44},\eta, +} \left(\left[\exp\left(h \sum \limits_{i \in \{1,2,4,5\}} c^{i}A_{i}\right)\right]\right)}{h} \\ &\leq \frac{2\eta-1}{2\eta}\; C^{\frac{2\eta}{2\eta-1}} \; |\sum \limits_{i \in \{1,2,4,5\}} c^{i}A_{i}|^{\frac{2\eta}{2\eta-1}}
\end{array}
\end{equation}
with $A_{i}=\left.\mathcal{A}_{i}\right|_{e} \in T_{e}(SE(3))$, $c^{i} \in \R$, $C\geq 1$ and $|\sum \limits_{i \in \{1,2,4,5\}} c^{i}A_{i}|= \sqrt{\sum \limits_{i \in \{1,2,4,5\}}\frac{1}{D^{ii}}|c^{i}|^{2}}$ with $D^{11}=D^{22}$ and $D^{44}=D^{55}$ and where {\small $\left[\exp\left(h \sum \limits_{i \in \{1,2,4,5\}} c^{i}A_{i}\right)\right]$} denotes the left coset in $\R^{3}\rtimes S^{2}=SE(3)/(\{0\}\times SO(2))$ associated to the group element {\small $\exp\left(h \sum \limits_{i \in \{1,2,4,5\}} c^{i}A_{i}\right) \in SE(3)$}.
\end{lemma}
\textbf{Proof } We first derive an explicit upper-bound. As the smooth exponential curve
\[
[0,h]\ni \tilde{h} \mapsto \left[\exp\left(\tilde{h} \sum \limits_{i \in \{1,2,4,5\}} c^{i}A_{i}\right)\right] \in \R^{3}\rtimes S^2
 \]
is horizontal and connects $[e]$ and {\small $\left[\exp\left(h \sum \limits_{i \in \{1,2,4,5\}} c^{i}A_{i}\right)\right]$} we have by the definition
of $k_{h}^{D^{11},D^{33},\eta, +}$ in Eq.\!~(\ref{distR3XS2}) that
\[
\frac{k_{h}^{D^{11},D^{44},\eta, +} \left(\left[\exp\left(h \sum \limits_{i \in \{1,2,4,5\}} c^{i}A_{i}\right)\right]\right)}{h} \\ \leq \frac{2\eta-1}{2\eta}\; \tilde{C}^{\frac{2\eta}{2\eta-1}} \; |\sum \limits_{i \in \{1,2,4,5\}} c^{i}A_{i}|^{\frac{2\eta}{2\eta-1}}
\]
with {\small $\tilde{C}^{2\eta}=\frac{\sqrt{\sum \limits_{i \in \{1,2,4,5\}}\frac{|c^i|^2}{D^{ii}}}}{\sqrt{\sum \limits_{i=1}^{3}(c^{1})^2}}$}
where we applied Lemma \ref{lemma:3} and used $\frac{dp}{ds}=\frac{p'(h)}{s'(h)}$ along the smooth horizontal exponential curve.
The lower bound follows by the theory of weighted subcoercive operators on Lie groups, \cite[ch:1,6]{TerElst3} where one should consider the special case $\mathcal{U}=\mathcal{R}$, $G=SE(3)$ and algebraic basis $\{\mathcal{A}_{1},\mathcal{A}_{2},\mathcal{A}_{4},\mathcal{A}_{5}\}$ as erosion/dilation only takes place in the directions $\mathcal{A}_{1}$, $\mathcal{A}_{2}$, $\mathcal{A}_{4}$ and $\mathcal{A}_{5}$ (i.e. $\mathbf{G}$-orthogonal to $\mathcal{A}_{3}$). As a result we must assign the following weights to the Lie-algebra elements $w_{1}=w_{2}=w_{5}=w_{4}=1$ and $w_{3}=w_{6}=2$.
$\hfill \Box$ \\
\\
 Clearly, not every element $(\ul{y},R_{\ul{e}_{z},\gamma}R_{\ul{e}_{y},\beta}) \in SE(3)$ is reached by an exponential curve of the type $h \mapsto e^{h \sum_{i\in \{1,2,4,5\}} c^{i}A_{i}}$.
 In fact, by the Campbell-Baker-Hausdorff formula and the commutator table (\ref{tabel}) one has
\begin{equation} \label{hintweights}
\begin{array}{ll}
(h (y^{1},y^{2},0), R_{\ul{e}_{z},\gamma h}R_{\ul{e}_{z},\beta h}) &= (h\, (y^{1},y^{2},0),I)(\ul{0},R_{\ul{e}_{z},\gamma h}R_{\ul{e}_{z},\beta h}) \\ &=
\exp\{h (y^{1}A_{1}+y^{2}A_{2})\} \exp\{h\gamma A_{4}\} \exp\{h\beta A_{5}\} \\
 &=\exp\{h (y^{1}A_{1}+y^{2}A_{2}) \} \exp\{h \gamma A_{4} +\beta h A_{5} + \frac{1}{2} \beta h^{2}[A_{4},A_{5}]+O(h^{3})\} \\
 &=\exp\{h (y^{1}A_{1}\!+\!y^{2}A_{2}) +h \gamma A_{4} +\beta h A_{5} + \frac{1}{2} h^{2} (y^{1}\beta \!-\!y^{2}\gamma)A_{3} + \frac{1}{2}h^{2}(\beta \gamma)A_{6} +O(h^{3}) \} \\
 &= \exp\{\sum \limits_{i=1}^{6} c^{i} A_{i} h^{w_{i}}+O(h^3)\}
\end{array}
\end{equation}
and consequently, we can extend the estimate in Lemma \ref{lemma:1};
\begin{lemma} \label{lemma:2}
Let $\Omega_{e}$ be a compact set around the unity element
$e=(\ul{0},\ul{e}_{z})\equiv [(\ul{0},I)]=\{(\ul{0},R_{\ul{e}_{z},\alpha})\; |\; \alpha \in [0,2\pi)\}$ within $\R^{3}\rtimes S^2$, then
there exists an $\epsilon>0$ and $C>1$ such that for all $h<\epsilon$ one has the following uniform estimate
on compact sets
\[
 \begin{array}{ll}
 \frac{2\eta-1}{2\eta}\, h^{\frac{-2\eta}{2\eta-1}}\,C^{-\frac{2\eta}{2\eta-1}} \left( \sum \limits_{i=1}^{6} (c^{i})^{\frac{2}{w_{i}}} \right)^{\frac{\eta}{2\eta-1}}
 &\leq \frac{k_{h}((\ul{y},\ul{n}(\beta,\gamma)))}{h}=\frac{k_{h}([\exp\{\sum \limits_{i=1}^{6}c^{i}A_{i}h^{w_{i}}\}])}{h} \\
 &
 \leq h^{\frac{-2\eta}{2\eta-1}}\,\frac{2\eta-1}{2\eta} C^{\frac{2\eta}{2\eta-1}} \left( \sum \limits_{i=1}^{6} (c^{i})^{\frac{2}{w_{i}}} \right)^{\frac{\eta}{2\eta-1}}
 \end{array}
\]
for all $(\ul{y},\ul{n}(\beta,\gamma)\equiv [\exp\{\sum \limits_{i=1}^{6}c^{i}A_{i}h^{w_{i}}\}] \in \Omega_{e}$,
with $w_{1}=w_{2}=w_{5}=w_{4}=1$ and $w_{3}=w_{6}=2$ and $c^{i}:=c^{i}(\ul{y},\alpha=0,\beta,\gamma)=\tilde{c}^{i}(\ul{y},\tilde{\alpha}=0,\tilde{\beta},\tilde{\gamma})$, given by Eq.~(\ref{log2ndchart}) and (\ref{logSE3}).
\end{lemma}
\textbf{Proof }
Regarding the assignment of the weights to the Lie-algebra elements we refer to Eq.~(\ref{hintweights}). Regarding the logarithmic mapping in $SE(3)$ we recall Section \ref{ch:exp}.
Akin to the metric (\ref{metrichor}), the formula (\ref{distR3XS2}) is well-defined on the partition $\R^{3}\rtimes S^{2}:=SE(3)/(\{0\}\times SO(2))$ of left cosets. However as we have the restriction $\langle {\rm d}\mathcal{A}^{6}, \dot{\gamma} \rangle= \langle {\rm d}\tilde{\alpha}, \dot{\gamma}\rangle=\langle {\rm d}\alpha, \dot{\gamma}\rangle=0$ we must apply a consistent cross-section in $\R^{3}\rtimes S^{2}$. For convenience, we take the unique
element from the left cosets with $\tilde{\alpha}=0$ in the Euler-angle parametrization in both the initial point $(\ul{0},R_{\ul{e}_{z}}=I)$ and endpoint
$(\ul{y},R_{\ul{n}}=R_{\ul{e}_{x},\tilde{\gamma}}R_{\ul{e}_{y},\tilde{\beta}})$ of the curve. For the rest the proof is the same as the proof of Lemma \ref{lemma:1}. 
$\hfill \Box$.
\begin{corollary} \label{corr:kernelapprox}
Let $\eta>\frac{1}{2}$, $D^{11}>0, D^{44}>0$. $w_{1}=w_{2}=w_{4}=w_{5}=1$, $w_{3}=w_{6}=2$.
For the morphological erosion (+) and dilation kernel (-) on $\R^{3}\rtimes S^{2}$ one can use the following asymptotical formula
\begin{equation}
k_{t}^{D^{11},D^{44}, \eta, \pm}(\ul{y},\tilde{\ul{n}}(\tilde{\beta},\tilde{\gamma})) \equiv \pm \, \frac{2\eta-1}{2\eta} C^{\frac{2\eta}{2\eta-1}} t^{-\frac{1}{2\eta-1}} \left(\sum \limits_{i=1}^{6}
\frac{|\tilde{c}^{i}(\ul{y},\tilde{\alpha}=0,\tilde{\beta},\tilde{\gamma})|^{\frac{2}{w_{i}}}}{D^{ii}} \right)^{\frac{\eta}{2\eta-1}}
\end{equation}
for sufficiently small time $t>0$.
\end{corollary}
\begin{lemma}\label{lemma:3}
Let $c^{i} \in \R$, $i=1,2,4,5$ and $\eta>\frac{1}{2}$ then
\[
\lim \limits_{h \downarrow 0}\frac{k^{D^{11},D^{44},\eta,+}_{h}\left(\left[\exp \left(h\sum \limits_{i \in \{1,2,4,5\}} c^{i}A_{i}\right)\right]\right)}{h} = \frac{2\eta-1}{2\eta} \; |\sum \limits_{i \in \{1,2,4,5\}} \frac{c^{i}}{\sqrt{D^{ii}}}A_{i}|^{\frac{2\eta}{2\eta-1}}= \frac{2\eta-1}{2\eta} \; \left|\sum \limits_{i \in \{1,2,4,5\}} \frac{|c^{i}|^2}{D^{ii}}\right|^{\frac{\eta}{2\eta-1}}\! \!.
\]
\end{lemma}
\textbf{Proof } Consider Eq.~(\ref{distR3XS2}), then to each $C^{2}$-curve with
$\langle \left. {\rm d}\mathcal{A}^{3}\right|_{\gamma}, \dot{\gamma} \rangle
=\langle \left. {\rm d}\mathcal{A}^{6}\right|_{\gamma}, \dot{\gamma} \rangle=0$, we have the following
$C^{1}$-functions
\[
[0,h] \ni p \mapsto \dot{\gamma}^{i}(p):=\langle \left. {\rm d}\mathcal{A}^{i}\right|_{\gamma(p)}, \dot{\gamma}(p) \rangle \in \R
\]
for $i\in \{1,2,4,5\}$. By applying a first order Taylor-approximation around $s=0$ we obtain
\[
\begin{array}{l}
\lim \limits_{h \downarrow 0} h^{-1}\int \limits_{0}^{h}
\mathcal{L}_{\eta}(\frac{1}{\sqrt{D^{11}}}\dot{\gamma}^{1}(p),\frac{1}{\sqrt{D^{11}}}\dot{\gamma}^{2}(p),
\frac{1}{\sqrt{D^{44}}}\dot{\gamma}^{4}(p),\frac{1}{\sqrt{D^{44}}}\dot{\gamma}^{5}(p) )\; \left( \frac{dp}{ds}(p)\right)^{\frac{1}{2\eta-1}}\;{\rm d}p
\\=
\lim \limits_{h \downarrow 0} h^{-1}\int \limits_{0}^{\psi^{-1}(h)}
\mathcal{L}_{\eta}(\frac{1}{\sqrt{D^{11}}}\dot{\gamma}^{1}(0)+O(s),\frac{1}{\sqrt{D^{11}}}\dot{\gamma}^{2}(0)+O(s),
\frac{1}{\sqrt{D^{44}}}\dot{\gamma}^{4}(0)+O(s),\frac{1}{\sqrt{D^{44}}}\dot{\gamma}^{5}(0)+O(s) )\; {\rm d}s \\
\end{array}
\]
and $\eta>\frac{1}{2}$ implies that this is of the order
\[
 O((\psi^{-1}(h))^2 h^{-1}) + \frac{2\eta-1}{2\eta} \left( \sum \limits_{i \in \{1,2,4,5\}} \frac{1}{D^{ii}}
(\dot{\gamma}^{i}(0))^{2}\right)^{\frac{\eta}{2\eta-1}}
\]
where $\psi^{-1}(0)=0$ and along a minimizing curve we have finite curvature (as we will see in Appendix
\ref{app:E}) so $\psi^{-1}(h)=O(h)$ and the first term vanishes as $h \downarrow 0$ and the result follows by definition of the morphological kernel, Eq.\!~(\ref{distR3XS2}).
$\hfill \Box$
\begin{theorem}\label{th:thetheorem}
The viscosity solutions of
\begin{equation}\label{Hamdi1}
\left\{
\begin{array}{l}
\frac{\partial W}{\partial t}(\ul{y},\ul{n},t) \mp \frac{1}{2\eta}\left(\ul{G}_{(\ul{y},\ul{n})}^{-1}\left(\left.{\rm d} W(\cdot,\cdot,t) \right|_{\ul{y},\ul{n}},\left. {\rm d} W(\cdot,\cdot,t) \right|_{\ul{y},\ul{n}}\right)\right)^{\eta}=0 \\
W(\ul{y},\ul{n},0)= U(\ul{y},\ul{n})
\end{array}
\right.
\end{equation}
with $\ul{G}_{(\ul{y},\ul{n})}= g_{ii}\, \left.{\rm d}\mathcal{A}^{i}\right|_{\ul{y},\ul{n}} \otimes \left.{\rm d}\mathcal{A}^{i}\right|_{\ul{y},\ul{n}}$ with $g_{ii}=\frac{1}{D^{ii}}$ with $D^{11}=D^{22}, D^{44}=D^{55}> 0$
are respectively given by ($+$ case) left-invariant erosion
\begin{equation}\label{Erosion2}
(k_{t}^{D^{11}, D^{44}, \eta, +} \ominus_{\R^{3} \rtimes S^2}U)(\ul{y},\ul{n})=
\inf \limits_{(\ul{y}',\ul{n}')\in \R^{3} \rtimes S^{2}} \left[U(\ul{y}',\ul{n}')+ k_{t}^{D^{11}, D^{44}, \eta, +}(R_{\ul{n}'}^{T}(\ul{y}-\ul{y}'),R_{\ul{n}'}^{T}\ul{n})\right]\ .
\end{equation}
and ($-$ case) left-invariant dilation
\begin{equation} \label{dilation2}
(k_{t}^{D^{11}, D^{44}, \eta, -} \oplus_{\R^{3} \rtimes S^2}U)(\ul{y},\ul{n})=
\sup \limits_{(\ul{y}',\ul{n}')\in \R^{3} \rtimes S^{2}} \left[k_{t}^{D^{11}, D^{44}, \eta, -}(R_{\ul{n}'}^{T}(\ul{y}-\ul{y}'),R_{\ul{n}'}^{T}\ul{n})+U(\ul{y}',\ul{n}')\right]
\end{equation}
\end{theorem}
\textbf{Proof }The proof consists of two parts, first we must show that they are indeed solutions and then we show that they are viscosity
solutions. We will only consider the erosion case since the dilation case can be treated analogously.

\emph{part I}. This part consists of two subparts. In part Ia we will show that if we set
\[
W(\ul{y},\ul{n},t):=(k_{t}^{D^{11}, D^{44}, \eta, +} \ominus_{\R^{3} \rtimes S^2}U)(\ul{y},\ul{n})
\]
with
morphological kernel given by (\ref{distR3XS2})
then
\begin{equation}\label{1a}
\begin{array}{l}
\frac{\partial W}{\partial t}(\ul{y},\ul{n},t) \mp \frac{1}{2\eta}\left(\ul{G}_{(\ul{y},\ul{n})}^{-1}\left(\left.{\rm d} W(\cdot,\cdot,t) \right|_{\ul{y},\ul{n}},\left. {\rm d} W(\cdot,\cdot,t) \right|_{\ul{y},\ul{n}}\right)\right)^{\eta} \leq 0\ .
\end{array}
\end{equation}
Subsequently, in part Ib we show that
\begin{equation}\label{1b}
\begin{array}{l}
\frac{\partial W}{\partial t}(\ul{y},\ul{n},t) \mp \frac{1}{2\eta}\left(\ul{G}_{(\ul{y},\ul{n})}^{-1}\left(\left.{\rm d} W(\cdot,\cdot,t) \right|_{\ul{y},\ul{n}},\left. {\rm d} W(\cdot,\cdot,t) \right|_{\ul{y},\ul{n}}\right)\right)^{\eta} \geq 0\ .
\end{array}
\end{equation}
so that part I is finished.
\begin{itemize}
\item[\emph{part Ia}] Again we resort to evolutions on the full group $SE(3)$ by setting
$\tilde{W}(\ul{y},R,t)=W(\ul{y},R\ul{e}_{z},t)$ for all $\ul{y}\in \R^{3}$, $R \in SO(3)$, $t>0$. Then for all $A=\sum \limits_{i \in \{1,2,4,5\}}c^{i}A_{i}$
one has by Lemma \ref{lemma:semigroup} that
\[
\begin{array}{ll}
\tilde{W}(g e^{hA},t+h) &= \min \limits_{v \in SE(3)}\tilde{k}_{h}^{D^{11},D^{44},\eta,+}(v^{-1}g e^{h A})+\tilde{W}(v,t) \\
 & \leq \tilde{k}_{h}^{D^{11},D^{44},\eta,+}(e^{hA})+\tilde{W}(g,t) \\
 \end{array}
\]
for all $g=(\ul{y},R)$
and consequently, one has
\[
\frac{\tilde{W}(ge^{hA},t+h)-\tilde{W}(g,t)}{h} \leq \frac{\tilde{k}_{h}^{D^{11},D^{44},\eta,+}(e^{h A})}{h}
\]
and thereby, by Lemma \ref{lemma:3} and the construction of the left-invariant vector fields
by means of the derivative of the right regular representation
$\mathcal{A}_{i}={\rm d}\mathcal{R}(A_{i})$ one gets
\begin{equation}
\frac{\partial \tilde{W}(g,t)}{\partial t}+ c^{i} \left.\mathcal{A}_{i}\right|_{g} \tilde{W}(g,t) \leq \mathcal{L}_{\eta}(\ul{c})
\end{equation}
for all $\ul{c}=(c^{1},c^{2},c^{4},c^{5}) \in \R^{4}$. So by subtracting the Lagrangian and taking the maximum over all $\ul{c}$ we apply the
Fenchel transform of the Lagrangian which yields the Hamiltonian, i.e.
\[
\begin{array}{l}
\frac{\partial \tilde{W}(g,t)}{\partial t}+ \sup \limits_{\ul{c}\in \R^{4}}\left\{
\sum \limits_{i \in \{1,2,4,5\}}
c^{i}\left.\mathcal{A}_{i}\right|_{g} \tilde{W}(g,t) - \mathcal{L}_{\eta}(\ul{c})\right\}\leq 0 \ \desda \\
\frac{\partial \tilde{W}(g,t)}{\partial t}+\frac{1}{2\eta}\left\{ \sum \limits_{i \in \{1,2,4,5\}}
\frac{|\left.\mathcal{A}_{i}\right|_{g} \tilde{W}(g,t)|^2}{D^{ii}}\right\}^{2\eta} \leq 0 .
\end{array}
\]
\item[\emph{part Ib}] Let $g^{*}=(\ul{y}^{*},R_{\ul{n}^{*}}) \in SE(3)$ be the minimizer in the erosion operator, i.e.
\[
\tilde{W}(g,t)=\tilde{W}(g^{*},0)+
\tilde{k}_{t}^{D^{11},D^{44},\eta,+}((g^{*})^{-1} g).
\]
Then we have
\begin{equation}
\begin{array}{ll}
\frac{\tilde{W}(g,t)-\tilde{W}(ge^{-\frac{h}{t}A}e^{+\frac{h}{t}A^*},t-h)}{h} & \geq
\frac{1}{h}\left(\tilde{W}(g,0)+ \tilde{k}_{t}^{D^{11},D^{44},\eta,+}((g^{*})^{-1}g) \right. \\ & \ \ \ \ \left. - \left(\tilde{W}(g,0)+ \tilde{k}_{t}^{D^{11},D^{44},\eta,+}
((g^{*})^{-1}g\,e^{-\frac{h}{t}A})e^{+\frac{h}{t}A^{*}} \right)\right)\\
 &=\frac{\tilde{k}_{t}^{D^{11},D^{44},\eta,+}((g^{*})^{-1}g)-\tilde{k}_{t-h}^{D^{11},D^{44},\eta,+}((g^{*})^{-1}
 g\,e^{-\frac{h}{t}A}e^{+\frac{h}{t}A^{*}})}{h} \\
 \end{array}
\end{equation}
where $g=e^{A}=e^{\sum \limits_{i=1}^{6}c^{i}A_{i}}$ and $g^{*}=e^{A^{*}}=e^{\sum \limits_{i=1}^{6}(c^*)^{i}A_{i}} \in SE(3)$. If we now let $h \downarrow 0$ then we obtain
\begin{equation} \label{limitingcase}
\begin{array}{cc}
\lim \limits_{h \downarrow 0}\frac{\tilde{W}(g,t)-\tilde{W}(ge^{-\frac{h}{t}(A-A^{*})},t-h)}{h} &\geq \lim \limits_{h \downarrow 0}
 \frac{\tilde{k}_{t}^{D^{11},D^{44},\eta,+}((g^{*})^{-1}g)-\tilde{k}_{t-h}^{D^{11},D^{44},\eta,+}((g^{*})^{-1}
 g\,e^{-\frac{h}{t}(A-A^{*})})}{h} \\
  & = \lim \limits_{h \downarrow 0}\frac{1}{h}\left(
 \sup \limits_{v \in SE(3)} \tilde{k}^{D^{11},D^{44}}_{t-h}(v)+ \tilde{k}_{h}^{D^{11},D^{44},\eta,+}(v^{-1}(g^{*})^{-1}g )\right.\\ & \left.\ \ -
 \tilde{k}_{t-h}^{D^{11},D^{44},\eta,+}((g^{*})^{-1}g e^{-\frac{h}{t}(A-A^{*})}) \right) \\
 & \geq \lim \limits_{h \downarrow 0}\frac{\tilde{k}_{h}^{D^{11},D^{44},\eta,+}(e^{\frac{h}{t}(A-A^{*})})}{h} 
\end{array}
\end{equation}
where we applied Lemma \ref{lemma:semigroup} and
where we note that by the Campbell-Baker-Hausdorff formula
\[
e^{-\frac{h}{t}A}\, e^{+\frac{h}{t}A^{*}}= e^{-\frac{h}{t}(A-A^{*})+ \frac{h^2}{2t^2}[A,A^{*}]+O(h^3)}
\]
so that
\[
\tilde{W}(g\, e^{-\frac{h}{t}(A-A^{*})},t-h)=\tilde{W}(g\, e^{-\frac{h}{t}A}\, e^{+\frac{h}{t}A^{*}},t-h)+O(h^2)
\]
for almost every
$g \in \R^{3}\rtimes S^{2}$ (where $\tilde{W}(\cdot,t-h)$ is differentiable).

Recall that the left-invariant vector
fields are obtained by the derivative of the right-regular representation
$\mathcal{A}_{i}={\rm d}\mathcal{R}(A_{i})$, $A_{i}=\left.\mathcal{A}_{i}\right|_{(0,I)}$ and as a result the limit in the left-hand side
of (\ref{limitingcase}) equals
\[
\lim \limits_{h \downarrow 0}\frac{\tilde{W}(g,t)-\tilde{W}(ge^{-\frac{h}{t}(A-A^{*})},t-h)}{h}= \frac{\partial }{\partial t}\tilde{W}(g,t)- \sum \limits_{i\in \{1,2,4,5\}}
\left(\frac{c^{i}-(c^{*})^{i}}{t} \right)\left.\mathcal{A}_{i}\right|_{g}\tilde{W}(g,t).
\]
Furthermore, we can assume without loss of generality that $c^{3}=c^{6}=(c^{*})^{3}=(c^{*})^{6}=0$, i.e.
\begin{equation} \label{assu}
A-A^{*} \in \textrm{span}\{A_{1},A_{2},A_{4},A_{5}\}.
\end{equation}
To this end we recall both $\R^{3}\rtimes S^{2}= SE(3)/(\{0\}\times SO(2))$ and $\tilde{W}(g,t)=W(\ul{y},R \ul{e}_{z},t)$, $g \in SE(3)$ so that
\[
\tilde{W}(g,t)=\tilde{W}(g e^{\alpha A_{6}},t) \textrm{ and } \tilde{k}_{t}^{D^{11},D^{44},\eta,+}(e^{-\beta A_{6}} g e^{\alpha A_{6}})=\tilde{k}_{t}^{D^{11},D^{44},\eta,+}(g )
\]
for all $\alpha, \beta \in \R$ and by the Campbell-Baker-Hausdorff formula one can always find $\alpha, \beta \in \R$ such that
\[
\log \left( e^{A}e^{\alpha A_{6}}\right)-\log \left( e^{A^{*}}e^{\beta A_{6}}\right) \subset \textrm{span}\{A_{1},A_{2},A_4,A_5\}.
\]

Assumption (\ref{assu}) allows us to compute the limit in the final right-hand side of inequality (\ref{limitingcase}) by means of Lemma \ref{lemma:3} and we obtain
\[
\begin{array}{c}
\frac{\partial }{\partial t}\tilde{W}(g,t)+ \sum \limits_{i\in \{1,2,4,5\}}
\left(\frac{c^{i}-(c^{*})^{i}}{t} \right)\left.\mathcal{A}_{i}\right|_{g}\tilde{W}(g,t) \\
 \geq  \lim \limits_{h \downarrow 0}\frac{\tilde{k}_{h}(e^{\frac{h}{t}(A-A^{*})})}{h}=
\mathcal{L}_{\eta}\left( \frac{c^{1}-(c^{*})^{1}}{t}, \frac{c^{2}-(c^{*})^{2}}{t}, \frac{c^{4}-(c^{*})^{4}}{t}),\frac{c^{5}-(c^{*})^{5}}{t}\right) \ .
\end{array}
\]
from which we conclude :
\[
\begin{array}{l}
\frac{\partial \tilde{W}(g,t)}{\partial t}+ \frac{1}{2\eta}\left( \sum \limits_{i \in \{1,2,4,5\}}\frac{(\mathcal{A}_{i}\tilde{W}(g,t))^2}{D^{ii}}\right)^{2\eta}
= \frac{\partial \tilde{W}(g,t)}{\partial t} + \sup \limits_{\tilde{\ul{c}}\in \R^{4}}\left(\sum \limits_{i \in \{1,2,4,5\}} \tilde{c}^{i}\left.\mathcal{A}_{i}\right|_{g}\tilde{W}(g,t)-\mathcal{L}_{\eta}(\tilde{c})\right) \\
\geq

\frac{\partial }{\partial t}\tilde{W}(g,t)+ \sum \limits_{i\in \{1,2,4,5\}}
\left(\frac{c^{i}-(c^{*})^{i}}{t} \right)\left.\mathcal{A}_{i}\right|_{g}\tilde{W}(g,t)- \mathcal{L}_{\eta}\left( \frac{c^{1}-(c^{*})^{1}}{t}, \frac{c^{2}-(c^{*})^{2}}{t}, \frac{c^{4}-(c^{*})^{4}}{t}),\frac{c^{5}-(c^{*})^{5}}{t}\right)
\end{array}
\]
from which the result (\ref{1b}) follows.
\end{itemize}
\emph{part II}.
Next we verify that erosion (\ref{dilation2}) with the Green's function $k_{t}^{D^{11},D^{44},\eta,+}$ indeed satisfies (\ref{visc1}) at a location $(\ul{y}_0,\ul{n}_{0}, t_{0}) \in (\R^{3}\rtimes S^{2}) \times R^{+}$
where  $W-V$ attains a local maximum. The other case (\ref{visc2}) can be shown analogously. Similarly one can show that dilation (\ref{Erosion}) with the kernel $k_{t}^{D^{11}, D^{44}, \eta, -}=-k_{t}^{D^{11}, D^{44}, \eta, +}$ is the viscosity solution of Eq.~(\ref{Hamdi2}).

First of all it follows by left-invariance (and the fact that $SE(3)$ acts transitively on $\R^{3}\rtimes S^2$ by means of Eq.~(\ref{actleft})) that without loss of generality we can restrict ourselves to the case $(\ul{y}_{0},\ul{n}_{0})=(\ul{0},\ul{e}_{z})$ and furthermore, by the semigroup property Lemma \ref{lemma:semigroup}, we can, again without loss of generality, restrict ourselves to the case $t_0=0$.
Since $W-V$ attains a maximum in $(\ul{0},\ul{e}_{z},0)$, there exists a small open set $\Omega$ around $(\ul{0},\ul{e}_{z},0)$ in $(\R^{3}\rtimes S^{2}) \times \R$ where
\begin{equation}\label{onelocal}
W(\ul{y},\ul{n},t)-V(\ul{y},\ul{n},t) \leq W(\ul{y},\ul{n},0)-V(\ul{y},\ul{n},0)
\end{equation}
Furthermore we have by Eq.~(\ref{Erosion2}) that
\begin{equation} \label{twomax}
W(\ul{y},\ul{n},t) \geq W(\ul{0},\ul{e}_{z},0) + k_{t}^{D^{11}, D^{44}, \eta, +}(\ul{y},\ul{n})\ ,
\end{equation}
for all $(\ul{y},\ul{n})\in \R^{3}\rtimes S^2$. Combining the estimates (\ref{onelocal}) with (\ref{twomax}) yields
\begin{equation}
-V(\ul{0},\ul{e}_{z},0)+V(\ul{y},\ul{n},t) \leq W(\ul{y},\ul{n},t) -W(\ul{0},\ul{e}_{z},0) \leq k_{t}^{D^{11}, D^{22}, \eta, +}(\ul{y},\ul{n})
\end{equation}
locally around $(\ul{0},\ul{e}_{z},0)$, that is within $\Omega$. Or equivalently for the corresponding function on $SE(3)$ given by $\tilde{V}(\ul{x},R)=V(\ul{x},R\ul{e}_{z})$ we have
\begin{equation}
\frac{-\tilde{V}(\ul{0},I,0)+\tilde{V}(\ul{y},R_{\ul{n}},t)}{t} \leq \frac{k_{t}^{D^{11}, D^{22}, \eta, +}(\ul{y},\ul{n})}{t}
\end{equation}
for $t \in [0,\epsilon_{1}), d_{SE(3)}(\ul{y},R_{\ul{n}},(0,I))<\epsilon_{2}$, for some $\epsilon_{1}>0$, $\epsilon_{2}>0$. Now set $h=t$, $(\ul{y},R_{\ul{n}})=e^{h A}$ with $A=\sum \limits_{i \in \{1,2,4,5\}}c^{i}A_{i}$ and take the limit $h \downarrow 0$ then
we have
\[
\frac{d \tilde{V}}{dt}(\ul{0},I,0)+ \sum \limits_{i \in \{1,2,4,5\}} c^{i} {\rm d}\mathcal{R}(A_{i}) \tilde{V}(\ul{0},I,0)
\leq \frac{2\eta-1}{2\eta}\; 
\; |\sum \limits_{i \in \{1,2,4,5\}} c^{i}A_{i}|^{\frac{2\eta}{2\eta-1}}
\]
where we applied Lemma \ref{lemma:3} (which tells us that the constant $C$ in Lemma \ref{lemma:1} can be set to $1$ in the limiting case $h\downarrow 0$). Now the left-invariant vector fields are given by the derivative of the right regular representation
$\mathcal{A}_{i}={\rm d}\mathcal{R}(A_{i})$ with $\left.\mathcal{A}_{i}\right|_{(\ul{0},I)}=A_{i}$ so we see that
\begin{equation}
\begin{array}{ll}
\frac{d \tilde{V}}{dt}(\ul{0},I,0) & + \sup \limits_{\sum \limits_{i\in \{1,2,4,5\}}c^{i}\mathcal{A}_{i}}\{
c^{i} A_{i} \tilde{V}(\ul{0},I,0)-\frac{2\eta-1}{2\eta} |\sum \limits_{i \in \{1,2,4,5\}} c^{i}A_{i}|^{\frac{2\eta}{2\eta-1}}\} \leq 0 \\
\frac{d \tilde{V}}{dt}(\ul{0},I,0) & + (\gothic{F}_{T_{e}(SE(3))}(A \mapsto \frac{2\eta-1}{2\eta}
|A|^{\frac{2\eta}{2\eta-1}}))(
{\rm d}\tilde{V}(\cdot,0))  \leq 0
\end{array}
\end{equation}
with $A=\sum \limits_{i\in \{1,2,4,5\}}c^{i}A_{i} \in T_{e}(SE(3)$ and ${\rm d}\tilde{V}(\cdot,t)=\sum \limits_{i=1}^{5} A_{i}\tilde{V}(\cdot,t){\rm d}A^{i} \in (T_{e}(SE(3)))^{*}$. Now by the final remark in Definition \ref{def:LT} and
Lemma \ref{lemma:LF} we see that
\begin{equation}
\begin{array}{l}
\frac{d \tilde{V}}{dt}(\ul{0},I,0)  + \left(\gothic{F}_{\R^4}((c^{1},c^{2},c^{4},c^{5})\mapsto \frac{2\eta-1}{2\eta} 
\cdot \right.\\
\left.
\left(\sum_{i \in \{1,2,4,5\}} |c^{i}|^2\right){\frac{\eta}{2\eta-1}}
)\right)(A_{1}\tilde{V}(\ul{0},I,0),(A_{2}\tilde{V}(\ul{0},I,0),(A_{4}\tilde{V}(\ul{0},I,0),(A_{5}\tilde{V}(\ul{0},I,0) ) \leq 0 \\
\desda \frac{d \tilde{V}}{dt}(\ul{0},I,0) + \frac{1}{2\eta} 
\left((A_{1}\tilde{V}(\ul{0},I,0))^2+(A_{2}\tilde{V}(\ul{0},I,0))^2 +(A_{4}\tilde{V}(\ul{0},I,0))^2 +(A_{5}\tilde{V}(\ul{0},I,0))^2\right)^{\eta} \leq 0
\end{array}
\end{equation}
As a result we conclude that the erosions Eq.~(\ref{Erosion2}) with morphological kernel (\ref{distR3XS2}) are the (unique) viscosity solutions of
the Hamilton-Jacobi equations (\ref{Hamdi1}) on the contact manifold $(\R^{3}\rtimes S^{2},{\rm d}\mathcal{A}^{3})$. $\hfill \Box$
\section{The Hamilton-Jacobi equation on $(\R^{3}\rtimes S^2, {\rm d}\mathcal{A}^{3})$ and the propagation of geodesically equidistant surfaces in $\R^{3}\rtimes S^2$ \label{app:meaming}}
\vspace{-0.6 cm} \mbox{}
Recall from Figure \ref{fig:intuition} that our left-invariant erosions on $\R^{3}\rtimes S^2$ take place on the contact/sub-Riemannian manifold $(SE(3), {\rm d}\mathcal{A}^{3}, {\rm d}\mathcal{A}^{6})$.
On $(SE(3), {\rm d}\mathcal{A}^{3}, {\rm d}\mathcal{A}^{6})$ we set the
Lagrangian (for $\eta>\frac{1}{2}$)
\begin{equation}
\overline{\mathcal{L}}_{\eta}(\tilde{\gamma}(p),\dot{\tilde{\gamma}}(p))= \frac{2\eta-1}{2\eta} \left(\frac{|\dot{\tilde{\gamma}}^{1}(p)|^2+
|\dot{\tilde{\gamma}}^{2}(p)|^2}{D^{11}} + \frac{|\dot{\tilde{\gamma}}^{4}(p)|^2 + |\dot{\tilde{\gamma}}^{5}(p)|^2}{D^{44}}
\right)^{\frac{\eta}{2\eta-1}} \cdot \left(\frac{dp}{ds}\right)^{\frac{1}{2\eta-1}},
\end{equation}
expressed in the $(SE(3), {\rm d}\mathcal{A}^{3}, {\rm d}\mathcal{A}^{6})$ arc-length parameter $p$ given by
\begin{equation}\label{arclengthR3XS2b1}
\begin{array}{ll}
p(\tau) &=  \int \limits_{0}^{\tau} \sqrt{\sum \limits_{i \in \{1,2,4,5\}}\frac{1}{D^{ii}}|\langle \left. {\rm d}\mathcal{A}^{i}\right|_{\gamma(\tilde{\tau})}, \dot{\gamma}(\tilde{\tau})\rangle |^{2}  }\, {\rm d}\tilde{\tau}.
\end{array}
\end{equation}
yielding (when $\eta \to \infty$) the well-defined metric on $(\R^{3}\rtimes S^2, {\rm d}\mathcal{A}^{3})$ given by
\begin{equation}\label{EN2}
\begin{array}{l}
d((\ul{y},\ul{n}),(\ul{y}',\ul{n}'))=d(((R_{\ul{n}'})^T(\ul{y}-\ul{y}'),(R_{\ul{n}'})^T\ul{n}),(\ul{0},\ul{e}_{z}))\ , \\
d((\ul{y},\ul{n}),(\ul{0},\ul{e}_{z}))=
\inf \limits_{{\tiny \begin{array}{c}
\tilde{\gamma}=(\ul{x}(\cdot),R(\cdot)) \in C^{\infty}((0,1), SE(3)),\\
\tilde{\gamma}(0)=(\ul{0},R_{\ul{e}_z}), \tilde{\gamma}(p_{\textrm{max}})=(\ul{y},R_{\ul{n}}),\\
\langle \left.{d \rm}\mathcal{A}^{3}\right|_{\tilde{\gamma}}, \dot{\tilde{\gamma}}\rangle=\langle \left.{d \rm}\mathcal{A}^{6}\right|_{\tilde{\gamma}},\dot{\tilde{\gamma}}\rangle=0
\end{array}}} \int \limits_{0}^{p_{\textrm{max}}} \overline{\mathcal{L}}_{\infty}(\tilde{\gamma}(p),\dot{\tilde{\gamma}}(p))\, {\rm d}p.
\end{array}
\end{equation}
Let us return to the general Lagrangian minimization problem on $(\R^{3}\rtimes S^2, {\rm d}\mathcal{A}^{3})$ for $\eta > \frac{1}{2}$
\begin{equation}
E_{\eta}(\ul{y},\ul{n},t):=
\begin{array}{l}
\inf \limits_{{\tiny \begin{array}{c}
\tilde{\gamma}=(\ul{x}(\cdot),R_{\ul{n}}(\cdot)) \in C^{\infty}((0,p_{\textrm{max}}), SE(3)),\\
\tilde{\gamma}(0)=g_{0}, \tilde{\gamma}(t)=(\ul{y},R_{\ul{n}}),\\
\langle \left.{d \rm}\mathcal{A}^{3}\right|_{\tilde{\gamma}}, \dot{\tilde{\gamma}}\rangle=\langle \left.{d \rm}\mathcal{A}^{6}\right|_{\tilde{\gamma}},\dot{\tilde{\gamma}}\rangle=0
\end{array}}} \int \limits_{0}^{p_{\textrm{max}}=t} \overline{\mathcal{L}}_{\eta}(\tilde{\gamma}(p),\dot{\tilde{\gamma}}(p))\, \cdot \left(\frac{dp}{ds}\right)^{\frac{1}{2\eta-1}}\, {\rm d}p.
\end{array}
\end{equation}
for some given $g_{0} \in SE(3)$ and
where the total arc-length $p_{\textrm{max}}$ in $\R^{3}\rtimes S^2$ is fixed.\footnote{Otherwise for fixed $\eta$ there does not exist a minimizer.}, say $p_{\textrm{max}}=t>0$.
Along an optimizing curve $\tilde{\gamma}^{*}=(\ul{x}^{*},R^{*})$ (``geodesic'') in $SE(3)$ with corresponding curve $s\mapsto \gamma^{*}(s)=(\ul{x}^{*}(s), \ul{n}^{*}(s):=R^{*}(s)\ul{e}_{z}) \in \R^{3}\rtimes S^2$ we have
\[
\frac{\partial }{\partial t}E_{\eta}(\gamma^{*}(t),t)= \overline{\mathcal{L}}_{\eta}(\gamma^{*}(t),\dot{\gamma}^{*}(t))\cdot \left(\frac{dp}{ds}(\psi^{-1}(t))\right)^{\frac{1}{2\eta-1}},
\]
where we recall that $\psi(s)=p$ and $s=\psi^{-1}(p)$, Eq.~(\ref{psis}).
Now consider a family of surfaces associated to a smooth function $W : \R^{3} \rtimes S^{2} \times \R^{+} \to \R^{+}$ given by
\[
\mathcal{S}_{t}:=\{(\ul{y},\ul{n}) \in \R^{3}\rtimes S^2 \; |\; W(\ul{y},\ul{n},t)=W_{0}\}
\]
parameterized by $t\geq 0$, where $W_{0}>0$ is some positive constant. Such a family of surfaces is called \emph{geodesically equidistant}, \cite{Rund}, if
\begin{equation}\label{geodesicallyequidistant}
\begin{array}{l}
\frac{dW}{dt}(\gamma^{*}(t),t)=\left(\frac{ds}{dp}(t)\right)^{\frac{1}{2\eta-1}}\cdot \frac{\partial E_{\eta}}{\partial t}(\gamma^{*}(t),t)= \overline{\mathcal{L}}_{\eta}(\gamma^{*}(t),\dot{\gamma}^{*}(t))\ , \\[8pt]
\langle \left. {\rm d}W(\cdot,t)\right|_{\gamma^{*}(t)}, \dot{\gamma}^{*}(t) \rangle=\left.\left(\nabla_{\ul{c}}\overline{\mathcal{L}}_{\eta}(\gamma^{*}(t),\ul{c})\cdot \ul{c}\right) \right|_{\ul{c}=\dot{\gamma}^{*}(t)}.
\end{array}
\end{equation}
The chain-law now gives
\[ 
\begin{array}{l}
\frac{\partial W}{\partial t}(\gamma^{*}(t),t)+ \langle \left. {\rm d}W(\cdot,t)\right|_{\gamma^{*}(t)}, \dot{\gamma}^{*}(t) \rangle = \overline{\mathcal{L}}_{\eta}(\gamma^{*}(t),\dot{\gamma}^{*}(t))  \desda \\
-\frac{\partial W}{\partial t}(\gamma^{*}(t),t)= -\overline{\mathcal{L}}_{\eta}(\gamma^{*}(t),\dot{\gamma}^{*}(t))+ \langle \left. {\rm d}W(\cdot,t)\right|_{\gamma^{*}(t)}, \dot{\gamma}^{*}(t) \rangle
\end{array}
\] 
where according to the 2nd equality in (\ref{geodesicallyequidistant}) we may rewrite the righthand side as
\[
\sup \limits_{\ul{c}=\! \!\sum \limits_{i\in \{1,2,4,5\}}\! \!c^{i}\mathcal{A}_{i} } \left\{-\overline{\mathcal{L}}_{\eta}(\gamma^{*}(t),\ul{c})+ \langle \left. {\rm d}W(\cdot,t)\right|_{\gamma^{*}(t)}, \ul{c} \rangle \right\}
\]
which equals the Legendre-Fenchel transform of $\mathcal{L}_{\eta}(\gamma^{*}(t),\cdot)$ on $T_{\gamma^{*}(t)}(SE(3))$ that can be computed by the results in Appendix \ref{app:viscosity}, recall Def. \ref{def:LT} and
Lemma \ref{lemma:LF} so that we obtain the Hamilton-Jacobi equation on $\R^{3}\rtimes S^2$:
\[
\begin{array}{l}
-\frac{\partial W}{\partial t}(\gamma^{*}(t),t)= \frac{1}{2\eta} \left( \mathbf{G}^{-1}_{\gamma^{*}(t)}({\rm d}W(\gamma^{*}(t),t),{\rm d}W(\gamma^{*}(t),t))\right)^{\eta} \desda \\[8pt]
-\frac{\partial W}{\partial t}(\gamma^{*}(t),t)= \frac{1}{2\eta} \left(
D^{11}(\left.\mathcal{A}_{1}\right|_{\gamma^{*}(t)}W(\gamma^{*}(t),t))^2+D^{11}(\left.\mathcal{A}_{2}\right|_{\gamma^{*}(t)}W(\gamma^{*}(t),t))^2 \right. \\
\left.
+
D^{44}(\left.\mathcal{A}_{4}\right|_{\gamma^{*}(t)}W(\gamma^{*}(t),t))^2+D^{44}(\left.\mathcal{A}_{5}\right|_{\gamma^{*}(t)}W(\gamma^{*}(t),t))^2
\right)^{\eta},
\end{array}
\]
along the characteristic curves.
So we conclude that iso-contours of the solutions of $W$ in Eq.~(\ref{Hamdi2}) and Eq.~(\ref{Hamdi1}) are (at least) locally geodesically equidistant.
This allows us to use the Hamilton-Jacobi equations (\ref{Hamdi1}) for
wavefront propagation methods \cite{Sethian,Osher} for finding geodesics in $\R^{3} \rtimes S^2$. We will consider this in future work.

\section{Asymptotical Expansions around the Origin of the $k$-step Time-integrated Heat Kernel on $\R^{3}\rtimes S^{2}$ \label{ch:asymptotics}}

Recall from Section \ref{ch:GE} the Gaussian estimates for the heat-kernels (contour enhancement);
\[
K_{t}^{D^{33},D^{44}}(\ul{y},\ul{n})= \frac{1}{16 \pi^2 (D^{33})^{2}  (D^{44})^{2} t^4}\, e^{-\frac{|g=(\ul{y},R_{\ul{n}})|_{D^{33},D^{44}}^{2}}{4t}}.
\]
with $|g|$ given by Eq.~(\ref{diffusionmodulus}).
Then we obtain the following relation (with again short notation $g=(\ul{y},R_{\ul{n}}) \in SE(3)$):
\[
\begin{array}{ll}
R_{\lambda, k}^{D^{33}, D^{44}}(\ul{y},\ul{n}) &:=\left((-D^{33}(\mathcal{A}_{3})^2 -D^{44}(\mathcal{A}_{4})^2-D^{44}(\mathcal{A}_{5})^2+\lambda I)^{-k} \delta_{e}\right)(\ul{y},\ul{n}) \\
&=
\int \limits_{0}^{\infty} K_{t}^{D^{33},D^{44}}(\ul{y},\ul{n}) \; \Gamma(t\, ;\, k,\lambda)\; {\rm d}t \\
 &=\frac{1}{16 \pi^2 (D^{33})^2 (D^{44})^2} \frac{\lambda^{k}}{\Gamma(k)} \int \limits_{0}^{\infty} t^{k-5} e^{-\frac{|g|^2}{4t}} e^{-\lambda t}\, {\rm d}t \\
 &=\frac{2^{1-k}}{\pi^2 (D^{33})^2 (D^{44})^2} \frac{\lambda^{\frac{k+4}{2}}}{(k-1)!}\, |g|^{k-4} \, K_{4-k}(|g|\sqrt{\lambda})
\end{array}
\]
with $\Delta_{S^{2}}=\mathcal{A}_{4}^{2}+\mathcal{A}_{5}^2$ and
now we have the following asymptotical formula for the Bessel functions
\[
K_{\nu}(z) \approx
\left\{
\begin{array}{ll}
 -\log(z/2)-\gamma_{ {\small EULER}} &\textrm{ if }\nu=0 \\
 \frac{1}{2} (|\nu|-1)! \left( \frac{z}{2}\right)^{-|\nu|} &\textrm{ if }\nu \in \mathbb{Z}, \nu \neq 0.
\end{array}
\right.
\]
with Euler's constant $\gamma_{ {\small EULER}}$,
which holds for $0< z << 1$,
so that for $|g|=|(\ul{y},R_{\ul{n}})|<<1$ we have
\[
R_{\lambda, k}^{D^{33}, D^{44}}(\ul{y},\ul{n}) \approx
\left\{
\begin{array}{ll}
\frac{2^{-k + |k-4|}}{\pi^2 (D^{33})^{2} (D^{44})^{2}} |g|^{-|k-4| +k-4} \lambda^{-\frac{|k-4|}{2}} &\textrm{ if }k\neq 4, \\
\frac{1}{\pi^2 (D^{33})^{2} (D^{44})^{2}} (-\log \frac{|g|\sqrt{\lambda}}{2}-\gamma_{{\small EULER}}) \frac{2^{-3}\lambda^{4}}{3! } &\textrm{ if }k=4
\end{array}
\right.
\]
Likewise in the contour completion kernel, recall (\ref{kresolvent}), we get rid of the singularity at the origin in
\[
R_{\lambda,k}=R_{\lambda,1}*^{k-1}_{\R^{3}\rtimes S^2}R_{\lambda,1}
\]
by setting $k \geq \textrm{dim}(\R^{3}\rtimes S^2)=5$. In other words to avoid a singularity in the Green's function one must use
at least $k=5$ iteration steps in the contour-enhancement process.

\section{Conditions on a Left-invariant Metric Tensor on $\R^{3}\rtimes S^{2}$ and the Choice of $D^{ii}$ and $g^{ij}$ \label{app:C}}

In this section we derive sufficient conditions on a metric tensor (\ref{metrictensor}) to be both left-invariant
and well-defined on the quotient $\R^{3}\rtimes S^2=SE(3)/(\{\ul{0}\}\times SO(2))$.
\begin{definition}
A metric tensor $\tilde{\mathbf{G}}: SE(3) \times T(SE(3)) \times T(SE(3)) \to \mathbb{C}$
\begin{itemize}
\item is left-invariant iff
\begin{equation}
\tilde{\mathbf{G}}_{gq}((L_{g})_{*}\tilde{X}_q,(L_{g})_{*}\tilde{Y}_q)= \tilde{\mathbf{G}}_{q}(\tilde{X}_q,\tilde{Y}_q)
\end{equation}
for all $g,q \in SE(3)$ and all $\tilde{X},\tilde{Y} \in T(SE(3))$.
\item provides a well-defined metric tensor on $\R^{3}\rtimes S^2=SE(3)/(\{0\}\times SO(2))$ iff
\begin{equation}\label{welldefined}
\tilde{\mathbf{G}}_{gh}((R_{h})_{*}\tilde{X}_g,(R_{h})_{*}\tilde{Y}_g)= \tilde{\mathbf{G}}_{g}(\tilde{X}_g,\tilde{Y}_g)
\end{equation}
for all $g \in SE(3)$ and all $h \in (\{0\}\times SO(2))$ and all $\tilde{X},\tilde{Y} \in T(SE(3))$, in which case the corresponding metric tensor on the quotient is then given by
\begin{equation}\label{metricshit}
\ul{G}_{(\ul{y},{\ul{n}})}(\sum_{i}c^{i}\left.\mathcal{A}_{i}\right|_{(\ul{y},\ul{n})},\sum_{j}d^{j}\left.\mathcal{A}_{j}\right|_{(\ul{y},\ul{n})}):=
\ul{G}_{(\ul{y},R_{\ul{n}})}(\sum_{i}c^{i}\left.\tilde{\mathcal{A}}_{i}\right|_{(\ul{y},R_{\ul{n}})},\sum_{j}d^{j}\left.\tilde{\mathcal{A}}_{j}\right|_{(\ul{y},R_{\ul{n}})})
\end{equation}
where vector fields are described by the differential operators on $C^{1}(\R^{3}\times S^2)$:
\[
\begin{array}{l}
(\left.\mathcal{A}_{j}\right|_{(\ul{y},\ul{n})}U)(\ul{y},\ul{n})=
\lim \limits_{h \to 0} \frac{U(\ul{y}+h R_{\ul{n}}\ul{e}_{j},\ul{n})-U(\ul{y}-h R_{\ul{n}}\ul{e}_{j},\ul{n})}{2h}, \\
(\left.\mathcal{A}_{3+j}\right|_{(\ul{y},\ul{n})}U)(\ul{y},\ul{n})=
\lim \limits_{h \to 0} \frac{U(\ul{y}, (R_{\ul{n}}R_{\ul{e}_{j},h})\ul{e}_z)-
U(\ul{y}, (R_{\ul{n}}R_{\ul{e}_{j},-h})\ul{e}_z)}{2h}, \ \ j=1,2,3,
\end{array}
\]
where $R_{\ul{e}_{j},h}$ denotes the counter-clockwise rotation around axis $\ul{e}_{j}$ by angle $h$,
with {\small $\ul{e}_{1}=(1,0,0)^{T}$, $\ul{e}_{2}=(0,1,0)^{T}$, $\ul{e}_{3}=(0,0,1)^{T}$}. Note that the choice of $R_{\ul{n}}$ satisfying Eq.~(\ref{Rn}) does not affect Eq.~(\ref{metricshit}) because of Eq.~(\ref{welldefined}).
\end{itemize}
\end{definition}
\begin{remark} In the other sections in this article, for the sake of simplicity,
we do not distinguish between $\tilde{X}$ and $X$. For example we write both
$\left.\mathcal{A}_{3}\right|_{(\ul{y},\ul{n})}U$ and $\left.\mathcal{A}_{3}\right|_{(\ul{y},R_{\ul{n}})}\tilde{U}$.
\end{remark}
Set $\underline{\tilde{\mathcal{A}}}=(\tilde{\mathcal{A}}_{1},\ldots,\tilde{\mathcal{A}}_{6})^{T}$ as a column vector of left-invariant vector fields on $SE(3)$, then one has the following identity
\[
\underline{\tilde{\mathcal{A}}}_{gh}= Z_{\alpha} (\mathcal{R}_{h})_{*} \underline{\tilde{\mathcal{A}}}_{g} (\equiv Z_{\alpha}\mathcal{A}_{g}),
\]
with $Z_{\alpha} \in SO(6)$ is given by Eq.~(\ref{Zalpha}),
which is straightforwardly verified by Eq.~(\ref{eq:LeftInvVFSEthree}) where we computed the left-invariant vector fields
explicitly in Euler-angles. As the functions $\tilde{U}$ are invariant under right-multiplication with elements $(\ul{0},R_{\ul{e}_{\alpha}}) \in (\{0\}\times SO(2))$ this yields the following identity
\[
\begin{array}{l}
(\left.\tilde{\mathcal{A}}_{i}\right|_{gh}\tilde{U})(gh)= \sum \limits_{j=1}^{6} (Z_{\alpha})^{i}_{j}
(\left.\tilde{\mathcal{A}}_{i}\right|_{g}\tilde{U})(g),
\end{array}
\]
for all $h=(\ul{0},R_{\ul{e}_{z},\alpha})$, where $Z_{\alpha}=R_{\ul{e}_{z},\alpha} \oplus R_{\ul{e}_{z},\alpha} \in SO(6)$, with
$R_{\ul{e}_{z},\alpha} \in SO(3)$ is given
by (\ref{Zalpha}). Recall that by right multiplication with
$(\ul{0}, R_{\ul{e}_{z},\alpha})$ one takes a different section in the partition of left-cosets in
$\R^{3}\rtimes S^{2}$, boiling down to a rotation over $\alpha$ simultaneously in the two grey planes depicted
in Figure~\ref{fig:intuition}.
\begin{theorem} \label{th:only}
A metric tensor $\mathbf{G}_{(\ul{y},\ul{n})}:(\R^{3}\rtimes S^{2}) \times T(\R^{3}\rtimes S^{2})\times T(\R^{3}\rtimes S^{2}) \to \mathbb{R}$  given by
\[
\mathbf{G}_{(\ul{y},\ul{n})}= \sum \limits_{i,j=1}^{6} g_{ij}(\ul{y},R_{\ul{n}})\; {\rm d}\mathcal{A}^{i} \otimes {\rm d}\mathcal{A}^{j}
\]
with $g_{6j}=g_{j6}=0$, $j=1,\ldots,6$
is well-defined and left-invariant iff
\begin{equation}\label{condit}
\begin{array}{l}
Z_{\alpha} \, [g_{ij}]\, Z_{\alpha}^{T}
=[g_{ij}]\ , \textrm{ for all }\alpha \in [0,2\pi), \textrm{ and } \\
g_{ij}(\ul{y},R_{\ul{n}})=g_{ij} \textrm{ are constant }
\end{array}
\end{equation}
which is satisfied iff
\begin{equation} \label{Fr}
[g_{ij}]=\textrm{diag}\{g_{11},g_{11},g_{33},g_{44},g_{44},0\}. 
\end{equation}
\end{theorem}
\textbf{proof }
The left-invariant vector fields satisfy $\left.\tilde{\mathcal{A}}_{i}\right|_{gq}=(L_{g})_{*} \left.\tilde{\mathcal{A}}_{i}\right|_{q}$, since $(L_{g})_{*}(L_{q})_{*}=(L_{gq})_{*}$.
So the left-invariance requirement reduces to
\[
L_{g^{*}}g_{ij}=g_{ij}
\]
for all $g \in SE(3)$. Now $SE(3)$ acts transitively onto $SE(3)$ (and onto $\R^{3}\rtimes S^2$) so $g_{ij}$ must be constant.
Regarding (\ref{welldefined}) we note that
\[
\begin{array}{l}
\mathbf{G}_{gh}((R_{h})_{*} \sum \limits_{i=1}^{6}c^{i} \left. \tilde{\mathcal{A}}_{i}\right|_{g},(R_{h})_{*} \sum \limits_{j=1}^{6}c^{j} \left. \tilde{\mathcal{A}}_{j}\right|_{g} )
= \sum \limits_{i,j,k,l=1}^{6}\mathbf{G}_{gh}( c^{k} (Z_{\alpha}^{T})^{k}_{i} \left.\tilde{\mathcal{A}}_{k}\right|_{gh}, Z_{\alpha}^{T})^{l}_{j} \left.\tilde{\mathcal{A}}_{l}\right|_{gh} \\
= \sum \limits_{i,j,k,l=1}^{6} c^{k}c^{l} g_{ij} Z_{\alpha}^{T})^{k}_{i} Z_{\alpha}^{T})^{k}_{j}
=\mathbf{G}_{g}( \sum \limits_{i=1}^{6}c^{i} \left. \tilde{\mathcal{A}}_{i}\right|_{g}, \sum_{j}c^{j} \left. \tilde{\mathcal{A}}_{j}\right|_{g} )
=\sum \limits_{i,j=1}^{6}c^{i}c^{j} g_{ij}
\end{array}
\]
for all $(c^{1},\ldots c^{6}) \in \R^{6}$ and all $\alpha \in [0,2\pi)$ so that
\[
\forall_{\alpha \in [0,2\pi)}\; : \;
Z_{\alpha}^{T} [g_{ij}]Z_{\alpha}=[g_{ij}].
\]
The final results follows by Schur's lemma on $SO(2)$ (operators commuting with irreducible representations are multiples of the identity).$\hfill \Box$

The actual (horizontal) metric (induced by the metric tensor given by Eq.~(\ref{metrictensorR3S2})
on the contact manifold $(\R^{3}\rtimes S^{2}, {\rm d}\mathcal{A}^{1},{\rm d}\mathcal{A}^{2}$)
is given by
\begin{equation} \label{dist}
d(\ul{y},\ul{n}\, , \,  \ul{y}', \ul{n}'):=
d(R_{\ul{n}'}^{T}(\ul{y}-\ul{y}'),R_{\ul{n}'}^{T}\ul{n}\, ,\, \ul{0},\ul{e}_z )
\end{equation}
with (horizontal) modulus given by
\begin{equation}\label{metrichor}
\begin{array}{l}
d_{\R^{3}\rtimes S^{2}}(\ul{y}, \ul{n} ,\, \ul{0},\ul{e}_z )= \\
\inf \limits_{{\tiny \begin{array}{c}
\tilde{\gamma}=(\ul{x}(\cdot),R(\cdot)) \in C^{\infty}((0,t), SE(3)),\\
\tilde{\gamma}(0)=(\ul{0},I=R_{\ul{e}_z}), \tilde{\gamma}(t)=(\ul{y},R_{\ul{n}}),\\
\langle \left.{d \rm}\mathcal{A}^{1}\right|_{\tilde{\gamma}}, \dot{\tilde{\gamma}}\rangle=
\langle \left.{d \rm}\mathcal{A}^{2}\right|_{\tilde{\gamma}}, \dot{\tilde{\gamma}}\rangle=
\langle \left.{d \rm}\mathcal{A}^{6}\right|_{\tilde{\gamma}},\dot{\tilde{\gamma}}\rangle=0
\end{array}}}
\int \limits_{0}^{t} \sqrt{\sum \limits_{i,j \in \{3,4,5\}} g_{ij} \,
\langle \left. {\rm d}\mathcal{A}^{i}\right|_{\tilde{\gamma}(p)}, \dot{\tilde{\gamma}}(p)
\rangle \,
\langle \left. {\rm d}\mathcal{A}^{j}\right|_{\tilde{\gamma}(p)}, \dot{\tilde{\gamma}}(p)
\rangle}
\, {\rm d}p\ ,
\end{array}
\end{equation}
\begin{remark}
When measuring the distance (Eq.~(\ref{dist}) and (\ref{metrichor})) between the cosets $[g_{1}]=g_{1}H$ and $[g_{2}]=g_{2}H$, with
$H=\{0\} \times SO(2)$ and $g_{i}=(\ul{x}_{i},R_{\ul{e}_{z},\gamma_{i}}R_{\ul{e}_{y},\beta_{i}}R_{\ul{e}_{z},\alpha_{i}})$, $i=1,2$ one first determines
the elements from both co-sets with vanishing first Euler angle, i.e
$(\ul{x}_{i},R_{\ul{e}_{z},\gamma_{i}}R_{\ul{e}_{y},\beta_{i}})$ and then compute the distance w.r.t. these elements in the full
group where the connecting curves are not allowed to use the ``illegal'' $\mathcal{A}_{6}$-directions. Note that this procedure leads to
a well-defined metric on $(\R^{3}\rtimes S^{2}, {\rm d}\mathcal{A}^{1},{\rm d}\mathcal{A}^{2}$) despite taking the section. The condition
$\langle \left. {\rm d}\mathcal{A}^{6} \right|_{\gamma}, \dot{\gamma}\rangle =0$ avoids possible short-cuts via the ``illegal''
$\mathcal{A}_{6}$-direction.
\end{remark}
\subsection{Data adaption and left invariance}

In the erosion algorithms one can include
adaptivity by making $D^{44}$ depend on the local laplace-Beltrami-operator, recall Subsection \ref{ch:adapt},
whereas in the diffusion algorithms one can include adaptivity by replacing $D^{33}$ in $\mathcal{A}_{3} D^{33} \mathcal{A}_{3}$ by a data adaptive conductivity
\[
D^{33}(U)= e^{-\frac{|\mathcal{A}_{3}U|^{2}}{K^2}}, \  K>0.
\]
In both cases this does \emph{not} correspond to making a left-invariant adaptive inverse metric tensor, as this would contradict Theorem \ref{th:only}. However, an adaptive erosion generator ($\eta=1$):
\[
U \mapsto \sum \limits_{i,j=1}^{5} g^{ij}(U)(y,n) \cdot (\mathcal{A}_{i} \otimes  \mathcal{A}_{j})({\rm d}U, {\rm d}U)=\sum \limits_{i,j=1}^{5} g^{ij}(U)(y,n) \cdot (\mathcal{A}_{i}U)(\mathcal{A}_{j}U)
\]
is left-invariant iff $g^{ij}$ is left-invariant, that is $g^{ij}$ commutes with left regular action of $SE(3)$ onto
$\mathbb{L}_{2}(\R^{3}\rtimes S^2)$. Similarly
an adaptive diffusion generator
\[
U \mapsto \sum \limits_{i,j=1}^{5}  \mathcal{A}_{i}( D^{ij}(U)  \mathcal{A}_{j}(U)),
\]
is left-invariant iff $D^{ij}$ is left-invariant.

\section{Putting the left-invariant vector fields and the diffusion generator in matrix form
in case of linear interpolation\label{app:G}}

We will derive the $N^{3} N_{o} \times N^{3} N_{o}$ matrix form of the \emph{forward} approximation (\ref{LIforward}) of the left-invariant vector fields and subsequently
for the (hypo)-elliptic left-invariant diffusion generator.
Here $N$ stands for the number of spatial pixels (where we assume a cubic domain) and where $N_{o}$ denotes the number of
discrete orientations on $S^{2}$. Note that the matrix-form of the backward and central differences can be derived analogously.

Recall that the forward approximations of the left-invariant vector fields are given by
\[
\begin{array}{rcl}
(\mathcal{A}_{p+3}U)(\ul{y},\ul{n}_{l}) &\approx& \frac{U(\ul{y},R_{\ul{n}_{l}} R_{\ul{e}_{p},h_{a}}\ul{e}_{z})-U(\ul{y},\ul{n}_{l})}{h_{a}}, \\
(\mathcal{A}_{p}U)(\ul{y},\ul{n}_{l}) &\approx& \frac{U(\ul{y}+ h R_{\ul{n}_{l}}\ul{e}_{p},\ul{n}_{l})-U(\ul{y},\ul{n}_{l})}{h},
\end{array}
\]
$p=1,2,3$ and
where $l \in \{1,\ldots, N_{o}\}$ enumerates the discrete orientations $\ul{n}_{l} \in S^{2}$ and where
$\ul{y}=(y^{1},y^{2},y^{3}) \in \{1,\ldots,N\}^{3}$, where $h$ denotes spatial step-size and $h_{a}$ denotes denotes the angular stepsize.

We approximate the angular vector fields by
\[
(\mathcal{A}_{p+3}U)(\ul{y},\ul{n}_{l}) \approx \frac{1}{h_{a}} \left(-U(\ul{y},\ul{n}_{l}) + \sum \limits_{l'=1}^{N_{0}} M_{ll'}^{f,h_{a},p+3} \; U(\ul{y},\ul{n}_{l'}),  \right)
\]
where the upper-index f stands for ``forward'' and with
\begin{equation} \label{angularentries}
M_{ll'}^{f,h_{a},p+3} = \left\{
\begin{array}{ll}
1-\sum \limits_{\ul{n}_{j} \in A_{p,l}} (\ul{n}_{p,l}-\ul{n}_{l'})\cdot (\ul{n}_{j}-\ul{n}_{l'}) &\textrm{if }\ul{n}_{l'} \in A_{p,l} \\
0 &\textrm{else}.
\end{array}
\right.
\end{equation}
where $A_{p,l}$ is the unique spherical triangle in the spherical triangularization containing $\ul{n}_{p,l}:=R_{\ul{n}_{l}} R_{\ul{e}_{p},h_{a}} \ul{e}_{z}$.

We approximate the spatial vector fields by
\[
(\mathcal{A}_{p}U)(\ul{y},\ul{n}_{l}) \approx -\frac{1}{h} U(\ul{y},\ul{n}_{l}) + \frac{1}{h}
\sum \limits_{y^{1'}, y^{2'}, y^{3'}} k_{l}^{p,h}[y^{1}-y^{1'},y^{2}-y^{2'}, y^{3}-y^{3'}] \; U(y^{1'},y^{2'},y^{3'},\ul{n}_{l})\ ,
\]
with $l$-indexed discrete spatial kernel given by
\[
k_{l}^{p,h}(y^{1},y^{2},y^{3}) = \prod \limits_{m=1}^{3} v_{(\ul{y}_{p,l})^{m}}(y^{m})\ ,
\]
with $(\ul{y}_{p,l})^{m}$ the m-th component of the vector $\ul{y}_{p,l}:= h R_{\ul{n}_{l}} \ul{e}_{p}$, $p=1,2,3$, and
with linear interpolation kernel $v_{a}: \mathbb{Z} \to [0,1]$ given by
\[
v_{a}(b)=
\left\{
\begin{array}{ll}
1- |a| &\textrm{if }b=0 \\
H(ab) |a| & \textrm{if }b \in \{-1,1\} \\
0 &\textrm{else }
\end{array}
\right.
\]
with heavyside function $H(u)$ while assuming $|a|<1$.

So if we store $((U(\ul{y},\ul{n}_{l}))_{l \in \{1, \ldots, N_{o}\}})_{\ul{y} \in \{1, \ldots,N\}^{3}}$ in one long column vector
$\ul{u} \in \R^{N^{3}N_{o}}$ we can represent the (forward) angular left-invariant vector fields by the matrix
\[
A_{p+3}^{f}:=\frac{1}{h_{a}}\left(\bigoplus \limits_{l=1}^{N_{3}} M^{f,h_{a}}_{p+3}-I_{N^{3} N_{o}} \right)=\frac{1}{h_{a}}\left(I_{N^{3}} \otimes
M^{f,h_{a}}_{p+3}-I_{N^{3} N_{o}} \right)
\]
with $\otimes$ the Kronecker product and with $M^{f,h_{a}}_{p+3} \in \R^{N_{0}\times N_{o}}$ the matrix with entries (\ref{angularentries}) and
where $I_{N^{3}N_{o}} \in \R^{N^{3}N_{o}\times N^{3}N_{o}}$ denotes the identity matrix.

The (forward) spatial left-invariant vector fields is now represented by the block matrix
\[
A_{p}^{f}:= \frac{1}{h}\left(\bigoplus \limits_{l=1,\ldots,N_{o}} M^{f,h,l}_{p}-I_{N^{3}N_{o}}\right)
\]
with $M^{f,h,l}_{p} \in \R^{N^{3} \times N^{3}}$ the matrix with entries $h^{-1}\, k_{l}^{p,h}[y^{1}-y^{1'},y^{2}-y^{2'}, y^{3}-y^{3'}]$.

Similarly, one defines central (index by $c$) and backward differences (index by $b$). The matrix-representation of the generator of
the (hypo)-elliptic diffusion, Eq.~\!(\ref{goal}),
is now
given by $J_{\R^{3}} + J_{S^{2}} \in \R^{N^{3}N_{o}\times N^{3}N_{o}}$ with
\begin{equation}\label{Js}
\begin{array}{l}
J_{\R^{3}}:= D^{33} A_{3}^{f} A_{3}^{b} + D^{11}\left( A_{1}^{f} A_{1}^{b} + A_{2}^{f} A_{2}^{b}\right), \\
J_{S^{2}}:= D^{44} \left( A_{4}^{f} A_{4}^{b}+  A_{5}^{f} A_{5}^{b}\right).
\end{array}
\end{equation}

\section{Solving (the Pfaffian system) for geodesics on $(SE(3),{\rm d}\mathcal{A}^{1},{\rm d}\mathcal{A}^{2}, {\rm d}\mathcal{A}^{6})$ \label{app:E}}

\begin{figure}[h] \vspace{-0.15cm}\mbox{}
\centerline{
\hfill
\includegraphics[width=0.7\hsize]{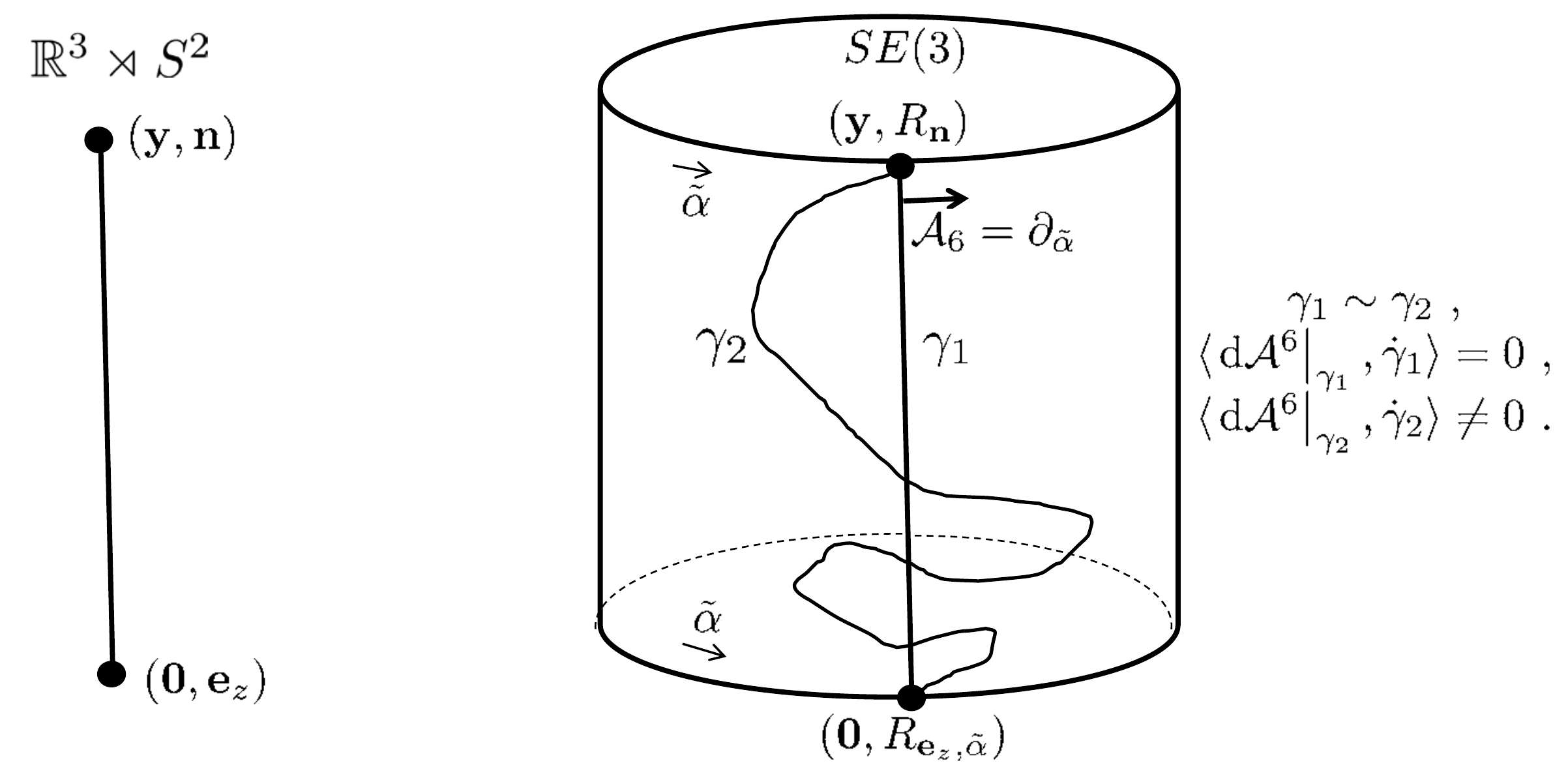}
\hfill
}
\vspace{-0.4cm}\mbox{}
\caption{An (optimal) curve in $\R^{3}\rtimes S^{2}$ corresponds to a class of equivalent curves in $SE(3)$. For convenience we impose the extra constraint $\langle \left. {\rm d}\mathcal{A}^{6}\right|_{\gamma},\dot{\gamma}\rangle=0$ to avoid spurious torsion and curvature.
In the figure we see that $\gamma_{1}$ and $\gamma_{2}$ are equivalent, i.e.
$\gamma_{1}(s)=\gamma_{2}(s) h(s)$ with $h(s)=(\ul{0},R_{\ul{e}_{z}, \alpha(s)})$, but only $\gamma_{1}$ satisfies the convenient constraint. As the constraint $\rm{d}\mathcal{A}^{6}=0$ does not affect the final geodesic in $\R^{3}\rtimes S^{2}$, we can set the corresponding Lagrange-multiplier $\lambda_{6}=0$ in Eq.~(\ref{Pfaff2}).
}
\label{fig:schematical}
\vspace{-0.3cm}\mbox{}
\end{figure}

Here we follow the same approach as in \cite[App.A]{DuitsAMS2} where we minimized on the contact manifold $(SE(2), {\rm d}\mathcal{A}^{3}_{\textrm{SE(2)}})$.
We extend the manifold $SE(3)$ into $Z=SE(3) \times (\R^{+})^2
\times \R^{+} \times (T(SE(3)))^{*}$ on which $SE(3)$ acts as follows
{\small
\[
\eta_{g=(x,y,z, R)}((y^{1},y^{1},y^{3},R'), \KKK, \sigma, (\lambda^{1},\lambda^{2}, \lambda^{3},\lambda^{4}, \lambda^{5},\lambda^{6}))=
(g(y^{1},y^{1},y^{3},R'), \KKK,\sigma,\textrm{CoAd}_{g}(\lambda^{1},\lambda^{2}, \lambda^{3},\lambda^{4}, \lambda^{5},\lambda^{6}))
\]
}
where $\textrm{CoAd}$ is the coadjoint action of $SE(3)$ onto the dual tangent space $(T(SE(3)))^{*}$, where
$\KKK=\ddot{\ul{x}}=-\kappa^{2}\mathcal{A}_{1}+\kappa^{1}\mathcal{A}_{2}+ 0 \mathcal{A}_{3}$ and ${\rm ds}=\sigma {\rm d}t=\|\ul{x}'(t)\| {\rm d}t$
represent curvature and arc-length $s$ of the spatial $\R^{3}$-part
of curves in $SE(3)$, as by definition on this extended contact manifold $(Z,\theta^{1},\ldots \theta^{6})$ the following Pfaffian forms vanish
\begin{equation} \label{vanish}
\begin{array}{l}
\theta^{1}:={\rm d}\mathcal{A}^{1}=0\ , \\
\theta^{2}:={\rm d}\mathcal{A}^{2}=0\ , \\
\theta^{3}:={\rm d}\mathcal{A}^{3}-\sigma {\rm d}t=0\ , \\
\theta^{4}:={\rm d}\mathcal{A}^{4}- \kappa^{1} \sigma {\rm d}t=0\ , \\
\theta^{5}:={\rm d}\mathcal{A}^{5}- \kappa^{2} \sigma {\rm d}t=0\ , \\
\theta^{6}:={\rm d}\mathcal{A}^{6}=0
\end{array}
\end{equation}
The first five Pfaffian forms have to vanish in order to ensure that
$\gamma:[0,L] \to SE(3)$ yields a horizontal curve $s \mapsto (\ul{x}(s),\ul{n}(s))$
in $\R^{3}\rtimes S^{2}$ where the spatial tangent corresponds to point on the sphere, i.e.
\begin{equation} \label{horR3S2}
\dot{\ul{x}}(s)=\ul{n}(s).
\end{equation}
The 6th Pfaffian form is imposed to avoid a spurious correspondence between horizontal curves
$\gamma=(\ul{x},R)$
in the sub-Riemannian manifold
$(SE(3), {\rm d}\mathcal{A}^{1}, {\rm d}\mathcal{A}^{2}, {\rm d}\mathcal{A}^{3})$
to curves $(\ul{x},\ul{n})$ in $\R^{3}\rtimes S^{2}$ satisfying (\ref{horR3S2}) via
\[
(\ul{x}(s),\ul{n}(s))= (\ul{x}(s), R(s) \ul{e}_{z}).
\]
To this end we recall that
$
\mathcal{A}_{6} \tilde{U} = 0 \textrm{ for all smooth }
\tilde{U}:SE(3) \to \R \textrm{ given by }\tilde{U}(x,R)=U(x,R\ul{e}_{z}) \textrm{ with }U:\R^{3}\rtimes S^{2} \to \R.
$
and $\langle \left.{\rm d}\mathcal{A}^{6}\right|_{\gamma}, \dot{\gamma} \rangle =0$
\emph{excludes} spurious curves in $SE(3)$ by means of
\[
(\ul{x}(s),R(s)) \mapsto (\ul{x}(s),R(s) R_{\ul{e}_{z},\alpha(s)})
\]
with $\dot{\alpha} \neq 0$, that would correspond to the
same horizontal curves in $\R^{3}\rtimes S^{2}$, see Figure \ref{fig:schematical}.
Our goal is to solve for the minimizing curves of the optimization problem (\ref{semimetric}) on $(SE(3),{\rm d}\mathcal{A}^{1},{\rm d}\mathcal{A}^{2}, {\rm d}\mathcal{A}^{6})$
\begin{equation} \label{semimetricfollowup}
\begin{array}{l}
\gothic{d}(g_{1},g_{2}):= \gothic{d}(g_{2}^{-1}g_{1},e=(\ul{0},I))= \\
\inf \limits_{h_{1}, h_2 \in \{0\} \times SO(2)}\! \! \! \! \! \!\inf \limits_{{\tiny \begin{array}{c}
\gamma=(\ul{x}(\cdot),R_{\ul{n}}(\cdot)) \in C^{\infty}((0,1), SE(3)),\\
\gamma(0)=e h_{2}, \gamma(1)=g_{2}^{-1}g_{1}h_{1},\\
\langle \left.{d \rm}\mathcal{A}^{1}\right|_{\gamma}, \dot{\gamma}\rangle =
\langle \left.{d \rm}\mathcal{A}^{2}\right|_{\gamma}, \dot{\gamma}\rangle =
\langle \left.{d \rm}\mathcal{A}^{6}\right|_{\gamma}, \dot{\gamma}\rangle =0
\end{array}
}}\! \! \! \! \! \! \!
\int_{0}^{1} \sqrt{\sum \limits_{i \in \{3,4,5\}}\frac{1}{D^{ii}}|{\rm d}\mathcal{A}^{i}(\dot{\gamma}(s))|^{2}}\, {\rm d}s
\end{array}
\end{equation}
where we set $D^{33}=\beta^{-2}$, $s$ spatial arc-length, and $D^{44}=D^{55}=1$, so that
the equivalent problem on $\R^3$ (recall Eq.~(\ref{metricred})) is:
\begin{equation}\label{metricredtwo}
d^{\textrm{hor}}((\ul{y},\ul{n}), (\ul{y}', \ul{n}'))=
\inf \limits_{{\scriptsize \begin{array}{c}
\ul{x}(\cdot) \in C^{\infty}((0,L), \R^3),\\
\ul{x}(0)=\ul{0}, \; \dot{\ul{x}}(0)=\ul{e}_{z}, \\
\ul{x}(L)=R_{\ul{n}}^{T}(\ul{y}'-\ul{y}), \\ \dot{\ul{x}}(L)=R_{\ul{n}}^{T} R_{\ul{n}'}\ul{e}_{z}=R_{\ul{n}}^{T}\ul{n}', \\
\end{array}
}}\!
\int \limits_{0}^{L} \sqrt{\|\KKK(s)\|^2 + \beta^2} {\rm d}s
\end{equation}
where $s$, $L>0$, and $\KKK(s)$ are respectively
spatial arclength, total length, and curvature of the
spatial part of the curve.
\begin{theorem}
The stationary curves $s \mapsto \ul{x}(s) \in \R^{3}$ 
that minimize
\[
\int \limits_{0}^{L} \sqrt{\kappa^2+\beta^{2}}\; {\rm d}s
\]
with arclength parameter $s \in [0,L]$, curvature $\kappa$ and $\beta>0$ and free length $L>0$,
subject to boundary conditions
$\ul{x}(0)=\ul{x}_{0}=\ul{0} \in \R^{3}$, $\dot{\ul{x}}(0)=\ul{n}_{0}=\ul{e}_{z} \in S^{2}$ and
$\ul{x}(L)=\ul{x}_{1} \in \R^{3}$, $\dot{\ul{x}}(L)=\ul{n}_{1} \in S^{2}$
coincide with the stationary curves
$s \mapsto (\ul{x}(s),\dot{\ul{x}}(s)) \in \R^{3} \rtimes S^{2}$ of the variational problem
(\ref{metricred}) on $\R^{3}\rtimes S^{2}$ and they satisfy
\[
\frac{d}{dt}
\int_{N_{t}}\psi=0 \textrm{ with }  \psi = \sqrt{\|\KKK\|^2 +\beta^{2}}\; \sigma {\rm d}t
+ \sum \limits_{i=1}^{6}\lambda_{i}\theta^{i}
\]
with ${\rm d}s=\sigma {\rm d}t$, $\sigma=\|\ul{x}'(t)\|$
for all horizontal curve $t \mapsto N_{t}$ perturbations
on the extended 15-dimensional manifold $Z=SE(3)\times (\R^{+})^{3} \times (T(SE(3)))'$. \\
Along the the minimizing curves
$\gamma=(g, \kappa^{1},\kappa^{2},\sigma, \lambda_{1},\ldots,\lambda_{6})$ in $Z$ the
Lagrange-multipliers are given by
\[
\begin{array}{l}
\lambda_{1}(s)=-\dot{\lambda_{5}}(s), \
\lambda_{2}(s)=\dot{\lambda_{4}}(s), \
\lambda_{3}(s)= \beta \sqrt{1-(\lambda_{4}(s))^2-(\lambda_{5}(s))^2}, \\
\lambda_{4}(s)= z^{1}(s), \
\lambda_{5}(s)= z^{2}(s), \
\lambda_{6}=0,
\end{array}
\]
expressed in the normalized curvature
\[
\ul{z}(s)=(z^{1}(s),z^{2}(s)):=
\frac{1}{\sqrt{(\kappa(s))^2+\beta^{2}}}(\kappa^{1}(s),\kappa^{2}(s))= \cosh(\beta s)\; \ul{z}(0) + \frac{\sinh (\beta s)}{\beta}\; \ul{z}'(0),
\]
with
$\kappa=\sqrt{(\kappa^{1})^2+(\kappa^{2})^2}$. They satisfy the
following preservation-laws
\[
\begin{array}{l}
(\lambda_{1})^{2} + (\lambda^{2})^2 +(\lambda_{3})^2=c^{2}\beta^{2}\ , \textrm{ (co-adjoint orbits) },\\
\beta^{-2} (\lambda_{3})^2 +(\lambda_{4})^2 +(\lambda_{5})^2=1\ ,
\end{array}
\]
with $c>0$ constant.
The $SE(3)$-part $s \mapsto g(s)$ of the stationary curves in $Z$ is obtained by means of 
\begin{equation}\label{mu}
\lambda(s)\; m(g^{-1})(s)=\lambda(0)\; m(g^{-1})(0)=\lambda(L)\; m(g^{-1})(L),
\end{equation}
with row-vector $\lambda=(\lambda_{1},\ldots,\lambda_{6})$ and
where $m: SE(3) \to \R^{6 \times 6}$ is the matrix group representation given by
\begin{equation}\label{mg}
m(g)=\left(
\begin{array}{cc}
R & \sigma_{\ul{x}}R \\
0 & R
\end{array}
\right)\ , \textrm{ for all }g=(\ul{x},R) \in SE(3),
\end{equation}
with $\sigma_{\ul{x}} \in \R^{3\times 3}$ such that $\sigma_{\ul{x}}\ul{y}=\ul{x} \times \ul{y}$.
The curvature and torsion magnitude of (the spatial part of) these stationary curves are given by
\[
\kappa(s) = \beta \frac{z(s)}{\sqrt{1-z^2(s)}} \textrm{ and }
\tau(s) = \frac{\det(\ul{z}_{0}\;|\; \ul{z}_{0}')}{z^{2}(s)}
\]
with $z(s)=\sqrt{(z^{1}(s))^2+(z^{2}(s))^2}=\|\ul{z}(s)\|$
with constant $\ul{z}_{0}, \ul{z}_{0}' \in \R^{2}, L>0$ such that
\[
\ul{x}(L)=\ul{x}_{1} \textrm{ and }\dot{\ul{x}}(L)=\ul{n}_{1}.
\]
The spatial part $s \mapsto \ul{x}(s)$ of the curves $s \mapsto  g(s)$ can be obtained
by integrating the Frenet formulas
\[
\left(
\begin{array}{l}
\dot{\ul{T}}(s) \\
\dot{\ul{N}}(s) \\
\dot{\ul{B}}(s)
\end{array}
\right)
=
\left(
\begin{array}{ccc}
0 & \kappa(s) &  0 \\
-\kappa(s) & 0 &\tau(s) \\
0 & -\tau(s) & 0
\end{array}
\right)
\left(
\begin{array}{l}
\ul{T}(s) \\
\ul{N}(s) \\
\ul{B}(s)
\end{array}
\right)
\]
where again $s \in [0,L]$ denotes arclength parametrization of the spatial part of the
curve with initial condition
\[
\left(
\begin{array}{l}
\ul{T}(0) \\
\ul{N}(0) \\
\ul{B}(0)
\end{array}
\right)
=
\left(
\begin{array}{c}
\ul{e}_{z} \\
\frac{1}{\|\ul{z}_0\|}\left( (z_0)_{2} \ul{e}_{x} - (z_0)_{1}
\ul{e}_{y}\right)
 \\
\frac{1}{\|\ul{z}_0\|}\left( (z_0)_{1} \ul{e}_{x} + z_0)_{2}
\ul{e}_{y}\right)

\end{array}
\right)
.
\]
If these curves are a global minimizer then they
coincide with the geodesics on the sub-Riemannian manifold \\
$(SE(3),{\rm d}\mathcal{A}^{1},{\rm d}\mathcal{A}^{2},{\rm d}\mathcal{A}^{6})$ subject to
left-invariant metric
\[
\mathcal{G}_{\beta}= \beta^{2}{\rm d}\mathcal{A}^{3} \otimes {\rm d}\mathcal{A}^{3} +
{\rm d}\mathcal{A}^{4} \otimes {\rm d}\mathcal{A}^{4}+{\rm d}\mathcal{A}^{5} \otimes {\rm d}\mathcal{A}^{5}.
\]
\end{theorem}
\textbf{Proof }
We follow the general theory on optimizing Lagrangians on contact manifolds, \cite{Bryant,Bryantbook} and
we set the 1-form on $Z$
\[
\psi = \sqrt{\|\KKK\|^2 +\beta^{2}}\; \sigma {\rm d}t + \sum \limits_{i=1}^{6}\lambda_{i}\theta^{i}
\]
where we must first find the Lagrange-multipliers.
Now suppose we have a 1-parameter family of Legendre sub-manifolds $\{N_t\}_{t \in \R}$
within $Z$, this corresponds to a one parameter family of horizontal vector fields on
$(SE(3))$.
Then compute the variation of the integrated Lagrangian-form $\psi$
along $N_t$:
\[
\frac{d}{dt} \int_{N_t} \psi = \int_{N_t} \mathcal{L}_{\frac{\partial }{\partial t}} \psi
= \int_{N_t} \frac{\partial}{\partial t} \rfloor {\rm d}\psi +
\int_{N_t} {\rm d}(\frac{\partial}{\partial t} \rfloor \psi)
= \int_{N_t} \frac{\partial}{\partial t} \rfloor {\rm d}\psi
\]
where we used the well-known Stokes Theorem
$\int_{N_t} {\rm d}(\frac{\partial}{\partial_t} \rfloor \psi)=
\oint_{\partial{N_t}} \frac{\partial}{\partial_t} \rfloor \psi=0$ and the formula
for Lie derivatives of volume forms along vector fields
$\mathcal{L}_{X}A=X \rfloor {\rm d}A + d(X \rfloor A)$.
Consequently, the optimal/characteristic curves are entirely determined by
\[
\gamma'(s) \rfloor \left. {\rm d\psi} \right|_{\gamma(s)}=0  \qquad \textrm{ for all }s>0.
\]
But this by definition of the insert operator $\rfloor$ means that
$\left. {\rm d}\psi \right|_{\gamma} (\gamma',v)=0 \textrm{ for all }v \in T(Z)$,
or equivalently formulated
\begin{equation}\label{Pfaff1}
(v \rfloor \left.{\rm d}\psi\right|_{\gamma})(\gamma') = 0 \textrm{ for all }v \in T(Z).
\end{equation}
For further details see \cite{Bryantbook}.
So in order to find the optimal trajectories
one has to integrate the Pfaffian system
\begin{equation}\label{Pfaff}
v \rfloor {\rm d}\psi = 0 \textrm{ for all }v \in T(Z),
\end{equation}
consisting of 15 vanishing dual forms.

\begin{figure}
\centerline{
\hfill
\includegraphics[width=0.4\hsize]{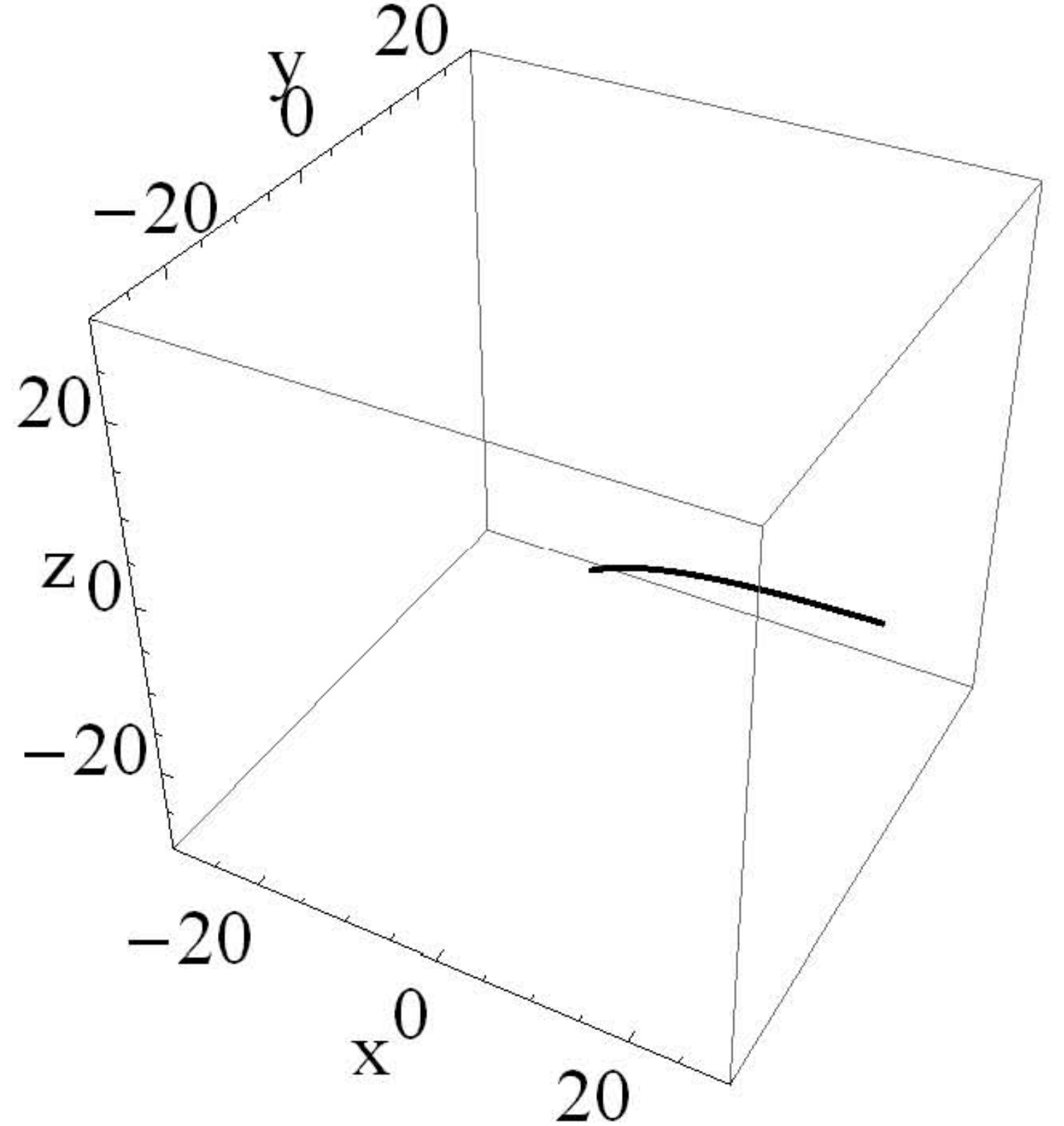}
\includegraphics[width=0.4\hsize]{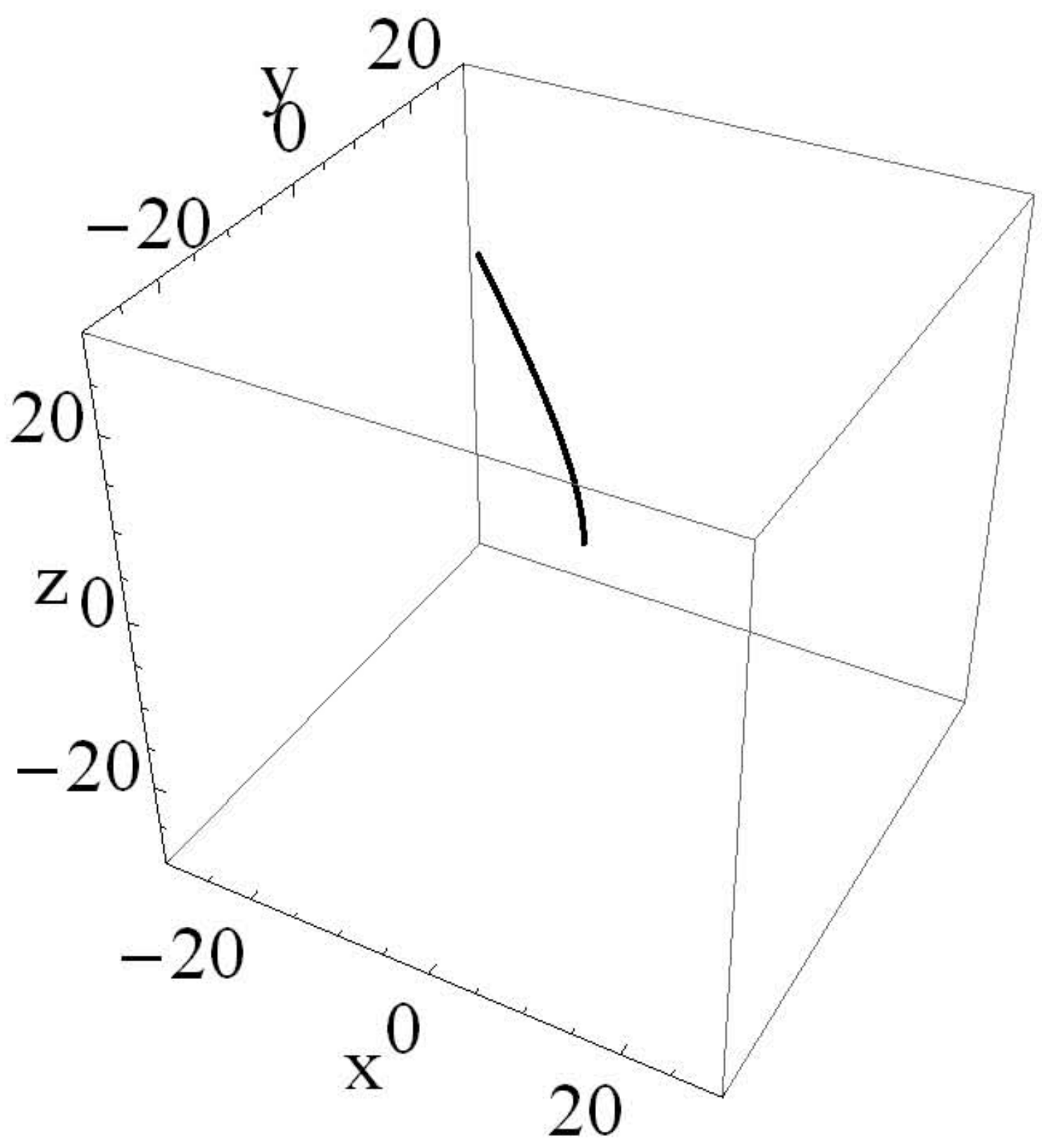}
\hfill}
\caption{Example of a geodesic in $\R^{3}\rtimes S^{2}$. Left: $s \mapsto \tilde{\ul{x}}(s)$. Right: $s \mapsto \ul{x}(s)$.
Parameter settings $\beta=\frac{1}{10}$, $\ul{z}_{0}=(\frac{1}{2},0)^{T}$, $\ul{z}_{0}'=-\beta \ul{z}_{0}$, (so $c=1$ and $W=0$),
length $L=30$, $\dot{\ul{x}}(0)=(0,0,1)^{T}$, $\ul{x}(0)=(0,0,0)$, the rotation and translation needed to map the curve
$\tilde{\ul{x}}$ onto the geodesic $\ul{x}$
is given by $h_{0}^{-1}=(m(\overline{\ul{x}}_{0},\overline{R}_{0}))^{-1}$.
Note that $\ul{x}(s)= \overline{R}_{0}^{T}(\tilde{\ul{x}}(s)-\overline{\ul{x}}_{0})$, with $\overline{\ul{x}}_{0}=(0,0,-5)$ and $\overline{R}_{0} \in SO(3)$
given by Eq.~(\ref{four4}).
In this example both the curvature $\KKK$ and $\ul{z}$ are always aligned with the
$y$-axis which explains that the curve stays within the $yz$-plane.
}\label{fig:examplenot}
\end{figure}

Now by Cartan's structural formula we have
\[
{\rm d}({\rm d}\mathcal{A}^{k})=
-\sum_{i,j=1}^{6} c^{k}_{ij} {\rm d}\mathcal{A}^{i} \otimes {\rm d}\mathcal{A}^{j}    ,
\]
where $c^{k}_{ij}$ are the structure constants of the Lie-algebra $\mathcal{L}(SE(3))$,
recall Eq.~\ref{structureconstants}) and using this crucial identity
one finds
after some straightforward (but intense) computation the following explicit form of (\ref{Pfaff}):
\begin{equation}\label{Pfaff2}
\begin{array}{ll}
\begin{array}{l}
\partial_{\sigma} \rfloor {\rm d}\psi =
(\sqrt{\|\KKK\|^2+\beta^2}-\lambda_{3} -\lambda_{4}\kappa^{1}-\lambda_{5}\kappa^{2}){\rm d}t =0, \\
\partial_{\kappa^{1}} \rfloor {\rm d}\psi = (\kappa^{1}/\sqrt{\|\KKK\|^2+\beta^2}-\lambda_{4})\sigma {\rm d}t =0, \\
\partial_{\kappa^{2}} \rfloor {\rm d}\psi = (\kappa^{2}/\sqrt{\|\KKK\|^2+\beta^2}-\lambda_{5})\sigma {\rm d}t =0, \\[7pt]
-\mathcal{A}_{1} \rfloor {\rm d}\psi = {\rm d}\lambda_{1} -\lambda_{2} {\rm d}\mathcal{A}^{6} +
\lambda_{3} {\rm d}\mathcal{A}^{5}=0, \\
-\mathcal{A}_{2} \rfloor {\rm d}\psi = {\rm d}\lambda_{2} +\lambda_{1} {\rm d}\mathcal{A}^{6} -
\lambda_{3} {\rm d}\mathcal{A}^{4}=0, \\
-\mathcal{A}_{3} \rfloor {\rm d}\psi = {\rm d}\lambda_{3} -\lambda_{1} {\rm d}\mathcal{A}^{5} +
\lambda_{2} {\rm d}\mathcal{A}^{4}=0, \\
-\mathcal{A}_{4} \rfloor {\rm d}\psi = {\rm d}\lambda_{4}-\lambda_{5}{\rm d}\mathcal{A}^{6} +
\lambda_{6} {\rm d}\mathcal{A}^{5} -
\lambda_{2} {\rm d}\mathcal{A}^{3} +\lambda_{3}{\rm d}\mathcal{A}^{2}=0, \\
-\mathcal{A}_{5} \rfloor {\rm d}\psi = {\rm d}\lambda_{5}-\lambda_{6}{\rm d}\mathcal{A}^{4} +
\lambda_{4} {\rm d}\mathcal{A}^{6} -
\lambda_{3} {\rm d}\mathcal{A}^{1} +\lambda_{1}{\rm d}\mathcal{A}^{3}=0, \\
-\mathcal{A}_{6} \rfloor {\rm d}\psi = {\rm d}\lambda_{6}-\lambda_{4}{\rm d}\mathcal{A}^{5} +
\lambda_{5} {\rm d}\mathcal{A}^{4} -
\lambda_{1} {\rm d}\mathcal{A}^{2} +\lambda_{2}{\rm d}\mathcal{A}^{1}=0, \\
\end{array}
&
\begin{array}{l}
\partial_{\lambda_{1}} \rfloor {\rm d}\psi = {\rm d}\mathcal{A}^{1}=0, \\
\partial_{\lambda_{2}} \rfloor {\rm d}\psi = {\rm d}\mathcal{A}^{2}=0,\\
\partial_{\lambda_{3}} \rfloor {\rm d}\psi = {\rm d}\mathcal{A}^{3}-\sigma {\rm d}t=0, \\
\partial_{\lambda_{4}} \rfloor {\rm d}\psi = {\rm d}\mathcal{A}^{4}-\kappa^{1} \sigma {\rm d}t=0, \\
\partial_{\lambda_{5}} \rfloor {\rm d}\psi = {\rm d}\mathcal{A}^{5}-\kappa^{2} \sigma {\rm d}t=0, \\
\partial_{\lambda_{6}} \rfloor {\rm d}\psi = {\rm d}\mathcal{A}^{6}=0, \\
\end{array}
\end{array}
\end{equation}
where in the last three Pfaffian forms in the left column one may set ${\rm d}\mathcal{A}^{1}={\rm d}\mathcal{A}^{2}=0$.
Now by Noether's  Theorem the characteristic curves in $Z$ are contained within the co-adjoint orbits, cf.\cite{Bryant}.
As the adjoint action is given by push-forward of conjugation $\textrm{Ad}(g)=(L_{g}R_{g^{-1}})_{*}$,
\[
\textrm{Ad}(g) A_{i} =
(R_{g^{-1}})_{*} \left.\mathcal{A}_{i}\right|_{g}=
(R_{g^{-1}})_{*} \sum \limits_{i=1}^{3} \left. c^{i}(g) \partial_{y^{i}}\right|_{g}=
\sum \limits_{i=1}^{3}  c^{i}(g) \left.\partial_{y^{i}}\right|_{e}
\]
with $(c^{1})^2+ (c^{2})^2 + (c^{3})^2=1$, where
we used
$\left.\mathcal{A}_{i}\right|_{g}=(L_{g})_{*}A_{i}$
and as the coadjoint action
given by
\[
\langle \textrm{CoAd}(g) \sum \limits_{i=1}^{6}
\lambda_{i}\,{\rm d}\mathcal{A}^{i}, \sum_{j} c^{j}\mathcal{A}_{j}  \rangle=
\langle \sum \limits_{i=1}^{6} \lambda_{i}\, {\rm d}\mathcal{A}^{i}, \textrm{Ad}(g^{-1}) \sum_{j} c^{j}\mathcal{A}_{j} \rangle
\]
acts transitively
on the span by the dual angular generators
we see that $SE(3)$ co-adjoint orbits are given by
\[
(\lambda_{1})^{2}+(\lambda_{2})^{2}+(\lambda_{3})^{2}=C^2 ,
\]
$C>0$.
Indeed it directly follows from (\ref{Pfaff}) that $\sum \limits_{i=1}^{3}\lambda_{i} {\rm d}\lambda_{i}=0$.
The first three Pfaffian forms in (\ref{Pfaff}) yield
\begin{equation}
\begin{array}{ll}
\lambda_{3}= \frac{\beta^2}{\sqrt{\|\KKK\|^2+\beta^2}},  & \\
\lambda_{4}= \frac{\kappa^{1}}{\sqrt{\|\KKK\|^2+\beta^2}}, & \lambda_{2}= \dot{\lambda}_{4}, \\
\lambda_{5}=\frac{\kappa^{2}}{\sqrt{\|\KKK\|^2+\beta^2}}, & \lambda_{1}= -\dot{\lambda}_{5},\\
\end{array}
\end{equation}
so that another preservation law is given by
\[
\beta^{-2} (\lambda_{3})^2 +(\lambda_{4})^2 +(\lambda_{5})^2=1\ ,
\]
and as $\dot{\lambda}_{6}=\lambda_{4} \kappa^{2} -\lambda_{5}\kappa^{1}=0$, $\lambda_{6}$ is constant. Now
since ${\rm d}\mathcal{A}^{6}=0$ is an extra constraint (recall Figure \ref{fig:schematical}) that allows us to choose convenient
representatives in the cosets, the minimizers on the quotient space do not depend on this constraint. Therefore\footnote{Comparison to an alternative
derivation generalizing the approach in \cite{Mumford} on $\R^{2}$ to $\R^{3}$ confirmed this choice.} we can set $\lambda_{6}=0$
and we find
\[
\dot{\lambda}_{5}+\lambda_{1}=0 \Rightarrow \ddot{\lambda}_{5}+\dot{\lambda}_{1}= \ddot{\lambda}_{5}-
\frac{\beta^2 \kappa^{2}}{\sqrt{\|\KKK\|^2+\beta^2}}= \ddot{\lambda}_{5}- \beta^{2}\lambda_{5}=0,
\]
so that $\ddot{\lambda}_{5}= \beta^{2}\lambda_{5}$. Analogously one finds $\ddot{\lambda}_{4}= \beta^{2}\lambda_{4}$,
so similar to the $(SE(2),{\rm d}\mathcal{A}^{3})$-case one has
\begin{equation} \label{ODElinear}
\ddot{\ul{z}}(s)=\beta^2 \ul{z}(s)\ , \textrm{ with }\ul{z}(s)=\frac{1}{\sqrt{\|\KKK(s)\|^2+\beta^{2}}} \KKK(s),
\end{equation}
for the normalized curvature $\ul{z}(s)$.
Furthermore we find
\[
\dot{\lambda}_{1} =-\lambda_{3} \kappa^{2}
\textrm{ and }  \dot{\lambda}_{2} =\lambda_{3} \kappa^{1}.
\]
So again in contrast to the common and well-known elastica the geodesics do not involve special functions and we have
\begin{equation}\label{zs}
\ul{z}(s)= \cosh(\beta s)\; \ul{z}(0) + \frac{\sinh (\beta s)}{\beta}\; \ul{z}'(0).
\end{equation}
Next we compute curvature magnitude and torsion magnitude. Curvature magnitude is given by
\[
\kappa= \sqrt{(\kappa^{1})^2+(\kappa^{2})^2}.
\]
In order to compute the torsion, we compute the normal and binormal from the tangent
\[
\ul{T}(s)=\left.\mathcal{A}_{3}\right|_{\gamma(s)}.
\]
When taking derivatives one should use covariant derivatives, recall Theorem \ref{th:CartanSE2}
to compute
\[
\frac{d}{ds} \left. \mathcal{A}_{i} \right|_{\gamma(s)}= \nabla_{\dot{\gamma}(s)} \mathcal{A}_{i}=
\sum \limits_{k,j=1}^{6} c^{k}_{ji} \langle \left.{\rm d}\mathcal{A}^{j}\right|_{\gamma(s)},
\dot{\gamma}(s)\rangle \, \left.\mathcal{A}_{k}\right|_{\gamma(s)}.
\]
More explicitly this yields
\begin{equation} \label{diffAi}
\begin{array}{l}
\frac{d}{ds} \left.\mathcal{A}_{1}\right|_{\gamma(s)}= c^{3}_{51}
\langle \left.{\rm d}\mathcal{A}^{5}\right|_{\gamma(s)},
\dot{\gamma}(s)\rangle \, \left.\mathcal{A}_{3}\right|_{\gamma(s)}+
c^{2}_{61}
\langle \left.{\rm d}\mathcal{A}^{6}\right|_{\gamma(s)},
\dot{\gamma}(s)\rangle \, \left.\mathcal{A}_{2}\right|_{\gamma(s)}=
-\kappa^{2}(s) \left.\mathcal{A}_{3} \right|_{\gamma(s)}, \\
\frac{d}{ds} \left.\mathcal{A}_{1}\right|_{\gamma(s)}= \kappa^{1}(s) \left.
\mathcal{A}_{3}\right|_{\gamma(s)}, \\
\frac{d}{ds} \left.\mathcal{A}_{3}\right|_{\gamma(s)}= \kappa^{2}(s) \left.
\mathcal{A}_{1}\right|_{\gamma(s)}-\kappa^{1}(s) \left.
\mathcal{A}_{2}\right|_{\gamma(s)}, \\
\frac{d}{ds} \left.\mathcal{A}_{4}\right|_{\gamma(s)}=  \left.
\mathcal{A}_{2}\right|_{\gamma(s)}-\kappa^{2}(s) \left.
\mathcal{A}_{6}\right|_{\gamma(s)}, \\
\frac{d}{ds} \left.\mathcal{A}_{5}\right|_{\gamma(s)}= - \left.
\mathcal{A}_{1}\right|_{\gamma(s)}+\kappa^{1}(s) \left.
\mathcal{A}_{6}\right|_{\gamma(s)} \\
\end{array}
\end{equation}
along the horizontal optimal curves in $SE(3)$, where we used (\ref{vanish}).
Using (\ref{diffAi}) we find
\[
\begin{array}{ll}
\ul{T}'(s)&= \kappa^{2}(s) \mathcal{A}_{1}-\kappa^{1} \mathcal{A}_{2}, \\
\ul{N}(s) &= \frac{1}{\sqrt{(\kappa^{1}(s))^2+ (\kappa^{2}(s))^2}}(\kappa^{2}(s) \mathcal{A}_{1}-\kappa^{1} \mathcal{A}_{2}), \\
\ul{N}'(s) &= \left( \frac{\kappa^{2}(s)}{\sqrt{(\kappa^{1}(s))^2+ (\kappa^{2}(s))^2}}\right)'
\left.\mathcal{A}_{1} \right|_{\gamma(s)} -
\left( \frac{\kappa^{1}(s)}{\sqrt{(\kappa^{1}(s))^2+ (\kappa^{2}(s))^2}}\right)'
\left.\mathcal{A}_{2} \right|_{\gamma(s)} \\
 &+
\frac{\kappa^{2}(s)}{\sqrt{(\kappa^{1}(s))^2+ (\kappa^{2}(s))^2}} \frac{d}{ds}\left.
\mathcal{A}_{1}\right|_{\gamma(s)}-
\frac{\kappa^{1}(s)}{\sqrt{(\kappa^{1}(s))^2+ (\kappa^{2}(s))^2}} \frac{d}{ds}\left.
\mathcal{A}_{2}\right|_{\gamma(s)}, \\
\ul{B}(s) &= \frac{1}{\sqrt{(\kappa^{1}(s))^2+ (\kappa^{2}(s))^2}}(\kappa^{1}(s) \mathcal{A}_{1}+\kappa^{2} \mathcal{A}_{2})
\end{array}
\]
Now the Frenet-frame satisfies $\dot{\ul{N}}=- \kappa \ul{T}+\tau \ul{B}$ so
curvature magnitude and torsion magnitude are given by
\begin{equation} \label{kt}
\kappa=\sqrt{(\kappa^{1})^2+(\kappa^{2})^2} \textrm{ and } \tau =
\frac{\kappa^{1}\dot{\kappa}^{2}-\kappa^{2}\dot{\kappa}^{1}}{\kappa^2}= \frac{z^{1} \dot{z}^{2} -z^{2} \dot{z}^{1}}{(z^{1})^2+(z^{2})^2}.
\end{equation}
So if we write $\|\ul{z}\|=\sqrt{(z^{1})^{2}+(z^{2})^2}$ we see that
\begin{equation} \label{kt2}
\kappa = \beta \frac{\|\ul{z}\|}{\sqrt{1-\|\ul{z}\|^2}} \textrm{ and }\tau = \frac{\det(\ul{z}_{0}\;|\; \ul{z}_{0}')}{\|\ul{z}\|^{2}}.
\end{equation}
where $\ul{z}_{0}=\ul{z}(0) \in \R^{2}$ and $\ul{z}_{0}'=\ul{z}'(0) \in \R^{2}$ and where we note
that $z^{1}$ and $z^{2}$ satisfy the same linear ODE with (constant) Wronskian
\[
z^{1}(s) \dot{z}^{2}(s) -z^{2}(s) \dot{z}^{1}(s)= z^{1}(0) \dot{z}^{2}(0) -z^{2}(0) \dot{z}^{1}(0)=\det(\ul{z}_{0}\;|\; \ul{z}_{0}').
\]
We verified our results by extending Mumford's approach to elastica on $\R^{2}$ to
curves minimizing $\int_{0}^{L} \sqrt{\kappa^{2}(s)+ \beta^{2}} {\rm d}s$ with $L$ free
on $\R^{3}$. The details will follow in future work. For now we mention that such approach
finally results in the following \emph{non-linear} ODE's for $\kappa$ and $\tau$:
\begin{equation} \label{Mumford}
\|\ul{z}\|''=(\beta^{2} +\tau^2) \|\ul{z}\| \textrm{ and } 2\tau \|\ul{z}\|' +\tau ' \|\ul{z}\| =0.
\end{equation}
It can be verified that the \emph{ linear }ODE's Eq.~(\ref{ODElinear})
and Eq.~(\ref{kt2}) indeed imply Eq.~(\ref{Mumford}).
So, besides of all the preservation-laws, we see another
advantage of the symplectic differential geometrical approach on the extended manifold
15 dimensional manifold $Z$ over the more basic approach on $\R^{3}$.
Another advantage of the symplectic differential geometrical approach is the
Frenet frame integration.
Basically, one has to integrate
\begin{equation}\label{integrate}
\left(
\begin{array}{l}
\dot{\ul{T}}(s) \\
\dot{\ul{N}}(s) \\
\dot{\ul{B}}(s)
\end{array}
\right)
=
\left(
\begin{array}{ccc}
0 & \kappa(s) &  0 \\
-\kappa(s) & 0 &\tau(s) \\
0 & -\tau(s) & 0
\end{array}
\right)
\left(
\begin{array}{l}
\ul{T}(s) \\
\ul{N}(s) \\
\ul{B}(s)
\end{array}
\right)
\end{equation}
where again $s \in [0,L]$ denotes arclength parametrization of the spatial part of the
curve with initial condition
\[
\left(
\begin{array}{l}
\ul{T}(0) \\
\ul{N}(0) \\
\ul{B}(0)
\end{array}
\right)
=
\left(
\begin{array}{c}
\ul{e}_{z} \\
\frac{1}{\|\ul{z}_0\|}\left( (z_0)_{2} \ul{e}_{x} - (z_0)_{1}
\ul{e}_{y}\right)
 \\
\frac{1}{\|\ul{z}_0\|}\left( (z_0)_{1} \ul{e}_{x} +(z_0)_{2}
\ul{e}_{y}\right)
\end{array}
\right)
.
\]
where we recall that the angles in the second chart $(\tilde{\beta}_{0},\tilde{\gamma}_{0})$
can be computed from $\ul{n}_{0}$ by means of Eq.~\!(\ref{anglesfromn}).
Now everything is expressed in the five constants
\[
L, \ul{z}_{0}=((z_{0})_{1},(z_{0})_{2})^{T} \textrm{ and }\ul{z}_{0}'=((z_{0}')_{1},(z_{0}')_{2})^{T},
\]
which yield five degrees of freedom needed
to meet the boundary condition
\[
\ul{x}(L)=\ul{x}_{1} \in \R^{3} \textrm{ and }\dot{\ul{x}}=\ul{n}_{1} \in S^{2}.
\]
Now in order to investigate on how the length $L$ and the integration constants
$\ul{z}_{0}$ and $\ul{z}_{0}'$ depend on the boundary condition, we follow a similar
approach as in \cite[App.A]{DuitsAMS2} where we managed to obtain full control over the
equivalent\footnote{In \cite[App.A]{DuitsAMS2} we had to choose $z_{0}$ and $z_{0}'$ and $L$
such that $\ul{x}(L)=\ul{x}_{1} \in \R^{2}$ and $\dot{\ul{x}}(L)=\ul{n}_{1} \in S^{1}$.}
problem in $SE(2)$.
To this end we re-express the last 6 Pfaffian forms in (\ref{Pfaff2}) as
\begin{equation} \label{6}
{\rm d}\lambda = - \lambda m(g)^{-1} {\rm d}m(g) ,
\end{equation}
where we recall (\ref{mg}), where $m : SE(3) \to \R^{6 \times 6}$ forms a
(uncommon)
group representation as we will show next. Taking the outer product $\ul{y} \mapsto \sigma_{\ul{x}}\ul{y}=\ul{x}\times \ul{y}$
is tensorial with respect to rotations
and $\sigma_{R \ul{x}}= R \sigma_{\ul{x}} R^{-1}$ indeed implies
\[
m(g_{1} g_{2})= m(g_{1}) m(g_{2}) \textrm{ and }m(e)=I
\]
for all $g_{1},g_{2} \in SE(3)$. The identity (\ref{6}) follows from Eq.~\!(\ref{Pfaff2}) by means of
\[
m(g)^{-1} {\rm d}m(g)
=\left(
\begin{array}{cc}
R^{-1} {\rm d}R &  R^{-1} {\rm d}\sigma_{\ul{x}} R \\
0 & R^{-1} {\rm d}R
\end{array}
\right)=\left(
\begin{array}{cc}
R^{-1} {\rm d}R & \sigma_{R^{-1} {\rm d}\ul{x}}  \\
0 & R^{-1} {\rm d}R
\end{array}
\right)
\]
and recall from Eq.~\!(\ref{duals2}) that $R^{-1} {\rm d}\ul{x}=({\rm d}\mathcal{A}^{1}, {\rm d}\mathcal{A}^{2}, {\rm d}\mathcal{A}^{3})^T$ and
where the matrix
$g^{-1}dg$ contains $\tilde{\omega}^{k}_{j}= \sum \limits_{i} c^{k}_{ij} {\rm d}\mathcal{A}^{j}$, recall Theorem \ref{th:CartanSE2},
as elements. 
Now from (\ref{6}) we deduce that
\[
{\rm d}(\lambda \, m(g)^{-1}) = \lambda {\rm d}(m(g)^{-1}) +({\rm d}\lambda) (m(g))^{-1}=(\lambda (m(g))^{-1} {\rm d}m(g) + {\rm d}\lambda )\,(m(g))^{-1}=0
\]
and the result (\ref{mu}) follows. $\hfill \Box$
\begin{theorem}
The Frenet system (\ref{integrate}) for the stationary curves can be integrated analytically. By means of Eq.~(\ref{mu}) we find
\begin{equation}\label{een1}
m(g)  =h_{0}^{-1} m(\tilde{g}), \quad g=(x,y,z,R),
\end{equation}
where $h_{0}^{-1} \in SE(3)$ is given by
\begin{equation} \label{twee2}
h_{0}^{-1}=\left(
\begin{array}{cc}
\overline{R}_{0}^{T} & \sigma_{-\overline{R}_{0}^{T} \overline{\ul{x}}_{0}} \overline{R}_{0}^{T} \\
0 & \overline{R}_{0}^{T}
\end{array}
\right)=(m(\overline{\ul{x}}_{0},\overline{R}_{0}))^{-1}=m((\overline{\ul{x}}_{0},\overline{R}_{0})^{-1})
\end{equation}
with $\overline{\ul{x}}_{0}=(\overline{x}_0,\overline{y}_0,\overline{z}_0) \in \R^{3}$ and $\overline{R}_{0} \in SO(3)$ given by
\begin{equation} \label{drie3}
\begin{array}{l}
\overline{x}_{0}=0  \ (\textrm{free}), \\
\overline{y}_{0}=\frac{W\sqrt{c^2 \beta^2-\|\ul{z}_{0}'\|^2}}{c^2 \beta^2 \|\ul{z}_{0}'\|},   \\
\overline{z}_{0}= \frac{\ul{z}_{0}\cdot \ul{z}_{0}'}{c \beta \|\ul{z}_{0}'\|},
\end{array}
\end{equation}
with Wronskian $W=z_{1}(s)\dot{z}_{2}(s)-z_{2}(s)\dot{z}_{1}(s)=z_{1}(0)\dot{z}_{2}(0)-z_{2}(0)\dot{z}_{1}(0)$ and
\begin{equation}\label{four4}
\overline{R}_{0}=
\frac{1}{c\beta}
\left(
\begin{array}{ccc}
-(z_{0}')_{2} & (z_{0}')_{1} & \sqrt{c^2 \beta^{2} -\|\ul{z}_{0}'\|^2}  \\
- c\beta \frac{(z_{0}')_1}{\|\ul{z}_{0}'\|}&  -c\beta \frac{(z_{0}')_2}{\|\ul{z}_{0}'\|}&  0 \\
\frac{(z_{0}')_2 \sqrt{c^{2}\beta^2 - \|\ul{z}_{0}'\|^2}}{\|\ul{z}_{0}'\|} & -\frac{(z_{0}')_{1} \sqrt{c^{2}\beta^2 - \|\ul{z}_{0}'\|^2}}{\|\ul{z}_{0}'\|} & \|\ul{z}_{0}'\|
\end{array}
\right)
\end{equation}
with $\ul{z}_{0}'=((z_{0}')_{1},(z_{0}')_{2})$ if $\ul{z}_{0}'\neq \ul{0}$. In case
$\ul{z}_{0}'=\ul{0}$ we have
\begin{equation} \label{five5}
\overline{\ul{x}}_{0}=-\frac{1}{c\beta}
\left(
\begin{array}{c}
0 \\
(z_{0})_{1} \\
(z_{0})_{2}
\end{array}
\right)
\textrm{ and }\overline{R}_{0}=R_{\ul{e}_{y},\frac{\pi}{2}}.
\end{equation}
Furthermore, we have
\[
\tilde{g}(s)=(\tilde{\ul{x}}(s),\tilde{R}(s))=(\tilde{x}(s),\tilde{y}(s),\tilde{z}(s),\tilde{R}(s))
\]
for all $s \in [0,L]$ with
\[
\begin{array}{l}
\tilde{x}(s)= \tilde{x}(0)-\frac{i\sqrt{1-\gamma}}{\beta}\frac{\sqrt{1+c^{2}}}{c\sqrt{2}}\, \left(E\left((\beta s+ \frac{\varphi}{2})i,m\right)-
E\left((\frac{\varphi}{2})i,m\right)\right), \\
(\tilde{y}(s), \tilde{z}(s))^T= e^{\int \limits_{0}^{s} A(s') {\rm d}s'} (\tilde{y}(0), \tilde{z}(0))^T=
\textrm{Re}\left\{(\tilde{z}_{0}-i \tilde{y}_{0}) e^{\int_{0}^{s} \frac{\|\dot{z}(\tau)\|^{2}}{\ul{z}(\tau)\cdot\dot{\ul{z}}(\tau)-
i W \sqrt{1- c^{-2}\beta^{-2}\|\dot{\ul{z}}(\tau)\|^2}}{\rm d}\tau }
\left(
\begin{array}{c}
i \\
1
\end{array}
\right) \right\}
\end{array}
\]
where the elliptic integral of the second kind is given by
$E(\phi,m)=\int \limits_{0}^{\phi} \sqrt{1- m \sin^{2} \theta} {\rm d}\theta$
with $m=\frac{2\gamma}{\gamma-1}$ and with $\gamma=\frac{\|\ul{z}_0 -\beta^{-1} \ul{z}_{0}'\|\|\ul{z}_0 +\beta^{-1} \ul{z}_{0}'\|}{1+c^2} \leq 1$, with
$\phi = \log \left(\frac{\|\ul{z}_{0}+\beta^{-1}\ul{z}_{0}'\|}{\|\ul{z}_{0}-\beta^{-1}\ul{z}_{0}'\|} \right)$ and 
with (commuting)
\[
A(s)= \frac{1}{ |\ul{z}(s)|^2- \frac{W^{2}}{c^2 \beta^{2}}}
\left(
\begin{array}{cc}
 \ul{z}(s) \cdot \dot{\ul{z}}(s) &  -\frac{W}{c\beta} \sqrt{c^{2}\beta^{2}- \|\dot{\ul{z}}(s)\|^2}  \\
\frac{W}{c\beta} \sqrt{c^{2}\beta^{2}- \|\dot{\ul{z}}(s)\|^2} &  \ul{z}(s) \cdot \dot{\ul{z}}(s)
\end{array}
\right)
\]
and with $\tilde{R}(s)$ in $SO(3)$ such that
\[
\tilde{R}(s)\ul{e}_{z}=\ul{n}(s)=
\left(\dot{\tilde{x}}(s), (A(s) (\tilde{y}(s),\tilde{z}(s))^{T})^T\right)^{T}=
\left(\frac{\lambda_{3}(s)}{c\beta}, (A(s) (\tilde{y}(s),\tilde{z}(s))^{T})^T\right)^{T}
.
\]
\end{theorem}
\textbf{Proof }
With respect to Eq.~(\ref{een1}),(\ref{twee2}), (\ref{drie3}), (\ref{four4}) and (\ref{five5}) we note
that by means of Eq.~(\ref{mu}) we have
\[
\lambda (m(g))^{-1} =
\mu \desda \lambda= \mu m(g) = \mu h_{0}^{-1} m(\tilde{g}) = \lambda(0) (m(g(0)))^{-1} h_{0}^{-1} m(\tilde{g})=\lambda(0) h_{0}^{-1} m(\tilde{g})
\]
as $g(0)=(\ul{0},I) \in SE(3)$ implies that its corresponding matrix representation equals the
$6\times6$-identity matrix $m(g(0))=I_{6}$.
Then we choose $h_{0}^{-1}$ such that
\[
\mu h_{0}^{-1} = \lambda(0) h_{0}^{-1} = (c\beta,0,0,-\frac{W}{c\beta},0,0)
\]
with $\mu_{i}=\lambda_{i}(0)$, $i=1,\ldots,6$, which follows from the previous theorem.
As a result $h_{0}$ now follows from
\begin{equation}\label{ess}
\begin{array}{l}
h_{0}^{-T} \left(
\begin{array}{c}
\mu_{1} \\
\mu_{2} \\
\mu_{3} \\
\mu_{4} \\
\mu_{5} \\
\mu_{6}
\end{array}
\right)
=\left(
\begin{array}{c}
c\beta \\
0 \\
0 \\
-\frac{W}{c\beta} \\
0 \\
0
\end{array}
\right) \desda
\left\{
\begin{array}{l}
\overline{R}_{0} \left( \begin{array}{c} (z_{0}')_{1} \\
(z_{0}')_{2} \\
\sqrt{c^{2}\beta^{2}-\|\ul{z}_{0}'\|^2}
\end{array}\right)=
\left(
\begin{array}{c}
c\beta \\
0 \\
0
\end{array}
\right) \\
\overline{\ul{x}}_{0} \times \left(\begin{array}{c}
c\beta \\
0 \\
0
\end{array} \right) + \overline{R}_{0}
\left(\begin{array}{c}
(z_{0})_{1} \\
(z_{0})_{2} \\
0
\end{array} \right)
=
\left(\begin{array}{c}
-\frac{W}{c\beta} \\
0 \\
0
\end{array} \right).
\end{array}
\right.
\end{array}
\end{equation}
Next we shall derive the ODE for $(\tilde{y},\tilde{z})$.
We parameterize by 
$
\tilde{R}=R_{\ul{e}_{x},\tilde{\gamma}}
R_{\ul{e}_{y},\tilde{\beta}} R_{\ul{e}_{z},\tilde{\alpha}} \in SO(3)$
and substitution of this parametrization in
\[
\lambda = \mu \, h_{0}^{-1} \, m(\tilde{g})= (c\beta,0,0,-\frac{W}{c\beta},0,0)\, m(\tilde{g}),
\]
yields the system
\begin{equation} \label{sys1}
\tilde{R} \left(
\begin{array}{c}
\lambda_{1}\\
\lambda_{2} \\
\lambda_{3}
\end{array}
\right)=
\left(
\begin{array}{c}
c\beta \\
0 \\
0
\end{array}
\right)
\end{equation}
and the system
{\small
\begin{equation}\label{sys2}
\begin{array}{l}
\frac{1}{c\beta} R_{\ul{e}_{x},\tilde{\gamma}}
\left( \begin{array}{ccc}
\sqrt{\lambda_{1}^2+\lambda_{2}^2} & 0 & \lambda_{3} \\
0 & c \beta & 0 \\
-\lambda_{3} & 0 & \sqrt{\lambda_{1}^{2}+\lambda_{2}^{2}}
\end{array}\right)
\left( \begin{array}{ccc}
\frac{\lambda_{1}}{\sqrt{\lambda_{1}^2+\lambda_{2}^2}} & \frac{\lambda_{2}}{\sqrt{\lambda_{1}^2+\lambda_{2}^2}} & 0 \\
-\frac{\lambda_{2}}{\sqrt{\lambda_{1}^2+\lambda_{2}^2}} & \frac{\lambda_{1}}{\sqrt{\lambda_{1}^2+\lambda_{2}^2}} & 0 \\
0 & 0 & 1
\end{array}\right)
\left(
\begin{array}{c}
\lambda_{4} \\
\lambda_{5} \\
0
\end{array}
\right)
=
\left(
\begin{array}{c}
-\frac{W}{c\beta} \\
-c\beta \tilde{z} \\
c\beta\tilde{y}
\end{array}
\right)
\end{array}
\end{equation}
}
where we used $\tilde{\beta}=\arg(\sqrt{\lambda_{1}^{2}+\lambda_{2}^{2}}+i\lambda_{3})$ and
$\tilde{\alpha}=\arg (\lambda_{1}-i \lambda_{2})$ that follow by Eq.~(\ref{sys1}).
Taking the derivative with respect to spatial arc-length $s$ in system (\ref{sys2}) yields
\[
\frac{d}{ds} \left( \tilde{R} (\lambda^{(2)})^{T}\right)=
\tilde{R} (\dot{\lambda}^{(2)})^{T} +
\dot{\tilde{R}} (\lambda^{(2)})^T =
\tilde{R} (\dot{\lambda}^{(2)})^{T} =c\beta
\left(
\begin{array}{c}
0 \\
-\dot{\tilde{z}} \\
\dot{\tilde{y}}
\end{array}
\right) , \quad \textrm{ with }\lambda^{(2)}=(\lambda_{4},\lambda_{5},0),
\]
where we stress that
\[
\begin{array}{l}
\tilde{R}^{T} {\rm d}\tilde{R}= \sigma_{({\rm d}\mathcal{A}^{4},{\rm d}\mathcal{A}^{5},0)}= \sigma_{(\kappa^{1}{\rm d}s,\kappa^{2}{\rm d}s,0)} \Rightarrow \\
\dot{\tilde{R}} (\lambda^{(2)})^{T} =\tilde{R} \, \sigma_{(\kappa^{1},\kappa^{2},0)^{T}}
\left(
\begin{array}{c}
\lambda_{4} \\
\lambda_{5} \\
0
\end{array}
\right)= \tilde{R} \left(
\left(
\begin{array}{c}
\kappa^{1} \\
\kappa^{2} \\
0
\end{array}
\right)
\times
\left(
\begin{array}{c}
\kappa^{1} \\
\kappa^{2} \\
0
\end{array}
\right)
\right) \frac{1}{\sqrt{(\kappa^{1})^{2}+(\kappa^{2})^{2} +\beta^2}}=\ul{0}.
\end{array}
\]
As a result, since $\sqrt{(\lambda_{1})^{2}+(\lambda_{2})^{2}}=\|\dot{\ul{z}}\|$, we obtain
\begin{equation}
\dot{\tilde{z}}= \frac{\cos \tilde{\gamma} \|\dot{\ul{z}}\|}{c\beta }
\textrm{ and }\dot{\tilde{y}}= \frac{-\sin \tilde{\gamma}\|\dot{\ul{z}}\| }{c\beta},
\end{equation}
whereas (\ref{sys2}) yields
\[
\begin{array}{l}
\frac{-\ul{z} \cdot \dot{\ul{z}}}{\|\dot{\ul{z}}\|}\, \cos \tilde{\gamma} -
\frac{W \sqrt{c^{2}\beta^{2}-\|\dot{\ul{z}}\|^2}}{c\beta\|\dot{\ul{z}}\|} \, \sin \tilde{\gamma} = -c\beta \tilde{z}, \\
\frac{-\ul{z} \cdot \dot{\ul{z}}}{\|\dot{\ul{z}}\|}\, \sin \tilde{\gamma} +
\frac{W \sqrt{c^{2}\beta^{2}-\|\dot{\ul{z}}\|^2}}{c\beta\|\dot{\ul{z}}\|} \, \cos \tilde{\gamma} = c\beta \tilde{y} \\
\end{array}
\]
from which we deduce that
\[
(A(s))^{-1}
\left(
\begin{array}{c}
\dot{\tilde{y}}(s)\\
\dot{\tilde{z}}(s)
\end{array}
\right)=
\left(
\begin{array}{c}
\tilde{y}(s)\\
\tilde{z}(s)
\end{array}
\right) \desda
\left(
\begin{array}{c}
\dot{\tilde{y}}(s)\\
\dot{\tilde{z}}(s)
\end{array}
\right)=A(s)
\left(
\begin{array}{c}
\tilde{y}(s)\\
\tilde{z}(s)
\end{array}
\right).
\]
As the matrices $A(s)$ share the same eigen-vectors ($(i,1)^{T}$ and $(-i,1)^{T}$) they commute and we have
\[
\begin{array}{l}
\left(
\begin{array}{c}
\tilde{y}(s)\\
\tilde{z}(s)
\end{array}
\right)= e^{\int_{0}^{s} A(\tau)\,{\rm d}\tau}
\left(
\begin{array}{c}
\tilde{y}(0)\\
\tilde{z}(0)
\end{array}
\right)=
\frac{1}{2i}(\tilde{y}(0)+ i\tilde{z}(0)) e^{\int_{0}^{s} \frac{\|\dot{\ul{z}}(\tau)\|^{2}}{\dot{\ul{z}}(\tau)\cdot\ul{z}(\tau)-i W \lambda_{3}(\tau)/(c\beta) }{\rm d}\tau
}\left(
\begin{array}{c}
i \\
1
\end{array}
\right) \\ +
\frac{1}{2i}(-\tilde{y}(0)+ i\tilde{z}(0)) e^{\int_{0}^{s} \frac{\|\dot{\ul{z}}(\tau)\|^{2}}{\dot{\ul{z}}(\tau)\cdot\ul{z}(\tau)+i W \lambda_{3}(\tau)/(c\beta)}{\rm d}\tau
}\left(
\begin{array}{r}
-i \\
1
\end{array}
\right)
\end{array}
\]
where the latter two vectors are each other's conjugate. Finally, $\tilde{x}(s)$ follows by integration
\[
\begin{array}{ll}
\tilde{x}(s)&=\tilde{x}(0)+ \int \limits_{0}^{s} \frac{\lambda_{3}(\tau)}{c\beta}\, {\rm d}\tau \\
            &=\tilde{x}(0)+ \frac{1}{c}\int \limits_{0}^{s} \sqrt{1-\|\ul{z}(\tau)\|^{2}}\, {\rm d}\tau \\
            &=\tilde{x}(0)+ \frac{\sqrt{1+c^2}}{c \sqrt{2}}\int \limits_{0}^{s} \sqrt{1- c_{1} \cosh(2\beta \tau)-c_{2} \sinh(2\beta \tau)}\, {\rm d}\tau
\end{array}
\]
with $c_{1}=\frac{1}{(1+c^2)}(\|\ul{z}_{0}\|^{2}-\|\dot{\ul{z}}_{0}\|^{2})$, $c_{2}=\frac{2}{\beta(1+c^2)} \ul{z}_{0} \cdot \dot{\ul{z}}_{0}$, $m=\frac{2\gamma}{\gamma-1}$,
which can be expressed as
\[
\begin{array}{ll}
\tilde{x}(s)&=\tilde{x}(0)+ \frac{\sqrt{1+c^2}}{c \sqrt{2}}\int \limits_{0}^{s} \sqrt{1- \gamma \cos(2i \beta \tau +i \varphi)}\, {\rm d}\tau \\
 &= \tilde{x}_{0}+\frac{\sqrt{1-\gamma} }{\beta} \frac{\sqrt{1+c^2}}{c \sqrt{2}}\,  \frac{1}{i}\, \int \limits_{\frac{i \varphi}{2}}^{\beta s i + \frac{i \varphi}{2}} \sqrt{1-m \sin^{2}(v)} {\rm d}v \\
 &= \tilde{x}(0)-\frac{i\sqrt{1-\gamma}}{\beta}\frac{\sqrt{1+c^{2}}}{c\sqrt{2}}\, \left(E\left((\beta s+ \frac{\varphi}{2})i,m\right)-
E\left((\frac{\varphi}{2})i,m\right)\right),
\end{array}
\]
where $\varphi=\frac{1}{2} \log \frac{c_{1}+c_{2}}{c_{1}-c_{2}}$, $v=i\beta\tau$, $\gamma=\sqrt{(c_{1})^{2}-(c_{2})^2}$ and where we note that
\[
\gamma=
\frac{1}{1+c^2}\|\ul{z}_{0}-\beta^{-1}\dot{\ul{z}}_{0}\|\|\ul{z}_{0}+\beta^{-1}\dot{\ul{z}}_{0}\| \leq \frac{1}{1+c^{2}}\frac{1}{2}\left(
\|\ul{z}_{0}-\beta^{-1}\dot{\ul{z}}_{0}\|^{2}+\|\ul{z}_{0}+\beta^{-1}\dot{\ul{z}}_{0}\|^{2}\right)=\frac{1-c^{2}}{1+c^{2}} \leq 1. \hfill \Box
\]
\begin{corollary}
The solution curves stay in a plane iff $\ul{z}_{0}$ and $\ul{z}_{0}'$ are linear dependent, i.e.
$W=\det (\ul{z}_{0}|\ul{z}_{0}')=0$.
In the special case $\ul{z}_{0}'=-\beta \ul{z}_{0}$, ($c=1 \textrm{ and }W=0$) we find
\[
\begin{array}{l}
\tilde{x}(s)= \tilde{x}(0) + s + \frac{1}{c\beta} \left( \sqrt{1-\|\ul{z}_{0}\|^2}-\sqrt{1-\|\ul{z}_{0}\|^2e^{-2s\beta}}
\; + \;
\log
\left( \frac{1+\sqrt{1-\|\ul{z}_{0}\|^2e^{-2s\beta}}e^{-2\beta s}}{1+\sqrt{1-\|\ul{z}_{0}\|^2}} \right)
\right)       , \\
(\tilde{y}(s), \tilde{z}(s))^T= e^{-\beta s} (\tilde{y}(0), \tilde{z}(0))^T
\end{array}
\]
In the case where solution curves stay in the plane, i.e. $W=0$, with general $c$, we find
\[
\begin{array}{l}
\tilde{x}(s)= \tilde{x}(0) +\frac{1}{c} \int \limits_{0}^{s} \sqrt{1- \|\ul{z}(\tau)\|^2}\, {\rm d}\tau , \\
\tilde{y}(s)= 0
\end{array}
\]
and
\[
\tilde{z}(s)= \left\{
\begin{array}{ll}
\tilde{z}(0) \frac{\|\ul{z}(s)\|}{\|\ul{z}(0)\|} &\textrm{ for } c\leq 1\\
sgn(s_0-s)\tilde{z}(0) \frac{\|\ul{z}(s)\|}{\|\ul{z}(0)\|} &\textrm{ for } c>1
\end{array}
\right.
\]
where $\ul{z}(s)$ is given by Eq.~(\ref{zs}) and $s_0=\frac{1}{2\beta}(ln(\|\ul{z}_{0}'\|^2-\beta^2\|\ul{z}_{0}\|^2)-ln(\|\ul{z}_{0}'+\beta\ul{z}_{0}\|^2))$ so that $\|\ul{z}(s_0)\|=0$. These solutions indeed coincide with the geodesics on the contact manifold
$(SE(2)\equiv \R^{2}\rtimes S^{1}, -\sin \theta {\rm d}x + \cos \theta {\rm d}y)$ where elements are parameterized by $(x,y, e^{i \theta})$,
cf.~\cite[App.A]{DuitsIJCV2010}.
\end{corollary}
See Figure \ref{fig:examplenot} for an example of a geodesic with $W=0$ and $c=1$.

\end{document}